
\ifx\shlhetal\undefinedcontrolsequence\let\shlhetal\relax\fi


\input amstex
\expandafter\ifx\csname mathdefs.tex\endcsname\relax
  \expandafter\gdef\csname mathdefs.tex\endcsname{}
\else \message{Hey!  Apparently you were trying to
  \string\input{mathdefs.tex} twice.   This does not make sense.} 
\errmessage{Please edit your file (probably \jobname.tex) and remove
any duplicate ``\string\input'' lines}\endinput\fi




\catcode`\X=12\catcode`\@=11

\def\n@wcount{\alloc@0\count\countdef\insc@unt}
\def\n@wwrite{\alloc@7\write\chardef\sixt@@n}
\def\n@wread{\alloc@6\read\chardef\sixt@@n}
\def\r@s@t{\relax}\def\v@idline{\par}\def\@mputate#1/{#1}
\def\l@c@l#1X{\firstpart.#1}\def\gl@b@l#1X{#1}\def\t@d@l#1X{{}}

\def\crossrefs#1{\ifx\all#1\let\tr@ce=\all\else\def\tr@ce{#1,}\fi
   \n@wwrite\cit@tionsout\openout\cit@tionsout=\jobname.cit 
   \write\cit@tionsout{\tr@ce}\expandafter\setfl@gs\tr@ce,}
\def\setfl@gs#1,{\def\@{#1}\ifx\@\empty\let\next=\relax
   \else\let\next=\setfl@gs\expandafter\xdef
   \csname#1tr@cetrue\endcsname{}\fi\next}
\def\m@ketag#1#2{\expandafter\n@wcount\csname#2tagno\endcsname
     \csname#2tagno\endcsname=0\let\tail=\all\xdef\all{\tail#2,}
   \ifx#1\l@c@l\let\tail=\r@s@t\xdef\r@s@t{\csname#2tagno\endcsname=0\tail}\fi
   \expandafter\gdef\csname#2cite\endcsname##1{\expandafter
     \ifx\csname#2tag##1\endcsname\relax?\else\csname#2tag##1\endcsname\fi
     \expandafter\ifx\csname#2tr@cetrue\endcsname\relax\else
     \write\cit@tionsout{#2tag ##1 cited on page \folio.}\fi}
   \expandafter\gdef\csname#2page\endcsname##1{\expandafter
     \ifx\csname#2page##1\endcsname\relax?\else\csname#2page##1\endcsname\fi
     \expandafter\ifx\csname#2tr@cetrue\endcsname\relax\else
     \write\cit@tionsout{#2tag ##1 cited on page \folio.}\fi}
   \expandafter\gdef\csname#2tag\endcsname##1{\expandafter
      \ifx\csname#2check##1\endcsname\relax
      \expandafter\xdef\csname#2check##1\endcsname{}%
      \else\immediate\write16{Warning: #2tag ##1 used more than once.}\fi
      \multit@g{#1}{#2}##1/X%
      \write\t@gsout{#2tag ##1 assigned number \csname#2tag##1\endcsname\space
      on page \number\count0.}%
   \csname#2tag##1\endcsname}}

\def\multit@g#1#2#3/#4X{\def\t@mp{#4}\ifx\t@mp\empty%
      \global\advance\csname#2tagno\endcsname by 1 
      \expandafter\xdef\csname#2tag#3\endcsname
      {#1\number\csname#2tagno\endcsnameX}%
   \else\expandafter\ifx\csname#2last#3\endcsname\relax
      \expandafter\n@wcount\csname#2last#3\endcsname
      \global\advance\csname#2tagno\endcsname by 1 
      \expandafter\xdef\csname#2tag#3\endcsname
      {#1\number\csname#2tagno\endcsnameX}
      \write\t@gsout{#2tag #3 assigned number \csname#2tag#3\endcsname\space
      on page \number\count0.}\fi
   \global\advance\csname#2last#3\endcsname by 1
   \def\t@mp{\expandafter\xdef\csname#2tag#3/}%
   \expandafter\t@mp\@mputate#4\endcsname
   {\csname#2tag#3\endcsname\lastpart{\csname#2last#3\endcsname}}\fi}
\def\t@gs#1{\def\all{}\m@ketag#1e\m@ketag#1s\m@ketag\t@d@l p
\let\realscite\scite
\let\realstag\stag
   \m@ketag\gl@b@l r \n@wread\t@gsin
   \openin\t@gsin=\jobname.tgs \re@der \closein\t@gsin
   \n@wwrite\t@gsout\openout\t@gsout=\jobname.tgs }
\outer\def\localtags{\t@gs\l@c@l}
\outer\def\globaltags{\t@gs\gl@b@l}
\outer\def\newlocaltag#1{\m@ketag\l@c@l{#1}}
\outer\def\newglobaltag#1{\m@ketag\gl@b@l{#1}}

\newif\ifpr@ 
\def\m@kecs #1tag #2 assigned number #3 on page #4.%
   {\expandafter\gdef\csname#1tag#2\endcsname{#3}
   \expandafter\gdef\csname#1page#2\endcsname{#4}
   \ifpr@\expandafter\xdef\csname#1check#2\endcsname{}\fi}
\def\re@der{\ifeof\t@gsin\let\next=\relax\else
   \read\t@gsin to\t@gline\ifx\t@gline\v@idline\else
   \expandafter\m@kecs \t@gline\fi\let \next=\re@der\fi\next}
\def\pretags#1{\pr@true\pret@gs#1,,}
\def\pret@gs#1,{\def\@{#1}\ifx\@\empty\let\n@xtfile=\relax
   \else\let\n@xtfile=\pret@gs \openin\t@gsin=#1.tgs \message{#1} \re@der 
   \closein\t@gsin\fi \n@xtfile}

\newcount\sectno\sectno=0\newcount\subsectno\subsectno=0
\newif\ifultr@local \def\ultralocal{\ultr@localtrue}
\def\firstpart{\number\sectno}
\def\lastpart#1{\ifcase#1 \or a\or b\or c\or d\or e\or f\or g\or h\or 
   i\or k\or l\or m\or n\or o\or p\or q\or r\or s\or t\or u\or v\or w\or 
   x\or y\or z \fi}

\def\resetall{\global\advance\sectno by 1\subsectno=0
   \gdef\firstpart{\number\sectno}\r@s@t}
\def\resetsub{\global\advance\subsectno by 1
   \gdef\firstpart{\number\sectno.\number\subsectno}\r@s@t}
\def\newsection#1\par{\resetall\vskip0pt plus.3\vsize\penalty-250
   \vskip0pt plus-.3\vsize\bigskip\bigskip
   \message{#1}\leftline{\bf#1}\nobreak\bigskip}
\def\subsection#1\par{\ifultr@local\resetsub\fi
   \vskip0pt plus.2\vsize\penalty-250\vskip0pt plus-.2\vsize
   \bigskip\smallskip\message{#1}\leftline{\bf#1}\nobreak\medskip}


\newdimen\marginshift

\newdimen\margindelta
\newdimen\marginmax
\newdimen\marginmin

\def\margininit{       
\marginmax=3 true cm                  
				      
\margindelta=0.1 true cm              
\marginmin=0.1true cm                 
\marginshift=\marginmin
}    

\def\t@gsjj#1,{\def\@{#1}\ifx\@\empty\let\next=\relax\else\let\next=\t@gsjj
   \def\@@{p}\ifx\@\@@\else
   \expandafter\gdef\csname#1cite\endcsname##1{\citejj{##1}}
   \expandafter\gdef\csname#1page\endcsname##1{?}
   \expandafter\gdef\csname#1tag\endcsname##1{\tagjj{##1}}\fi\fi\next}
\newif\ifshowstuffinmargin
\showstuffinmarginfalse
\def\jjtags{\ifx\shlhetal\relax 
  \else
\ifx\shlhetal\undefinedcontrolseq
\else
\showstuffinmargintrue
\ifx\all\relax\else\expandafter\t@gsjj\all,\fi\fi \fi
}

\def\tagjj#1{\realstag{#1}\oldmginpar{\zeigen{#1}}}
\def\citejj#1{\rechnen{#1}\mginpar{\zeigen{#1}}}     

\def\rechnen#1{\expandafter\ifx\csname stag#1\endcsname\relax ??\else
                           \csname stag#1\endcsname\fi}

\newdimen\theight

\def\marginfont{\sevenrm}

\def\trymarginbox#1{\setbox0=\hbox{\marginfont\hskip\marginshift #1}%
		\global\marginshift\wd0 
		\global\advance\marginshift\margindelta}

\def \oldmginpar#1{%
\ifvmode\setbox0\hbox to \hsize{\hfill\rlap{\marginfont\quad#1}}%
\ht0 0cm
\dp0 0cm
\box0\vskip-\baselineskip
\else 
             \vadjust{\trymarginbox{#1}%
		\ifdim\marginshift>\marginmax \global\marginshift\marginmin
			\trymarginbox{#1}%
                \fi
             \theight=\ht0
             \advance\theight by \dp0    \advance\theight by \lineskip
             \kern -\theight \vbox to \theight{\rightline{\rlap{\box0}}%
\vss}}\fi}

\newdimen\upordown
\global\upordown=8pt
\font\tinyfont=cmtt8 
\def\mginpar#1{\smash{\hbox to 0cm{\kern-10pt\raise7pt\hbox{\tinyfont #1}\hss}}}
\def\mginpar#1{{\hbox to 0cm{\kern-10pt\raise\upordown\hbox{\tinyfont #1}\hss}}\global\upordown-\upordown}


\def\t@gsoff#1,{\def\@{#1}\ifx\@\empty\let\next=\relax\else\let\next=\t@gsoff
   \def\@@{p}\ifx\@\@@\else
   \expandafter\gdef\csname#1cite\endcsname##1{\zeigen{##1}}
   \expandafter\gdef\csname#1page\endcsname##1{?}
   \expandafter\gdef\csname#1tag\endcsname##1{\zeigen{##1}}\fi\fi\next}
\def\verbatimtags{\showstuffinmarginfalse
\ifx\all\relax\else\expandafter\t@gsoff\all,\fi}
\def\zeigen#1{\hbox{$\scriptstyle\langle$}#1\hbox{$\scriptstyle\rangle$}}


\def\margintag#1{\ifshowstuffinmargin\oldmginpar{\zeigen{#1}}\fi}

\def\marginplain#1{\ifshowstuffinmargin\mginpar{{#1}}\fi}
\def\marginbf#1{\marginplain{{\bf \ \ #1}}}

\def\(#1){\edef\dot@g{\ifmmode\ifinner(\hbox{\noexpand\etag{#1}})
   \else\noexpand\eqno(\hbox{\noexpand\etag{#1}})\fi
   \else(\noexpand\ecite{#1})\fi}\dot@g}

\newif\ifbr@ck
\def\eat#1{}
\def\[#1]{\br@cktrue[\br@cket#1'X]}
\def\br@cket#1'#2X{\def\temp{#2}\ifx\temp\empty\let\next\eat
   \else\let\next\br@cket\fi
   \ifbr@ck\br@ckfalse\br@ck@t#1,X\else\br@cktrue#1\fi\next#2X}
\def\br@ck@t#1,#2X{\def\temp{#2}\ifx\temp\empty\let\neext\eat
   \else\let\neext\br@ck@t\def\temp{,}\fi
   \def\teemp{#1}\ifx\teemp\empty\else\rcite{#1}\fi\temp\neext#2X}
\def\resetbr@cket{\gdef\[##1]{[\rtag{##1}]}}
\def\references{\resetbr@cket\newsection References\par}

\newtoks\symb@ls\newtoks\s@mb@ls\newtoks\p@gelist\n@wcount\ftn@mber
    \ftn@mber=1\newif\ifftn@mbers\ftn@mbersfalse\newif\ifbyp@ge\byp@gefalse
\def\defm@rk{\ifftn@mbers\n@mberm@rk\else\symb@lm@rk\fi}
\def\n@mberm@rk{\xdef\m@rk{{\the\ftn@mber}}%
    \global\advance\ftn@mber by 1 }
\def\rot@te#1{\let\temp=#1\global#1=\expandafter\r@t@te\the\temp,X}
\def\r@t@te#1,#2X{{#2#1}\xdef\m@rk{{#1}}}
\def\b@@st#1{{$^{#1}$}}\def\str@p#1{#1}
\def\symb@lm@rk{\ifbyp@ge\rot@te\p@gelist\ifnum\expandafter\str@p\m@rk=1 
    \s@mb@ls=\symb@ls\fi\write\f@nsout{\number\count0}\fi \rot@te\s@mb@ls}
\def\byp@ge{\byp@getrue\n@wwrite\f@nsin\openin\f@nsin=\jobname.fns 
    \n@wcount\currentp@ge\currentp@ge=0\p@gelist={0}
    \re@dfns\closein\f@nsin\rot@te\p@gelist
    \n@wread\f@nsout\openout\f@nsout=\jobname.fns }
\def\m@kelist#1X#2{{#1,#2}}
\def\re@dfns{\ifeof\f@nsin\let\next=\relax\else\read\f@nsin to \f@nline
    \ifx\f@nline\v@idline\else\let\t@mplist=\p@gelist
    \ifnum\currentp@ge=\f@nline
    \global\p@gelist=\expandafter\m@kelist\the\t@mplistX0
    \else\currentp@ge=\f@nline
    \global\p@gelist=\expandafter\m@kelist\the\t@mplistX1\fi\fi
    \let\next=\re@dfns\fi\next}
\def\symbols#1{\symb@ls={#1}\s@mb@ls=\symb@ls} 
\def\bigsymbol{\textstyle}
\symbols{\bigsymbol\ast,\dagger,\ddagger,\sharp,\flat,\natural,\star}
\def\ftnumbers{\ftn@mberstrue} \def\ftsymbols{\ftn@mbersfalse}
\def\paginal{\byp@ge} \def\resetftnumbers{\ftn@mber=1}
\def\ftnote#1{\defm@rk\expandafter\expandafter\expandafter\footnote
    \expandafter\b@@st\m@rk{#1}}

\long\def\jump#1\endjump{}
\def\ssum{\mathop{\lower .1em\hbox{$\textstyle\Sigma$}}\nolimits}

\def\qed{\nobreak\kern 1em \vrule height .5em width .5em depth 0em}
\def\newneq{\hbox{\rlap{\hbox to 1\wd9{\hss$=$\hss}}\raise .1em 
   \hbox to 1\wd9{\hss$\scriptscriptstyle/$\hss}}}
\def\subsetne{\setbox9 = \hbox{$\subset$}\mathrel{\hbox{\rlap
   {\lower .4em \newneq}\raise .13em \hbox{$\subset$}}}}
\def\supsetne{\setbox9 = \hbox{$\subset$}\mathrel{\hbox{\rlap
   {\lower .4em \newneq}\raise .13em \hbox{$\supset$}}}}

\def\vbar{\mathchoice{\vrule height6.3ptdepth-.5ptwidth.8pt\kern-.8pt}
   {\vrule height6.3ptdepth-.5ptwidth.8pt\kern-.8pt}
   {\vrule height4.1ptdepth-.35ptwidth.6pt\kern-.6pt}
   {\vrule height3.1ptdepth-.25ptwidth.5pt\kern-.5pt}}
\def\f@dge{\mathchoice{}{}{\mkern.5mu}{\mkern.8mu}}
\def\b@c#1#2{{\rm \mkern#2mu\vbar\mkern-#2mu#1}}
\def\b@b#1{{\rm I\mkern-3.5mu #1}}
\def\b@a#1#2{{\rm #1\mkern-#2mu\f@dge #1}}
\def\bb#1{{\count4=`#1 \advance\count4by-64 \ifcase\count4\or\b@a A{11.5}\or
   \b@b B\or\b@c C{5}\or\b@b D\or\b@b E\or\b@b F \or\b@c G{5}\or\b@b H\or
   \b@b I\or\b@c J{3}\or\b@b K\or\b@b L \or\b@b M\or\b@b N\or\b@c O{5} \or
   \b@b P\or\b@c Q{5}\or\b@b R\or\b@a S{8}\or\b@a T{10.5}\or\b@c U{5}\or
   \b@a V{12}\or\b@a W{16.5}\or\b@a X{11}\or\b@a Y{11.7}\or\b@a Z{7.5}\fi}}

\catcode`\X=11 \catcode`\@=12




\let\thischap\jobname

\def\partof#1{\csname returnthe#1part\endcsname}
\def\chapof#1{\csname returnthe#1chap\endcsname}

\def\setchapter#1,#2,#3;{%
  \expandafter\def\csname returnthe#1part\endcsname{#2}%
  \expandafter\def\csname returnthe#1chap\endcsname{#3}%
}

\setchapter 300a,A,II.A;
\setchapter 300b,A,II.B;
\setchapter 300c,A,II.C;
\setchapter 300d,A,II.D;
\setchapter 300e,A,II.E;
\setchapter 300f,A,II.F;
\setchapter 300g,A,II.G;
\setchapter  E53,B,N;
\setchapter  88r,B,I;
\setchapter  600,B,III;
\setchapter  705,B,IV;
\setchapter  734,B,V;

\def\cprefix#1{
\edef\theotherpart{\partof{#1}}\edef\theotherchap{\chapof{#1}}%
\ifx\theotherpart\thispart
   \ifx\theotherchap\thischap 
    \else 
     \theotherchap%
    \fi
   \else 
     \theotherchap\fi}

\def\sectioncite[#1]#2{%
     \cprefix{#2}#1}

\def\chaptercite#1{Chapter \cprefix{#1}}

\edef\thispart{\partof{\thischap}}
\edef\thischap{\chapof{\thischap}}

\def\lastpage of '#1' is #2.{\expandafter\def\csname lastpage#1\endcsname{#2}}


\def\spuriousreset{}


\expandafter\ifx\csname citeadd.tex\endcsname\relax
\expandafter\gdef\csname citeadd.tex\endcsname{}
\else \message{Hey!  Apparently you were trying to
\string\input{citeadd.tex} twice.   This does not make sense.} 
\errmessage{Please edit your file (probably \jobname.tex) and remove
any duplicate ``\string\input'' lines}\endinput\fi

\sectno=-1   
\ifx\epsfannounce\undefined \def\epsfannounce{\immediate\write16}\fi
 \epsfannounce{This is `epsf.tex' v2.7k <10 July 1997>}%
\newread\epsffilein    
\newif\ifepsfatend     
\newif\ifepsfbbfound   
\newif\ifepsfdraft     
\newif\ifepsffileok    
\newif\ifepsfframe     
\newif\ifepsfshow      
\epsfshowtrue          
\newif\ifepsfshowfilename 
\newif\ifepsfverbose   
\newdimen\epsfframemargin 
\newdimen\epsfframethickness 
\newdimen\epsfrsize    
\newdimen\epsftmp      
\newdimen\epsftsize    
\newdimen\epsfxsize    
\newdimen\epsfysize    
\newdimen\pspoints     
\pspoints = 1bp        
\epsfxsize = 0pt       
\epsfysize = 0pt       
\epsfframemargin = 0pt 
\epsfframethickness = 0.4pt 
\def\epsfbox#1{\global\def\epsfllx{72}\global\def\epsflly{72}%
   \global\def\epsfurx{540}\global\def\epsfury{720}%
   \def\lbracket{[}\def\testit{#1}\ifx\testit\lbracket
   \let\next=\epsfgetlitbb\else\let\next=\epsfnormal\fi\next{#1}}%
%
%
\def\epsfgetlitbb#1#2 #3 #4 #5]#6{%
   \epsfgrab #2 #3 #4 #5 .\\%
   \epsfsetsize
   \epsfstatus{#6}%
   \epsfsetgraph{#6}%
}%
\def\epsfnormal#1{%
    \epsfgetbb{#1}%
    \epsfsetgraph{#1}%
}%
\newhelp\epsfnoopenhelp{The PostScript image file must be findable by
TeX, i.e., somewhere in the TEXINPUTS (or equivalent) path.}%
\def\epsfgetbb#1{%
%
%
    \openin\epsffilein=#1
    \ifeof\epsffilein
        \errhelp = \epsfnoopenhelp
        \errmessage{Could not open file #1, ignoring it}%
    \else                       
        {
            \chardef\other=12
            \def\do##1{\catcode`##1=\other}%
            \dospecials
            \catcode`\ =10
            \epsffileoktrue         
            \epsfatendfalse     
            \loop               
                \read\epsffilein to \epsffileline
                \ifeof\epsffilein 
                \epsffileokfalse 
            \else                
                \expandafter\epsfaux\epsffileline:. \\%
            \fi
            \ifepsffileok
            \repeat
            \ifepsfbbfound
            \else
                \ifepsfverbose
                    \immediate\write16{No BoundingBox comment found in %
                                    file #1; using defaults}%
                \fi
            \fi
        }
        \closein\epsffilein
    \fi                         
    \epsfsetsize                
    \epsfstatus{#1}%
}%
%
%
\def\epsfclipoff{\def\epsfclipstring{\ifepsfdraft\space clip\fi}}%
\epsfclipoff 
%
%
\def\epsfspecial#1{%
     \epsftmp=10\epsfxsize
     \divide\epsftmp\pspoints
     \ifnum\epsfrsize=0\relax
       \includegraphics{\ifepsfdraft}%
     \else
       \epsfrsize=10\epsfysize
       \divide\epsfrsize\pspoints
       \includegraphics{\ifepsfdraft}%
     \fi
}%
%
\def\epsfframe#1%
{%
  \leavevmode                   
  \setbox0 = \hbox{#1}%
  \dimen0 = \wd0                                
  \advance \dimen0 by 2\epsfframemargin         
  \advance \dimen0 by 2\epsfframethickness      
  \vbox
  {%
    \hrule height \epsfframethickness depth 0pt
    \hbox to \dimen0
    {%
      \hss
      \vrule width \epsfframethickness
      \kern \epsfframemargin
      \vbox {\kern \epsfframemargin \box0 \kern \epsfframemargin }%
      \kern \epsfframemargin
      \vrule width \epsfframethickness
      \hss
    }
    \hrule height 0pt depth \epsfframethickness
  }
}%
\def\epsfsetgraph#1%
{%
   %
   %
   \leavevmode
   \hbox{
     \ifepsfframe\expandafter\epsfframe\fi
     {\vbox to\epsfysize
     {%
        \ifepsfshow
            \vfil
            \hbox to \epsfxsize{\epsfspecial{#1}\hfil}%
        \else
            \vfil
            \hbox to\epsfxsize{%
               \hss
               \ifepsfshowfilename
               {%
                  \epsfframemargin=3pt 
                  \epsfframe{{\tt #1}}%
               }%
               \fi
               \hss
            }%
            \vfil
        \fi
     }%
   }}%
   %
   %
   \global\epsfxsize=0pt
   \global\epsfysize=0pt
}%
%
%
\def\epsfsetsize
{%
   \epsfrsize=\epsfury\pspoints
   \advance\epsfrsize by-\epsflly\pspoints
   \epsftsize=\epsfurx\pspoints
   \advance\epsftsize by-\epsfllx\pspoints
%
%
   \epsfxsize=\epsfsize{\epsftsize}{\epsfrsize}%
   \ifnum \epsfxsize=0
      \ifnum \epsfysize=0
        \epsfxsize=\epsftsize
        \epsfysize=\epsfrsize
        \epsfrsize=0pt
%
%
      \else
        \epsftmp=\epsftsize \divide\epsftmp\epsfrsize
        \epsfxsize=\epsfysize \multiply\epsfxsize\epsftmp
        \multiply\epsftmp\epsfrsize \advance\epsftsize-\epsftmp
        \epsftmp=\epsfysize
        \loop \advance\epsftsize\epsftsize \divide\epsftmp 2
        \ifnum \epsftmp>0
           \ifnum \epsftsize<\epsfrsize
           \else
              \advance\epsftsize-\epsfrsize \advance\epsfxsize\epsftmp
           \fi
        \repeat
        \epsfrsize=0pt
      \fi
   \else
     \ifnum \epsfysize=0
       \epsftmp=\epsfrsize \divide\epsftmp\epsftsize
       \epsfysize=\epsfxsize \multiply\epsfysize\epsftmp
       \multiply\epsftmp\epsftsize \advance\epsfrsize-\epsftmp
       \epsftmp=\epsfxsize
       \loop \advance\epsfrsize\epsfrsize \divide\epsftmp 2
       \ifnum \epsftmp>0
          \ifnum \epsfrsize<\epsftsize
          \else
             \advance\epsfrsize-\epsftsize \advance\epsfysize\epsftmp
          \fi
       \repeat
       \epsfrsize=0pt
     \else
       \epsfrsize=\epsfysize
     \fi
   \fi
}%
%
%
\def\epsfstatus#1{
   \ifepsfverbose
     \immediate\write16{#1: BoundingBox:
                  llx = \epsfllx\space lly = \epsflly\space
                  urx = \epsfurx\space ury = \epsfury\space}%
     \immediate\write16{#1: scaled width = \the\epsfxsize\space
                  scaled height = \the\epsfysize}%
   \fi
}%
%
%
{\catcode`\%=12 \global\let\epsfpercent=
\global\def\epsfatend{(atend)}%
%
%
%
%
%
%
%
\long\def\epsfaux#1#2:#3\\%
{%
   \def\testit{#2}
   \ifx#1\epsfpercent           
       \ifx\testit\epsfbblit    
            \epsfgrab #3 . . . \\%
            \ifx\epsfllx\epsfatend 
                \global\epsfatendtrue
            \else               
                \ifepsfatend    
                \else           
                    \epsffileokfalse
                \fi
                \global\epsfbbfoundtrue
            \fi
       \fi
   \fi
}%
%
%
\def\epsfempty{}%
\def\epsfgrab #1 #2 #3 #4 #5\\{%
   \global\def\epsfllx{#1}\ifx\epsfllx\epsfempty
      \epsfgrab #2 #3 #4 #5 .\\\else
   \global\def\epsflly{#2}%
   \global\def\epsfurx{#3}\global\def\epsfury{#4}\fi
}%
%
%
\def\epsfsize#1#2{\epsfxsize}%
%
%

\localtags
\jjtags
\newbox\noforkbox \newdimen\forklinewidth
\forklinewidth=0.3pt   
\setbox0\hbox{$\textstyle\bigcup$}
\setbox1\hbox to \wd0{\hfil\vrule width \forklinewidth depth \dp0
                        height \ht0 \hfil}
\wd1=0 cm
\setbox\noforkbox\hbox{\box1\box0\relax}
\def\unionstick{\mathop{\copy\noforkbox}\limits}
\def\nonfork#1#2_#3{#1\unionstick_{\textstyle #3}#2}
\def\nonforkin#1#2_#3^#4{#1\unionstick_{\textstyle #3}^{\textstyle #4}#2}     
%
\setbox0\hbox{$\textstyle\bigcup$}
\setbox1\hbox to \wd0{\hfil{\sl /\/}\hfil}
\setbox2\hbox to \wd0{\hfil\vrule height \ht0 depth \dp0 width
                                \forklinewidth\hfil}
\wd1=0cm
\wd2=0cm
\newbox\doesforkbox
\setbox\doesforkbox\hbox{\box1\box0\relax}
\def\nunionstick{\mathop{\copy\doesforkbox}\limits}

\def\fork#1#2_#3{#1\nunionstick_{\textstyle #3}#2}
\def\forkin#1#2_#3^#4{#1\nunionstick_{\textstyle #3}^{\textstyle #4}#2}     
\NoBlackBoxes
\define\ortp{\text{\bf tp}}

\define\sftp{\text{\rm tp}}

\define\mr{\medskip\roster}
\define\sn{\smallskip\noindent}
\define\mn{\medskip\noindent}
\define\bn{\bigskip\noindent}
\define\ub{\underbar}
\define\wilog{\text{without loss of generality}}
\define\ermn{\endroster\medskip\noindent}

\define\dbcu{\dsize\bigcup}
\define \nl{\newline}
\magnification=\magstep 1
\documentstyle{amsppt}

{    
\catcode`@11

\ifx\alicetwothousandloaded@\relax
  \endinput\else\global\let\alicetwothousandloaded@\relax\fi

\gdef\subjclass{\let\savedef@\subjclass
 \def\subjclass##1\endsubjclass{\let\subjclass\savedef@
   \toks@{\def\usualspace{{\rm\enspace}}\eightpoint}%
   \toks@@{##1\unskip.}%
   \edef\thesubjclass@{\the\toks@
     \frills@{{\noexpand\rm2000 {\noexpand\it Mathematics Subject
       Classification}.\noexpand\enspace}}%
     \the\toks@@}}%
  \nofrillscheck\subjclass}
} 


\expandafter\ifx\csname alice2jlem.tex\endcsname\relax
  \expandafter\xdef\csname alice2jlem.tex\endcsname{\the\catcode`@}
\else \message{Hey!  Apparently you were trying to
\string\input{alice2jlem.tex}  twice.   This does not make sense.}
\errmessage{Please edit your file (probably \jobname.tex) and remove
any duplicate ``\string\input'' lines}\endinput\fi

\expandafter\ifx\csname bib4plain.tex\endcsname\relax
  \expandafter\gdef\csname bib4plain.tex\endcsname{}
\else \message{Hey!  Apparently you were trying to \string\input
  bib4plain.tex twice.   This does not make sense.}
\errmessage{Please edit your file (probably \jobname.tex) and remove
any duplicate ``\string\input'' lines}\endinput\fi

\def\renewcommand{\newcommand}	       
\edef\cite{\the\catcode`@}%
\catcode`@ = 11
\let\@oldatcatcode = \cite
\chardef\@letter = 11
\chardef\@other = 12
%
%
%
%
\def\@innerdef#1#2{\edef#1{\expandafter\noexpand\csname #2\endcsname}}%
%
%
\@innerdef\@innernewcount{newcount}%
\@innerdef\@innernewdimen{newdimen}%
\@innerdef\@innernewif{newif}%
\@innerdef\@innernewwrite{newwrite}%
%
%
%
\def\@gobble#1{}%
%
%
%
\ifx\inputlineno\@undefined
   \let\@linenumber = \empty 
\else
   \def\@linenumber{\the\inputlineno:\space}%
\fi
%
%
%
\def\@futurenonspacelet#1{\def\cs{#1}%
   \afterassignment\@stepone\let\@nexttoken=
}%
\begingroup 
\def\\{\global\let\@stoken= }%
\\ 
\endgroup
\def\@stepone{\expandafter\futurelet\cs\@steptwo}%
\def\@steptwo{\expandafter\ifx\cs\@stoken\let\@@next=\@stepthree
   \else\let\@@next=\@nexttoken\fi \@@next}%
\def\@stepthree{\afterassignment\@stepone\let\@@next= }%
%
%
%
\def\@getoptionalarg#1{%
   \let\@optionaltemp = #1%
   \let\@optionalnext = \relax
   \@futurenonspacelet\@optionalnext\@bracketcheck
}%
%
%
\def\@bracketcheck{%
   \ifx [\@optionalnext
      \expandafter\@@getoptionalarg
   \else
      \let\@optionalarg = \empty
      \expandafter\@optionaltemp
   \fi
}%
\def\@@getoptionalarg[#1]{%
   \def\@optionalarg{#1}%
   \@optionaltemp
}%
%
%
%
\def\@nnil{\@nil}%
\def\@fornoop#1\@@#2#3{}%
\def\@for#1:=#2\do#3{%
   \edef\@fortmp{#2}%
   \ifx\@fortmp\empty \else
      \expandafter\@forloop#2,\@nil,\@nil\@@#1{#3}%
   \fi
}%
\def\@forloop#1,#2,#3\@@#4#5{\def#4{#1}\ifx #4\@nnil \else
       #5\def#4{#2}\ifx #4\@nnil \else#5\@iforloop #3\@@#4{#5}\fi\fi
}%
\def\@iforloop#1,#2\@@#3#4{\def#3{#1}\ifx #3\@nnil
       \let\@nextwhile=\@fornoop \else
      #4\relax\let\@nextwhile=\@iforloop\fi\@nextwhile#2\@@#3{#4}%
}%
%
%
%
\@innernewif\if@fileexists
\def\@testfileexistence{\@getoptionalarg\@finishtestfileexistence}%
\def\@finishtestfileexistence#1{%
   \begingroup
      \def\extension{#1}%
      \immediate\openin0 =
         \ifx\@optionalarg\empty\jobname\else\@optionalarg\fi
         \ifx\extension\empty \else .#1\fi
         \space
      \ifeof 0
         \global\@fileexistsfalse
      \else
         \global\@fileexiststrue
      \fi
      \immediate\closein0
   \endgroup
}%
%
%
%
%
\def\bibliographystyle#1{%
   \@readauxfile
   \@writeaux{\string\bibstyle{#1}}%
}%
\let\bibstyle = \@gobble
%
%
\let\bblfilebasename = \jobname
\def\bibliography#1{%
   \@readauxfile
   \@writeaux{\string\bibdata{#1}}%
   \@testfileexistence[\bblfilebasename]{bbl}%
   \if@fileexists
      \nobreak
      \@readbblfile
   \fi
}%
\let\bibdata = \@gobble
%
%
\def\nocite#1{%
   \@readauxfile
   \@writeaux{\string\citation{#1}}%
}%
\@innernewif\if@notfirstcitation
%
%
\def\cite{\@getoptionalarg\@cite}%
%
%
\def\@cite#1{%
   \let\@citenotetext = \@optionalarg
   \printcitestart
   \nocite{#1}%
   \@notfirstcitationfalse
   \@for \@citation :=#1\do
   {%
      \expandafter\@onecitation\@citation\@@
   }%
   \ifx\empty\@citenotetext\else
      \printcitenote{\@citenotetext}%
   \fi
   \printcitefinish
}%
\newif\ifweareinprivate
\weareinprivatetrue
\ifx\shlhetal\undefinedcontrolseq\weareinprivatefalse\fi
\ifx\shlhetal\relax\weareinprivatefalse\fi
\def\@onecitation#1\@@{%
   \if@notfirstcitation
      \printbetweencitations
   \fi
   \expandafter \ifx \csname\@citelabel{#1}\endcsname \relax
      \if@citewarning
         \message{\@linenumber Undefined citation `#1'.}%
      \fi
     \ifweareinprivate
      \expandafter\gdef\csname\@citelabel{#1}\endcsname{%
\strut 
\vadjust{\vskip-\dp\strutbox
\vbox to 0pt{\vss\parindent0cm \leftskip=\hsize 
\advance\leftskip3mm
\advance\hsize 4cm\strut\openup-4pt 
\rightskip 0cm plus 1cm minus 0.5cm ?  #1 ?\strut}}
         {\tt
            \escapechar = -1
            \nobreak\hskip0pt\pfeilsw
            \expandafter\string\csname#1\endcsname
             \pfeilso
            \nobreak\hskip0pt
         }%
      }%
     \else  
      \expandafter\gdef\csname\@citelabel{#1}\endcsname{%
            {\tt\expandafter\string\csname#1\endcsname}
      }%
     \fi  
   \fi
   \csname\@citelabel{#1}\endcsname
   \@notfirstcitationtrue
}%
%
%
\def\@citelabel#1{b@#1}%
%
%
\def\@citedef#1#2{\expandafter\gdef\csname\@citelabel{#1}\endcsname{#2}}%
%
%
%
\def\@readbblfile{%
   \ifx\@itemnum\@undefined
      \@innernewcount\@itemnum
   \fi
   \begingroup
      \def\begin##1##2{%
         \setbox0 = \hbox{\biblabelcontents{##2}}%
         \biblabelwidth = \wd0
      }%
      \def\end##1{}
      %
      %
      \@itemnum = 0
      \def\bibitem{\@getoptionalarg\@bibitem}%
      \def\@bibitem{%
         \ifx\@optionalarg\empty
            \expandafter\@numberedbibitem
         \else
            \expandafter\@alphabibitem
         \fi
      }%
      \def\@alphabibitem##1{%
         \expandafter \xdef\csname\@citelabel{##1}\endcsname {\@optionalarg}%
         \ifx\biblabelprecontents\@undefined
            \let\biblabelprecontents = \relax
         \fi
         \ifx\biblabelpostcontents\@undefined
            \let\biblabelpostcontents = \hss
         \fi
         \@finishbibitem{##1}%
      }%
      \def\@numberedbibitem##1{%
         \advance\@itemnum by 1
         \expandafter \xdef\csname\@citelabel{##1}\endcsname{\number\@itemnum}%
         \ifx\biblabelprecontents\@undefined
            \let\biblabelprecontents = \hss
         \fi
         \ifx\biblabelpostcontents\@undefined
            \let\biblabelpostcontents = \relax
         \fi
         \@finishbibitem{##1}%
      }%
      \def\@finishbibitem##1{%
         \biblabelprint{\csname\@citelabel{##1}\endcsname}%
         \@writeaux{\string\@citedef{##1}{\csname\@citelabel{##1}\endcsname}}%
         \ignorespaces
      }%
      %
      %
      \let\em = \bblem
      \let\newblock = \bblnewblock
      \let\sc = \bblsc
      \frenchspacing
      \clubpenalty = 4000 \widowpenalty = 4000
      \tolerance = 10000 \hfuzz = .5pt
      \everypar = {\hangindent = \biblabelwidth
                      \advance\hangindent by \biblabelextraspace}%
      \bblrm
      \parskip = 1.5ex plus .5ex minus .5ex
      \biblabelextraspace = .5em
      \bblhook
      \input \bblfilebasename.bbl
   \endgroup
}%
%
%
\@innernewdimen\biblabelwidth
\@innernewdimen\biblabelextraspace
%
%
%
\def\biblabelprint#1{%
   \noindent
   \hbox to \biblabelwidth{%
      \biblabelprecontents
      \biblabelcontents{#1}%
      \biblabelpostcontents
   }%
   \kern\biblabelextraspace
}%
%
%
%
\def\biblabelcontents#1{{\bblrm [#1]}}%
%
%
\def\bblrm{\rm}%
%
%
\def\bblem{\it}%
%
%
\def\bblsc{\ifx\@scfont\@undefined
              \font\@scfont = cmcsc10
           \fi
           \@scfont
}%
%
%
\def\bblnewblock{\hskip .11em plus .33em minus .07em }%
%
%
\let\bblhook = \empty
%
%
%
\def\printcitestart{[}
\def\printcitefinish{]}
\def\printbetweencitations{, }
\def\printcitenote#1{, #1}
%
%
%
\let\citation = \@gobble
%
%
%
\@innernewcount\@numparams
%
%
\def\newcommand#1{%
   \def\@commandname{#1}%
   \@getoptionalarg\@continuenewcommand
}%
%
%
\def\@continuenewcommand{%
   \@numparams = \ifx\@optionalarg\empty 0\else\@optionalarg \fi \relax
   \@newcommand
}%
%
%
\def\@newcommand#1{%
   \def\@startdef{\expandafter\edef\@commandname}%
   \ifnum\@numparams=0
      \let\@paramdef = \empty
   \else
      \ifnum\@numparams>9
         \errmessage{\the\@numparams\space is too many parameters}%
      \else
         \ifnum\@numparams<0
            \errmessage{\the\@numparams\space is too few parameters}%
         \else
            \edef\@paramdef{%
               \ifcase\@numparams
                  \empty  No arguments.
               \or ####1%
               \or ####1####2%
               \or ####1####2####3%
               \or ####1####2####3####4%
               \or ####1####2####3####4####5%
               \or ####1####2####3####4####5####6%
               \or ####1####2####3####4####5####6####7%
               \or ####1####2####3####4####5####6####7####8%
               \or ####1####2####3####4####5####6####7####8####9%
               \fi
            }%
         \fi
      \fi
   \fi
   \expandafter\@startdef\@paramdef{#1}%
}%
%
%
%
%
\def\@readauxfile{%
   \if@auxfiledone \else 
      \global\@auxfiledonetrue
      \@testfileexistence{aux}%
      \if@fileexists
         \begingroup
            \endlinechar = -1
            \catcode`@ = 11
            \input \jobname.aux
         \endgroup
      \else
         \message{\@undefinedmessage}%
         \global\@citewarningfalse
      \fi
      \immediate\openout\@auxfile = \jobname.aux
   \fi
}%
%
%
\newif\if@auxfiledone
\ifx\noauxfile\@undefined \else \@auxfiledonetrue\fi
%
%
%
%
\@innernewwrite\@auxfile
\def\@writeaux#1{\ifx\noauxfile\@undefined \write\@auxfile{#1}\fi}%
%
%
%
\ifx\@undefinedmessage\@undefined
   \def\@undefinedmessage{No .aux file; I won't give you warnings about
                          undefined citations.}%
\fi
%
%
\@innernewif\if@citewarning
\ifx\noauxfile\@undefined \@citewarningtrue\fi
%
%
%
\catcode`@ = \@oldatcatcode

\def\pfeilso{\leavevmode
            \vrule width 1pt height9pt depth 0pt\relax
           \vrule width 1pt height8.7pt depth 0pt\relax
           \vrule width 1pt height8.3pt depth 0pt\relax
           \vrule width 1pt height8.0pt depth 0pt\relax
           \vrule width 1pt height7.7pt depth 0pt\relax
            \vrule width 1pt height7.3pt depth 0pt\relax
            \vrule width 1pt height7.0pt depth 0pt\relax
            \vrule width 1pt height6.7pt depth 0pt\relax
            \vrule width 1pt height6.3pt depth 0pt\relax
            \vrule width 1pt height6.0pt depth 0pt\relax
            \vrule width 1pt height5.7pt depth 0pt\relax
            \vrule width 1pt height5.3pt depth 0pt\relax
            \vrule width 1pt height5.0pt depth 0pt\relax
            \vrule width 1pt height4.7pt depth 0pt\relax
            \vrule width 1pt height4.3pt depth 0pt\relax
            \vrule width 1pt height4.0pt depth 0pt\relax
            \vrule width 1pt height3.7pt depth 0pt\relax
            \vrule width 1pt height3.3pt depth 0pt\relax
            \vrule width 1pt height3.0pt depth 0pt\relax
            \vrule width 1pt height2.7pt depth 0pt\relax
            \vrule width 1pt height2.3pt depth 0pt\relax
            \vrule width 1pt height2.0pt depth 0pt\relax
            \vrule width 1pt height1.7pt depth 0pt\relax
            \vrule width 1pt height1.3pt depth 0pt\relax
            \vrule width 1pt height1.0pt depth 0pt\relax
            \vrule width 1pt height0.7pt depth 0pt\relax
            \vrule width 1pt height0.3pt depth 0pt\relax}

\def\pfeilsw{ \leavevmode 
            \vrule width 1pt height0.3pt depth 0pt\relax
            \vrule width 1pt height0.7pt depth 0pt\relax
            \vrule width 1pt height1.0pt depth 0pt\relax
            \vrule width 1pt height1.3pt depth 0pt\relax
            \vrule width 1pt height1.7pt depth 0pt\relax
            \vrule width 1pt height2.0pt depth 0pt\relax
            \vrule width 1pt height2.3pt depth 0pt\relax
            \vrule width 1pt height2.7pt depth 0pt\relax
            \vrule width 1pt height3.0pt depth 0pt\relax
            \vrule width 1pt height3.3pt depth 0pt\relax
            \vrule width 1pt height3.7pt depth 0pt\relax
            \vrule width 1pt height4.0pt depth 0pt\relax
            \vrule width 1pt height4.3pt depth 0pt\relax
            \vrule width 1pt height4.7pt depth 0pt\relax
            \vrule width 1pt height5.0pt depth 0pt\relax
            \vrule width 1pt height5.3pt depth 0pt\relax
            \vrule width 1pt height5.7pt depth 0pt\relax
            \vrule width 1pt height6.0pt depth 0pt\relax
            \vrule width 1pt height6.3pt depth 0pt\relax
            \vrule width 1pt height6.7pt depth 0pt\relax
            \vrule width 1pt height7.0pt depth 0pt\relax
            \vrule width 1pt height7.3pt depth 0pt\relax
            \vrule width 1pt height7.7pt depth 0pt\relax
            \vrule width 1pt height8.0pt depth 0pt\relax
            \vrule width 1pt height8.3pt depth 0pt\relax
            \vrule width 1pt height8.7pt depth 0pt\relax
            \vrule width 1pt height9pt depth 0pt\relax
      }


\def\widestnumber#1#2{}

\def\citewarning#1{\ifx\shlhetal\relax 
    \else
    \par{#1}\par
    \fi
}

\def\rm{\fam0 \tenrm}

\def\fakesubhead#1\endsubhead{\bigskip\noindent{\bf#1}\par}



%
%
%

%

\font\textrsfs=rsfs10
\font\scriptrsfs=rsfs7
\font\scriptscriptrsfs=rsfs5

\newfam\rsfsfam
\textfont\rsfsfam=\textrsfs
\scriptfont\rsfsfam=\scriptrsfs
\scriptscriptfont\rsfsfam=\scriptscriptrsfs

\edef\oldcatcodeofat{\the\catcode`\@}
\catcode`\@11

\def\Cal@@#1{\noaccents@ \fam \rsfsfam #1}

\catcode`\@\oldcatcodeofat


\expandafter\ifx \csname margininit\endcsname \relax\else\margininit\fi

\long\def\red#1\endred{}
\long\def\green#1\endgreen{}
\long\def\blue#1\endblue{}
\long\def\private#1\endprivate{}

\def\endred{ \unmatched endred! }
\def\endgreen{ \unmatched endgreen! }
\def\endblue{ \unmatched endblue! }
\def\endprivate{ \unmatched endprivate! }

\ifx\latexcolors\undefinedcs\def\latexcolors{}\fi

\def\emptycs{}
\def\evaluatelatexcolors{%
        \ifx\latexcolors\emptycs\else
        \expandafter\xxevaluate\latexcolors\xxfertig\evaluatelatexcolors\fi}
\def\xxevaluate#1,#2\xxfertig{\setupthiscolor{#1}%
        \def\latexcolors{#2}}


\font\smallfont=cmsl7
\def\rutgerscolor{\ifmmode\else\endgraf\fi\smallfont
\advance\leftskip0.5cm\relax}
\def\setupthiscolor#1{\edef\tmptmpcs{\noexpand\bgroup\noexpand\rutgerscolor
\noexpand\def\noexpand\currentcolor{#1}%
\noexpand}%
\expandafter\let\csname#1\endcsname\tmptmpcs
\def\tmptmpcs{\checkColorUnmatched{#1}\popthecolor}
\expandafter\let\csname end#1\endcsname\tmptmpcs}

\def\checkColorUnmatched#1{\def\expectcolor{#1}%
    \ifx\expectcolor\currentcolor   
    \else \edef\failhere{\noexpand\tryingToClose '\currentcolor' with end\expectcolor}\failhere\fi}

\def\currentcolor{???}

\def\popthecolor{\ifmmode\else\endgraf\fi\egroup}

\expandafter\def\csname#1\endcsname{}

\evaluatelatexcolors

 \let\outerhead\head
 \def\head{\innerhead}
 \let\innerhead\outerhead

 \let\outersubhead\subhead
 \def\subhead{\innersubhead}
 \let\innersubhead\outersubhead

 \let\outersubsubhead\subsubhead
 \def\subsubhead{\innersubsubhead}
 \let\innersubsubhead\outersubsubhead

 \let\outerproclaim\proclaim
 \def\proclaim{\innerproclaim}
 \let\innerproclaim\outerproclaim

 %
 %
 %
 %

\def\demo#1{\medskip\noindent{\it #1.\/}}
\def\enddemo{\smallskip}

\def\remark#1{\medskip\noindent{\it #1.\/}}
\def\endremark{\smallskip}

\def\beginaside{\endgraf\leftskip2cm \vrule width 0pt\relax}
\def\endaside{\endgraf\leftskip0cm \vrule width 0pt\relax}

\pageheight{8.5truein}
\topmatter
\title {Categoricity in abstract elementary classes: going up
inductive step} \endtitle
\rightheadtext{Categoricity in abstract elementary classes}
\author {Saharon Shelah \thanks{\null\newline
I thank Alice Leonhardt for the beautiful typing \null\newline
This research was supported by the United States-Israel Binational
Science Foundation and in its final stages also by the Israel Science
Foundation (Grant no. 242/03). Publication 600. } \endthanks} \endauthor

\affil{The Hebrew University of Jerusalem \\
Einstein Institute of Mathematics \\
Edmond J. Safra Campus, Givat Ram
Jerusalem 91904, Israel
\medskip
Department of Mathematics \\
Hill Center-Busch Campus \\
Rutgers, The State University of New Jersey\\
110 Frelinghuysen Road
Piscataway, NJ 08854-8019 USA} \endaffil

\subjclass  03C45, 03C75, 03C95, 03C50  \endsubjclass 
\mn

\abstract We deal with beginning stability theory for ``reasonable" 
non-elementary classes without
any remnants of compactness like dealing with models 
above Hanf number or by the class being definable by
$\Bbb L_{\omega_1,\omega}$. \endabstract
\endtopmatter
\document  
 
\pretags{300a,300b,300c,300d,300e,300f,300g,300x,300y,300z,88r,705,E53,734}
\newpage

\head {\S0 Introduction} \endhead  \resetall \sectno=0
 \spuriousreset
\bn
The paper's main explicit result is proving Theorem \scite{600-0.A} below.  
It is done axiomatically, in a ``superstable" abstract framework with the 
set of ``axioms" of the frame, verified by applying earlier works, so
it suggests this frame as the, or at least a major,
non-elementary parallel of superstable. \nl
A major case to which this is applied, is the one from 
\cite{Sh:576}; we continue this work in several ways but 
the use of \cite{Sh:576}
is only in verifying the basic framework;  we refer the reader
to the book's introduction or 
\cite[\S0]{Sh:576} for background and some further claims but all 
the definitions and basic properties appear here.
Otherwise, the heavy use of earlier works is in 
proving that our abstract framework applies in those contexts.  
If $\lambda = \aleph_0$ is O.K. for you, you may use  \chaptercite{88r} or 
\cite{Sh:48} instead of \cite{Sh:576}.

Naturally, our deeper aim is to develop stability theory (actually a parallel
of the theory of superstable elementary classes) for non-elementary
classes.  We use the number of non-isomorphic models as test problem.
Our main conclusion is \scite{600-0.A} below.
As a concession to supposedly general opinion, we restrict ourselves here
to the $\lambda$-good framework and delay dealing with weak relatives
(see \cite{Sh:849}, \cite{Sh:838}, Jarden-Shelah \cite{JrSh:875},
\cite{Sh:E46}; in
the first we start from \cite[\S1-\S5]{Sh:576}, i.e., before proving that
$K^{3,\text{uq}}_\lambda$ is dense).
Also, we assume the (normal) weak-diamond ideal on the $\lambda^{+
\ell}$ is not saturated (for $\ell=1,\dotsc,n-1$) and delay dealing with the 
elimination of this extra assumption (to \cite{Sh:838} there we waive
the ``not $\lambda^{+\ell+1}$-saturated$^1$ (ideal)", the price is
that we replace most $< 2^{\lambda^{+\ell+1}}$ by $<
\mu_{\text{unif}}(\lambda^{+\ell+1},2^{\lambda^{+ \ell}})$, see on it
\marginbf{!!}{\cprefix{88r}.\scite{88r-0.wD}}(3)). 
\bigskip

\proclaim{\stag{600-0.A} Theorem}  Assume $2^\lambda < 2^{\lambda^{+1}} < \dots
< 2^{\lambda^{+n+1}}$ and the (so called weak diamond) 
normal \footnote{recall that as $2^{\lambda_{\ell -1}} <
2^{\lambda_\ell}$ this ideal is not trivial, i.e., $\lambda^{+ \ell}$
is not in the ideal}
ideal ${\text{\rm WDmId\/}}(\lambda^{+ \ell})$ is
not $\lambda^{+ \ell +1}$-saturated
\footnote{actually the statement ``some normal ideal on $\mu^+$ is
$\mu^{++}$-saturated" is ``expensive", i.e., of large consistency
 strength, etc., so it is ``hard" for this assumption to fail}
 for $\ell = 1,\dotsc,n$. \nl
1) Let ${\frak K}$ be an abstract elementary class (see \S1 below) 
categorical in $\lambda$ and $\lambda^+$ with 
${\text{\rm LS\/}}({\frak K}) \le \lambda$ (e.g. the class of models of
$\psi \in 
\Bbb L_{\lambda^+,\omega}$ with $\le_{\frak K}$ defined naturally).  
If $1 \le \dot I(\lambda^{+2},{\frak K})$ 
and $2 \le \ell \le n \Rightarrow \dot I(\lambda^{+ \ell},
{\frak K}) < 2^{\lambda^{+ \ell}}$, \ub{then}
${\frak K}$ has a model of cardinality $\lambda^{+n+1}$. \nl
2)  Assume $\lambda = \aleph_0$, and $\psi \in
\Bbb L_{\omega_1,\omega}(\bold Q)$.  
If $1 \le \dot I(\lambda^{+ \ell},\psi) < 2^{\lambda^{+ \ell}}$
for $\ell = 1, \dotsc,n-1$ \ub{then} $\psi$ has a model in 
$\lambda^{+ n}$  (see \cite{Sh:48}).
\endproclaim
\bn
Note that if $n = 3$, then \scite{600-0.A}(1) is already 
proved in \cite{Sh:576}.  If ${\frak K}$ is the class
of models of some $\psi \in \Bbb L_{\omega_1,\omega}$ this is proved in
\cite{Sh:87a}, \cite{Sh:87b}, but the proof here \ub{does not generalize} the
proofs there.  It is a different one (of course, they are related).  
There, for proving the theorem for $n$, we have to consider
a few statements on $(\aleph_m,{\Cal P}^-(n-m))$-systems for all $m \le n$,
(going up and down).
A major point (there) is 
that for $n=0$, as $\lambda = \aleph_0$ we have the omitting
type theorem and the types are ``classical", that is, are 
sets of formulas.  This helps in
proving strong dichotomies; so the analysis of what occurs in $\lambda^{+n}
= \aleph_n$ is helped by those dichotomies.  Whereas here we deal with 
$\lambda,\lambda^+,\lambda^{+2},
\lambda^{+3}$ and then ``forget" $\lambda$ and deal with $\lambda^+,
\lambda^{+2},\lambda^{+3},\lambda^{+4}$, etc.  However, there are
some further theorems proved in \cite{Sh:87a}, \cite{Sh:87b}, 
whose parallels are not proved 
here, mainly that if for every $n$, in $\lambda^{+n}$ we get the
``structure" side, then the class has models in every $\mu \ge
\lambda$, and theorems about categoricity.  
We shall deal with them in subsequent works, mainly \chaptercite{705}.  
Also in \cite{Sh:48}, \cite{Sh:88} = \chaptercite{88r} 
we started to deal with $\psi \in 
\Bbb L_{\omega_1,\omega}(\bold Q)$ 
dealing with $\aleph_1,\aleph_2$; we put it in our present 
framework.  Of course, also the framework
of \chaptercite{88r} is integrated into our present context.
In the axiomatic framework (introduced in \S2) we are able to 
present a lemma, speaking
on only 4 cardinals, and which implies the theorem \scite{600-0.A}.  (Why?
Because in \S3 by \cite{Sh:576} we can get a so-called good
$\lambda^+$-frame ${\frak s}$ with $K^{\frak s} \subseteq {\frak K}$,
and then we prove a similar theorem on good frames by induction on
$n$, with the induction step done by the lemma mentioned above).
For this, parts of the proof are a generalization of the proof 
of \cite[\S8,\S9,\S10]{Sh:576}. 
\bn
A major theme here (and even more so in \chaptercite{705}) is \nl
\margintag{600-0.B}\ub{\stag{600-0.B} Thesis}:  It is worthwhile to develop model theory (and
superstability in particular) in the context of ${\frak K}_\lambda$ or
$K_{\lambda^{+ \ell}},\ell \in \{0,\dotsc,n\}$, i.e., restrict
ourselves to one, few, an interval of cardinals.  We may have
good understanding of the class in this context, while in general
cardinals we are lost.
\bn
As in \cite{Sh:c} for first order classes \nl 
\margintag{600-0.C}\ub{\stag{600-0.C} Thesis}:  It is reasonable first to develop the theory for the
class of (quite) saturated enough models as it is smoother and even if you
prefer to investigate the non-restricted case, the saturated case will
clarify it.  In our case this will mean investigating ${\frak s}^{+n}$
for each $n$ and then $\cap\{{\frak K}^{{\frak s}^{+n}}:n < \omega\}$.
\bn
\margintag{600-0.D}\ub{\stag{600-0.D} The Better to be poor Thesis}:  Better to know what 
is essential. e.g., you may have better closure properties (here a
major point of poverty is having no formulas, this is even more
noticeable in \chaptercite{705}).
\sn
I thank John Baldwin, Alex Usvyatsov, Andres Villaveces and Adi Yarden 
for many complaints and corrections.

\S1 gives a self-contained introduction to a.e.c. (abstract elementary
classes), including definitions of types, $M_2$ is
$(\lambda,\kappa)$-brimmed over $M_1$ and saturativity = universality
+ model homogeneity.  An interesting point is observing that any
$\lambda$-a.e.c. ${\frak K}_\lambda$ can be lifted to ${\frak K}_{\ge
\lambda}$, uniquely; so it does not matter if we deal with ${\frak
K}_\lambda$ or ${\frak K}_{\ge \lambda}$ (unlike the situation for
good $\lambda$-frames, which if we lift, we in general, lose some
essential properties).

The good $\lambda$-frames introduced in \S2 are a very central notion
here.  It concentrates on one cardinal $\lambda$, in ${\frak K}_\lambda$
we have amalgamation and more, hence types, in the orbital, 
not in the classical sense of set of formulas, for models of cardinality
$\lambda$ can be reasonably defined (we concentrate on so-called basic
types) and we axiomatically have a non-forking relation for them.

In \S3 we show that starting with classes belonging to reasonably large
families, from assumptions on categoricity (or few models), good
$\lambda$-frames arise.  In \S4 we deduce some things on good
$\lambda$-frames; mainly: stability in $\lambda$, existence and (full)
uniqueness of $(\lambda,*)$-brimmed extensions of $M \in K_\lambda$.

Concerning \S5 we know that if $M \in K_\lambda$ and $p \in {\Cal
S}^{\text{bs}}(M)$ then there is $(M,N,a) \in K^{3,\text{bs}}_\lambda$
such that \ortp$(a,M,N)=p$.  But can we find a special (``minimal" or
``prime") triple in some sense?  Note that if $(M_1,N_1,a)
\le_{\text{bs}} (M_2,N_2,a)$ then $N_2$ is an amalgamation of
$N_1,M_2$ over $M_1$ (restricting ourselves to the case
``\ortp$(a,M_2,N_2)$ does not fork over $M_1$") and we may wonder is this
amalgamation unique (i.e., allowing to increase or decrease $N_2$).
If this holds for any such $(M_2,N_2,a)$ we say $(M_1,N_1,a)$ has
uniqueness (= belongs to $K^{3,\text{uq}}_\lambda =
K^{3,\text{uq}}_{\frak s}$).  Specifically we ask: is
$K^{3,\text{uq}}_\lambda$ dense in
$(K^{3,\text{bs}}_\lambda,\le_{\text{bs}})$?  If no, we get a
non-structure result; if yes, we shall (assuming categoricity) deduce
the ``existence for $K^{3,\text{uq}}_{\frak s}$" and this is used
later as a building block for non-forking amalgamation of models.

So our next aim is to find ``non-forking" amalgamation of models (in
\S6).  We first note that there is at most one such notion which
fulfills our expectations (and ``respect" ${\frak s}$).  Now if
$\nonfork{}{}_{}(M_0,M_1,a,M_3),M_0 \le_{\frak K} M_2 \le_{\frak K}
M_3$ and $(M_0,M_2,a) \in K^{3,\text{uq}}_\lambda$ by our demands we
have to say that $M_1,M_2$ are in non-forking amalgamation over $M_0$
inside $M_3$.   Closing this family under the closure demands we expect
to arrive to a notion NF$_\lambda = \text{\rm NF}_{\frak s}$ 
which should be the
right one (if a solution exists at all).  But then we have to work on
proving that it has all the properties it hopefully has.

A major aim in advancing to $\lambda^+$ is having a superlimit model.
So in \S7 we find out who he should be: the saturated model of ${\frak
K}_{\lambda^+}$, but is it superlimit?  We use our NF$_\lambda$ to
define a ``nice" order $\le^*_{\lambda^+}$ on ${\frak K}_{\lambda^+}$,
investigate it and prove the existence of a superlimit model under
this partial order.  To
complete the move to $\lambda^+$ we would like to have that the class
of $\lambda^+$-saturated model with the partial order $\le^*_{\lambda^+}$ is a
$\lambda^+$-a.e.c.  Well, we do not prove it but rather use it as a
dividing line: if it fails we eventually get many models in ${\frak
K}_{\lambda^{++}}$ (coding a stationary subset of $\lambda^{++}$
(really any $S \subseteq \{\delta < \lambda^{++}:\text{cf}(\delta) =
\lambda^+\}$)), see \S8.  

Lastly, we pay our debts: prove the theorems which were the motivation
of this work, in \S9.
\bn
\centerline {$* \qquad * \qquad *$}
\bn
\centerline {Reading Instructions} 
\bn
As usual these are instructions on what you can avoid reading.

Note that \S3 contains the examples, i.e., it shows how ``good
$\lambda$-frame", our main object of study here, arise in previous
works.  This, on the one hand, may help the reader to understand what
is a good frame and, on the other hand, helps us in the end to draw 
conclusions continuing those works.  However, it is \ub{not} necessary
here otherwise, so you may ignore it.

Note that we treat the subject axiomatically, in a 
general enough way to treat the
cases which exist without trying too much to eliminate axioms as long as the
cases are covered (and probably most potential readers will feel they are
more than general enough). \nl
We shall assume
\mr
\item "{$(*)_0$}"  $2^\lambda < 2^{\lambda^+} < 2^{\lambda^{+2}} < \ldots
< 2^{\lambda^{+n}}$ and $n \ge 2$.
\ermn
In the end of \S1 there are some basic definitions.
\bn
\ub{Reading Plan 0}:  We accept the good frames as interesting per se
so ignore \S3 (which gives ``examples") and: \S1 tells you all you
need to know on abstract elementary classes; \S2 presents frames, etc.
\bn
\ub{Reading Plan 1}:  The reader decides to understand why we reprove
the main theorem of \cite{Sh:87a}, \cite{Sh:87b} so
\mr
\item "{$(*)_1$}"  $K$ is the class of models of some $\psi \in
\Bbb L_{\lambda^+,\omega}$ (with a natural notion of elementary embedding
$\prec_{\Cal L}$ for ${\Cal L}$ a fragment of $\Bbb L_{\lambda^+,\omega}$ 
of cardinality $\le \lambda$ to which $\psi$ belongs). 
\ermn
So in fact (as we can replace, for this result, $K$ by any class with
fewer models still satisfying the assumptions)  \wilog \,
\mr
\item "{$(*)'_1$}"  if $\lambda = \aleph_0$ then $K$ 
is the class of atomic models of some complete first order theory,
$\le_{\frak K}$ is being elementary submodel.
\ermn
The theorems we are seeking are of the form
\mr
\item "{$(*)_2$}"  if $K$ has few models in $\lambda + \aleph_1,\lambda^+,
\dotsc,\lambda^{+n}$ then it has a model in $\lambda^{+n+1}$. \nl
[Why  ``$\lambda + \aleph_1$"?  If $\lambda > \aleph_0$ this means
$\lambda$ whereas if $\lambda = \aleph_0$ this means that we do not
require ``few model in $\lambda = \aleph_0$".  The reason is that for
the class or models of $\psi \in \Bbb L_{\omega_1,\omega}$ (or $\in
\Bbb L_{\omega_1,\omega}(\bold Q)$ or an a.e.c. which is
PC$_{\aleph_0}$, see Definition \scite{600-0.4}) we 
have considerable knowledge of general methods of
building models of cardinality $\aleph_1$, for general $\lambda$ we
are very poor in such knowedge (probably as there is much less).]
\ermn
But, of course, what we would really like to have are rudiments of stability
theory (non-forking amalgamation, superlimit models, etc.).  Now
reading plan 1 is to follow reading plan 2 below \ub{but}
replacing the use of Claim \scite{600-Ex.4} and \cite{Sh:576} by the 
use of a simplified version of \scite{600-Ex.1} and \cite{Sh:87a}.
\bn
\ub{Reading Plan 2}:  The reader wants to understand the proof of
$(*)_2$ for arbitrary ${\frak K}$ and $\lambda$.
The reader 
\mr
\item "{$(a)$}"   knows at least the main definitions and results 
of \cite{Sh:576},
\nl
or just
\sn
\item "{$(b)$}"   reads the main definitions of \S1 here (in
\scite{600-0.1} - \scite{600-0.6}) and is 
willing to believe some quotations of results of \cite{Sh:576}.
\ermn
We start assuming 
${\frak K}$ is an abstract elementary class, LS$({\frak K}) \le
\lambda$ (or read \S1 here until \scite{600-0.22}) and
${\frak K}$ is categorical in $\lambda,
\lambda^+$ and $1 \le \dot I(\lambda^{++},K) 
< 2^{\lambda^{++}}$ and moreover, $1 \le \dot I(\lambda^{++},K) < 
2^{\lambda^{++}}$.  As an appetizer and to understand types and the
definition of types and saturated (in the present context) and
brimmed, read from \S1 until \scite{600-4a.10}.
\nl
He should read in \S2 Definition \scite{600-1.1} of $\lambda$-good frame, an 
axiomatic framework and then read the following two Definitions
\scite{600-1.6}, \scite{600-1.7} and Claim \scite{600-1.8}.  In \S3, \scite{600-Ex.4}
show how by \cite{Sh:576} the context there gives a $\lambda^+$-good frame;
of course the reader may just believe instead of reading proofs, and
he may remember that our basic types are minimal in this case. \nl
In \S4 he should read some consequences of the axioms proved with weaker
axioms, understanding here and later $<^*$ as $\le_{\frak K} \restriction
K_\lambda$. \nl
Then in \S5 we show some
amount of unique amalgamation.  Then \S6,\S7,\S8 do a parallel to 
\cite[\S8,\S9,\S10]{Sh:576} in our context; still there are
differences, in particular our context is not necessarily
uni-dimensional which complicates matters.  But if we restrict ourselves to
continuing \cite{Sh:576}, our frame is ``uni-dimensional", we could have
simplified the proofs 
by using ${\Cal S}^{\text{bs}}(M)$ as the set of minimal types.
\bn
\ub{Reading Plan 3}:  $\psi \in \Bbb L_{\omega_1,\omega}(\bold Q)$ so 
$\lambda= \aleph_0,1 \le \dot I(\aleph_1,\psi) < 2^{\aleph_1}$.  

For this,
\cite{Sh:576} is irrelevant (except quoting the ``\ub{black box}" use of
the combinatorial section \S3
of \cite{Sh:576} when using the weak diamond to get many non-isomorphic
models in \S5).

Now reading plan 3 is to follow reading plan 2 but 
\scite{600-Ex.4} is replaced by \scite{600-Ex.1A} which relies on
\cite{Sh:48}, i.e., it proves that we get an $\aleph_1$-good frame
investigating $\psi \in \Bbb L_{\omega_1,\omega}(\bold Q)$.

Note that our class may well be such that ${\frak K}$ is the parallel
of ``superstable non-multidimension complete first order theory";
e.g., $\psi_1 = (\bold Q x)(P(x)) \wedge(\bold Q x)(\neg P(x)),
\tau_\psi = \{P\},P$ a unary predicate; this is
categorical in $\aleph_1$ and has no models in $\aleph_0$ and $\psi_1$
has 3 models in $\aleph_2$.  But if we
use $\psi_0 = (\forall x)(P(x) \equiv P(x))$ we have
$\dot I(\aleph_1,\psi_0) = \aleph_0$; however, even starting with $\psi_1$, the
derived a.e.c. ${\frak K}$ has exactly three non-isomorphic models in
$\aleph_1$.  In general we derived an a.e.c. ${\frak K}$ from $\psi$
such that: ${\frak K}$ is an a.e.c. with LS number $\aleph_0$,
categorical in $\aleph_0$, and the number of somewhat ``saturated"
models of ${\frak K}$ in $\lambda$ is $\le \dot I(\lambda,\psi)$ for
$\lambda \ge \aleph_1$.  The relationship of $\psi$ and ${\frak K}$ is
not comfortable; as it means that, for  general results to be applied,
they have to be somewhat stronger, e.g. ``${\frak K}$ has
$2^{\lambda^{++}}$ non-isomorphic \ub{$\lambda^+$-saturated} models of
cardinality $\lambda^{++}$".  The reason is
LS$({\frak K}) = \lambda = \aleph_0$; we have to find many 
somewhat $\lambda^+$-saturated models as we have first in a sense
eliminate the quantifier $\bold Q = 
\exists^{\ge \aleph_1}$, (i.e., the choice of
the class of models and of the order guaranteed that what has to be
countable is countable, and $\lambda^+$-saturation guarantees that what
should be uncountable is uncountable).
\bn
\ub{Reading Plan 4}:  ${\frak K}$ an abstract elementary class which is 
PC$_\omega$ ($= \aleph_0$-presentable, see Definition \scite{600-0.4}); 
see \chaptercite{88r} or \cite{Mw85a} which includes a 
friendly presentation of \cite[\S1-\S3]{Sh:88} so of \sectioncite[\S1-\S3]{88r}).

Like plan 3 but we have to use \scite{600-Ex.1} instead of 
\scite{600-Ex.1A} and fortunately
the reader is encouraged to read \sectioncite[\S4,\S5]{88r} to understand why we
get a $\lambda$-good quadruple.
\newpage

\head {\S1 Abstract elementary classes} \endhead  \resetall \sectno=1
 \spuriousreset
\bn

First we present the basic material on a.e.c. ${\frak K}$, so types,
saturativity and brimmness (so most is repeating some things from
\sectioncite[\S1]{88r} and from \chaptercite{300b}).  
\nl
Second we show that the situation in $\lambda = \text{ LS}({\frak K})$
determine the situation above $\lambda$, moreover such lifting always
exists; so a $\lambda$-a.e.c. can be lifted to a $(\ge
\lambda)$-a.e.c. in one and only one way.
\bigskip

\demo{\stag{600-0.1} Conventions}  Here ${\frak K} = (K,\le_{\frak K})$, where
$K$ is a class of $\tau$-models for a fixed vocabulary $\tau = \tau_K =
\tau_{\frak K}$ and $\le_{\frak K}$ is a two-place relation on the models in 
$K$.  We do not always 
strictly distinguish between 
${\frak K},K$ and $(K,\le_{\frak K})$.  We shall 
assume 
that $K,\le_{\frak K}$ are fixed, and $M \le_{\frak K} N \Rightarrow 
M,N \in K$; and we assume that it is an abstract elementary class, see
Definition \scite{600-0.2} below.  When we use $<_{\frak K}$ in the $\prec$ 
sense (elementary submodel for first order logic), we write
$\prec_{\Bbb L}$.
\enddemo
\bigskip

\definition{\stag{600-0.1A} Definition}  For a class of $\tau_K$-models we
let $\dot I(\lambda,K) = |\{M/\cong:M \in K,
\|M\| = \lambda\}|$.
\enddefinition
\bigskip

\definition{\stag{600-0.1B} Definition}  1) We say $\bar M =  \langle M_i:i <
\mu \rangle$ is a representation or filtration of a model $M$ of
cardinality $\mu$ \ub{if} $\tau_{M_i} = \tau_M,M_i$
is $\subseteq$-increasing continuous, $\|M_i\| < \|M\|$ and $M =
\cup\{M_i:i < \mu\}$ and 
$\mu = \chi^+ \Rightarrow \|M_i\| = \chi$. \nl
2) We say $\bar M$ is a $\le_{\frak K}$-representation or 
$\le_{\frak K}$-filtration
\ub{if} in addition $M_i \le_{\frak K} M$ for $i < \|M\|$
(hence $M_i,M \in K$ and $\langle M_i:i < \mu \rangle$ is $\le_{\frak
K}$-increasing continuous, by Av V from Definition \scite{600-0.2}).
\enddefinition
\bigskip

\definition{\stag{600-0.2} Definition}  We say ${\frak K} = 
(K,\le_{\frak K})$ is an
abstract elementary class, a.e.c. in short, if 
($\tau$ is as in \scite{600-0.1}, $Ax 0$ holds and) AxI-VI hold where: \nl
$Ax 0$: The holding of 
$M \in K,N \le_{\frak K} M$ depends on $N,M$ only up to isomorphism, i.e., 
$[M \in K,M \cong N \Rightarrow N \in K]$, and [if $N \le_{\frak K} M$ and 
$f$ is an isomorphism from $M$ onto the $\tau$-model $M'$ mapping $N$ onto 
$N'$ \underbar{then} $N' \le_{\frak K} M'$].
\medskip

$Ax I$: If $M \le_{\frak K} N$ then $M \subseteq N$ (i.e. $M$ is a submodel
of $N$).
\medskip

$Ax II$: $M_0 \le_{\frak K} M_1 \le_{\frak K} M_2$ implies $M_0 \le_{\frak K}
M_2$ and $M \le_{\frak K} M$ for $M \in K$.
\medskip

$Ax III$: If $\lambda$ is a regular cardinal, $M_i$ (for $i < \lambda$) is
$\le_{\frak K}$-increasing (i.e. $i < j < \lambda$ implies $M_i \le_{\frak K}
M_j$) and continuous (i.e. for limit ordinal $\delta < \lambda$ we have
\newline
$M_\delta = \dsize \bigcup_{i < \delta} M_i$) \underbar{then}
$M_0 \le_{\frak K} \dsize \bigcup_{i < \lambda} M_i$.
\medskip

$Ax IV$: If $\lambda$ is a regular cardinal, $M_i$ (for $i < \lambda)$ is
$\le_{\frak K}$-increasing continuous and 
$M_i \le_{\frak K} N$ for $i < \lambda$ \underbar{then}
$\dsize \bigcup_{i < \lambda} M_i \le_{\frak K} N$.
\medskip

$Ax V$:  If $M_0 \subseteq M_1$ and $M_\ell \le_{\frak K} N$ for $\ell =
0,1$, \ub{then} $M_0 \le_{\frak K} M_1$.
\medskip

$Ax VI$:  LS$({\frak K})$ exists 
\footnote{We normally assume $M \in {\frak K}
\Rightarrow \|M\| \ge \text{ LS}({\frak K})$ so may forget to write
$\|M\| ``+\text{ LS}({\frak K})"$ instead $\|M\|$, here 
there is no loss in it.  It is also
natural to assume $|\tau({\frak K})| \le \text{ LS}
({\frak K})$ which means just
increasing LS$({\frak K})$, but no real need here; dealing with Hanf
numbers it is natural.}, where LS$({\frak K})$ is 
the minimal cardinal $\lambda$ such that: if \newline
$A \subseteq N$ and $|A| \le \lambda$ \underbar{then} for some 
$M \le_{\frak K} N$ we have $A \subseteq |M|$ and $\|M\| \le \lambda$.
\enddefinition
\bn
\margintag{600-0.2B}\underbar{\stag{600-0.2B} Notation}:  1) $K_\lambda = 
\{ M \in K:\|M\| = \lambda\}$ and
$K_{< \lambda} = \dsize \bigcup_{\mu < \lambda} K_\mu$. 
\bigskip

\definition{\stag{600-0.3} Definition}  1) The embedding $f:N \rightarrow M$ 
is $\le_{\frak K}$-embedding when its range is 
the universe of a model $N' \le_{\frak K} M$, 
(so $f:N \rightarrow N'$ is an isomorphism onto).
\nl
2) We say $f$ is a $\le_{\frak K}$-embedding of $M_1$ into $M_2$ over
$M_0$ when for some $M'_1$ we have: $M_0 \le_{\frak K} M_1,M_0 
\le_{\frak K} M'_1 \le_{\frak K} M_2$ and $f$ is an isomorphism from
$M_1$ onto $M'_1$ extending the mapping id$_{M_0}$.
\enddefinition
\bn
Recall
\demo{\stag{600-0.6} Observation}  Let $I$ be a directed set (i.e., $I$ is 
partially ordered by $\le = \le^I$, such that 
any two elements have a common upper bound). \newline
1) If $M_t$ is defined for $t \in I$, and $t \le s \in I$ implies $M_t
\le_{\frak K} M_s$ \underbar{then} for every $t \in I$ we have
$M_t \le_{\frak K} \dsize \bigcup_{s \in I} M_s$. \newline
2) If in addition $t \in I$ implies $M_t \le_{\frak K} N$ \underbar{then}
$\dsize \bigcup_{s \in I} M_s \le_{\frak K} N$.
\enddemo
\bigskip

\demo{Proof}  Easy or see \marginbf{!!}{\cprefix{88r}.\scite{88r-1.6}} which does not rely on
anything else.  \hfill$\square_{\scite{600-0.6}}$
\enddemo
\bigskip

\proclaim{\stag{600-0.7} Claim}  1) For every $N \in K$ there is a directed 
partial order $I$ of cardinality $\le \|N\|$ and sequence
$\bar M = \langle M_t:t \in
I \rangle$ such that $t \in I \Rightarrow M_t \le_{\frak K} N,\|M_t\| \le
{ \text{\rm LS\/}}({\frak K}),
I \models s < t \Rightarrow M_s \le_{\frak K} M_t$ and $N =
\dbcu_{t \in I} M_t$.  If $\|N\| \ge \text{\rm LS}({\frak K})$ we can add
$\|M_t\| = \text{\rm LS}({\frak K})$ for $t \in I$.
\nl
2) For every $N_1 \le_{\frak K} N_2$ we can find $\langle M^\ell_t:t \in
I_\ell \rangle$ as in part (1) for $\ell=1,2$ 
such that $I_1 \subseteq I_2$ and
$t \in I_1 \Rightarrow M^2_t = M^1_t$. \nl
3) Any 
$\lambda \ge {\text{\rm LS\/}}({\frak K})$ 
satisfies the requirement in the definition of ${\text{\rm LS\/}}
({\frak K})$.
\endproclaim
\bigskip

\demo{Proof}  Easy or see \marginbf{!!}{\cprefix{88r}.\scite{88r-1.7}} which does not require
anything else.  \hfill$\square_{\scite{600-0.7}}$
\enddemo
\bn
We now (in \scite{600-0.12}) recall the (non-classical) definition  of type
(note that it is natural to look at types only over models which are
amalgamation bases, see part (4) of \scite{600-0.12} below and 
consider only extensions of
the models of the same cardinality). 
Note that though the choice of the name indicates that they are
supposed to behave like complete types over models as in classical
model theory (on which we are not relying), this does not
guarantee most of the basic properties.  E.g., when cf$(\delta) =
\aleph_0$, uniqueness of $p_\delta \in {\Cal S}(M_\delta)$ such that
$i < \delta \Rightarrow p_\delta \restriction M_i = p_i$ is not 
guaranteed even if $p_i \in {\Cal S}(M_i),M_i$ is
$\le_{\frak K}$-increasing continuous for $i \le \delta$ and $i < j
< \delta \Rightarrow p_i = p_j \restriction M_i$.  Still 
we have existence: if for 
$i < \delta,p_i \in {\Cal S}(M_i)$
increasing with $i$, then there is $p_\delta \in {\Cal S}(\cup \{M_i:i
< \delta\})$ such that $i < \delta \Rightarrow p_i = p_\delta
\restriction M_i$.  But when cf$(\delta) > \aleph_0$ even existence
is not guaranteed.
\definition{\stag{600-0.12} Definition}   1) For $M \in K_\mu$ we 
define ${\Cal S}(M) = {\Cal S}_{\frak K}(M)$ as 
$\{\text{\ortp}(a,M,N):M \le_{\frak K} N \in K_\mu 
\text{ and } a \in N\}$ where
$\text{\ortp}(a,M,N) = \text{ \ortp}_{\frak K}(a,M,N) = 
(M,N,a)/{\Cal E}_M$ where ${\Cal E}_M$ is the transitive closure
of ${\Cal E}^{\text{at}}_M$, and the two-place relation 
${\Cal E}^{\text{at}}_M$  is defined by:

$$
\align
(M,N_1,a_1){\Cal E}^{\text{at}}_M &(M,N_2,a_2) 
\text{ \ub{iff} } M \le_{\frak K}
N_\ell,a_\ell \in N_\ell,\|N_\ell\| = \mu = \|M\| \\
  &\text{ for } \ell=1,2 \text{ and } \text{ there is }
N \in K_\mu \text{ and } \le_{\frak K} \text{-embeddings} \\
  &f_\ell:N_\ell \rightarrow N \text{ for } \ell = 1,2 \text{ such that:} \\
  &f_1 \restriction M = \text{ id}_M = f_2 \restriction M \text{ and }
f_1(a_1) = f_2(a_2).
\endalign
$$
\mn
We may say $p = \text{ \ortp}(a,M,N)$ is the type which $a$ realizes over
$M$ in $N$.  Of course, all those notions depend on ${\frak K}$ so we
may write \ortp$_{\frak K}(a,M,N)$ and ${\Cal E}_M[{\frak K}],
{\Cal E}^{\text{at}}_M[{\frak K}]$. 
\nl
(If in Definition \scite{600-0.2} we do not require $M \in K \Rightarrow
\|M\| \ge \text{\rm LS}({\frak K})$, here we should allow any $N$ such
that $\|M\| \le \|N\| \le M + \text{\rm LS}({\frak K})$.)
\nl
1A) For $M \in {\frak K}_\mu$ let\footnote{we can insist that $N \in
K_\mu$, the difference is not serious} ${\Cal S}_{\frak K}(M) =
\{\text{\rm \ortp}(a,M,N):M \le_{\frak K} N$ and $N \in K_\mu$ or just
$N \in K_{\le(\mu +\text{LS}({\frak K}))}$ and $a \in N\}$ and ${\Cal
S}^{\text{na}}_{\frak K}(M) = \{\ortp(a,M,N):M \le_{\frak K} N$ and
$N \in K_{\le(\mu + \text{LS}({\frak K}))}$ and $a \in N
\backslash M\}$, na stands for non-algebraic.  We may write ${\Cal
S}^{\text{na}}(M)$ omitting ${\frak K}$ when ${\frak K}$ is clear from
the context; so omitting na means $a \in N$ rather than $a \in N
\backslash M$.
\nl
2) Let $M \in K_\mu$ and $M \le_{\frak K} N$.  
We say ``$a$ realizes $p$ in $N$" and ``$p = \text{ \ortp}(a,M,N)$" when:
if $a \in N,p \in {\Cal S}(M)$ and 
$N' \in K_{\le(\mu+\text{LS}({\frak K}))}$ satisfies 
$M \le_{\frak K} N' \le_{\frak K} N$ and $a \in N'$ then 
$p = \text{ \ortp}(a,M,N')$ and there is at least one such $N'$; 
so $M,N' \in K_\mu$  (or just 
$M \le \|N'\| \le \mu + \text{\rm LS}({\frak K}))$
but possibly $N \notin K_\mu$. 
\nl
3) We say ``$a_2$ strongly \footnote{note that 
${\Cal E}^{\text{at}}_M$ is not an
equivalence relation and hence in general is not ${\Cal E}_M$}
realizes $(M,N^1,a_1)/{\Cal E}^{\text{at}}_M$ in $N$" \ub{when} 
for some $N^2$ of cardinality $\le \|M\| + \text{\rm LS}({\frak K})$ 
we have $M \le_{\frak K} N^2 \le_{\frak K} N$ and
$a_2 \in N^2$ and $(M,N^1,a_1)\,{\Cal E}^{\text{at}}_M\,(M,N^2,a_2)$. 
\nl
4) We say $M_0 \in K_\lambda$ is an amalgamation base (in ${\frak K}$,
but normally ${\frak K}$ is understood from the context) \ub{if}: for every
$M_1,M_2 \in K_\lambda$ and $\le_{\frak K}$-embeddings
$f_\ell:M_0 \rightarrow M_\ell$ (for $\ell = 1,2$) there is $M_3 \in
K_\lambda$ and 
$\le_{\frak K}$-embeddings $g_\ell:M_\ell \rightarrow
M_3$ (for $\ell=1,2$) such that $g_1 \circ f_1 = g_2 \circ f_2$.
Similarly for ${\frak K}_{\le \lambda}$.
\nl
4A) ${\frak K}$ has amalgamation in $\lambda$ (or
$\lambda$-amalgamation or ${\frak K}_\lambda$ has amalgamation) \ub{when}
every $M \in K_\lambda$ is an amalgamation base.
\nl
4B) ${\frak K}$ has the $\lambda$-JEP or JEP$_\lambda$ (or ${\frak
K}_\lambda$ has the JEP) when any $M_1,M_2 \in K_\lambda$ can be
$\le_{\frak K}$-embedded into some $M \in K_\lambda$.
\nl
5) We say ${\frak K}$ is stable in $\lambda$ \underbar{if} (LS$({\frak K})
\le \lambda$ and) $M \in K_\lambda \Rightarrow |{\Cal S}(M)| \le \lambda$.
\nl
6) We say $p=q \restriction M$ if $p \in {\Cal S}(M),q \in {\Cal S}(N),
M \le_{\frak K} N$ and for some $N^+,N \le_{\frak K} N^+$ and $a \in N^+$ we
have $p = \text{ \ortp}(a,M,N^+)$ and $q = \text{ \ortp}(a,N,N^+)$; see
\scite{600-0.12A}(1),(2).
We may express this also as ``$q$ extends $p$ or $p$ is
the restriction of $q$ to $M$".  
\nl
7) For finite $m$, for $M \le_{\frak K} N,\bar a \in {}^m N$ we can define
\ortp$(\bar a,M,N)$ and ${\Cal S}^m_{\frak K}(M)$ 
similarly and ${\Cal S}^{< \omega}_{\frak K}(M) = 
\dbcu_{m < \omega} {\Cal S}^m_{\frak K}(M)$; similarly for 
${\Cal S}^\alpha(M)$ (but we shall not use this in any
essential way, so we agree ${\Cal S}(M) = {\Cal S}^1(M)$.)  Again we
may omit ${\frak K}$ when clear from the context.
\nl
8) We say that $p \in {\Cal S}_{\frak K}(M)$ is algebraic when some $a
\in M$ realizes it.
\nl
9) We say that $p \in {\Cal S}_{\frak K}(M)$ is minimal \ub{when} 
it is not algebraic
and for every $N \in K$ of cardinality $\le \|M\| + \text{ LS}({\frak K})$
which $\le_{\frak K}$-extend $M$, the type $p$ has at most one
non-algebraic extension in ${\Cal S}_{\frak K}(M)$.
\enddefinition
\bigskip

\remark{\stag{600-0.12Y} Remark}   1) Note that here ``amalgamation base"
means only for extensions of the same cardinality!
\nl
2) The notion ``minimal type" is important (for categoricity)
but not used much in this chapter.
\endremark
\bigskip

\demo{\stag{600-0.12A} Observation}  
0) Assume 
$M \in K_\mu$ and $M \le_{\frak K} N,a \in N$ \ub{then} \ortp$(a,M,N)$ is
well defined and is $p$ \ub{if} for some $M' \in K_\mu$ we have $M
\cup \{a\} \subseteq M' \le_{\frak K} N$ and 
$p = \text{ \ortp}(a,M,M')$. 
\nl
1) If $M \le_{\frak K} N_1 \le_{\frak K}
N_2,M \in K_\mu$ and $a \in N_1$ \ub{then} 
\ortp$(a,M,N_1)$ is well defined and equal to \ortp$(a,M,N_2)$, (more
transparent if ${\frak K}$ has the $\mu$-amalgamation which is the
real case anyhow).
\nl
2) If $M \le_{\frak K} N$ and $q \in {\Cal S}(N)$ \ub{then} for one and
only one $p$ we have $p = q \restriction M$.
\nl
3) If $M_0 \le_{\frak K} M_1 \le_{\frak K} M_2$ and $p \in {\Cal S}
(M_2)$ then $p \restriction M_0 = (p \restriction M_1) \restriction
M_0$.
\nl
4) If $M \in {\frak K}_\mu$ is an amalgamation base then
${\Cal E}^{\text{at}}_M$ is a transitive relation hence is equal to
${\Cal E}_M$.
\nl
5) If $M \le_{\frak K} N$ are from ${\frak K}_\lambda,M$ is an
amalgamation base and $p \in {\Cal S}(M)$ then there is $q \in {\Cal
S}(N)$ extending $p$, so the mapping $q \mapsto q \restriction M$ is a
function from ${\Cal S}(N)$ onto ${\Cal S}(M)$. 
\enddemo
\bigskip

\demo{Proof}  Easy.  \hfill$\square_{\scite{600-0.12A}}$
\enddemo
\bigskip

\definition{\stag{600-0.13} Definition} 1) We say 
\ub{$N$ is $\lambda$-universal over $M$} when $\lambda \ge \|N\|$ and 
for every $M',M 
\le_{\frak K} M' \in K_\lambda$, there is a $\le_{\frak K}$-embedding of $M'$
into $N$ over $M$.  If we omit $\lambda$ we mean $\|N\|$;
clearly if $N$ is universal over $M$ then $M$ is an amalgamation base. \nl
2)  $K^3_\lambda = \{(M,N,a):M \le_{\frak K} N,a \in N \backslash M$ and
$M,N \in {\frak K}_\lambda\}$, with the partial order $\le$ defined by
$(M,N,a) \le (M',N',a')$ iff $a = a',M \le_{\frak K} M'$ and $N
\le_{\frak K} N'$.  \nl
3) We say $(M,N,a) \in K^3_\lambda$ is minimal \ub{when}: if 
$(M,N,a) \le (M',N_\ell,a)
\in K^3_\lambda$ for $\ell =1,2$ implies \ortp$(a,M',N_1) = 
\text{ \ortp}(a,M',N_2)$ moreover, $(M',N_1,a) 
{\Cal E}^{\text{at}}_\lambda (M',N_2,a)$
(this strengthening is not needed 
if every $M' \in K_\lambda$ is an amalgamation bases).
\nl
4) $N \in {\frak K}$ is $\lambda$-universal if every $M \in 
{\frak K}_\lambda$ can be $\le_{\frak K}$-embedded into it.
\enddefinition
\bigskip

\remark{Remark}  Why do we use 
$\le$ on $K^3_\lambda$?  Because those triples serve us
as a representation of types for which direct limit exists.
\endremark
\bigskip

\definition{\stag{600-0.15} Definition}   1) $M^* \in K_\lambda$ is 
\ub{superlimit} if: clauses (a) + (b) + (c) below hold (and locally
superlimit if clauses (a)$^- + (b) + (c)$ below hold and is pseudo
superlimit if clauses (b) + (c) below hold) where:
\mr
\item "{$(a)$}"  it is universal, (i.e. every $M \in K_\lambda$ 
can be $\le_{\frak K}$-embedded into $M^*$), and
\sn
\item "{$(b)$}"   if $\langle M_i:i \le \delta \rangle$ is 
$\le_{\frak K}$-increasing continuous, 
$\delta < \lambda^+$ and $i < \delta \Rightarrow M_i \cong M^*$ then 
$M_\delta \cong M^*$
\sn
\item "{$(a)^-$}"   if $M^* \le M_1 \in K_\lambda$ then there is $M_2
\in K_2$ which $\le_{\frak K}$-extend $M_1$ and is isomorphic to $M^*$
\sn
\item "{$(c)$}"   there is $M^{**}$ isomorphic to $M^*$ such that $M^*
<_{\frak K} M^{**}$.
\ermn
2) $M$ is $\lambda$-saturated above $\mu$ \ub{when} 
$\|M\| \ge \lambda > \mu \ge
\text{ LS}({\frak K})$ and: $N \le_{\frak K} M,\mu \le \|N\| <
\lambda,N \le_{\frak K} N_1,\|N_1\| \le \|N\| + \text{ LS}({\frak K})$
and $a \in N_1$ then some $b \in M$ strongly realizes $(N,N_1,a)/{\Cal
E}^{\text{at}}_N$ in $M$, see Definition \scite{600-0.12}(3).  
Omitting ``above $\mu$" means ``for some $\mu < \lambda$" hence
``$M$ is $\lambda^+$-saturated" mean that ``$M$ is $\lambda^+$-saturated above
$\lambda$" and $K(\lambda^+$-saturated) 
$= \{M \in K:M$ is $\lambda^+$-saturated$\}$ and ``$M$ is saturated" 
means ``$M$ is $\|M\|$-saturated". 
\enddefinition
\bigskip

\proclaim{\stag{600-0.19}  The Model-homogeneity = Saturativity Lemma}  Let 
$\lambda > \mu + { \text{\rm LS\/}}({\frak K})$ and $M \in K$. \newline
1) $M$ is $\lambda$-saturated above $\mu$ \underbar{iff} $M$ is 
$({\Bbb D}_{{\frak K}_{\ge\mu}},\lambda)$-homogeneous above $\mu$, which 
means: for every
$N_1 \le_{\frak K} N_2 \in K$ such that $\mu \le \|N_1\| \le \|N_2\| < 
\lambda$ and $N_1 \le_{\frak K} M$, there is a $\le_{\frak K}$-embedding $f$ 
of $N_2$ into $M$ over $N_1$. \newline
2)  If $M_1,M_2 \in K_\lambda$ are $\lambda$-saturated above $\mu < \lambda$ 
and for some $N_1 \le_{\frak K} M_1,N_2 \le_{\frak K} M_2$, both 
of cardinality $\in [\mu,\lambda)$, we have $N_1 \cong N_2$ \ub{then} 
$M_1 \cong M_2$; in fact, any isomorphism $f$ from $N_1$ onto $N_2$ 
can be extended to an isomorphism from $M_1$ onto $M_2$. \newline
3) If in (2) we demand only ``$M_2$ is $\lambda$-saturated" and $M_1 \in
K_{\le \lambda}$ \underbar{then} $f$ can be extended to a 
$\le_{\frak K}$-embedding from $M_1$ into $M_2$. \newline
4) In part (2) instead of $N_1 \cong N_2$ it suffices to assume that
$N_1$ and $N_2$ can be $\le_{\frak K}$-embedded into some $N \in K$, 
which holds if ${\frak K}$ has the {\rm JEP} or just {\rm JEP}$_\mu$.
\nl
5) If $N$ is $\lambda$-universal over $M \in K_\mu$ and ${\frak K}$
has $\mu$-JEP then $N$ is $\lambda$-universal (where $\lambda \ge
\text{\rm LS}({\frak K})$ for simplicity).
\nl
6) Assume $M$ is $\lambda$-saturated above $\mu$.  If $N \le_{\frak K}
M$ and $\mu \le \|N\| < \lambda$ \ub{then} $N$ is an amalgamation base
(in $K_{\le(\|N\|+\text{LS}({\frak K}))}$ and 
even in ${\frak K}_{\le \lambda}$) and $|{\Cal S}(N)| \le \|M\|$.
So if every $N \in K_\mu$ can be $\le_{\frak K}$-embedded into $M$
then ${\frak K}$ has $\mu$-amalgamation.
\endproclaim
\bigskip

\demo{Proof}  1)  The ``if" direction is easy as $\lambda > \mu \ge
\text{ LS}({\frak K})$.  
Let us prove the other direction.

We prove this by induction on $\|N_2\|$.  Now first consider the case 
$\|N_2\| > \|N_1\| +
\text{\rm LS}({\frak K})$ then we can find a $\le_{\frak
K}$-increasing continuous sequence $\langle
N_{1,\varepsilon}:\varepsilon < \|N_2\|\rangle$ with union $N_2$ with
$N_{1,0} = N_1$ and $\|N_{1,\varepsilon}\| \le \|N_1\| +
|\varepsilon|$.  Now we choose $f_\varepsilon$, a $\le_{\frak
K}$-embedding of $N_{1,\varepsilon}$ into $M$, increasing continuous
with $\varepsilon$ such that $f_0 = \text{\rm id}_{N_1}$.  For
$\varepsilon = 0$ this is trivial for $\varepsilon$ limit take unions
and for $\varepsilon$ successor use the induction hypothesis.  So
\wilog \, $\|N_2\| \le \|N_1\| + \text{\rm LS}({\frak K})$.

Let $|N_2| = \{ a_i:i < \kappa \}$, and we know $\mu \le \kappa'' :=
\|N_1\| \le \kappa := \|N_2\| \le \kappa' := 
\|N_1\| + \text{\rm LS}({\frak K}) 
< \lambda$; so if, as usual, $\|N_1\| \ge \text{\rm LS}({\frak K})$
then $\kappa' = \kappa$.  We define by induction on $i \le \kappa,
N^i_1,N^i_2,f_i$ such that:
\mr
\item "{(a)}"  $N^i_1 \le_{\frak K} N^i_2$ and $\|N^i_1\| \le
\|N^i_2\| \le \kappa'$
\sn
\item "{(b)}" $N^i_1$ is $\le_{\frak K}$-increasing continuous with $i$
\sn
\item "{(c)}" $N^i_2$ is $\le_{\frak K}$-increasing continuous with $i$
\sn
\item "{(d)}" $f_i$ is a $\le_{\frak K}$-embedding of $N^i_1$ into $M$
\sn
\item "{(e)}" $f_i$ is increasing continuous with $i$
\sn
\item "{(f)}" $a_i \in N^{i+1}_1$
\sn
\item "{(g)}" $N^0_1 = N_1,N^0_2 = N_2,f_0 = \text{\rm id}_{N_1}$.
\ermn
For $i = 0$, clause $(g)$ gives the definition.  For $i$ limit let: \newline
\sn
$N^i_1 = \dsize \bigcup_{j < i} N^j_1$ and 
\sn
$N^i_2 = \dsize \bigcup_{j < i} N^j_2$ and 
\sn
$f_i = \dsize \bigcup_{j < i} f_j$. \newline
\sn
Now (a)-(f) continues to hold by continuity (and $\|N^i_2\| \le 
\kappa'$ easily).
\medskip

For $i$ successor we use our assumption; more elaborately, let
$M^{i-1}_1 \le_{\frak K} M$ be $f_{i-1}(N^{i-1}_1)$ and let $M^{i-1}_2,
g_{i-1}$ be such that $g_{i-1}$ is an isomorphism from $N^{i-1}_2$ \newline
onto $M^{i-1}_2$ extending
$f_{i-1}$, so $M^{i-1}_1 \le_{\frak K} M^{i-1}_2$ (but \wilog \nl
$M^{i-1}_2 \cap M = M^{i-1}_1$).  Now apply the saturation
assumption with $M,M^{i-1}_1$, \nl
$\text{\ortp}(g_{i-1}(a_{i-1}),M^{i-1}_1,M^{i-1}_2)$ 
here standing for $N,M,p$ there
(note: $a_{i-1} \in N_2 = N^0_2 \subseteq N^{i-1}_2$ and 
$\lambda > \kappa' \ge \|N^{i-1}_2\| = \|M^{i-1}_2\| \ge \|M^{i-1}_1\| =
\|N^{i-1}_1\| \ge \|N^0_1\| = \|N_1\| = \kappa'' \ge \mu$ so the requirements
including the requirements on the cardinalities 
in Definition \scite{600-0.15}(2) holds).
So there is $b \in M$ such that \ortp$(b,M^{i-1}_1,M) = \text{ \ortp}(g_{i-1}
(a_{i-1}),M^{i-1}_1,M^{i-1}_2)$.  Moreover (if ${\frak K}$ has amalgamation
in $\mu$ the proof is slightly shorter) remembering the end of the
first sentence in
\scite{600-0.15}(2) which speaks about ``strongly realizes" and recalling
Definition \scite{600-0.12}(3) there is $b \in M$
such that $b$ strongly realizes $(M^{i-1}_1,M^{i-1}_3,g_{i-1}(a_{i-1}))/
{\Cal E}^{\text{at}}_{M^{i-1}_1}$ in $M$.  
This means (see Definition \scite{600-0.12}(3)) that for some 
$M^{i,*}_1$ we have $b \in M^{i,*}_1$ and
$M^{i-1}_1 \le_{\frak K} M^{i,*}_1
\le_{\frak K} M$ and $(M^{i-1}_1,M^{i-1}_2,g_{i-1}(a_{i-1}))
{\Cal E}^{\text{at}}_{M^{i-1}_1}(M^{i-1}_1,M^{i,*}_1,b)$.  This means (see
Definition \scite{600-0.12}(1)) that $M^{i,*}_1$ too has cardinality 
$\le \kappa'$ and there is $M^{i,*}_2 \in K_{\le \kappa'}$ 
such that $M^{i-1}_1 \le_{\frak K} M^{i,*}_2$ and
there are $\le_{\frak K}$-embeddings $h^i_2,h^i_1$ of $M^{i-1}_2,M^{i,*}_1$
into $M^{i,*}_2$ over $M^{i-1}_1$ respectively, such that $h^i_2(g_{i-1}
(a_{i-1})) = h^i_1(b)$. \nl
Now changing names, \wilog \, $h^i_1$ is the identity.
\nl
Let $N^i_2,h_i$ be such that $N^{i-1}_2 \le_{\frak K} N^i_2$ and $h_i$ an 
isomorphism from $N^i_2$ onto $M^{i,*}_2$ extending $g_{i-1}$.  Let
$N^i_1 = h^{-1}_i (M^{i,*}_1)$ and $f_i = (h_i \restriction N^i_1)$.
\medskip

We have carried the induction.  Now $f_\kappa$ is a 
$\le_{\frak K}$-embedding of $N^\kappa_1$ into $M$ over $N_1$, 
but $|N_2| = \{ a_i:i <
\kappa \} \subseteq N^\kappa_1$ hence by AxV of Definition
\scite{600-0.2}, $N_2 \le_{\frak K} N^\kappa_1$, so 
$f_\kappa \restriction N_2:N_2 \rightarrow M$ is as required. 
\newline
2), 3)  By the hence and forth argument (or see \marginbf{!!}{\cprefix{88r}.\scite{88r-2.3}},\marginbf{!!}{\cprefix{88r}.\scite{88r-2.4}}
or see \cite[II,\S3]{Sh:300} = \sectioncite[\S3]{300b}). \nl
4),5),6)  Easy, too.  \hfill$\square_{\scite{600-0.19}}$
\enddemo
\bigskip

\definition{\stag{600-0.21} Definition}  1) For $\sigma = \text{
cf}(\sigma) \le \lambda^+$, we say \ub{$N$ is
$(\lambda,\sigma)$-brimmed over $M$} if ($M \le_{\frak K} N$ are in
$K_\lambda$ and) we can find a sequence $\langle
M_i:i < \sigma \rangle$ which is $\le_{\frak K}$-increasing
\footnote{we have not asked continuity; because in the direction we are
going, it makes no difference if we add ``continuous".  Then we have in
general fewer cases of existence, uniqueness (of being
$(\lambda,\sigma)$-brimmed over $M \in K_\lambda$) does not need extra
assumptions and existence is harder}, $M_i \in
K_\lambda,M_0 = M,M_{i+1}$ is $\le_{\frak K}$-universal
\footnote{hence $M_i$ is an amalgamation base}  over $M_i$ and
$\dbcu_{i < \sigma} M_i = N$.  We say $N$ is $(\lambda,\sigma)$-brimmed
over $A$ if $A \subseteq N \in K_\lambda$ and we can find $\langle M_i:i <
\sigma \rangle$ as above such that $A \subseteq M_0$ but $M_0
\restriction A \le_{\frak K} M_0 \Rightarrow M_0 = A$; if $A =
\emptyset$ we may omit ``over $A$".  We say continuously
$(\lambda,\sigma)$-brimmed (over $M$) \ub{when} the sequence $\langle M_i:i
< \sigma\rangle$ is $\le_{\frak K}$-increasing continuous; if ${\frak
K}_\lambda$ has amalgamation, the two notions coincide.
\nl
2) We say 
$N$ is $(\lambda,*)$-brimmed over $M$ \ub{if} for some $\sigma \le
\lambda,N$ is $(\lambda,\sigma)$-brimmed over $M$.  We say $N$ is
$(\lambda,*)$-brimmed if for some $M,N$ is $(\lambda,*)$-brimmed
over $M$. 
\nl
3) If $\alpha < \lambda^+$ let ``$N$ is $(\lambda,\alpha)$-brimmed over
$M$" mean $M \le_{\frak K} N$ are from $K_\lambda$ and cf$(\alpha) \ge
\aleph_0 \Rightarrow N$ is $(\lambda,\text{cf}(\alpha))$-brimmed over $M$.
\enddefinition
\bn
On the meaning of $(\lambda,\sigma)$-brimmed for elementary classes,
see \scite{600-Ex.0}(2) below.  Recall
\proclaim{\stag{600-0.22} Claim}  Assume $\lambda \ge \text{\rm LS}({\frak K})$.
\nl
1) If ${\frak K}$ has amalgamation in $\lambda$, is stable in $\lambda$
and $\sigma = {\text{\rm cf\/}}(\sigma) \le \lambda$, \ub{then}
\mr
\item "{$(a)$}"   for every $M \in {\frak K}_\lambda$ there is $N,M
\le_{\frak K} N \in K_\lambda$, universal over $M$
\sn
\item "{$(b)$}"  for every $M \in {\frak K}_\lambda$ there is 
$N \in {\frak K}_\lambda$ which is $(\lambda,\sigma)$-brimmed 
over $M$ 
\sn
\item "{$(c)$}"  if $N$ is $(\lambda,\sigma)$-brimmed over $M$
\ub{then} $N$ is universal over $M$. 
\ermn
2) If $N_\ell$ is $(\lambda,\aleph_0)$-brimmed over $M$ for $\ell =1,2$,
\ub{then} $N_1,N_2$ are isomorphic over $M$. 
\nl
3) Assume $\sigma = { \text{\rm cf\/}}(\sigma) \le \lambda^+$, and for
every $\aleph_0 \le \theta = { \text{\rm cf\/}}(\theta) 
< \sigma$ any $(\lambda,\theta)$-brimmed model
is an amalgamation base (in ${\frak K}$).  \ub{Then}:
\mr
\item "{$(a)$}"  if $N_\ell$ is $(\lambda,\sigma)$-brimmed 
over $M$ for $\ell=1,2$ \ub{then} $N_1,N_2$ are isomorphic over $M$
\sn
\item "{$(b)$}"  if ${\frak K}$ has JEP$_\lambda$ (i.e., the 
joint embedding property in $\lambda$) and 
$N_1,N_2$ are $(\lambda,\sigma)$-brimmed \ub{then} $N_1,N_2$ are isomorphic.
\ermn
3A) There is a $(\lambda,\sigma)$-brimmed model $N$ over $M \in
K_\lambda$ \ub{when}: for every $\le_{{\frak K}_\lambda}$-extension
$M_1$ of $M$ there is a $\le_{{\frak K}_\lambda}$-extension $M_2$ of
$M_1$ which is an amalgamation base and there is a $\lambda$-universal
extension $M_3 \in K_\lambda$ of $M_2$.
\nl
4) Assume ${\frak K}$ has $\lambda$-amalgamation and the $\lambda$-JEP
and $\bar M = \langle M_i:i \le \lambda \rangle$ is $\le_{\frak K}$-increasing
continuous and $M_i \in K_\lambda$ for $i \le \lambda$.
\mr
\item "{$(a)$}"  If $\lambda$ is regular and for every $i < \lambda,p
\in {\Cal S}(M_i)$ for some $j \in (i,\lambda)$, some $a \in M_j$
realizes $p$, \ub{then} $M_\lambda$ is universal over $M_0$ and is
$(\lambda,\lambda)$-brimmed over $M_0$
\sn
\item "{$(b)$}"  if for every $i < \lambda$ every $p \in {\Cal
S}(M_i)$ is realized in $M_{i+1}$ \ub{then} $M_\lambda$ is
$(\lambda,\text{\rm cf}(\lambda))$-brimmed over $M_0$.
\ermn 
5)  Assume $\sigma = {\text{\rm cf\/}}(\sigma) \le \lambda$ and $M \in
{\frak K}$ is continuous $(\lambda,\sigma)$-brimmed.  \ub{Then} $M$ is
locally a
$(\lambda,\{\sigma\})$-strongly limit model in ${\frak K}_\lambda$
(see Definition \marginbf{!!}{\cprefix{88r}.\scite{88r-3.1}}(2),(7), not used).
\endproclaim
\bigskip

\demo{Proof}  1) Clause (c) holds by Definition \scite{600-0.21}.  

As for clause (a), for any given $M \in K_\lambda$, easily there 
is an $\le_{\frak K}$-increasing continuous sequence 
$\langle M_i:i \le \delta \rangle$ of models from $K_\lambda,M_0 = M$
such that $p \in {\Cal S}(M_i) \Rightarrow p$ is realized in $M_{i+1}$,
this by stability + amalgamation.  So $\langle M_i:i \le \lambda
\rangle$ is as in part (4) hence by clause (b) of part (4) we get that
$M_\delta$ is $\le_{\frak K}$-universal over $M_0=M$ so we are done.
Clause (b) follows by (a) as ${\frak K}_\lambda$ has $\lambda$-amalgamation.
\nl
2) By (3)(a).
\nl
3) \ub{Clause (a)} holds by the hence and forth argument, that is
assume $\langle N_{\ell,i}:i \le \sigma \rangle$ be $\le_{\frak K}$-increasing
continuous, $N_{\ell,0} =M,N_{\ell,i+1}$ is universal over
$N_{\ell,i}$ and $N_\ell = N_{\ell,\sigma}$ so $N_{\ell,i} \in
{\frak K}_\lambda$.  We now choose $f_i$ by induction on $i \le
\sigma$ such that:
\mr
\widestnumber\item{$(iii)$}
\item "{$(i)$}"   if $i$ is odd, $f_i$ is a $\le_{\frak K}$-embedding
of $N_{1,i}$ into $N_{2,i}$
\sn
\item "{$(ii)$}"   if $i$ is even, $f^{-1}_i$ is a 
$\le_{\frak K}$-embedding of $N_{2,i}$ into $N_{1,i}$
\sn
\item "{$(iii)$}"   if $i$ is limit then $f_i$ is an
isomorphism from $N_{1,i}$ onto $N_{2,i}$
\sn
\item "{$(iv)$}"   $f_i$ is increasing continuous with $i$
\sn
\item "{$(v)$}"   if $i=0$ then $f_0 = \text{\rm id}_M$.
\ermn
For $i=0$ let $f_0 = \text{ id}_M$.  If $i=2j+2$ use ``$N_{1,i}$ is
a universal extension of $N_{1,2j+1}$ (in ${\frak K}_\lambda$) and
$f_{2j+1}$ is a $\le_{\frak K}$-embedding of $N_{1,2j+1}$ into
$N_{2,2j+1}$ (by clause (i) applied to $2j+1$)
 and $N_{1,2j+1}$ is an amalgamation base".  That is, $N_{2,i}$ is a
$\le_{\frak K}$-extension of $f_{2j+1}(N_{2j+1})$ which is an
amalgamation base so $f^{-1}_{2j+1}$ can be extended to a $\le_{\frak
K}$-embedding of $f^{-1}_i$ of $N_{2,i}$ into $N_{1,i}$.
For $i=2j+1$ use ``$N_{2,i}$ is a universal extension 
(in ${\frak K}_\lambda$) of
$N_{2,2j}$ and $f^{-1}_{2j}$ is a $\le_{\frak K}$-embedding of
$N_{2,2j}$ into $N_{1,2j}$ and $N_{2,2j}$ is an amalgamation base (in
${\frak K}_\lambda$)".
\nl
For $i$ limit let $f_i = \cup\{f_j:j < i\}$.  Clearly $f_\sigma$ is an
isomorphism from $N_1 = N_{1,\sigma}$ onto $N_{2,\sigma} = N_2$ so we
are done.

As for clause (b), we can assume that $\langle N_{\ell,i}:i \le \sigma
\rangle$ exemplifies ``$N_\ell$ is $(\lambda,\sigma)$-brimmed" for
$\ell=1,2$.  By the JEP$_\lambda$ there is a pair $(g_1,N)$ such that
$N_{1,0} \le_{\frak K} N \in K_\lambda$ and $g_1$ is a $\le_{\frak
K}$-embedding of $N_{2,0}$ into $N$.  As above there is a $\le_{\frak
K}$-embedding $g_2$ of $N$ into $N_{1,1}$ over $N_{1,0}$.  Let $f_0 =
g_2$ and continue as in the proof of clause (a).
\nl
3A) Easy, too.
\nl
4) We first proved weaker versions of (a) and of (b) called
(a)$^-$,(b)$^-$ respectively.
\bn
\ub{Clause $(a)^-$}:  Like (a) but we conclude only: $M_\lambda$ is
universal over $M_0$.

Let $M_0 \le_{\frak K} N \in K_\lambda$ and we
let $\langle S_i:i < \lambda \rangle$ be a partition of $\lambda$ to
$\lambda$ sets each with $\lambda$ members, $i \le \text{\rm
Min}(S_i)$.  Let $M_{1,i} = M_i$ for $i \le \lambda$ and we choose
$\langle M_{2,i}:i \le \delta \rangle$ which is $\le_{\frak
K}$-increasing such that $M_{2,i} \in {\frak K},
M_{2,0} = M_{1,0},N \le_{\frak K} M_{2,1}$ and every type
$p \in {\Cal S}(M_{2,i})$ is realized in $M_{2,i+1}$.  We shall prove
that $M_{1,\lambda},M_{2,\lambda}$ are isomorphic over 
$M_0 = M_{1,0}$,  this clearly suffices.  We choose a quintuple 
$(j_i,M_{3,i},f_{1,i},f_{2,i},\bar{\bold a}_i)$ by induction
on $i < \lambda$ such that
\mr
\item "{$\circledast$}"  $(a) \quad j_i < \lambda$ is increasing
continuous
\sn
\item "{${{}}$}"  $(b) \quad M_{3,i} \in K_\lambda$ is $\le_{\frak
K}$-increasing continuous
\sn
\item "{${{}}$}"  $(c) \quad f_{\ell,i}$ is a 
$\le_{\frak K}$-embedding of $M_{\ell,j_i}$ into $M$ for $\ell=1,2$
\sn
\item "{${{}}$}"  $(d) \quad f_{\ell,i}$ is increasing continuous with
$i$ for $\ell=1,2$
\sn
\item "{${{}}$}"  $(e) \quad \bar{\bold a}_i = \langle
a^i_\varepsilon:\varepsilon \in S_i\rangle$ lists the members of
$M_{3,i}$
\sn
\item "{${{}}$}"  $(f) \quad$ if $\varepsilon \in S_i$ then
$a^i_\varepsilon \in \text{\rm Rang}(f_{1,2 \varepsilon +1})$ and
$a^i_\varepsilon \in \text{\rm Rang}(f_{2,2 \varepsilon +2})$.
\ermn
If we succeed then $f_\ell := \cup\{f_{\ell,i}:i < \lambda\}$ is a
$\le_{\frak K}$-embedding of $M_{1,\lambda}$ into $M_{3,\lambda} :=
M_3 := \cup\{M_{3,i}:i < \lambda\}$ and this embedding is onto because
$a \in M_3 \Rightarrow$ for some $i < \lambda,a \in M_{3,i}
\Rightarrow$ for some $i < \lambda$ and $\varepsilon \in S_i,a =
a^i_\varepsilon \Rightarrow a = a^i_\varepsilon \in \text{\rm
Rang}(f_{\ell,\varepsilon +1}) \Rightarrow a \in \text{\rm
Rang}(f_\ell)$.  So $f^{-1}_1 \circ f_2$ is an isomorphism from
$M_{2,\lambda}$ onto $M_{1,\lambda} = M_\lambda$ so as said above we are
done.

Carrying the induction; for $i=0$ use ``${\frak K}$ has the
$\lambda$-JEP" for $M_{1,0},M_{2,0}$.

For $i$ limit take unions.

For $i= 2 \varepsilon +1$ let $j_i = \text{ min}\{j < \lambda_i:j >
j_{2 \varepsilon}$ and $(f^1_{2 \varepsilon})^{-1} \,
(\text{tp}(a^i_\varepsilon,f^1_{2 \varepsilon}(M_{1,i}),M_{3,i}) \in
{\Cal S}_{\frak K}(M_{1,i})$ is realized in $M_j$ and continue as in
the proof of \scite{600-0.19}(1), so can avoid using ``$(f^1_i)^{-2}$ of a
type.  

For $i = 2 \varepsilon +2$, the proof is similar.
\bn
\ub{Clause $(b)^-$}:  Like clause (b) but we conclude only:
$M_\lambda$ is universal over $M_0$.

Similar to the proof of $(a)^-$ except that we weaken clause (f) to
\mr
\item "{$(f)^-$}"  if $\varepsilon \in S_i$ then $a^i_\varepsilon \in
\text{ Rang}(f_{1,2 \varepsilon +1})$.
\endroster
\bn
\ub{Clauses $(a),(b)$}:

As above but now in the choice of $\langle M_{2,i}:i \le \lambda
\rangle$ we demand that ``$M_{2,i+1}$ is $\le_{\frak K}$-universal
over $M_{2,i}$" for $i < \lambda$.  This is permissible as by $(b)^-$
which we have already proved every member of $K_\lambda$ has a
$<_{{\frak K}_\lambda}$-extension
from ${\frak K}_\lambda$ which is universal over it.  So
$M_{2,\lambda}$ is $(\lambda,\text{\rm cf}(\lambda))$-brimmed over
$M_{2,0} = M_0$ hence also $M_\lambda$ being isomorphic to
$M_{2,\lambda}$ over $M_0$ is $(\lambda,\text{\rm
cf}(\lambda))$-brimmed over $M_0$, as required.
\nl
5) Easy and not used.  (Let $\langle M_i:i \le \sigma\rangle$ witness
``$M$ is $(\lambda,\sigma)$-brimmed", so $M$ can be $\le_{\frak
K}$-embedded into $M_i$, hence \wilog \, $M_0 \cong M_1$.)  Now use
$\bold F$ such that $\bold F(M')$ is a 
$\le_{{\frak K}_\lambda}$-extension of $M'$ which is 
$\le_{{\frak K}_\lambda}$-universal over it and is an amalgamation base.
  \hfill$\square_{\scite{600-0.22}}$ 
\enddemo
\bigskip

\proclaim{\stag{600-4a.10} Claim}  1) Assume that ${\frak K}$ is an {\rm
a.e.c., LS}$({\frak K}) \le \lambda$ and ${\frak K}$ has
$\lambda$-amalgamation and is stable in $\lambda$.  \ub{Then} there is a
saturated $N \in K_{\lambda^+}$.  Also for every saturated $N \in
K_{\lambda^+}$ (in ${\frak K}$, above $\lambda$ of course) we can find
a $\le_{\frak K}$-representation $\bar N = \langle N_i:i < \lambda^+
\rangle$, with $N_{i+1}$ being $(\lambda,{\text{\rm
cf\/}}(\lambda))$-brimmed over $N_i$ and $N_0$ being
$(\lambda,\lambda)$-brimmed. 
\nl
2) If for $\ell=1,2$ we have $\bar N^\ell = \langle N^\ell_i:i < \lambda^+
\rangle$ as in part (1), \ub{then} there is an isomorphism $f$ from
$N^1$ onto $N^2$ mapping $N^1_i$ onto $N^2_i$ for each $i <
\lambda^+$.  Moreover, for any $i < \lambda^+$ and isomorphism $g$
from $N^1_i$ onto $N^2_i$ we can find an isomorphism $f$ from $N^1$
onto $N^2$ extending $g$ and mappng $N^1_j$ onto $N^2_j$ for each $j
\in [i,\lambda^+)$. 
\nl
3) If $N^0 \le_{\frak K} N^1$ are both saturated (above $\lambda$) and
are in $K_{\lambda^+}$ (hence {\rm LS}$({\frak K}) \le \lambda$), 
\ub{then} we can find $\le_{\frak K}$-representation $\bar
N^\ell$ of $N^\ell$ as in (1) for $\ell=1,2$ 
with $N^0_i = N^0 \cap N^1_i$, (so $N^0_i \le_{\frak K} N^1_i$) for $i
< \lambda^+$. 
\nl
4) If $M \in K_{\lambda^+}$ and ${\frak K}$ has $\lambda$-amalgamation
and is stable in $\lambda$ (and {\rm LS}$({\frak K}) \le \lambda$), 
\ub{then} for some $N \in K_{\lambda^+}$ saturated (above 
$\lambda$) we have $M \le_{\frak K} N$.
\endproclaim
\bigskip

\demo{Proof}  Easy (for (2),(3) using \scite{600-0.19}(6)), e.g.
\nl
4)  There is a $\le_{\frak K}$-increasing continuous sequence
$\langle M_i:i < \lambda^+\rangle$ with union $M$ such that $M_i \in
K_\lambda$.  Now we choose $N_i$ by induction on $i < \lambda$
\mr
\item "{$(*)$}"  $(a) \quad N_i \in K_\lambda$ is $\le_{\frak
K}$-increasing continuous
\sn
\item "{${{}}$}"  $(b) \quad N_{i+1}$ is 
$(\lambda,\text{\rm cf}(\lambda))$-brimmed over $N_i$
\sn
\item "{${{}}$}"  $(c) \quad N_0 = M_0$.
\ermn
This is possible by \scite{600-0.22}(1).  Then by induction on $i \le
\lambda^+$ we choose a $\le_{\frak K}$-embedding $f_i$ of $M_i$ into $N_i$,
increasing continuous with $i$.  For $i=0$ let $f_i = \text{\rm
id}_{M_0}$.   For $i$ limit use union.

Lastly, for $i=j+1$ use ``${\frak K}$ has $\lambda$-amalgamation" and
``$N_j$ is universal over $N_i$".  Now by renaming \wilog \,
$f_{\lambda^+} = \text{\rm id}_{N_{\lambda^+}}$ and we are done.  (Of
course, we hae assumed less).    \hfill$\square_{\scite{600-4a.10}}$
\enddemo
\bn
You may wonder why in this work we have not restricted our ${\frak K}$ to 
``abstract elementary class in $\lambda$" say in \S2 below (or in
\cite{Sh:576}); by the following facts (mainly \scite{600-0.31}) 
this is immaterial.
\bigskip

\definition{\stag{600-0.30A} Definition}  1) We say that ${\frak K}_\lambda$ is
a $\lambda$-abstract elementary class or $\lambda$-a.e.c.
in short, \ub{when}:
\mr
\item "{$(a)$}"  ${\frak K}_\lambda = (K_\lambda,\le_{{\frak K}_\lambda})$,
\sn
\item "{$(b)$}"  $K_\lambda$ is a class of $\tau$-models of cardinality
$\lambda$ closed under isomorphism for some vocabulary $\tau = 
\tau_{{\frak K}_\lambda}$,
\sn
\item "{$(c)$}"  $\le_{{\frak K}_\lambda}$ a partial order of $K_\lambda$,
closed under isomorphisms
\sn
\item "{$(d)$}"  axioms (0 and) I,II,III,IV,V of abstract elementary classes
(see \scite{600-0.2}) hold except that in AxIII we demand $\delta <
\lambda^+$ (you can demand this also in AxIV).
\ermn
2) For an abstract elementary class ${\frak K}$ let ${\frak K}_\lambda = 
(K_\lambda,\le_{\frak K} \restriction K_\lambda)$ and similarly
${\frak K}_{\ge \lambda},{\frak K}_{\le \lambda},{\frak
K}_{[\lambda,\mu]}$ and define $(\le \lambda)$-a.e.c. and
$[\lambda,\mu]$-a.e.c., etc.
\nl
3) Definitions \scite{600-0.12}, \scite{600-0.13}, \scite{600-0.15}, \scite{600-0.21}
apply to $\lambda$-a.e.c. ${\frak K}_\lambda$.
\enddefinition
\bigskip

\demo{\stag{600-0.30B} Observation}  If ${\frak K}^1$ is an a.e.c. with
$K^1_\lambda \ne \emptyset$ then
\mr
\item "{$(a)$}"   ${\frak K}^1_\lambda$ is a $\lambda$-a.e.c.
\sn
\item "{$(b)$}"  if ${\frak K}^2_\lambda$ is a $\lambda$-a.e.c., and
${\frak K}^1_\lambda = {\frak K}^2_\lambda$ \ub{then} Definitions
\scite{600-0.12}, \scite{600-0.13}, \scite{600-0.15}, \scite{600-0.21} when applied to
${\frak K}^1$ but restricting ourselves to models of cardinality
$\lambda$ and when applied to ${\frak K}^2_\lambda$ are equivalent.
\endroster
\enddemo
\bigskip

\demo{Proof}  Just read the definitions.
\enddemo
\bn
We may wonder 
\nl
\margintag{600-0.30C}\ub{\stag{600-0.30C} Problem}:  Suppose ${\frak K}^1,{\frak K}^2$ are
a.e.c. such that for some $\lambda > \mu \ge \text{ LS}({\frak K}^1)$,
$\text{LS}({\frak K}^2),{\frak K}^1_\lambda = {\frak K}^2_\lambda$.
Can we bound the first such $\lambda$ above $\mu$?  (Well, better
bound than the Lowenheim number of $\Bbb L_{\mu^+,\mu^+}$(second order)).
\bigskip

\demo{\stag{600-0.30D} Observation}  1) Let ${\frak K}$ be an a.e.c. with
$\lambda = \text{ LS}({\frak K})$ and $\mu \ge \lambda$ and we define
${\frak K}_{\ge \mu}$ by: $M \in {\frak K}_{\ge \mu}$ iff $M \in K
\and \|M\| \ge\mu$ and $M \le_{{\frak K}_{\ge \mu}} N$ if $M
\le_{\frak K} N$ and $\|M\|,\|N\| \ge \mu$.  \ub{Then} ${\frak K}_{\ge
\mu}$ is an a.e.c. with LS$({\frak K}_{\ge \mu}) = \mu$.
\nl
2) If ${\frak K}_\lambda$ is a $\lambda$-a.e.c. then observation
\scite{600-0.6} holds when $|I| \le \lambda$.
\nl
3) Claims \scite{600-0.12A}(2)-(5), \scite{600-0.22} apply to $\lambda$-a.e.c.
\enddemo
\bigskip

\remark{\stag{600-0.30D.1} Remark}  If 
${\frak K}$ is an a.e.c. with Lowenheim-Skolem number $\lambda$, then
every model of ${\frak K}$ can be written as a direct limit
(by $\le_{\frak K}$) of members of ${\frak K}_\lambda$ (see
\scite{600-0.7}(1)).  Alternating we prove below that
given a $\lambda$-abstract elementary class ${\frak K}_\lambda$, the
class of direct limits of members of ${\frak K}_\lambda$ is an
a.e.c. ${\frak K}^{\text{up}}$.  
We show below $({\frak K}_\lambda)^{\text{up}} = {\frak K}$, hence
${\frak K}_\lambda$ determines ${\frak K}_{\ge \lambda}$.
\endremark
\bigskip

\proclaim{\stag{600-0.31} Lemma}  Suppose 
${\frak K}_\lambda$ is a $\lambda$-abstract elementary class. \nl
1) The pair $(K',\le_{{\frak K}'})$ is an abstract elementary class with
Lowenheim-Skolem number $\lambda$ which we denote also by 
${\frak K}^{\text{up}}$ where we define

$$
\align
K' = \biggl\{ M:&M \text{ is a } \tau_{{\frak K}_\lambda} \text{-model, and
for some directed partial order} \tag"{$(a)$}" \\
  &I \text{ and } \bar M = \langle M_s:s \in I \rangle \text{ we have} \\
  &\qquad \qquad M = \dsize \bigcup_{s \in I} M_s \\
  &\qquad \qquad s \in I \Rightarrow M_s \in K_\lambda \\
  &I \models s < t \Rightarrow M_s \le_{{\frak K}_\lambda} M_t \biggr\}.
\endalign
$$
\mn
We call such $\langle M_s:s \in I \rangle$ a witness for $M \in K'$,
we call it reasonable if $|I| \le \|M\|$

$$
\align
M \le_{{\frak K}'} N \text{ \underbar{iff} } &\text{for some directed partial
order } J,\text{ and} \tag"{$(b)$}" \\
  &\text{directed } I \subseteq J \text{ and } \langle M_s:s \in J \rangle
\text{ we have} \\
  &M = \dsize \bigcup_{s \in I} M_s,N = \dsize \bigcup_{t \in J} M_t,
M_s \in K_\lambda \text{ and} \\
  &J \models s < t \Rightarrow M_s \le_{{\frak K}_\lambda} M_t.
\endalign
$$
\mn
We call such $I, \langle M_s:s \in J \rangle$ witnesses for $M 
\le_{{\frak K}'} N$ or say $(I,J,\langle M_s:s \in J \rangle)$ witness
$M \le_{{\frak K}'} N$. \nl
2) Moreover, $K'_\lambda = K_\lambda$ and $\le_{{\frak K}'_\lambda}$
(which means $\le_{{\frak K}'} \restriction K'_\lambda$) is equal to 
$\le_{{\frak K}_\lambda}$ so $({\frak K}')_\lambda = {\frak K}_\lambda$. 
\nl
3) If ${\frak K}''$ is an abstract elementary class satisfying
(see \scite{600-0.30D})
$K''_\lambda = K_\lambda,<_{{\frak K}''} \restriction K_\lambda = 
\le_{{\frak K}_\lambda}$ and ${\text{\rm LS\/}}({\frak K}'') \le 
\lambda$ \ub{then}
\footnote{if we assume in addition that 
$M \in {\frak K}'' \Rightarrow \|M\| \ge \lambda$ then we can show that
equality holds} ${\frak K}''_{\ge \lambda} = {\frak K}'$. \nl
4)  If ${\frak K}''$ is an a.e.c., $K_\lambda \subseteq K''_\lambda$ 
and $\le_{{\frak K}_\lambda} = \le_{{\frak K}''} \restriction
K_\lambda$,   
\ub{then} $K' \subseteq K''$ and $\le_{{\frak K}'} = \le_{{\frak K}''}
\restriction K'$.
\endproclaim
\bigskip

\demo{Proof}  The proof of part (2) is straightforward (recalling
\scite{600-0.6}) and part (3) follows
from \scite{600-0.7} and part (4) is easier than part (3) (both statements
being proved by induction on the cardinality of the relevant models)
hence we concentrate on part (1).  So let us check the axioms one by one.
\enddemo
\bigskip
\noindent
\underbar{Ax 0}:  $K'$ is a class of $\tau$ models, $\le_{{\frak K}'}$
a two-place relation on $K'$, both closed under isomorphisms. \newline
[Why?  trivially.]
\bn
\underbar{Ax I}:  If $M \le_{{\frak K}'} N$ then $M \subseteq N$. \newline
[Why?  trivial.]
\bn
\underbar{Ax II}:  $M_0 \le_{{\frak K}'} M_1 \le_{{\frak K}'} M_2$ 
implies $M_0 \le_{{\frak K}'} M_2$ and $M \in K' \Rightarrow M 
\le_{{\frak K}'} M$. \newline
[Why?  The second phrase is trivial (as if $\bar M = \langle M_t:t \in
I\rangle$ witness $M \in K'$ then $(I,I,\bar M)$ witness $M
\le_{{\frak K}'} M$ above).  
For the first phrase let for $\ell \in
\{1,2\}$ the directed partial orders $I_\ell \subseteq J_\ell$ and
$\bar M^\ell = \langle M^\ell_s:s \in J_\ell \rangle$ witness 
$M_{\ell-1} \le_{{\frak K}'}M_\ell$ and let $\bar M^0 = 
\langle M^0_s:s \in I_0 \rangle$ witness $M_0 \in K'$.  
Now \wilog \, $\bar M^0$ is reasonable, i.e.
$|I_0| \le \|M_0\|$, why? by
\mr
\item "{$\boxtimes_1$}"   every $M \in K'$ has a reasonable witness,
in fact, if $\bar M = \langle M_t:t \in I \rangle$ is a witness for
$M$ then for some $I' \subseteq I$ of cardinality $\le \|M\|$ we have
$\bar M \restriction I'$ is a reasonable witness for $M$.
\nl
[Why?  If $\bar M = \langle M_t:t \in I \rangle$ is a witness, for each
$a \in M$ choose $t_a \in I$ such that $a \in M_{t_a}$ and let
$F:[I]^{< \aleph_0} \rightarrow I$ be such that $F(\{t_1,\dotsc,t_n\})$ is an
upper bound of $\{t_1,\dotsc,t_n\}$ and let $J$ be the closure of
$\{t_a:a \in M\}$ under $F$; now $\bar M \restriction J$ is a reasonable
witness of $M \in K'$.]
\ermn
Similarly
\mr
\item "{$\boxtimes_2$}"   if $(I,J,\langle M_s:s \in J\rangle$ witness 
$M \le_{{\frak K}'} N$ then for some directed $I' \subseteq I,|I'| \le
\|M\|$ we have $(I',J,\langle M_s:s \in J\rangle)$ witness $M \le_{K'} N$ 
\sn
\item "{$\boxtimes_3$}"  if $I,\bar M = \langle M_t:t \in J \rangle$ 
witness $M \le_{{\frak K}'} N$ \ub{then} for some directed 
$J' \subseteq J$ we have $\|J'\| \le |I| + \|N\|,I \subseteq J'$ and 
$I,\bar M \restriction J'$ witness $M \le_{{\frak K}'} N$.
\ermn
Clearly $\boxtimes_1$ (and $\boxtimes_2,\boxtimes_3$) are 
cases of the LS-argument.
We shall find a witness $(I,J,\langle M_s:s \in J \rangle$) for $M_0 
\le_{{\frak K}'} M_2$ such that $\langle M_s:s \in I \rangle = \langle M^0_s:
s \in I_0 \rangle$ so $I=I_0$ and $|J| \le \|M_2\|$.
This is needed for the proof of Ax III below.  
Without loss of generality $I_1,I_2$ has cardinality $\le \|M_0\|,
\|M_1\|$ respectively, by the proof of $\boxtimes_2$.  Also \wilog
\, $\bar M^1,\bar M^1 \restriction I_1,\bar M^2,\bar M^2 \restriction
I_2$ are reasonable as by the same argument we can have 
$|J_1| \le \|M_1\|,|J_2| \le \|M_2\|$ by $\boxtimes_3$.

As $\langle M^0_s:s \in I_0 \rangle$ 
is reasonable, there is a one-to-one 
function $h$ from $I_0$ into $M_2$ (and even $M_0$); the function $h$
will be used to get that $J$ defined below is directed.
We 
choose by induction on $m < \omega$, for every $\bar c \in {}^m(M_2)$, sets
$I_{0,\bar c},I_{1,\bar c},I_{2,\bar c},J_{1,\bar c},J_{2,\bar c}$ such
that:
\mr
\item "{$\otimes_1(a)$}"  $I_{\ell,\bar c}$ 
is a directed subset of $I_\ell$ of
cardinality $\le \lambda$ for $\ell \in \{0,1,2\}$
\sn
\item "{$(b)$}"  $J_{\ell,\bar c}$ is a directed subset of $J_\ell$ of
cardinality $\le \lambda$ for $\ell \in \{1,2\}$
\sn
\item "{$(c)$}"  $\dsize \bigcup_{s \in I_{\ell +1,\bar c}} M^{\ell + 1}_s =
\bigl( \dsize \bigcup_{s \in J_{\ell + 1,\bar c}} M^{\ell + 1}_s \bigr) \cap
M_\ell$ for $\ell = 0,1$
\sn
\item "{$(d)$}"  $\dsize \bigcup_{s \in I_{0,\bar c}} M^0_s =
(\dsize \bigcup_{s \in I_{1,\bar c}} M^1_s) \cap M_0$
\sn
\item "{$(e)$}"  $\dsize \bigcup_{s \in J_{1,\bar c}} M^1_s =
\dsize \bigcup_{s \in I_{2,\bar c}} M^2_s$
\sn
\item "{$(f)$}"  $\bar c \subseteq \dsize \bigcup_{s \in J_{2,\bar c}} 
M^2_s$
\sn
\item "{$(g)$}"  if $\bar d$ is a permutation of $\bar c$ 
(i.e., letting $m = \ell g(\bar c)$ for some one
to one $g:\{0,\dotsc,m-1\} \rightarrow \{0,\dotsc,m-1\}$ we have 
$d_\ell = c_{g(\ell)}$) \underbar{then} $I_{\ell,\bar c} = 
I_{\ell,\bar d},J_{m,\bar c} = J_{m,\bar d}$ \nl
(for $\ell \in \{0,1,2\},m \in \{1,2\}$)
\sn
\item "{$(h)$}"  if $\bar d$ is a subsequence of $\bar c$ (equivalently:
an initial segment of some permutation of $\bar c$)
\underbar{then} $I_{\ell,\bar d} \subseteq I_{\ell,\bar c},
J_{m,\bar d} \subseteq J_{m,\bar c}$ for $\ell \in \{0,1,2\},m \in \{1,2\}$
\sn
\item "{$(i)$}"  if $h(s) = c$ so $s \in I_0$ then $s \in I_{0,<c>}$.
\ermn
There is no problem to carry the definition by LS-argument recalling
clauses (a) + (b) and $\|M^\ell_s\| = \lambda$ when $\ell=0 \wedge s \in
I_0$ or $\ell=1 \wedge s \in J_1$ or $\ell=2 \wedge s \in J_2$.  
Without loss of generality 
$I_\ell \cap {}^{\omega >}(M_2) = \emptyset$.

Now let $J$ have as set of elements $I_0 \cup \{\bar c:\bar c \text{ a finite
sequence from }M_2\}$ ordered by: $J \models x \le y$ iff $I_0 \models x \le
y$ \underbar{or} $x \in I_0,y \in J \backslash I_0,\exists z \in
I_{0,y}[x \le_{I_0} z]$ \ub{or} 
$x,y \in J \backslash I_0$ and $x$ is an initial segment of a
permutation of $y$ (or you may identify $\bar c$ with its set of 
permutations).
\mn
Let $I = I_0$. \newline
Let $M_x$ be $M^0_x$ if $x \in I_0$ and $\dbcu_{s \in J_{2,x}}
M^2_s$ if $x \in J \backslash I_0$. \nl
Now
\mr
\item "{$(*)_1$}"  $J$ is a partial order
\nl
[Clearly $x \le_J y \le_J x \Rightarrow x=y$, hence it is enough to
prove transitivity.
Assume $x \le_J y \le_J z$; if all three are in $I_0$ use ``$I_0$ is a
partial order", if all three are not in $J \backslash I_0$, use 
the definition of
the order.  As $x' \le_J y' \in I_0 \Rightarrow x' \in I_0$ \wilog \,
$x \in I_0,z \in J \backslash I_0$.  If $y \in I_0$ then (as $y \le_J
z$) for some $y',y \le_{I_0} y' \in I_{0,z}$ but $x \le_{I_0} y$ (as
$x,y \in I_0,x \le_J y$) hence $x \le_{I_0} y' \in I_{0,z}$ so $x \le_J
z$.  If $y \notin I_0$ then $I_{0,y} \subseteq I_{0,z}$ (by clause
(h)) so we can finish similarly.  So we have covered all cases.]
\sn
\item "{$(*)_2$}"  $J$ is directed
\nl
[Let $x,y \in J$ and we shall find a common upper bound.  If $x,y
\notin I_0$ their concatanation $x \char 94 y$ can serve.  If $x,y
\in I_0$ use ``$I_0$ is directed".  If $x \in I_0,y \in J \backslash
I_0$, then $\langle h(x) \rangle \in J \backslash I_0$ and $z = y
\char 94 \langle h(x) \rangle \in J \backslash I_0$ is $<_J$ above $y$
(by the choice of $\le_J$) and is $\le_J$-above $x$ as $x \in
I_{0,\langle h(x)\rangle} \subseteq I_{0,z}$ by clause (i) of
$\otimes_1$ so we are done.   If $x \in J \backslash
I_0,y \in J_0$ the dual proof works.]
\sn
\item "{$(*)_3$}"  if $x \in J \backslash I_0$ then $M_x \cap M_\ell
\le_{{\frak K}_x} M_x$ for $\ell=0,1$
\nl
[Why?  Clearly $M_x \cap M_0 = \cup\{M^1_t:t \in I_{1,x}\}$ by clause
(c) and $M_x \cap M_1 = \cup\{M^1_t:t \in J_{1,x}\}$ by clause (c),
too.  Now the sets $I_{1,x} \subseteq J_{1,x} (\subseteq J_1)$
are directed by $\le_{J_1}$ so by the assumption on $\langle
M^1_t:t \in J_1 \rangle$ and Lemma \scite{600-0.6} we have $M_x \cap M_0
\le_{{\frak K}_\lambda} M_x \cap M_1$.  Using $J_2$ we can similarly
prove $M_x \cap M_1 \le_{{\frak K}_\lambda} M_x \cap M_2$ and trivially
$M_x \cap M_2 = M_x$.  As $\le_{{\frak K}_\lambda}$ is transitive
we are done.]
\sn
\item "{$(*)_4$}"  if $x \le_J y$ then $M_x \le_{{\frak K}_\lambda} M_y$
\nl
[Why?  If $x,y \in I_0$ use the choice of $\langle M^0_s:s \in I_0
\rangle$.  If $x,y \in J \backslash I_0$ the proof is similar to that
of $(*)_3$ using $J_2$.  If $x \in I_0,y \in J \backslash I_0$ there
is $s \in I_{0,y}$ such that $x \le_{I_0} s$, hence $M_x = M^0_x
\le_{{\frak K}_\lambda} M^0_s$ and as $\langle M^0_t:t \in
I_{0,y}\rangle$ is $\le_{{\frak K}_\lambda}$-directed clearly
$M^0_s \le_{{\frak K}_\lambda} \cup
\{M^0_t:t \in I_{0,y}\} = M_y \cap M_0$ and
$M_y \cap M_0 \le_{{\frak K}_\lambda} M_y$ by $(*)_3$.  By the
transitivity of $\le_{{\frak K}_\lambda}$ we are done.]
\sn
\item "{$(*)_5$}"  $\cup\{M_x:x \in I\} = \cup\{M^0_x:x \in I_0\} = M_0$
\nl
[Why?  Trivially recalling $I_0 = I$.]
\sn  
\item "{$(*)_6$}"  $M_2 = \cup\{M_x:x \in J\}$
\nl
[Why?  Trivially.]
\ermn
As $I_0 \subseteq J$ is directed by $(*)_1 + (*)_2 + (*)_4 + (*)_5 + (*)_6$ we
have checked that $I,\langle M_x:x \in J \rangle$ witness $M_0
\le_{{\frak K}'} M_2$.  This completes the proof of AxII, but we also 
have proved
\mr
\item "{$\otimes_2$}"  if $\bar M = \langle M_t:t \in I \rangle$ is a 
reasonable witness to $M \in K'$ and $M \le_{{\frak K}'} N \in K'$, \ub{then}
there is a witness $I',\bar M' = \langle M'_t:t \in J' \rangle$ to $M
\le_{{\frak K}'} N$ such that $I' = I,\bar M' \restriction I = \bar M$
and $\bar M'$ is reasonable and $x \le_{J'} y \wedge y \in I'
\Rightarrow x \in I'$.]
\ermn
\underbar{Ax III}:  If $\theta$ is a regular cardinal, $M_i$ (for 
$i < \theta)$ is $\le_{{\frak K}'}$-increasing and continuous, 
\underbar{then} $M_0 \le_{{\frak K}'} \dsize \bigcup_{i < \theta} M_i$ 
(in particular $\dsize \bigcup_{i < \theta} M_i \in {\frak K}'$). \newline
[Why?  Let $M_\theta = \dsize \bigcup_{i < \theta} M_i$, \wilog \, $\langle
M_i:i < \theta \rangle$ is not eventually constant and so \wilog \, $i <
\theta \Rightarrow M_i \ne M_{i+1}$ hence $\|M_i\| \ge |i|$; (this
helps below to get ``reasonable", i.e. $|I_\ell| = \|M_i\|$ for limit $i$).
We choose by
induction on $i \le \theta$, a directed partial order $I_i$ and $M_s$ for
$s \in I_i$ such that:
\mr
\item "{$\otimes_3(a)$}"  $\langle M_s:
s \in I_i \rangle$ witness $M_i \in K'$
\smallskip
\noindent
\item "{$(b)$}"  for $j < i,I_j \subseteq I_i$ and ($I_j,I_i,\langle M_s:s
\in I_i \rangle$) witness $M_j \le_{{\frak K}'} M_i$
\sn
\item "{$(c)$}"  $I_i$ is of cardinality $\le \|M_i\|$
\sn
\item "{$(d)$}"  if $I_i \models s \le t$ and $j < i,t \in I_j$ then
$s \in I_j$
\ermn
For $i=0$ use the definition of $M_0 \in K'$. \newline
For $i$ limit let 
$I_i := \dsize \bigcup_{j < i} I_j$ (and the already defined
$M_s$'s) are as required because $M_i = \dsize \bigcup_{j <i} M_j$ and the
induction hypothesis (and $|I_i| \le \|M_i\|$ as we have assumed above that
$j < i \Rightarrow M_j \ne M_{j+1}$) . \newline
For $i=j+1$ use the proof of Ax.II above with $M_j,M_i,M_i,\langle M_s:s \in 
I_j \rangle$ here serving as $M_0,M_1,M_2,\langle M^0_j:s \in I_0 \rangle$
there, that is, we use $\otimes_2$ from there.  
Now for $i = \theta,\langle M_s:s \in I_\theta \rangle$ witness
$M_\theta \in K'$ and $(I_i,I_\theta,\langle M_s:s \in I_\theta \rangle)$
witness $M_i \le_{{\frak K}'} M_\theta$ for each $i < \theta$.]
\bn
\ub{Axiom IV}:  Assume $\theta$ is regular and $\langle M_i:i < \theta
\rangle$ is $\le_{\frak K}$-increasingly continuous, $M \in K'$ and $i <
\theta \Rightarrow M_i \le_{{\frak K}'} M$ and $M_\theta = \dbcu_{i < \theta}
M_i$ (so $M_\theta \subseteq M$).  \ub{Then} $M_\theta \le_{{\frak K}'} M$.
\nl
[Why?  By the proof of Ax.III there are $\langle M_s:s \in I_i \rangle$ for
$i < \theta$ satisfying clauses (a),(b),(c) and (d) of $\otimes_2$ 
there and \wilog \, $I_i \cap \theta = \emptyset$.   For each $i < \theta$
as $M_i \le_{{\frak K}'} M$ there are $J_i$ and $M_s$ for $s \in J_i
\backslash I_i$ such that $(I_i,J_i,\langle M_s:s \in J_i 
\rangle)$ witnesses it; \wilog \, with $\langle \dbcu_{i < \theta} I_i \rangle
\char 94 \langle J_i \backslash I_i:i < \theta
\rangle$ a sequence of pairwise disjoint sets;  exist by $\otimes_2$ above.  
Let $I := \dbcu_{i < \theta} I_i$, let
$\bold i:I \rightarrow \theta$ be $\bold i(s) = \text{ Min}\{i:s \in
I_i\}$ and recall $|I| \le \|M_\theta\|$ hence by clause (d) of
$\otimes_2$ we have $s \le_I t \Rightarrow \bold i(s) \le \bold i(t)$ and
let $h$ be a one-to-one function from $I$ into $M_\theta$.
Without loss of generality the union below is disjoint and let
\mr
\item "{$(*)_7$}"  $J := I \cup 
\bigl\{(A,S):A \text{ a finite subset of } M \text{ and }
S \text{ a finite } \text{ subset of } I 
\text{ with a maximal element} \bigr\}$.
\ermn
ordered by: $J \models x \le y$ \ub{iff} $x,y \in I,I \models x \le y$ or
$x \in I,y = (A,S) \in J \backslash I$ and $x \in S$ or $x = (A^1,S^1) \in J
\backslash I,y = (A^2,S^2) \in J \backslash I,A^1 \subseteq A^2,S^1 \subseteq
S^2$.
\sn
We choose $N_y$ for $y \in J$ as follows:
\sn
If $y \in I$ we let $N_y = M_y$.
\sn
By induction on $n <  \omega$, if 
$y = (A,S) \in J \backslash I$ saisfies $n = |A| + |S|$,
we choose the objects $N_y,I_{y,s},J_{y,s}$ for $s \in S$ such that:
\mr
\item "{$\otimes_3(a)$}"  $I_{y,s}$ 
is a directed subset of $I_{\bold i(s)}$ of
cardinality $\le \lambda$
\sn
\item "{$(b)$}"  $J_{y,s}$ is a directed subset of $J_{\bold i(s)}$ of
cardinality $\le \lambda$
\sn
\item "{$(c)$}"  $s \in I_{\bold i(s)}$ for $s \in S$ (follows from
the definition of $\bold i(s)$)
\sn
\item "{$(d)$}"  $I_{y,s} \subseteq J_{y,s}$ for $s \in S$ and for $s <_I t$ 
from $S$ we have $I_{y,s} \subseteq I_{y,t} \and J_{y,s} \subseteq J_{y,t}$
\sn
\item "{$(e)$}"  if $y_1 = (A_1,S_1) \in J \backslash I,(A_1,S_1) <_J
(A,S)$ and $s \in S_1$ then \nl
$I_{y_1,s} \subseteq I_{y,s},J_{y_1,s} \subseteq J_{y,s}$
\sn
\item "{$(f)$}"  $N_y = \dbcu_{t \in J_{y,s}} M_t$ for any $s \in S$
\sn
\item "{$(g)$}"  $A \subseteq M_t$ for some $t \in J_{y,s}$ for any $s
\in S$, hence $A \subseteq N_y$.
\ermn
No problem to carry the induction and check that $(I,J,\langle N_y:y \in
J \rangle)$ witness $M_\theta \le_{{\frak K}'} M$.
\bn
\ub{Axiom V}:  Assume $N_0 \le_{{\frak K}'} M$ and $N_1 \le_{{\frak K}'} M$.
\nl
If $N_0 \subseteq N_1$, then $N_0 \le_{{\frak K}'} N_1$. \nl
[Why?  Let $(I_0,J_0,\langle M^0_s:s \in J_0 \rangle)$ witness $N_0
\le_{{\frak K}'} M$.  Let $\langle M^1_s:s \in I_1 \rangle$ witness $N_1
\in {\frak K}'$ and  \wilog \, $I_1$ is isomorphic to
$([N_1]^{<\aleph_0},\subseteq)$ and let $h_1$ be an isomorphism from $I_1$
onto $([N_1]^{< \aleph_0},\subseteq)$.  
Now by induction on $n$, for $s \in I_1$ satisfying $n = 
|\{t:t <_{I_1} s\}|$ we choose directed
subsets $F_0(s),F_1(s)$ of $I_0,I_1$ respectively, each of cardinality
$\le \lambda$ such that:
\mr
\item "{$(i)$}"  $h(s) \in F_1(s)$ and 
$t <_{I_1} s \Rightarrow F_0(t) \subseteq F_0(s) \and
F_1(t) \subseteq F_1(s)$
\sn
\item "{$(ii)$}"   $\bigcup\{M^0_t:t \in F_0(s)\} = \bigcup\{M^1_t:t
\in F_1(s)\} \cap N_0$.
\ermn
Now replacing $M^1_s$ by $\cup\{M^1_t:t \in F_1(s)\}$ and letting
$F=F_0$ we get: 
\mr
\item "{$(iii)$}"  $t \in I_1 \wedge s \in F(t)(\subseteq I_0)
\Rightarrow M^0_s \subseteq M^1_t$
\sn
\item "{$(iv)$}"  $F$ is a function from $I_1$ to $[I_0]^{\le \lambda}$
\sn
\item "{$(v)$}"   for $s \in I_1,F(s)$ is a directed subset of
$I_0$ of cardinality 
$\le \lambda,M^1_s \cap N_0 = \cup\{M^0_t:t \in F(s)\}$
and $I_1 \models s \le t \Rightarrow F(s) \subseteq F(t)$. 
\ermn
As $N_1 \le_{\frak K} M$ by the 
proof of Ax.II, i.e., by $\otimes_2$ above we can 
find $J_1$ and $M^1_s$ 
for $s \in J_1 \backslash I_1$ such that
$(I_1,J_1,\langle M^1_s:s \in J_1 \rangle)$ witnesses $N_1 \le_{{\frak K}'}
M$.  We now prove
\mr
\item "{$\boxtimes_4$}"  if $r \in I_1,s \in I_0$ and $s \in F(r)$ then
$M^0_s \le_{{\frak K}_\lambda} M^1_r$. 
\ermn
[Why?  As $\langle M^0_s:s \in J_0 \rangle,\langle M^1_s:s \in J_1
\rangle$ are both witnesses for $M \in K'$, clearly for 
$r \in I_1$ we can find $J'_0(r) \subseteq J_0$ 
directed of cardinality $\le \lambda$ and $J'_1(r) \subseteq J_1$ directed 
of cardinality $\le \lambda$ such that $r \in J'_1(r),F(r) \subseteq 
J'_0(r)$ and $\dbcu_{t \in J'_0(r)} M^0_t = 
\dbcu_{t \in J'_1(r)} M^1_t$, call it $M^*_r$.
\nl
Now $M^*_r \in K'_\lambda = 
K_\lambda$ (by part (2) and \scite{600-0.6}) and $t \in J'_1(r) \Rightarrow
M^1_t \le_{{\frak K}_\lambda} M^*_r$ (as ${\frak K}_\lambda$ is 
a $\lambda$-abstract elementary class applying the parallel to
observation \scite{600-0.6}, i.e., \scite{600-0.30D}(2)) and 
similarly $t \in J'_0(r) \Rightarrow M^0_t 
\le_{{\frak K}_\lambda} M^*_r$.  Now the $s$ from $\boxtimes_4$ satisfied
$s \in F(r) \subseteq J'_0(r)$ hence 
$M^0_s \subseteq M^1_r$ (why?  by clause (iii) above $s \in F(r)$ is
as required in $\boxtimes_4$).  But above we got    
$M^0_s \le_{\frak K} M^*_r,M^1_r \le_{\frak K} M^*_r$, so 
by Ax.V for ${\frak K}_\lambda$ we have $M^0_s \le_{\frak K} M^1_r$ as
required in $\boxtimes_4$.]  

Without loss of generality $I_0 \cap I_1 = \emptyset$ and define the
partial order $J$ with set of elements 
$I_0 \cup I_1$ by $J \models x \le y$ iff $x,y \in I_0,I_0 \models
x \le y$ or 
$x \in I_0,y \in I_1$ and $x \in F(y)$ or $x,y \in I_1,I_1 \models
x \le y$.  Recalling clause (c) above, it is a partial order.
Define $M_s$ for $s \in J$ as $M^0_s$ if $s \in I_0$ and as
$M^1_s$ if 
$s \in I_1$.  Now check that $(I_0,J,\langle M_s:s \in J \rangle)$
witnesses $N_0 \le_{{\frak K}'} N_1$ as required.]
\bn
\ub{Axiom VI}:  LS$({\frak K}') = \lambda$. \nl
[Why?  Let $M \in K',A \subseteq M,|A| + \lambda \le \mu < \|M\|$ and let
$\langle M_s:s \in J \rangle$ witness $M \in K'$.  As $\|M\| > \mu$ we
can choose a directed $I \subseteq J$ of cardinality $\le \mu$ 
such that $A \subseteq M' := \dbcu_{s \in I} M_s$ and so $(I,J,
\langle M_s:s \in J \rangle)$
witnesses $M' \le_{{\frak K}'} M$, so as $A \subseteq M'$ and 
$\|M'\| \le |A| + \mu$
we are done.]  \hfill$\square_{\scite{600-0.31}}$
\bn
We may like to use ${\frak K}_{\le \lambda}$ instead of 
${\frak K}_\lambda$; no need as essentially ${\frak K}$ consists of
two parts ${\frak K}_{\le \lambda}$ and ${\frak K}_{\ge \lambda}$
which have just to agree in $\lambda$.  That is
\proclaim{\stag{600-0.32} Claim}  Assume
\mr
\item "{$(a)$}"  ${\frak K}^1$ is an abstract elementary class with
$\lambda = {\text{\rm LS\/}}({\frak K}^1),K^1 = K^1_{\ge \lambda}$
\sn
\item "{$(b)$}"  ${\frak K}^2_{\le \lambda}$ is a $(\le \lambda)$-abstract 
elementary class (defined as in \scite{600-0.30A}(1) with the obvious changes so
$M \in {\frak K}^2_{\le \lambda} \Rightarrow \|M\| \le \lambda$ and 
in Axiom III, $\|\dbcu_i M_i\| \le \lambda$ is required)
\sn
\item "{$(c)$}"  $K^2_\lambda = K^1_\lambda$ and $\le_{{\frak K}^2}
\restriction K^2_\lambda = \le_{{\frak K}^1} \restriction K^1_\lambda$
\sn
\item "{$(d)$}"  we define ${\frak K}$ as follows: 
$K = K^1 \cup K^2,M \le_{\frak K} N$
iff $M \le_{{\frak K}^1} N$ or $M \le_{{\frak K}^2} N$ or for some $M',M 
\le_{{\frak K}^2} M' \le_{{\frak K}^1} N$.
\ermn
\ub{Then} ${\frak K}$ is an abstract elementary class and ${\text{\rm
LS\/}}({\frak K}) = {\text{\rm LS\/}}({\frak K}^2)$ which trivially is
$\le \lambda$.
\endproclaim
\bigskip

\demo{Proof}  Straight.  E.g. \nl
\ub{Axiom V}:  We shall use freely
\mr
\item "{$(*)$}"  ${\frak K}_{\le \lambda} = {\frak K}^2$ and ${\frak K}
_{\ge \lambda} = {\frak K}^1$.
\ermn
So assume $N_0 \le_{\frak K} M,N_1 \le_{\frak K} M,N_0
\subseteq N_1$. \nl
Now if $\|N_0\| \ge \lambda$ use assumption (a), so we can assume
$\|N_0\| < \lambda$. If $\|M\| \le \lambda$ we can use assumption (b) so we
can assume $\|M\| > \lambda$ and by the definition of $\le_{\frak K}$ there
is $M'_0 \in K^1_\lambda = K^2_\lambda$ such that $N_0 \le_{{\frak K}^2} M'_0
\le_{{\frak K}^1} M$.  First assume 
$\|N_1\| \le \lambda$, so we can find $M'_1 \in K^1_\lambda$ such that
$N_1 \le_{{\frak K}^2} M'_1 \le_{{\frak K}^1} M$ 
(why?  if $N_1 \in 
K_{< \lambda}$, by the definition of $\le_{\frak K}$ and if
$N_1 \in K_\lambda$ just choose $M'_1 = N_1$).  Now we can by assumption (a) 
find $M'' \in K^1_\lambda$ such that $M'_0 \cup M'_1 \subseteq M'' 
\le_{{\frak K}^1} M$, hence by assumption (a) (i.e. AxV for ${\frak K}^1$)
we have $M'_0 \le_{{\frak K}^1} M'',M'_1 \le_{{\frak K}^1} M''$, so by
assumption (c) we have $M'_0 \le_{{\frak K}^2} M'',M'_1 
\le_{{\frak K}^2} M''$.  As $N_0 \le_{{\frak K}^2} M'_0 
\le_{{\frak K}^2} M'' \in K_{\le \lambda}$ 
by assumption (b) we have $N_0
\le_{{\frak K}^2} M''$, and similarly we have $N_1 \le_{{\frak K}^2}
M''$.   So $N_0 \subseteq N_1,N_0 \le_{{\frak K}^2} M'',
N_1 \le_{{\frak K}^2} M'$ so by 
assumption (b) we have $N_0 \le_{{\frak K}^2} N_1$ hence
$N_0 \le_{\frak K} N_1$.

We are left with the case $\|N_1\| > \lambda$; by assumption (a) there is
$N'_1 \in K_\lambda$ such that $N_0 \subseteq N'_1 \le_{{\frak K}^1} N_1$.  By
assumption (a) we have $N'_1 \le_{{\frak K}^1} M$, so by the previous paragraph
we get $N_0 \le_{{\frak K}^2} N'_1$, 
together with the previous sentence we have
$N_0 \le_{{\frak K}^2} N'_1 \le_{{\frak K}^1} N_1$ so by the definition of
$\le_{\frak K}$ we are done.  \hfill$\square_{\scite{600-0.32}}$
\enddemo
\bn
Recall
\definition{\stag{600-0.33} Definition}  If $M \in K_\lambda$ is locally
superlimit or just pseudo superlimit
let $K_{[M]} = K^{[M]}_\lambda = \{N \in K_\lambda:N \cong M\},{\frak
K}_{[M]} = {\frak K}^{[M]}_\lambda =
(K_{[M]},\le_{\frak K} \restriction K^{[M]}_\lambda)$ and let 
${\frak K}^{[M]}$ be the 
${\frak K}'$ we get in \scite{600-0.31}(1) for ${\frak K} =
{\frak K}_{[M]} = {\frak K}^{[M]}_\lambda$.  
We may write ${\frak K}_\lambda[M],{\frak K}[M]$.
\enddefinition
\bn
Trivially
\proclaim{\stag{600-0.34} Claim}  1) If ${\frak K}$ is an $\lambda$-{\rm a.e.c.}, 
$M \in K_\lambda$ is locally superlimit or just pseudo superlimit
\ub{then} ${\frak K}_{[M]}$ 
is a $\lambda$-{\rm a.e.c.} which
is categorical (i.e. categorical in $\lambda$).
\nl
2) Assume ${\frak K}$ is an a.e.c. and $M \in {\frak K}_\lambda$ is
not $\le_{\frak K}$- maximal.  $M$ is pseudo superlimit (in ${\frak
K}$, i.e., in ${\frak K}_\lambda$) \ub{iff} ${\frak K}_{[M]}$
is a $\lambda$-a.e.c. which is categorical \ub{iff} ${\frak K}^{[M]}$
is an a.e.c., categorical in $\lambda$ and $\le_{{\frak K}^{[M]}} =
\le_{\frak K} \restriction K^{[M]}$.
\nl
3) In (1) and (2), {\rm LS}$({\frak K}^{[M]}) = \lambda = \text{\rm Min}
\{\|N\|:N \in {\frak K}^{[M]}\}$.
\endproclaim
\bigskip

\demo{Proof}  Straightforward.  \hfill$\square_{\scite{600-0.34}}$
\enddemo
\bn
\margintag{600-0.37.3}\ub{\stag{600-0.37.3} Exercise}:  Assume ${\frak K}$ is a
$\lambda$-a.e.c. with amalgamation and stability in $\lambda$.
\ub{Then} for every $M_1 \in K_\lambda,p_1 \in {\Cal S}_{\frak K}(M_1)$ 
we can find $M_2 \in K$ and minimal $p_2 \in {\Cal S}_{\frak
K}(M_2)$ such that $M_1 \le_{\frak K} M_2,p_1 = p_2 \restriction M_1$.
\bn
\margintag{600-1c.29}\ub{\stag{600-1c.29} Exercise}:  1) Any $\le_{{\frak K}_\lambda}$-embedding $f_0$ of
$M^1_0$ into $M^2_0$ can be extended to an isomorphism $f$ from
$M^1_\delta$ onto $M^2_\delta$ such that $f(M^1_{2 \alpha})
\le_{{\frak K}_\lambda} M^2_{2 \alpha},f^{-1}(M^2_{2 \alpha +1})
\le_{{\frak K}_\lambda} M^1_{2 \alpha +1}$ for every $\alpha <
\delta$, \ub{provided that}
\mr
\item "{$\circledast$}"  $(a) \quad {\frak K}_\lambda$ is a
$\lambda$-a.e.c. and $\delta$ is a limit ordinal $\le \lambda^+$
\sn
\item "{${{}}$}"  $(b) \quad \langle M^\ell_\alpha:\alpha \le
\delta\rangle$ is $\le_{{\frak K}_\lambda}$-increasing continuous for
$\ell=1,2$
\sn
\item "{${{}}$}"  $(c) \quad M^\ell_\alpha$ is an amalgamation base in
${\frak K}_\lambda$ (for $\alpha < \delta$ and $\ell=1,2$)
\sn
\item "{${{}}$}"  $(d) \quad M^\ell_{\alpha +1}$ is $\le_{{\frak
K}_\lambda}$-universal extension of $M^\ell_\alpha$ for $\alpha <
\delta,\ell=1,2$.
\ermn
2) Write the axioms of ``a $\lambda$-a.e.c." which are used.
[Hint: Should be clear, and the argument appear.] 
\newpage

\head {\S2 Good Frames} \endhead  \resetall \sectno=2
 \spuriousreset
\bigskip

We first present our central definition: good $\lambda$-frame (in
Definition \scite{600-1.1}).  We are given the relation 
``$p \in {\Cal S}(N)$ does not fork over $M \le_{\frak K} N$ when $p$
is basic" (by the basic relations
and axioms) so it is natural to look at how well we can ``lift" the
definition of non-forking to models of cardinality $\lambda$ and later
to non-forking of models (and types over them) in
cardinalities $> \lambda$.  Unlike the lifting of $\lambda$-a.e.c. in
Lemma \scite{600-0.31}, life is not so easy.  
We define in \scite{600-1.6}, \scite{600-1.7},
\scite{600-1.9} and we prove basic properties in \scite{600-1.8},
\scite{600-1.10}, \scite{600-1.12} and less obvious ones in \scite{600-1.11},
\scite{600-1.13}, \scite{600-1.13B}.  This should serve as a reasonable exercise in the
meaning of good frames; however, the lifting, in general, does not
give good $\mu$-frames for $\mu > \lambda$.  There may be no $M \in
K_\mu$ at all and/or amalgamation may fail.  Also the existence and
uniqueness of non-forking types is problematic.  We do not give up and
will return to the lifting problem, under additional assumptions in
\sectioncite[\S12]{705} and \cite{Sh:842}. 

In \scite{600-1.14} (recalling \scite{600-0.34}) we show that the case
``${\frak K}^{\frak s}$ categorical in $\lambda$" is not so rare among
good $\lambda$-frames; in fact if there is a superlimit model in
$\lambda$ we can restrict ${\frak K}_\lambda$ to it.  So in a sense
superstability and categoricity are close.  
For elementary classes they are less close and note that if $T$ is a
complete first order superstable theory and $\lambda \ge 2^{|T|}$,
\ub{then} the class ${\frak K} = {\frak K}_{T,\lambda}$ of
$\lambda$-saturated model of $T$ is in general not an elementary class
(though is a PC$_\lambda$ class) but is an a.e.c. categorical in
$\lambda$ and for some good $\lambda$-frame ${\frak s},K_{\frak s} =
{\frak K}_{T,\lambda}$.  How justified is our restriction here to
something like ``the $\lambda$-saturated model"?  It is O.K. for our
test problems but more so it is justified as
our approach is to first analyze the quite saturated models.

Last but not least in \scite{600-1.15} we show that 
one of the axioms from \scite{600-1.1},
i.e., (E)(i), follows from the rest in our present definition;
additional implications are in Claims \scite{600-1.16A}, \scite{600-1.16B}.  Later
``Ax(X)(y)" will mean (X)(y) from Definition \scite{600-1.1}. 
\bn
Recall that good $\lambda$-frame is intended to be a parallel to (bare
bones) superstable
elementary class stable in $\lambda$;  here we restrict
ourselves to models of cardinality $\lambda$.
\definition{\stag{600-1.1} Definition}  We say ${\frak s} = 
({\frak K},\nonfork{}{}_\lambda,
{\Cal S}^{\text{bs}}_\lambda) = ({\frak K}^{\frak s},\nonfork{}{}_{\frak s},
{\Cal S}^{\text{bs}}_{\frak s})$ is a good frame in $\lambda$ or a 
good $\lambda$-frame ($\lambda$ may be omitted when clear, 
note that $\lambda = \lambda_{\frak s} =
\lambda({\frak s})$ is determined by ${\frak s}$ and we may write
${\Cal S}_{\frak s}(M)$ instead of ${\Cal S}_{{\frak K}^{\frak s}}(M)$
and \ortp$_{\frak s}(a,M,N)$ instead of \ortp$_{{\frak K}^{\frak s}}(a,M,N)$
when $M \in K^{\frak s}_\lambda,N \in K^{\frak s}$; we may write \ortp$(a,M,N)$
for \ortp$_{{\frak K}^{\frak s}}(a,M,N)$) 
\ub{when} the following conditions hold:
\mr
\widestnumber\item{$(D)(a)$}
\item "{$(A)$}"  ${\frak K} = (K,\le_{\frak K})$ is an abstract elementary
class also denoted by ${\frak K}[{\frak s}]$, 
the 
L\"owenheim Skolem number of ${\frak K}$, being $\le \lambda$ 
(see Definition \scite{600-0.2}); 
there is no harm in assuming $M \in K \Rightarrow \|M\| \ge \lambda$;
let ${\frak K}_{\frak s} = {\frak K}^{\frak s}_\lambda$ and
$\le_{\frak s} = \le_{\frak K} \restriction K_\lambda$, and let 
${\frak K}_{\frak s} = (K_\lambda,\le_{\frak s})$ and ${\frak K}[{\frak s}] =
{\frak K}^{\frak s}$ so we may write ${\frak s} = ({\frak K}_{\frak
s},\nonfork{}{}_{\frak s},{\Cal S}^{\text{bs}}_{\frak s})$
\sn
\item "{$(B)$}"  ${\frak K}$ has a superlimit model in $\lambda$ which
\footnote{in fact, the ``is not $<_{\frak K}$-maximal" follows by (C)} is
not $<_{\frak K}$-maximal. 
\sn
\item "{$(C)$}"  ${\frak K}_\lambda$ has the amalgamation property, the
JEP (joint embedding property), and has no $\le_{\frak K}$-maximal 
member.
\sn
\item "{$(D)(a)$}"  ${\Cal S}^{\text{bs}} = {\Cal S}^{\text{bs}}_\lambda$ 
(the class of basic types for ${\frak K}_\lambda$) is included in \newline
$\bigcup\{{\Cal S}(M):M \in K_\lambda\}$ and is closed under isomorphisms
including automorphisms; 
for $M \in K_\lambda$ let ${\Cal S}^{\text{bs}}(M) = 
{\Cal S}^{\text{bs}} \cap {\Cal S}(M)$; no harm in 
allowing types of finite sequences, i.e., replacing ${\Cal S}(M)$ by
${\Cal S}^{< \omega}(M)$, (${\Cal S}^\omega(M))$ is different as being
new (= non-algebraic) is not preserved under increasing unions).
\sn
\item "{${}(b)$}"   if $p \in {\Cal S}^{\text{bs}}(M)$, \ub{then} $p$ is
non-algebraic (i.e. not realized by any $a \in M$).
\sn
\item "{${}(c)$}"  \underbar{(density)} \newline
if $M \le_{\frak K} N$ are from $K_\lambda$ and $M \ne N$, \underbar{then}
for some $a \in N \backslash M$ we have $\text{\ortp}(a,M,N) \in 
{\Cal S}^{\text{bs}}$
\newline
\beginaside [intention: examples are: 
minimal types in \cite{Sh:576}, regular types for superstable first
order (= elementary) classes].
\sn
\endaside 
\item "{${}(d)$}"  \ub{bs-stability} \nl
${\Cal S}^{\text{bs}}(M)$ has cardinality $\le \lambda$ for $M \in
K_\lambda$.
\sn
\item "{$(E)(a)$}"  $\nonfork{}{}_{\lambda}$ denoted also by
$\nonfork{}{}_{\frak s}$ or just $\nonfork{}{}_{}$, is a four 
place relation
\footnote{we tend to forget to write the $\lambda$, this is justified by
\scite{600-1.8}(2), and see Definition \scite{600-1.7}} called non-forking with 
$\nonfork{}{}_{}(M_0,M_1,a,M_3)$
implying $M_0 \le_{\frak K} M_1 \le_{\frak K} M_3$ are from $K_\lambda,
a \in M_3 \backslash M_1$ and $\text{\ortp}(a,M_0,M_3) \in 
{\Cal S}^{\text{bs}}(M_0)$ and \newline
$\text{\ortp}(a,M_1,M_3) \in {\Cal S}^{\text{bs}}(M_1)$.  
Also $\nonfork{}{}_{}$ is preserved under isomorphisms and we
demand: if $M_0 = M_1 
\le_{\frak K} M_3$ both in $K_\lambda$ and $a \in M_3$, then: \newline
$\nonfork{}{}_{}(M_0,M_1,a,M_3)$ is equivalent to 
``\ortp$(a,M_0,M_3) \in {\Cal S}^{\text{bs} }(M_0)$".  The assertion
$\nonfork{}{}_{}(M_0,M_1,a,M_3)$ is also written as  
$\nonforkin{M_1}{a}_{M_0}^{M_3}$  and also 
as ``\ortp$(a,M_1,M_3)$ does not fork 
over $M_0$ (inside $M_3$)" (this is justified by clause (b) below).
So \ortp$(a,M_1,M_3)$ forks over $M_0$ (where $M_0 \le_{\frak s} M_1
\le_{\frak s} M_3,a \in M_3$) is just the negation
\nl
\beginaside
[Explanation: The 
intention is to axiomatize non-forking of types, but we already
commit ourselves to dealing with basic types only.  Note that in
\cite{Sh:576} we know something on minimal types but other types are
something else.]
\sn
\endaside
\item "{${}(b)$}"  \underbar{(monotonicity)}: \newline
if $M_0 \le_{\frak K} M'_0 \le_{\frak K} M'_1 \le_{\frak K} 
M_1 \le_{\frak K} M_3 \le_{\frak K} M'_3,M_1 \cup \{a\} 
\subseteq M''_3 \le_{\frak K} M'_3$ 
all of them in $K_\lambda$, \ub{then}
$\nonfork{}{}_{}(M_0,M_1,a,M_3) \Rightarrow 
\nonfork{}{}_{}(M'_0,M'_1,a,M'_3)$ and $\nonfork{}{}_{}(M'_0,M'_1,a,M'_3)
\Rightarrow \nonfork{}{}_{}
(M'_0,M'_1,a,M''_3)$, \ub{so} it is legitimate to just 
say ``$\text{\ortp}(a,M_1,M_3)$ 
does not fork over $M_0$". \nl
\beginaside
[Explanation: non-forking is preserved by decreasing the type,
increasing the basis (= the set over which it does not fork) and
increasing or decreasing 
the model inside which all this occurs.  The same holds for
stable theories only here we restrict ourselves to ``legitimate",
i.e., basic types.  But note that here the ``restriction of
\ortp$(a,M_1,M_3)$ to $M'_1$ is basic" is worthwhile information.]
\endaside
\sn
\item "{${}(c)$}"  \underbar{(local character)}: \nl
if $\langle M_i:i \le \delta + 1 \rangle$ is $\le_{\frak K}$-increasing
continuous in ${\frak K}_\lambda,a \in M_{\delta + 1}$ and \newline
\ortp$(a,M_\delta,M_{\delta +1}) \in {\Cal S}^{\text{bs}}
(M_\delta)$ \underbar{then} 
for every $i < \delta$ large enough \ortp$(a,M_\delta,M_{\delta +1})$ 
does not fork over $M_i$. \nl
\beginaside
[Explanation: This is a replacement for superstability which says
that: if $p \in {\Cal S}(A)$ then there is a finite 
$B \subseteq A$ such that $p$ does not fork over $B$.]
\endaside
\sn
\item "{${}(d)$}"  \underbar{(transitivity)}: \newline
if $M_0 \le_{\frak s} M'_0 \le_{\frak s} M''_0 \le_{\frak s} M_3$ are
from $K_\lambda$ and
$a \in M_3$ and \ortp$(a,M''_0,M_3)$ does not fork over $M'_0$ and
\ortp$(a,M'_0,M_3)$ does not fork over $M_0$ (all models are in $K_\lambda$, of 
course, and necessarily the three relevant types are in 
${\Cal S}^{\text{bs}}$), \ub{then} 
\ortp$(a,M''_0,M_3)$ does not fork over $M_0$
\smallskip
\noindent
\item "{${}(e)$}"  \underbar{uniqueness}: \newline
if $p,q \in {\Cal S}^{\text{bs}}(M_1)$ do not fork over 
$M_0 \le_{\frak K} M_1$ (all in $K_\lambda$) and \newline
$p \restriction M_0 = q \restriction M_0$ \underbar{then} $p = q$
\sn
\item "{${}(f)$}"  \underbar{symmetry}:  \newline
if $M_0 \le_{\frak K} M_3$ are in ${\frak K}_\lambda$ and for $\ell = 1,2$
we have \newline
$a_\ell \in M_3$ and $\text{\ortp}(a_\ell,M_0,M_3) \in 
{\Cal S}^{\text{bs}}(M_0)$, \ub{then} the following are equivalent:
\smallskip
\noindent
{\roster
\itemitem{ $(\alpha)$ }  there are $M_1,M'_3$ in $K_\lambda$ such that
$M_0 \le_{\frak K} M_1 \le_{\frak K} M'_3$, \newline
$a_1 \in M_1,M_3 \le_{\frak K} M'_3$ and 
$\text{\ortp}(a_2,M_1,M'_3)$ does not fork over $M_0$
\smallskip
\noindent
\itemitem{ $(\beta)$ }  there are $M_2,M'_3$ in $K_\lambda$ such that
$M_0 \le_{\frak K} M_2 \le_{\frak K} M'_3$, \newline
$a_2 \in M_2,M_3 \le_{\frak K} M'_3$ and $\text{\ortp}(a_1,M_2,M'_3)$ 
does not fork over $M_0$. 
\endroster}
\beginaside
[Explanation: this is a replacement to ``\ortp$(a_1,M_0 \cup
\{a_2\},M_3)$ forks over $M_0$ iff \ortp$(a_2,M_0 \cup \{a_1\},M_3)$
forks over $M_0$" which is not well defined in our context.]
\endaside
\sn
\item "{${}(g)$}"  \underbar{extension existence}: \newline
if $M \le_{\frak K} N$ are from $K_\lambda$ and $p \in 
{\Cal S}^{\text{bs}}(M)$ \ub{then} some 
$q \in {\Cal S}^{\text{bs}}(N)$ does not fork over $M$ and extends
$p$
\sn
\item "{${}(h)$}"  \ub{continuity}: \newline
if $\langle M_i:i \le \delta \rangle$ is $\le_{\frak K}$-increasing
continuous, all in $K_\lambda$ (recall $\delta$ is always a limit
ordinal), $p \in {\Cal S}(M_\delta)$ and $i < \delta
\Rightarrow p \restriction M_i \in {\Cal S}^{\text{bs}}(M_i)$ 
does not fork over
$M_0$ \underbar{then} $p \in {\Cal S}^{\text{bs}}(M_\delta)$ and moreover $p$
does not fork over $M_0$. \nl
\beginaside
[Explanation: This is a replacement to: for an increasing sequence of types
which do not fork over $A$, the union does not fork over $A$; 
equivalently if $p$ forks over $A$ then some finite subtype does.]
\endaside
\sn
\item "{${}(i)$}"  \ub{non-forking amalgamation}: \newline
if for $\ell = 1,2,M_0 \le_{\frak K} M_\ell$ are from $K_\lambda,a_\ell \in
M_\ell \backslash M_0$, \ortp$(a_\ell,M_0,M_\ell) \in {\Cal S}^{bs}(M_0)$,
\ub{then} we can find $f_1,f_2,M_3$ satisfying $M_0 \le_{\frak K} M_3 \in
K_\lambda$ such that for $\ell =1,2$ we have $f_\ell$ is a 
$\le_{\frak K}$-embedding of $M_\ell$ into $M_3$ over
$M_0$ and \ortp$(f_\ell(a_\ell),f_{3-\ell}(M_{3-\ell}),M_3)$ 
does not fork over $M_0$ for $\ell=1,2$.
\beginaside
[Explanation: This strengthens clause (g), (existence) saying we can
do it twice so close to (f), symmetry, see \scite{600-1.15}.]
\endaside
\endroster
\enddefinition
\bigskip
\centerline {$* \qquad * \qquad *$}
\bn
\margintag{600-1.1B}\ub{\stag{600-1.1B} Discussion}:  0) On connections between the axioms see
\scite{600-1.15}, \scite{600-1.16A}, \scite{600-1.16B}.
\nl
1) What 
justifies the choice of the good $\lambda$-frame as a parallel to
(bare bones) superstability?  Mostly starting from assumptions 
on few models around $\lambda$ in the a.e.c. ${\frak K}$ and reasonable,
``semi ZFC" set theoretic assumptions (e.g. involving categoricity and
weak cases of G.C.H., see \S3) 
we can prove that, essentially, for some $\nonfork{}{}_{},
{\Cal S}$ the demands in Definition \scite{600-1.1} hold.  
So here we get (i.e., applying our general theorem to the case of
\scite{600-Ex.1}) an alternative proof of the main theorem of \cite{Sh:87a},
\cite{Sh:87b} in a local version, i.e., dealing 
with few cardinals rather than having to deal with all the cardinals
$\lambda,\lambda^{+1},\lambda^{+2},\dotsc,\lambda^{+n}$ as in
\cite{Sh:87a}, \cite{Sh:87b} in an inductive proof.
That is, in \cite{Sh:87b}, we get dichotomies by the omitting type
theorem for countable models (and theories).  So problems on
$\aleph_n$ are ``translated" down to $\aleph_{n-1}$ (increasing the
complexity) till we arrive to $\aleph_0$ and then ``translated" back.
Hence it is important there to deal with $\aleph_0,\dotsc,\aleph_n$
together.  Here our $\lambda$ may not have special helpful properties,
so if we succeed to prove the relevant claims then they apply to
$\lambda^+$, too.  There are advantages to being poor.
\nl
2) Of course, we may just point out that the axioms seem reasonable
and that eventually we can say much more.
\nl
3) We may consider weakening bs-stability 
(i.e., Ax$(D)(d)$ in Definition \scite{600-1.1}) to
$M \in K_\lambda \Rightarrow |{\Cal S}^{\text{bs}}(M)| \le
\lambda^+$, we just have not looked into it here; Jarden-Shelah
\cite{JrSh:875} will; actually \chaptercite{88r} deals in a limited way
with this in a more restricted framework.  \nl
4) On stability in $\lambda$ and existence of $(\lambda,\sigma)$-brimmed
extensions see \scite{600-4a.1}. \nl
\bn
\relax From the rest of this section we shall use mainly the defintion of
$K^{3,\text{bs}}_\lambda$ in Definition \scite{600-1.6}(3), also
\scite{600-1.16} (restricting ourselves to a superlimit).  We sometimes
use implications among the axioms (in \scite{600-1.15} - \scite{600-1.16B}).
The rest is, for now  
an exercise to familiarize the reader with $\lambda$-frames, in
particular (\scite{600-1.5}-\scite{600-1.14}) to see what occurs to 
non-forking and basic types in cardinals $> \lambda$.  This is 
easy (but see below).  
For this we first present the basic definitions.  
\medskip

\demo{\stag{600-1.5} Convention}  1) We fix ${\frak s}$, a good
$\lambda$-frame so $K = K^{\frak s},{\Cal S}^{\text{bs}} = {\Cal
S}^{\text{bs}}_{\frak s}$.
\nl
2) By $M \in K$ we mean $M \in K_{\ge \lambda}$ if not said
otherwise. 
\enddemo
\bigskip

We lift the properties to ${\frak K}_{\ge \lambda}$ by reflecting to the
situation in $K_\lambda$.  But do not be too excited: the good properties
do not lift automatically, we shall be working on that later (under
additional assumptions).  Of course, from the definition below later
we shall use mainly $K^{3,\text{bs}}_{\frak s} = K^{3,\text{bs}}_\lambda$.
\definition{\stag{600-1.6} Definition}   1) 

$$
\align
K^{3,\text{bs}} = K^{3,\text{bs}}_{\ge {\frak s}} :=
\biggl\{ (M,N,a):&M \le_{\frak K} N,a \in N \backslash M 
\text{ and there is } M' \le_{\frak K} M \\
  &\text{satisfying } M' \in K_\lambda,\text{ such that for every } 
M'' \in K_\lambda\text{ we have:}\\
  &[M' \le_{\frak K} M'' \le_{\frak K} M \Rightarrow \text{ \ortp}
(a,M'',N) \in {\Cal S}^{\text{bs}}(M'') \\
  &\text{does not fork over } M']; \text{ equivalently } 
[M' \le_{\frak K} M'' \le_{\frak K} M \\
  &\and M'' \le_{\frak K} N'' \le_{\frak K} N \and N'' \in K_\lambda
\and a \in N'' \\
  &\Rightarrow \nonfork{}{}_{\lambda}(M',M'',a,N'')] \biggr\}.
\endalign
$$
\mn
2) $K^{3,\text{bs}}_{= \mu} = K^{3,\text{bs}}_{{\frak s},\mu} :=
\{(M,N,a) \in K^{3,\text{bs}}_{\ge {\frak s}}:M,N \in {\frak K}^{\frak
s}_\mu\}$. 
\nl
3) $K^{3,\text{bs}}_{\frak s} := K^{3,\text{bs}}_{= \lambda,{\frak
s}}$; and let $K^{3,\text{bs}}_\mu = K^{3,\text{bs}}_{= \mu}$, used
mainly for $\mu = \lambda_{\frak s}$ 
and $K^{3,\text{bs}}_{{\frak s},\ge \mu}$ is defined naturally.
\enddefinition
\bigskip

\definition{\stag{600-1.7} Definition}  We define 
$\nonfork{}{}_{< \infty}(M_0,M_1,a,M_3)$ 
(rather than $\nonfork{}{}_{\lambda}$) as follows: it holds \ub{iff}
$M_0 \le_{\frak K} M_1 \le_{\frak K}
M_3$ are from $K$ (not necessarily $K_\lambda$), $a \in M_3 \backslash M_1$ 
and there is $M'_0 \le_{\frak K} M_0$ which belongs to $K_\lambda$ 
satisfying: if
$M'_0 \le_{\frak K} M'_1 \le_{\frak K} M_1,M'_1 \in K_\lambda$, \newline
$M'_1 \cup \{a\}
\subseteq M'_3 \le_{\frak K} M_3$ and $M'_3 \in K_\lambda$ \underbar{then}
$\nonfork{}{}_{\lambda}(M'_0,M'_1,a,M'_3)$. \nl
We now check that $\nonfork{}{}_{< \infty}$ behaves correctly when
restricted to $K_\lambda$.
\enddefinition
\bigskip

\proclaim{\stag{600-1.8} Claim}  1) Assume $M \le_{\frak K} N$ are from
$K_\lambda$ and $a \in N$.  \underbar{Then} $(M,N,a) \in 
K^{3,\text{bs}}_{\frak s}$ 
\ub{iff} ${\text{\rm \ortp\/}}(a,M,N) \in {\Cal S}^{\text{bs}}_{\frak
s}(M)$. 
 \newline
2) Assume $M_0,M_1,M_3 \in K_\lambda$ and $a \in M_3$.  \underbar{Then}
$\nonfork{}{}_{< \infty}(M_0,M_1,a,M_3)$ \underbar{iff} \nl
$\nonfork{}{}_{\lambda}(M_0,M_1,a,M_3)$. \newline
3) Assume $M \le_{\frak K} N_1 \le_{\frak K} N_2$ and $a \in N_1$.  
\underbar{Then} \newline
$(M,N_1,a) \in K^{3,\text{bs}}_{\ge {\frak s}} \Leftrightarrow 
(M,N_2,a) \in K^{3,\text{bs}}_{\ge {\frak s}}$. \newline
4)  Assume 
$M_0 \le_{\frak K} M_1 \le_{\frak K} M_3 \le_{\frak K} M^*_3$ and
$a \in M_3$ \ub{then}: $\nonfork{}{}_{< \infty}
(M_0,M_1,a,M_3)$ \underbar{iff} \newline
$\nonfork{}{}_{< \infty}(M_0,M_1,a,M^*_3)$. 
\endproclaim
\bigskip

\demo{Proof}  1) First assume \ortp$(a,M,N) \in 
 {\Cal S}^{\text{bs}}_{\frak s}(M)$ 
and check the definition of $(M,N,a) \in K^{3,\text{bs}}$.  
Clearly $M \le_{\frak K} N,
a \in N$ and $a \in N \backslash M$; we have to find $M'$ as required 
in Definition \scite{600-1.6}(1); 
we let $M' = M$, so $M' \le_{\frak K} M,M' \in K_\lambda$ and

$$
\align
M' \le_{\frak K} M'' 
\le_{\frak K} M \and M'' \in K_\lambda &\Rightarrow M'' = M\\
  &\Rightarrow \text{ \ortp}_{{\frak K}_\lambda}(a,M'',N) 
= \text{ \ortp}_{{\frak K}_\lambda}(a,M,N) \in {\Cal
S}^{\text{bs}}_{\frak s}(M) = {\Cal S}^{\text{bs}}_{\frak s}(M'')
\endalign
$$
\mn
so we are done.

Second assume $(M,N,a) \in K^{3,\text{bs}}$ 
so there is $M' \le_{\frak K} M$ as
asserted in the definition \scite{600-1.6}(1) of 
$K^{3,\text{bs}}$ so $(\forall M'')[M' \le_{\frak K}
M'' \le_{\frak K} M \and M'' \in K_\lambda 
\Rightarrow \text{ \ortp}(a,M'',N) \in {\Cal S}^{\text{bs}}_{\frak s}(M'')]$ 
in particular this holds for $M'' = M$ and we get \ortp$(a,M,N) \in
{\Cal S}^{\text{bs}}_{\frak s}(M)$ as required. \newline
2) First assume $\nonfork{}{}_{< \infty} (M_0,M_1,a,M_3)$. \newline
So there is $M'_0$ as required in Definition \scite{600-1.7}; this means

$$
M'_0 \in K_\lambda,M'_0 \le_{\frak K} M_0 \text{ and}
$$

$$
\align
(\forall M'_1 \in K_\lambda)(\forall M'_3 \in K_\lambda)[M'_0
\le_{\frak K} M'_1 \le M_1 &\and M'_1 \cup \{a\} \subseteq M'_3 \le_{\frak K}
M_3 \\
  &\rightarrow \nonfork{}{}_{\lambda}(M'_0,M'_1,a,M'_3)].
\endalign
$$
\mn
In particular, we can choose $M'_1 = M_1,M'_3 = M_3$ so the antecedent holds
hence $\nonfork{}{}_{\lambda}(M'_0,M'_1,a,M'_3)$ which means
$\nonfork{}{}_{\lambda}(M'_0,M_1,a,M_3)$ and by clause $(E)(b)$ of
Definition \scite{600-1.1}, $\nonfork{}{}_{\lambda}(M_0,M_1,a,M_3)$ holds,
as required.

Second assume $\nonfork{}{}_{\lambda}(M_0,M_1,a,M_3)$.  So in Definition
\scite{600-1.7} the demands $M_0 \le_{\frak K} M_1 \le_{\frak K} M_3,a \in
M_3 \backslash M_1$ hold by clause $(E)(a)$ of Definition \scite{600-1.1}; and
we choose $M'_0$ as $M_0$; clearly $M'_0 \in K_\lambda \and M'_0 \le_{\frak K}
M_0$.  Now suppose $M'_0 \le_{\frak K} M'_1 \le_{\frak K} M_1 \and M'_1 \in
K_\lambda,M'_1 \cup \{a\} \le_{\frak K} M'_3 \le M_3$; by clause 
$(E)(b)$ of Definition 
\scite{600-1.1} we have $\nonfork{}{}_{\lambda}(M'_0,M'_1,a,M'_3)$; so $M'_0$
is as required so really $\nonfork{}{}_{< \infty}(M_0,M_1,a,M_3)$. \newline
3) We prove something stronger: for any $M' \in {\frak K}_{\frak s}$
which is $\le_{{\frak K}[{\frak s}]} M,M'$ witnesses 
$(M,N_1,a) \in K^{3,\text{bs}}$ iff $M'$ witnesses
$(M,N_2,a) \in K^{3,\text{bs}}$ (witness means: as required in
Definition \scite{600-1.6}).  So we have to check the statement there for
every $M'' \in K_\lambda$ such that $M' \le_{\frak s} M'' \le_{\frak
K} M$.  The equivalence holds 
because for every $M'' \le_{\frak K} M,M'' \in K_\lambda$ we
have \ortp$(a,M'',N_1) = \text{\ortp}(a,M'',N_2)$, by
\scite{600-0.12A}(2), more transparent as
${\frak K}_\lambda$ has the amalgamation property (by clause (C) of
Definition \scite{600-1.1}) and so one is ``basic" iff the other is 
by clause $(E)(b)$ of Definition \scite{600-1.1}. \newline
4) The direction $\Leftarrow$ is because if $M'_0$ witness
$\nonfork{}{}_{< \infty} (M_0,M_1,a,M^*_3)$ (see Definition
\scite{600-1.7}), \ub{then} it witnesses $\nonfork{}{}_{< \infty}
(M_0,M_1,a,M_3)$ as there are just fewer pairs $(M'_1,M'_3)$ to consider.
For the direction $\Rightarrow$ the demands $M_0 \le_{\frak K} M_1 
\le_{\frak K} M_3,a \in M_3 \backslash M_1$, of course, hold and let $M'_0$
be as required in the definition of $\nonfork{}{}_{< \infty}
(M_0,M_1,a,M_3)$; 
let $M'_0 \le_{\frak K} M'_1 \le_{\frak K} M_1,M'_1 \cup \{a\} \subseteq
M'_3 \le_{\frak K} M^*_3,M'_3 \in K_\lambda$.  
As $\lambda \ge \text{ LS}({\frak K})$ we can find $M''_3 \le_{\frak K} M_3$ 
such that $M'_1 \cup \{a\} \subseteq M''_3 \in K_\lambda$ and then
find $M'''_3 \le_{\frak s} M^*_3$ 
such that $M'_3 \cup M''_3 \subseteq M'''_3 \in K_\lambda$.  
So by the choice of $M'_0$ and $M''_3$ clearly $\nonfork{}{}_{\lambda}
(M'_0,M'_1,a,M''_3)$ and by clause $(E)(b)$ of Definition \scite{600-1.1}
we have

$$
\nonfork{}{}_{\lambda}(M'_0,M'_1,a,M''_3) \Leftrightarrow
\nonfork{}{}_{\lambda}(M'_0,M'_1,a,M'''_3) \Leftrightarrow
\nonfork{}{}_{\lambda}(M'_0,M'_1,a,M'_3)
$$
\medskip
\noindent
(note that we know the left statement and need the right statement)
so $M'_1$ is as required to complete the checking of
$\nonfork{}{}_{< \infty}(M_0,M_1,a,M^*_3)$. 
 \hfill$\square_{\scite{600-1.8}}$
\enddemo
\bn
We extend the definition of ${\Cal S}^{\text{bs}}_{\frak s}(M)$ 
from $M \in K_\lambda$ to arbitrary $M \in K$.

\definition{\stag{600-1.9} Definition}  1) For $M \in K$ we let

$$
\align
{\Cal S}^{\text{bs}}(M) = {\Cal S}^{\text{bs}}_{\ge {\frak s}}(M) = 
\biggl\{ p \in {\Cal S}(M):&\text{ for some } 
N \text{ and } a,\\
  &\,\,p = \text{ \ortp}(a,M,N) \text{ and } (M,N,a) 
\in K^{3,\text{bs}}_{\ge {\frak s}}\biggr\}
\endalign
$$
\mn
(for $M \in K_\lambda$ we get the old definition by \scite{600-1.8}(1); note that
as we do not have amalgamation (in general) the meaning of types is 
more delicate.  Not so in ${\frak K}_\lambda$ as in a good $\lambda$-frame we
have amalgamation in ${\frak K}_\lambda$ but not 
necessarily in ${\frak K}_{\ge \lambda}$).
\newline
2) We say that $p \in {\Cal S}^{\text{bs}}_{\ge {\frak s}}(M_1)$ 
does not fork over $M_0 
\le_{\frak K} M_1$ \ub{if} for some $M_3,a$ we have \newline
$p = \text{ \ortp}_{{\frak K}[{\frak s}]}(a,M_1,M_3)$ 
and $\nonfork{}{}_{< \infty}(M_0,M_1,a,M_3)$.
(Again, for $M \in K_\lambda$ this is equivalent to the old definition by
\scite{600-1.8}). \nl
3) For $M \in K$ let ${\Cal E}^\lambda_M$ be the following two-place relation
on ${\Cal S}(M):p_1 {\Cal E}^\lambda_M p_2$ iff $p_1,p_2 \in {\Cal
S}^{\text{bs}}(M)$ and if $p_\ell = \text{ \ortp}(a_\ell,M,M^*),N
\le_{\frak K} M,N \in K_\lambda$ then $p_1 \restriction N = p_2
\restriction N$.  Let ${\Cal E}^{\frak s}_M = {\Cal E}^{\lambda({\frak
s})}_M \restriction {\Cal S}^{\text{bs}}(M)$.
\nl
4) ${\frak K}$ is $(\lambda,\mu)$-local if every $M \in 
{\frak K}_\mu$ is $\lambda$-local which means that 
${\Cal E}^\lambda_M$ is equality; let $({\frak s},\mu)$-local
means $(\lambda_{\frak s},\mu)$-local.
\mn
Though we will prove 
below some nice things, having the extension property is 
more problematic.
We may define ``the extension" in a 
formal way, for $M \in K_{> \lambda}$ but then it is not
clear if it is realized in any $\le_{\frak K}$-extension of $M$.  Similarly 
for the uniqueness property.  
That is, assume $M_0 \le_{\frak K} M \le_{\frak K} N_\ell$ and $a_\ell
\in N_\ell \backslash M$, and $M_0 \in {\frak K}_{\frak s}$ and
\ortp$(a_\ell,M,N_\ell)$ does not fork over $M_0$ for $\ell=1,2$ and
\ortp$(a_1,M_0,N_1) = \text{\rm \ortp}(a_2,M_0,N_1)$.  Now does it
follow that $\ortp(a_1,M,N_1) = \text{\rm \ortp}(a_2,M,N_2)$?  
This requires the existence of some form of
amalgamation in ${\frak K}$, which we are not justified in assuming.
So we may prefer to define ${\Cal S}^{\text{bs}}(M)$ ``formally", the
set of stationarization of $p \in {\Cal S}^{\text{bs}}
(M_0),M_0 \in {\frak K}_{\frak s}$, see \cite{Sh:842}.
We now note that in definition \scite{600-1.9} 
``some" can be replaced by ``every".
\enddefinition
\bigskip

\demo{\stag{600-1.10} Fact}  1) For $M \in K$

$$
\align
{\Cal S}^{\text{bs}}_{\ge {\frak s}}(M) = 
\biggl\{ p \in {\Cal S}_{{\frak K}[{\frak s}]}(M):&\text{ for every } N,a \\
  &\text{ we have: if } M \le_{\frak K} N,a \in N \backslash M \text{ and} \\
  &\,p = \text{ \ortp}_{\frak K}(a,M,N) 
\text{ then } (M,N,a) \in K^{3,\text{bs}}_{\ge{\frak s}}\biggr\}.
\endalign
$$
\mn
2) The type $p \in {\Cal S}_{{\frak K}[{\frak s}]}(M_1)$ does not 
fork over $M_0 \le_{\frak K} M_1$ iff for every $a,M_3$ satisfying
$M_1 \le_{\frak K} M_3 \in K,a \in M_3 \backslash M_1$ and
$p = \text{ \ortp}_{{\frak K}[{\frak s}]}(a,M_1,M_3)$ we have 
$\nonfork{}{}_{< \infty}(M_0,M_1,a,M_3)$. \nl
3) $(M,N,a) \in K^{3,\text{bs}}_{\ge{\frak s}}$ is preserved by
isomorphisms.
\nl
4) If $M \le_{\frak K} N_\ell,a_\ell \in N_\ell \backslash M$ for
$\ell=1,2$ and $\ortp(a_1,M,N_1) {\Cal E}^{\frak s}_M
\ortp(a_2,M,N_2)$ then $(M,N_1,a_1) \in
K^{3,\text{bs}}_{\ge{\frak s}} \Leftrightarrow (M,N_2,a_2) \in
K^{3,\text{bs}}_{\ge{\frak s}}$.
\nl
5) ${\Cal E}^{\frak s}_M$ is an equivalence relation on ${\Cal
S}^{\text{bs}}_{\ge {\frak s}}(M)$ and if $p,q \in
{\Cal S}^{\text{bs}}_{\ge{\frak s}}(M)$ do not fork over $N \in 
K_\lambda$ so $N
\le_{\frak K} M$ then $p {\Cal E}^{\frak s}_M q \Leftrightarrow (p \restriction
N = q \restriction N)$.
\enddemo
\bigskip

\demo{Proof}  1) By \scite{600-1.8}(3) and the definition of type. \newline
2)  By \scite{600-1.8}(4) and the definition of type.  \nl
3) Easy.\nl
4) Enough to deal with the case $(M,N_1,a_1)
E^{\text{at}}_M,(M,N_2,a_2)$ or (by (3)) even $a_1 = a_2,N_1
\le_{\frak K} N_2$.  This is easy.
\nl
5) Easy, too.  \hfill$\square_{\scite{600-1.10}}$
\enddemo
\bn
We can also get that there are enough basic types, as follows:
\proclaim{\stag{600-1.11} Claim}  If $M \le_{\frak K} N$ and $M \ne N$, 
\underbar{then} for some $a \in N
\backslash M$ we have ${\text{\rm \ortp\/}}_{\frak K}(a,M,N) \in 
{\Cal S}^{\text{bs}}(M)$.
\endproclaim
\bigskip

\demo{Proof}  Suppose not, so by clause $(D)(c)$ of Definition \scite{600-1.1}
necessarily $\|N\| > \lambda$.  If $\|M\| = \lambda < \|N\|$ choose
$N'$ satisfying 
$M <_{\frak K} N' \le_{\frak K} N,N' \in K_\lambda$ and by clause $(D)(c)$ of
Definition \scite{600-1.1} choose $a^* \in N' \backslash M$ such 
that \ortp$_{\frak s}(a^*,M,N') \in {\Cal S}^{\text{bs}}_{\frak s}(M)$.  
So we can assume $\|M\| > \lambda$; choose 
$a^* \in N \backslash M$.  We choose by induction on 
$i < \omega,M_i,N_i,M_{i,c}$ (for $c \in N_i \backslash M_i)$
such that:
\mr
\item "{$(a)$}"  $M_i \le_{\frak K} M$ is $\le_{\frak K}$-increasing
\sn
\item "{$(b)$}"  $M_i \in K_\lambda$
\sn
\item "{$(c)$}"  $N_i \le_{\frak K} N$ is $\le_{\frak K}$-increasing
\sn
\item "{$(d)$}"  $N_i \in K_\lambda$
\sn
\item "{$(e)$}"  $a^* \in N_0$
\sn
\item "{$(f)$}"  $M_i \le_{\frak K} N_i$
\sn
\item "{$(g)$}"  if $c \in N_i \backslash M$, \ortp$_{\frak s}(c,M_i,N) \in
{\Cal S}^{\text{bs}}_{\frak s}(M_i)$ and there is
$M' \in K_\lambda$ such that $M_i \le_{\frak K} M' \le_{\frak K} M$ and
\ortp$_{\frak s}(c,M',N)$ forks over $M_i$ then 
$M_{i,c}$ satisfies this, otherwise $M_{i,c} = M_i$
\sn
\item "{$(h)$}"  $M_{i+1}$ includes the set 
$\dsize \bigcup_{c \in N_i \backslash M} M_{i,c} \cup(N_i \cap M)$.
\ermn
There is no problem to carry the definition.  Let $M^* = \dsize \bigcup_{i <
\omega} M_i$ and $N^* = \dsize \bigcup_{i < \omega} N_i$.  It is easy to
check that:
\mr
\widestnumber\item{(viii)}
\item "{$(i)$}"  $M_i \le_{\frak K} M^* \le_{\frak K} M$ for $i <
\omega$ 
\nl
(by clause (a))
\sn
\item "{$(ii)$}"  $M^* \in K_\lambda$ \newline
(by clause (i) we have $M^* \in K$ and $\|M^*\| = \lambda$ by the
choice of $M^*$ and clause (b))
\sn
\item "{$(iii)$}"  $N_i \le_{\frak K} N^* \le_{\frak K} N$ \newline
(by clause (c))
\sn
\item "{$(iv)$}"  $N^* \in K_\lambda$ \newline
(by clause (iii) we have $N^* \in K$ and $\|N^*\| = \lambda$ by the
choice of $N^*$ and clause (d))
\sn
\item "{$(v)$}"  $M_i \le_{\frak K} M^* 
\le_{\frak K} N^* \le_{\frak K} N$ \newline
(by clauses (a) + (f) + (iii) we have $M_i \le_{\frak K} N^*$ hence by
clause (a) and the choice of $M^*$ we have $M^* \le_{\frak K} N^*$,
and $N^* \le_{\frak K} N$ by clause (iii))
\sn
\item "{$(vi)$}"  $M^* = N^* \cap M$ \newline
(by clauses (f) + (h) and the choices of $M^*,N^*$)
\sn
\item "{$(vii)$}"  $M^* \ne N^*$ \nl
(as $a^* \in N \backslash M$ and $a^* \in N_0 \le_{\frak K} N^*
\le_{\frak K} N$ and $M^* = N^* \cap M$; \newline
they hold by the choice of $a^*$, clause (e), choice of $N^*$, clause
(iii) and clause (vi) respectively)
\sn
\item "{$(viii)$}"  there is $b^* \in N^* \backslash M^*$ such that
\ortp$(b^*,M^*,N^*) \in {\Cal S}^{\text{bs}}(M^*)$ \newline
[why?  by Definition \scite{600-1.1} clause (D)(c) (density)]
\sn
\item "{$(ix)$}"  for some $i < \omega$ 
we have $\nonfork{}{}_{}(M_i,M^*,b^*,N^*)$, so \nl
\ortp$(b^*,M^*,N^*) \in {\Cal S}^{\text{bs}}_{\frak s}(M^*)$ 
and \ortp$_{\frak s}(b^*,M_j,N^*) \in
{\Cal S}^{\text{bs}}_{\frak s}(M_j)$ for $j \in [i,\omega)$ \nl
[why?  by Definition \scite{600-1.1} clause $(E)(c)$ (local character)
applied to the sequence $\langle M_n:n < \omega \rangle \char 94
\langle M^*,N^* \rangle$ and the element $b^*$, using of course (E)(a)
of Definition \scite{600-1.1} and clause (viii)]
\sn
\item "{$(x)$}"  $\nonfork{}{}_{}(M_i,M_{i,b^*},b^*,N^*)$ \newline
[why?  by clause (ix) and Definition \scite{600-1.1}$(E)(b)$ (monotonicity) as 
\newline
$M_i \le_{\frak K} M_{i,b^*} \le_{\frak K} M_{i+1} \le_{\frak K} M^*$ by
clause (g) in the construction]
\sn
\item "{$(xi)$}"  if $M_i \le_{\frak K} M' \le_{\frak K} M$ and 
$M' \cup \{ b^*\} \subseteq N' \le_{\frak K} N$ 
and $M' \in K_\lambda,N' \in K_\lambda$ then
$\nonfork{}{}_{}(M_i,M',b^*,N')$ \newline
[why? by clause (x) and clause (g) in the construction.]
\ermn
So we are done. \hfill$\square_{\scite{600-1.11}}$
\enddemo
\bigskip

\proclaim{\stag{600-1.12} Claim}  If $M \le_{\frak K} N,a \in N \backslash M$, and
${\text{\rm \ortp\/}}(a,M,N) \in {\Cal S}^{\text{bs}}_{\ge {\frak s}}(M)$ 
\underbar{then} for some $M_0 \le_{\frak K} M$ we have
\mr
\item "{$(a)$}"  $M_0 \in K_\lambda$
\sn
\item "{$(b)$}"  ${\text{\rm \ortp\/}}(a,M_0,N) \in 
{\Cal S}^{\text{bs}}_{\frak s}(M_0)$
\sn
\item "{$(c)$}"  if $M_0 \le_{\frak K} M' \le_{\frak K} M$, 
\ub{then} ${\text{\rm \ortp\/}}(a,M',N) 
\in {\Cal S}^{\text{bs}}_{\frak s}(M')$ does not fork over $M_0$.
\endroster
\endproclaim
\bigskip

\demo{Proof}  Easy by now.
\enddemo
\bigskip

\proclaim{\stag{600-1.13} Claim}  1) Assume $M_1 \le_{\frak K} M_2$ and $p \in
{\Cal S}_{\frak K}(M_2)$.  \ub{Then} 
$p \in {\Cal S}^{\text{bs}}_{\ge {\frak s}}(M_2)$ and $p$ does not fork
over $M_1$ \underbar{iff} for some $N_1 \le_{\frak K} M_1,N_1 \in
K_\lambda$ and $p$ does not fork over $N_1$ \ub{iff} for some $N_1
\le_{\frak K} M_1,N_1 \in K_\lambda$ and we have 
$(\forall N)[N_1 \le_{\frak K} N \le_{\frak
K} M_2 \and N \in K_\lambda \Rightarrow p \restriction N \in {\Cal
S}^{\text{bs}}_{\frak s}(N) 
\and (p \restriction N$ does not fork over $N_1)]$; we
call such $N_1$ a witness, so every $N'_1 \in K_\lambda,N_1 \le_{\frak
K} N'_1 \le M_1$ is a witness, too. \newline 
2) Assume $M^* \in K$ and $p \in {\Cal S}_{\frak K}(M^*)$. \newline 
\ub{Then}: $p \in {\Cal S}^{\text{bs}}_{\ge {\frak s}}(M^*)$ \ub{iff} 
for some $N^* \le_{\frak K} M^*$ we have 
$N^* \in K_\lambda,p \restriction N^* \in 
{\Cal S}^{\text{bs}}(N^*)$ and $(\forall N \in K_\lambda)
(N^* \le_{\frak K} N \le_{\frak K} M^* \Rightarrow p
\restriction N \in {\Cal S}^{\text{bs}}(N)$ and does not fork over $N^*)$
(we say such $N^*$ is a witness, so any $N' \in K_\lambda,N^*
\le_{\frak K} N' \le_{\frak K} M$ is a witness, too). \newline 
3) (Monotonicity) \newline If $M_1 \le_{\frak K} M'_1 \le_{\frak K} M'_2
\le_{\frak K} M_2$ and $p \in {\Cal S}^{\text{bs}}_{\ge {\frak
s}}(M_2)$ does not fork over
$M_1$, \underbar{then} \newline $p \restriction M'_2 \in 
{\Cal S}^{\text{bs}}_{\ge{\frak s}}(M'_2)$ and it does not 
fork over $M'_1$. \newline 
4) (Transitivity) \newline If $M_0 \le_{\frak K} M_1 \le_{\frak K} M_2$
and $p \in {\Cal S}^{\text{bs}}_{\ge{\frak s}}(M_2)$ 
does not fork over $M_1$ and $p
\restriction M_1$ does not fork over $M_0$, \underbar{then} $p$ does
not fork over $M_0$. \newline 
5) (Local character) 
If $\langle M_i:i \le \delta + 1 \rangle$ is 
$\le_{\frak K}$-increasing continuous and $a \in M_{\delta + 1}$ and 
${\text{\rm \ortp\/}}_{\frak K}(a,M_\delta,M_{\delta + 1}) \in {\Cal S}^{\text{bs}}_{\ge{\frak s}}(M_\delta)$ \underbar{then} for some $i < \delta$ we have
${\text{\rm \ortp\/}}_{\frak K}(a,M_\delta,M_{\delta + 1})$ does not fork 
over $M_i$. \newline 
6) Assume that $\langle M_i:i \le \delta + 1 \rangle$ is 
$\le_{\frak K}$-increasing and 
$p \in {\Cal S}(M_\delta)$ and for every $i < \delta$ we have 
$p \restriction M_i \in {\Cal S}^{\text{bs}}_{\ge{\frak s}}(M_i)$ 
does not fork over $M_0$.  \ub{Then} 
$p \in {\Cal S}^{\text{bs}}_{\ge{\frak s}}(M_\delta)$ and 
$p$ does not fork over $M_0$.
\endproclaim
\bigskip

\demo{Proof}  1), 2)  Check the definitions. 
\newline
3) As $p \in {\Cal S}^{\text{bs}}_{\ge{\frak s}}(M_2)$ does not 
fork over $M_1$, there is $N_1
\in K_\lambda$ which witnesses it.

This same $N_1$ witnesses that $p \restriction M'_2$ does not fork over
$M'_1$. \newline
4) Let $N_0 \le_{\frak K} M_0$ witness that $p \restriction M_1$ does not
fork over $M_0$ (in particular $N_0 \in K_\lambda$); 
let $N_1 \le_{\frak K} M_1$ witness that $p$ does not fork over $M_1$ 
(so in particular $N_1 \in K_\lambda$).
Let us show that $N_0$ witnesses $p$ does not fork over $M_0$, so let
$N \in K_\lambda$ be such that $N_0 \le_{\frak K} N \le_{\frak K} M_2$
and we should just prove that $p \restriction N$ does not fork over
$N_0$.  We can find $N' \le_{\frak K} M_1,
N' \in K_\lambda$ such that $N_0 \cup N_1 \subseteq N'$, we can also find 
$N'' \le_{\frak K} M_2$ satisfying $N'' \in K_\lambda$ 
such that $N' \cup N \subseteq N''$.
As $N_1$ witnesses 
that $p$ does not fork over $M_1$, clearly $p \restriction 
N'' \in {\Cal S}^{\text{bs}}_{\frak s}(N'')$ 
does not fork over $N_1$, hence by 
monotonicity does not fork over $N'$.  As $N_0$ witnesses
$p \restriction M_1$ does not fork over
$M_0$, clearly $p \restriction N'$ belongs to ${\Cal
S}^{\text{bs}}(N')$ and does not fork over $N_0$, so by
transitivity (in ${\frak K}_{\frak s}$) we know that
$p \restriction N''$ does not fork over $N_0$; hence by
monotonicity $p \restriction N$ does not fork over $N_0$. \newline
5) Let $p = \text{ \ortp}_{\frak K}(a,M_\delta,M_{\delta +1})$ and let $N^*
\le_{\frak K} M_\delta$
witness $p \in {\Cal S}^{\text{bs}}
(M_\delta)$.  Assume toward contradiction that
the conclusion fails.  Without loss of generality cf$(\delta) = \delta$.
\enddemo
\bn
\ub{Case 0}:  $\|M_\delta\| \le \lambda (=\lambda_{\frak s})$.

Trivial.
\bn
\underbar{Case 1}:  $\delta < \lambda^+,\|M_\delta\| > \lambda$.

As $\|M_\delta\| > \lambda$, for some $i,\|M_i\| > \lambda$ so without loss
of generality $i < \delta \Rightarrow \|M_i\| > \lambda$.  We choose by
induction on $i < \delta$, models $N_i,N'_i$ such that:
\medskip
\roster
\item "{$(\alpha)$}"  $N_i \in K_\lambda$
\sn
\item "{$(\beta)$}"  $N_i \le_{\frak K} M_i$ (hence $N_i \le_{\frak K}
M_j$ for $j \in [i,\delta))$
\sn
\item "{$(\gamma)$}"  $N_i$ is $\le_{\frak K}$-increasing continuous
\sn
\item "{$(\delta)$}"  $N'_i \in K_\lambda,N^* \le_{\frak K} N'_0$
\sn
\item "{$(\varepsilon)$}"  $N_i \le_{\frak K} N'_i \le_{\frak K} 
M_\delta$,
\sn
\item "{$(\zeta)$}"  $N'_i$ is $\le_{\frak K}$-increasing continuous
\sn
\item "{$(\eta)$}"  $p \restriction N'_i$ forks over $N_i$ when $i \ne
0$ for simplicity
\sn
\item "{$(\theta)$}"   $N_i \cup 
\bigcup_{j \le i}(N'_j \cap M_{i+1}) \subseteq N_{i+1}$.
\ermn
No problem to carry the induction, but we give details.
\sn
First, if $i = 0$ trivial.  Second let $i$ be a limit ordinal.

Let $N_i = \cup\{N_j:j < i\}$, now $N_i \le_{\frak K} M_i$ by clauses
$(\beta) + (\gamma)$ and ${\frak K}$ being a.e.c. and 
$\|N_i\| = \lambda$ by clause $(\alpha)$, as
$i \le \delta < \lambda^+$; so clauses $(\alpha), (\beta), (\gamma)$
hold.  Next, let $N'_i = \cup\{N'_j:j <i\}$ and similarly clauses
$(\delta), (\varepsilon), (\zeta)$ hold.  Lastly, we shall prove 
clause $(\eta)$ and assume toward contradiction that it fails; so
$p \restriction N'_i$ does not fork over $N_i$  in
particular $p \restriction N_i \in {\Cal S}^{\text{bs}}_{\frak
s}(N_i)$ hence for some $j < i$ the type $p \restriction N'_i$ does not fork
over $N_j \le_{\frak K} N_i$, (by (E)(c) of Definition \scite{600-1.1})
hence by transitivity (for ${\frak K}_{\frak s}$), $p \restriction
N'_i$ does not fork over $N_j$ hence by monotonicity
$p \restriction N'_j$ does not fork over $N_j$ (see (E)(b) of
Definition \scite{600-1.1}) contradicting the induction hypothesis.

Lastly, clause $(\theta)$ is vacuous.  

Third assume $i=j+1$, so first choose $N_i$ satisfying clause
$(\theta)$ (with $j,i$ here standing for $i,i+1$ there), and
$(\alpha), (\beta), (\gamma)$; this is possible by the L.S. property.
Now $N_i$ cannot witness ``$p$ does not fork over $M_i$" hence for
some $N^*_i \in K_\lambda$ we have $N_i \le_{\frak K} N^*_i \le_{\frak
K} M_\delta$ and $p \restriction N^*_i$ forks over $N_i$; again by
L.S. choose $N'_i \in K_\lambda$ such that $N'_i \le_{\frak K} M_\delta$ and
$N^* \cup N_i \cup N'_j \cup N^*_i \subseteq N'_i$, easily
$(N_i,N'_i)$ are as required.

Let $N_\delta = \dsize \bigcup_{i < \delta} N_i$, so
by clause $(\beta), (\gamma)$ we have 
$N_\delta \le_{\frak K} M_\delta$ and by clause $(\alpha)$, as
$\delta < \lambda^+$ we have $N_\delta \in K_\lambda$ and by 
clauses $(\delta) + (\theta)$ in the construction we have $i < \delta
\Rightarrow N'_i = \cup\{N'_i \cap M_{j+1}:j \in [i,\delta)\}
\subseteq N$ so by clause $(\delta),N^* \le_{\frak K} N'_0 \le_{\frak
K} N_\delta$.  Hence by the choice of $N^*,
p \restriction N_\delta \in {\Cal S}^{\text{bs}}_{\frak s}(N_\delta)$
and it does not fork
over $N^*$.  Now as $p \restriction N_\delta \in {\Cal S}^{\text{bs}}_{\frak s}
(N_\delta)$ by local character, i.e., 
clause $(E)(c)$ of Definition \scite{600-1.1}, for some $i < \delta,
p \restriction
N_\delta$ does not fork over $N_i$ (so 
$p \restriction N_i \in {\Cal S}^{\text{bs}}_{\frak s}(N_i)$).  Now $N_i
\le_{\frak K} N'_i \le_{\frak K} 
M_\delta$ and by clause $(\theta)$ of the construction $N'_i \subseteq
N_\delta$ hence $N_i \le_{\frak K} N'_i \le_{\frak K} N_\delta$ hence by
monotonicity of non-forking (i.e. clause (E)(b) of Definition \scite{600-1.1}),
$p \restriction N'_i \in {\Cal S}^{\text{bs}}
(N_i)$ does not fork over $N_i$.  But this contradicts 
the choice of $N'_i$ (i.e., clause $(\eta)$ of the construction).
\bn
\underbar{Case 2}:  $\delta = \text{ cf}(\delta) > \lambda$.

Recall that $N^* \le_{\frak K} M_\delta,N^*$ is from $K_\lambda$ and $N^* 
\le_{\frak K} N \le_{\frak K} M_\delta \and N \in K_\lambda \Rightarrow$ \nl
$p \restriction N \in {\Cal S}^{\text{bs}}_{\frak s}(N)$.
Now as $\delta = \text{ cf}(\delta) > \lambda \ge \|N^*\|$ clearly
for some $i < \delta$ we have $N^* \subseteq M_i$ hence $N^*
\le_{\frak K} M_i$
(hence $i \le j < \delta \Rightarrow
p \restriction M_j \in {\Cal S}^{\text{bs}}_{\ge{\frak s}}(M_j)$), and $N^*$
witnesses that $p \in {\Cal S}^{\text{bs}}_{\ge{\frak s}}(M_\delta)$
does not fork over $M_i$ so we are clearly done.
\nl
6)  Let $N_0 \in
K_\lambda,N_0 \le_{\frak K} M_0$ witness $p \restriction M_0 \in
{\Cal S}^{\text{bs}}_{\ge{\frak s}}(M_0)$.  
By the proof of part (4) clearly $i < \delta \and N_0
\le_{\frak K} N \in K_\lambda \and N \le_{\frak K} M_i \Rightarrow p \restriction N$
does not fork over $N_0$.  If cf$(\delta) > \lambda$ we are done, so assume
cf$(\delta) \le \lambda$.  Let $N_0 \le_{\frak K} N^* \in K_\lambda
\and N^* \le_{\frak K} M_\delta$, and we shall prove that $p
\restriction N^*$ does not fork over $N_0$, this clearly suffices.  As
in Case 1 in the proof of part (5) we can find
$N_i \le_{\frak K} M_i$ for $i \in (0,\delta)$ such that $\langle
N_i:i \le \delta \rangle$ is $\le_{\frak K}$-increasing with $i$,
each $N_i$ belongs to ${\frak K}_\lambda$ and $N^* \cap M_i \subseteq
N_{i+1}$, hence $N^* \subseteq N_\delta := \dbcu_{i < \delta} N_i$.
Now $N_\delta \le_{\frak K} M_\delta$ and as said as $i < \delta
\Rightarrow p \restriction N_i \in {\Cal S}^{\text{bs}}_{\ge{\frak s}}
(N_i)$ does not
fork over $N_0$ hence $p \restriction N_\delta$ does not fork over
$N_0$ and by monotonicity 
$p \restriction N$ does not fork over $N_0$, as required.
\nl
${{}}$  \hfill$\square_{\scite{600-1.13}}$
\bigskip

\proclaim{\stag{600-1.13B} Lemma}  If $\mu = {\text{\rm cf\/}}
(\mu) > \lambda$ and $M \le_{\frak K} N$ are in $K_\mu$, \ub{then} 
we can find $\le_{\frak K}$-representations 
$\bar M,\bar N$ of $M,N$ respectively such that:
\mr
\item "{$(i)$}"  $N_i \cap M = M_i$ for $i < \mu$
\sn
\item "{$(ii)$}"   if $i < j < \mu \and a \in N_i$ \underbar{then}
\ermn
$$
\align
(a) \qquad 
{\text{\rm \ortp\/}}(a,M_i,N) \in {\Cal S}^{\text{bs}}_{\ge{\frak s}}
(M_i) &\Leftrightarrow {\text{\rm \ortp\/}}(a,M_j,
N) \in {\Cal S}^{\text{bs}}_{\ge{\frak s}}(M_j)  \\
  &\Leftrightarrow {\text{\rm \ortp\/}}(a,M,N) \text{ does not fork over } M_i \\
  &\Leftrightarrow { \text{\rm \ortp\/}}(a,M_j,N) 
\text{ is a non-forking extension of {\rm \ortp\/}} (a,M_i,N) 
\endalign
$$
$(b) \qquad M_i \le_{\frak K} N_i \le_{\frak K} N_j$ and $M_i
\le_{\frak K} M_j \le_{\frak K} N_j$ \nl

\hskip20pt (and obviously $M_i \le_{\frak K} N_j$ and $M_i \le_{\frak K} M,M_i
\le_{\frak K} N,N_i \le_{\frak K} N$).
\endproclaim
\bigskip

\remark{\stag{600-1.13C} Remark}  In fact for any representations $\bar
M,\bar N$ of $M,N$ respectively, for some club $E$ of $\mu$ the
sequences $\bar M \restriction E,\bar N \restriction E$ are as above.
\endremark
\bigskip

\demo{Proof}  Let $\bar M$ be a $\le_{\frak K}$-representation of $M$.
For $a \in N$ we define 
$S_a = \{\alpha < \mu:\text{\ortp}(a,M_\alpha,N) \in
{\Cal S}^{\text{bs}}_{\ge{\frak s}}(M_\alpha)\}$.  
Clearly if $\delta \in S_a$ is a limit ordinal 
then for some $i(a,\delta) < \delta$ we
have $i(a,\delta) \le i < \delta \Rightarrow i \in S_a \and
(\text{\ortp}(a,M_i,N)$ does not fork over $M_{i_{(a,\delta)}})$ by
\scite{600-1.13}(5).
So if $S_a$ is stationary, then for some $i(a) < \mu$ the set $S'_a =
\{\delta \in S_a:i(a,\delta)=i(a)\}$ is a stationary subset of
$\lambda$ hence by monotonicity we have $i(a)
\le i \le \mu \Rightarrow \text{ \ortp}(a,M_i,N)$ does not fork over
$M_{i(a)}$.  Let $E_a$ be a club of $\mu$ such that: if $S_a$ is not
stationary (subset of $\mu$) then $E_a \cap S_a = \emptyset$ and if
$S_a$ is not stationary then $S_a \cap E_a = \emptyset$.
\nl
Let $\bar N$ be a representation of $N$, and let

$$
\align
E^* = \{\delta < \mu:&N_\delta \cap M = M_\delta \text{ and }
M_\delta \le_{\frak K} M,N_\delta \le_{\frak K} N \\
  &\text{ and for every } a \in N_\delta \text{ we have } \delta \in E_a\}.
\endalign
$$
\medskip
\noindent
Clearly it is a club of $\mu$ and $\bar M \restriction E^*,\bar N
\restriction E^*$ are as required.  \hfill$\square_{\scite{600-1.13B}}$
\enddemo
\bn
\centerline{$* \qquad * \qquad *$}
\bn
We may treat the lifting of $K^{3,\text{bs}}_\lambda$ as a special
case of the ``lifting" of ${\frak K}_\lambda$ to ${\frak K}_{\ge \lambda} =
({\frak K}_\lambda)^{\text{up}}$ in Claim \scite{600-0.31}; this may be
considered a good exercise.
\proclaim{\stag{600-1.13E} Claim}  1) $(K^{3,\text{bs}}_\lambda,
\le_{\text{bs}})$ is a $\lambda$-a.e.c.
\nl
2) $(K^{3,\text{bs}}_{\ge {\frak s}},
\le_{\text{bs}})$ is $(K^{3,\text{bs}}_\lambda,\le_{\text{bs}})^{\text{up}}$.
\endproclaim
\bigskip

\remark{Remark}  What is the class in \scite{600-1.13E}(1)?  Formally let
$\tau^+ = \{R_{[\ell]}:R$ a predicate of $\tau_K,\ell=1,2\} \cup
\{F_{[\ell]}:F$ a function symbolf rom $\tau_K\} \cup \{c\}$ where
$R_{[\ell]}$ is an $n$-place predicate when $R \in \tau$ is an
$n$-place predicate and similarly $F_{[\ell]}$ and $c$ is an
individual constant.  A triple $(M,N,a)$ is identified with the
following $\tau^+$-model $N^+$ defined as follows:
\mr
\item "{$(a)$}"  its universe is the universe of $N$
\sn
\item "{$(b)$}"  $c^{N^+} = a$
\sn
\item "{$(c)$}"  $R^{N^+}_{[2]} = R^N$
\sn
\item "{$(d)$}"  $F^{N^+}_{[2]} = F^N$ 
\sn
\item "{$(e)$}"  $R^{N^+}_{[1]} = R^M$
\sn
\item "{$(f)$}"  $F^{N^+}_{[1]} = F^M$
\ermn
(if you do not like partial functions, extend them to functions with
full domain by $F(a_0,\ldots) = a_0$ when not defined if
$F$ has arity $>0$, if $F$ has arity zero it is an individual constant,
$F^{N^+} = F^N$ so no problem).
\endremark
\bigskip

\demo{Proof}  Left to the reader (in particular this means that
$K^{3,\text{bs}}_\lambda$ is closed under $\le_{\text{bs}}$-increasing
chains of length $< \lambda^+$).
\enddemo
\bn
Continuing 
\scite{600-0.31}, \scite{600-0.34} note that (and see more in \scite{600-1.16}):
\proclaim{\stag{600-1.14} Lemma}  Assume
\mr
\item "{$(a)$}"  ${\frak K}$ is an abstract elementary class with
${\text{\rm LS\/}}({\frak K}) \le \mu$
\sn
\item "{$(b)$}"  $K'_{\le \mu}$ is a class of $\tau_K$-model,
$K'_{\le \mu} \subseteq K_{\le \mu}$ is non-empty and
closed under $\le_{\frak K}$-increasing unions of length $< \mu^+$ and
isomorphisms  (e.g. the class of $\mu$-superlimit models of ${\frak
K}_\mu$, if there is one)
\sn
\item "{$(c)$}"  define $K' := 
\{M \in K:M \text{ is a $\le_{\frak K}$-directed
union of members of } K'_\mu\}$
\sn
\item "{$(d)$}"  let ${\frak K}' = (K',\le_{\frak K} \restriction K')$
so $\le_{{\frak K}'}$ is $\le_{\frak K} \restriction K'$, so
${\frak K}'_{\le \mu} := (K'_{\le \mu},\le_{\frak K} \restriction
K'_{\le \mu})$; or $\le_{\frak K}$ is as in
\scite{600-0.31}(1), see \scite{600-0.31}(4).
\ermn
\ub{Then}
\mr
\item "{$(A)$}"  ${\frak K}'$ is an abstract elementary class, 
${\text{\rm LS\/}}({\frak K}) \le { \text{\rm LS\/}}({\frak K}') \le \mu$
\sn
\item "{$(B)$}"  If $\mu \le \lambda$ and 
$({\frak K},\nonfork{}{}_{},{\Cal S}^{\text{bs}})$ 
is a good $\lambda$-frame and ${\frak K}'_\lambda$ has amalgamation
and JEP and $M \in {\frak K}'_\lambda \Rightarrow 
{\Cal S}_{{\frak K}'}(M) = {\Cal S}_{\frak K}(M)$, \ub{then} 
$({\frak K}',\nonfork{}{}_{},{\Cal S}^{\text{bs}})$ (with $\nonfork{}{}_{},
{\Cal S}^{\text{bs}}$ restricted to ${\frak K}'$) is a good $\lambda$-frame
\sn
\item "{$(C)$}"  in clause $(B)$, instead ``$M \in {\frak K}'_\lambda
\Rightarrow {\Cal S}_{{\frak K}'}(M)= {\Cal S}_{\frak K}(M)$, it
suffices to require: if $M \in {\frak K}'_\lambda,M \le_{\frak K} N
\in {\frak K}'_\lambda,p \in {\Cal S}^{\text{bs}}_{\frak s}(N),p$ does not
fork over $M$ and $p \restriction M$ is realized in some $M',M
\le_{{\frak K}'} M'$ \ub{then} $p$ is realized in some $N',N \le_{\frak
K} N' \in {\frak K}'_\lambda$.
\endroster
\endproclaim
\bigskip

\remark{Remark}  If in \scite{600-1.14}, $K'_\mu$ is not closed under
$\le_{\frak K}$-increasing unions, we can close it but then the ``so
${\frak K}'_{\le \mu} = \ldots$" in clause (d) may fail.
\endremark
\bigskip

\demo{Proof}  \ub{Clause (A)}:  As in \scite{600-0.31}.
\mn
\ub{Clauses (B),(C)}:  Check.  \hfill$\square_{\scite{600-1.14}}$
\enddemo
\bn
\centerline{$* \qquad * \qquad *$}
\bn
Next we deal with some implications between the axioms in \scite{600-1.1}.
\proclaim{\stag{600-1.15} Claim}  1) In Definition \scite{600-1.1} clause (E)(i) is
redundant, i.e., follows from the rest, recalling
\mr
\item "{$(E)(i)$}"  \ub{non-forking amalgamation}: \newline
if for $\ell = 1,2,M_0 \le_{\frak K} M_\ell$ are in $K_\lambda,a_\ell \in
M_\ell \backslash M_0,{\text{\rm \ortp\/}}
(a_\ell,M_0,M_\ell) \in {\Cal S}^{\text{bs}}(M_0)$,
\ub{then} we can find $f_1,f_2,M_3$ satisfying $M_0 \le_{\frak K} M_3 \in
K_\lambda$ such that for $\ell =1,2$ we have $f_\ell$ is a 
$\le_{\frak K}$-embedding of $M_\ell$ into $M_3$ over
$M_0$ and ${\text{\rm \ortp\/}}
(f_\ell(a_\ell),f_{3-\ell}(M_{3-\ell}),M_3)$ does not fork over $M_0$.
\ermn
2) In fact we use only Axioms (A),(C),(E)(b),(d),(f),(g) only.
\bn
\epsfbox{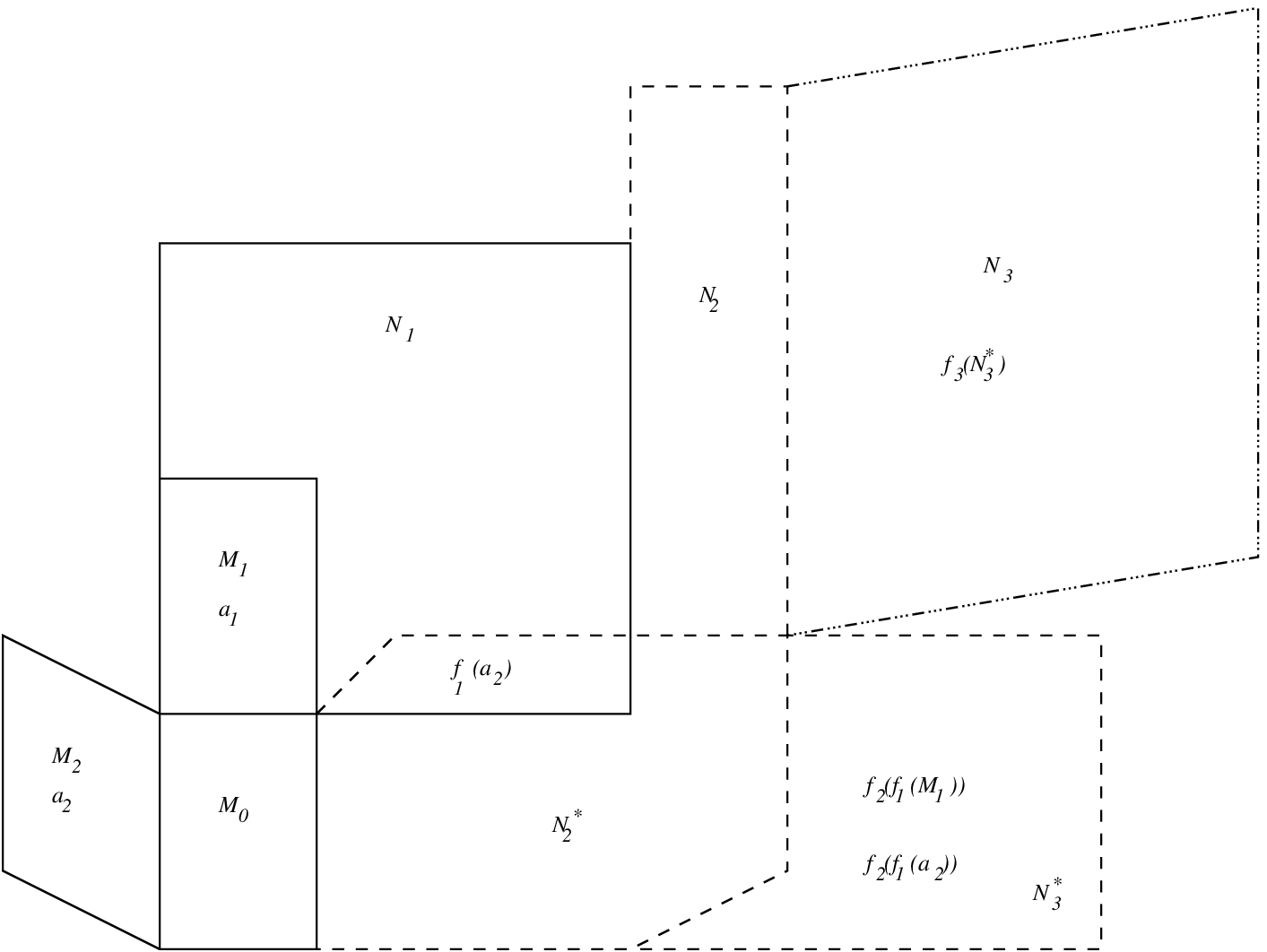}

\endproclaim
\bigskip

\demo{Proof}  By Axiom (E)(g) (existence) applied with \ortp$(a_2,M_0,M_2),
M_0,M_1$ here standing for $p,M,N$ there; there is $q_1$ such that:
\mr
\item "{$(a)$}"  $q_1 \in {\Cal S}^{\text{bs}}(M_1)$
\sn
\item "{$(b)$}"  $q_1$ does not fork over $M_0$
\sn
\item "{$(c)$}"  $q_1 \restriction M_0 = \text{ \ortp}(a_2,M_0,M_2)$.
\ermn
By the definition of types and as ${\frak K}_\lambda$ has amalgamation
(by Axiom (C)) there are $N_1,f_1$ such that
\mr
\item "{$(d)$}"  $M_1 \le_{\frak K} N_1 \in K_\lambda$
\sn
\item "{$(e)$}"  $f_1$ is a $\le_{\frak K}$-embedding of $M_2$ into $N_1$
over $M_0$
\sn
\item "{$(f)$}"  $f_1(a_2)$ realizes $q_1$ inside $N_1$.
\ermn
Now consider Axiom (E)(f) (symmetry) applied with $M_0,N_1,
a_1,f_1(a_2)$ here standing
for $M_0,M_3,a_1,a_2$ there; now as clause $(\alpha)$ of (E)(f) holds (use
$M_1,N_1$ for $M_1,M'_3$) we get that clause $(\beta)$ of (E)(f) holds which
means that there are $N_2,N^*_2$ (standing for $M'_3,M_2$ in clause $(\beta)$
of (E)(f)) such that:
\mr
\item "{$(g)$}" $N_1 \le_{\frak K} N_2 \in K_\lambda$
\sn
\item "{$(h)$}"  $M_0 \cup \{f_1(a_2)\} \subseteq N^*_2 \le_{\frak K} N_2$
\sn
\item "{$(i)$}"  \ortp$(a_1,N^*_2,N_2) \in 
{\Cal S}^{\text{bs}}(N^*_2)$ does not fork over $M_0$.
\ermn
As ${\frak K}_\lambda$ has amalgamation (see Axiom (C)) and the
definition of type and as
\nl
$\ortp(f_1(a_2),M_0,f_1(M_2)) = 
\text{ \ortp}(f_1(a_2),M_0,N_2)= 
\text{ \ortp}(f_1(a_2),M_0,N^*_2)$, we can find $N^*_3,f_2$ such that
\mr
\item "{$(j)$}"  $N^*_2 \le_{\frak K} N^*_3 \in K_\lambda$
\sn
\item "{$(k)$}"  $f_2$ is a $\le_{\frak K}$-embedding 
\footnote{we could have chosen $N^*_3 = N_2,f_2 = \text{ id}_{f_1(M_2)}$}
of $f_1(M_2)$ into
$N^*_3$ over $M_0 \cup \{f_1(a_2)\}$.
\ermn
As by clause (i) above \ortp$(a_1,N^*_2,N_2) \in {\Cal S}^{\text{bs}}(N^*_2)$,
so by Axiom (E)(g) (extension existence) there are $N_3,f_3$ such that
\mr
\item "{$(l)$}"  $N_2 \le_{\frak K} N_3 \in K_\lambda$
\sn
\item "{$(m)$}"  $f_3$ is a $\le_{\frak K}$-embedding of $N^*_3$ into $N_3$
over $N^*_2$
\sn
\item "{$(n)$}"  \ortp$(a_1,f_3(N^*_3),N_3) \in {\Cal S}^{\text{bs}}
(N^*_3)$ does not fork over $N^*_2$.
\ermn
By Axiom (E)(d) (transitivity) using clauses (i) + (n) above we have
\mr
\item "{$(o)$}"  \ortp$(a_1,f_3(N^*_3),N_3) \in {\Cal S}^{\text{bs}}
(N^*_3)$ does not fork over $M_0$.
\ermn
Letting $f = f_3 \circ f_2 \circ f_1$ as $f(M_2) \subseteq f_3(N^*_3)$
by clauses $(e),(k),(m)$ we have
\mr
\item "{$(p)$}"  $f$ is a $\le_{\frak K}$-embedding of $M_2$ into $N_3$
over $M_0$.
\ermn
By (E)(b) (monotonicity) and clause (o) and clause (p)
\mr
\item "{$(q)$}"  \ortp$(a_1,f(M_2),N_3) \in 
{\Cal S}^{\text{bs}}(f(M_2))$ does not
fork over $M_0$.
\ermn
As \ortp$(f_1(a_2),M_1,N_3) = \text{ \ortp}(f_1(a_2),M_1,N_1) = q_1$ does not fork
over $M_0$ by clauses (b) + (f), and $f_2(f_1(a_2)) = f_1(a_2)$ by
clause (k) and $f_3(f_1(a_2)) = f_1(a_2)$ by clauses (m) + (h), we get
\mr
\item "{$(r)$}"  \ortp$(f(a_2),M_1,N_3) \in {\Cal S}^{\text{bs}}(M_1)$ does not
fork over $M_0$.
\ermn
So by clauses (o) and (r) we have 
id$_{M_1},f,N_3$ are as required on $f_1,f_2,M_3$ in our desired
conclusion.  \hfill$\square_{\scite{600-1.15}}$
\enddemo
\bigskip

\proclaim{\stag{600-1.16A} Claim}  1) In the local character 
Axiom (E)(c) of Definition
\scite{600-1.1} if ${\Cal S}^{\text{bs}}_{\frak s} = {\Cal
S}^{\text{na}}_{{\frak K}_{\frak s}}$
where ${\Cal S}^{\text{na}}_{{\frak K}_s}(M) = \{\text{\rm
\ortp}(a,M,N):M \le_{\frak s} N$ and $a \in N \backslash M\}$
\ub{then} it
suffices to restrict ourselves to the case that $\delta$ has
cofinality $\aleph_0$ (i.e., the general case follows from this
special case and the other axioms). \nl
2) In fact in part (1) we need only Axioms (E)(b),(h) 
and you may say (A),(D)(a),(E)(a). \nl
3) If ${\Cal S}^{\text{bs}} = {\Cal S}^{\text{na}}$ \ub{then} the
continuity Axiom (E)(h) follows from the rest. \nl
4) In (3) actually we need only Axioms (E)(c), (local character) 
(d), (transitivity) and you may say (A),(D)(a),(E)(a). 
\endproclaim
\bigskip

\demo{Proof}  1), 2)  Let $\langle M_i:i \le \delta +1 \rangle$ be
$\le_{{\frak K}_\lambda}$-increasing, $a \in M_{\delta +1} \backslash
M_\delta$ and without loss of generality 
$\aleph_0 < \delta = \text{ cf}(\delta)$, so for every
$\alpha \in S := \{\alpha < \delta:\text{cf}(\alpha) = \aleph_0\}$,
\ortp$(a,M_\alpha,M_{\delta +1}) \in {\Cal S}^{\text{bs}}(M_\alpha)$ by
the assumption ``$S^{\text{bs}}_{\frak s} = {\Cal
S}^{\text{na}}_{{\frak K}_{\frak s}}$ hence 
there is $\beta_\alpha < \alpha$ such that \ortp$(a,M_\alpha,M_{\delta
+1})$ does not fork over $M_{\beta_\alpha}$, so for some
$\beta < \delta$ the set $S_1 = \{\alpha \in S:\beta_\alpha = \beta$) 
is a stationary subset of $\delta$.  
By Axiom $(E)(b)$ (monotonicity) it follows that for any
$\gamma_1 \le \gamma_2$ from $[\beta,\delta)$ the type
$\ortp(a,M_{\gamma_2},M_{\delta +1}) \in {\Cal
S}^{\text{bs}}(M_{\gamma_2})$ does not fork over $M_{\gamma_1}$.  Now
for any $\gamma \in [\beta,\delta)$ the type $\ortp(a,M_\delta,M_{\delta +1})$
does not fork over $M_\gamma$ by applying  $(E)(h)$ (continuity) to
$\langle M_\alpha:\alpha \in [\gamma,\delta +1]$ so we have finished. 
\nl
3),4) So assume $\langle M_i:i \le \delta \rangle$ is 
$\le_{\frak K}$-increasing all in $K_\lambda$ and $\delta$ is a limit ordinal,
$p \in {\Cal S}(M_\delta)$ and $p_i := p \restriction M_i \in 
{\Cal S}^{\text{bs}}(M_i)$ does not fork over $M_0$ for each $i <
\delta$;  we should prove
that $p \in {\Cal S}^{\text{bs}}(M_\delta)$ and $p$ does not fork over
$M_0$.

First, for each $i < \delta,p_i \in {\Cal S}^{\text{bs}}(M_i)$ hence
$p_i$ is not realized in $M_i$.  As
$M_\delta = \cup\{M_i:i < \delta\}$ clearly $p$ is not realized in
$M_\delta$ so $p \in {\Cal S}^{\text{na}}(M_\delta) = {\Cal
S}^{\text{bs}}(M_\delta)$.

Second, by Ax(E)(c) the type $p$ does not fork over $M_j$
for some $j < \delta$.  As $p_j = p \restriction M_j$ does not fork
over $M_0$ (by assumption) by the transitivity Axiom $(E)(d)$,  we get that $p$
does not fork over $M_0$, as required.   \hfill$\square_{\scite{600-1.16A}}$
\enddemo
\bigskip

\remark{Remark}  So in some sense by \scite{600-1.16A} we can omit in
\scite{600-1.1}, the local character Axiom $(E)(c)$ \ub{or} the continuity
Axiom $(E)(h)$ but \ub{not} both.  In fact (under reasonable
assumptions) they are equivalent.
\endremark
\bigskip

\proclaim{\stag{600-1.16B} Claim}  In Definition \scite{600-1.1}, Clause
(E)(d), i.e., transitivity of non-forking follows from
(A),(C),(D)(a),(b),(E)(a),(b),(e),(g). 
\endproclaim
\bigskip

\demo{Proof}  As ${\frak K}_\lambda$ is an $\lambda$-a.e.c. with amalgamation,
types as well as restriction of types are not only well defined but
are ``rasonable".

So assume $M_0 \le_{\frak s} M'_0 \le_{\frak s} M''_0 \le_{\frak s}
M_3,a \in M_3$ and $p'' := \text{ \ortp}_{\frak s}(a,M''_0,M_3)$ does not
fork over $M'_0$ and $p' := \text{ \ortp}_{\frak s}(a,M'_0,M_3)$ does not
fork over $M_0$.  Let $p = p' \restriction M_0$.  As $p'$ does not fork
over $M_0$, by Axiom $(E)(a)$ we have 
$p' \in {\Cal S}^{\text{bs}}(M'_0)$ and $p =
\text{ \ortp}(a,M_0,M_3) = p' \restriction M_0$ belongs to 
${\Cal S}^{\text{bs}}(M_0)$.  As $p''$ does not fork over $M'_0$ clearly $p''
\in {\Cal S}^{\text{bs}}(M''_0)$ and recall $p'' \restriction M'_0 = p'$.  By
the existence axiom $(E)(g)$ the type $p$ has an extension $q'' \in {\Cal
S}^{\text{bs}}(M''_0)$ which does not fork over $M_0$.  By the
monotonicity Axiom (E)(b) the type 
$q''$ does not fork over $M'_0$ and $q' = q''
\restriction M'_0$ does not fork over $M_0$.  As $p',q' \in {\Cal
S}^{\text{bs}}(M'_0)$ do not fork over $M_0$ and $p' \restriction
M_0 = p = q'' \restriction M_0 = q' \restriction M_0$, by the
uniqueness Axiom Ax(E)(e), we have $p'=q'$.
Similarly $p'' = q''$, but $q''$ does not fork over $M_0$ hence $p''$
does not fork over $M_0$ as required. \hfill$\square_{\scite{600-1.16B}}$
\enddemo
\bigskip

\proclaim{\stag{600-1.16E} Claim}  1) The symmetry axiom (E)(f) is
equivalent to (E)(f)$'$ and to (E)(f)$'$ below if we assume
(A),(B),(C),(D)(a0,(b),(E)(a),(b),(g) in Definition \scite{600-1.1}
\mr
\item "{$(E)(f)'$}"  there are no 
$M_\ell(\ell \le 3)$ and $a_\ell(\ell \le 2)$ such that
{\roster
\itemitem{ $(a)$ }  $M_0 \le_{\frak s} M_1 \le_{\frak s} M_2
\le_{\frak s} M_3$
\sn
\itemitem{ $(b)$ }  $\ortp(a_\ell,M_\ell,M_{\ell +1})$ does not fork over
$M_0$ for $\ell = 0,1,2$
\sn
\itemitem{ $(c)$ }  $\ortp_{\frak s}(a_0,M_0,M_1) = 
\ortp_{\frak s}(a_2,M_0,M_3)$
\sn
\itemitem{ $(d)$ }  $\ortp_{\frak s}(\langle a_0,a_1\rangle,M_0,M_1) \ne
\ortp_{\frak s}(\langle a_2,a_1\rangle,M_0,M_1)$.
\endroster}
\endroster
\endproclaim
\bigskip

\demo{Proof}  Easy.
\enddemo
\bn
\centerline{$* \qquad * \qquad *$}
\bn
A most interesting case of \scite{600-1.14} is the following.  In
particular it tells us that the categoricity assumption is not so rare
and it will have essential uses here.
\proclaim{\stag{600-1.16} Claim}  If 
${\frak s} = ({\frak K},\nonfork{}{}_{\lambda},
{\Cal S}^{\text{bs}})$ is a good $\lambda$-frame and $M \in K_\lambda$ is a
superlimit model in ${\frak K}_\lambda$ and we define ${\frak s}' =
{\frak s}^{[M]} = {\frak s}[M] = ({\frak K}[{\frak s}^{[M]}],
\nonfork{}{}_{\lambda}[{\frak s}^{[M]}],{\Cal S}^{\text{bs}}
[{\frak s}^{[M]}])$ by

$$
{\frak K}[{\frak s}^{[M]}] = {\frak K}^{[M]}, \text{ see Definition
\scite{600-0.33} so } {\frak K}_{{\frak s}[M]} = {\frak K} \restriction
\{N:N \cong M\}
$$

$$
\nonfork{}{}_{\lambda}[{\frak s}^{[M]}] = \{(M_0,M_1,a,M_3) \in
\nonfork{}{}_{\lambda}:M_0,M_1,M_3 \in K^{[M]}_\lambda\}
$$

$$
\align
{\Cal S}^{\text{bs}}[{\frak s}^{[M]}] = 
\bigl\{{\text{\rm \ortp\/}}_{{\frak K}[M]}(a,M_0,M_1):&M_0
\le_{\frak K} M_1,M_0 \in K^{[M]}_\lambda,N \in K^{[M]}_\lambda \\
  &\text{and {\rm \ortp}}_{\frak K} (a,M_0,M_1) \in 
{\Cal S}^{\text{bs}}(M_0) \bigr\}.
\endalign
$$
\mn
\ub{Then}
\mr
\item "{$(a)$}"  ${\frak s}'$ is a good $\lambda$-frame
\sn
\item "{$(b)$}"  ${\frak K}[{\frak s}'] \subseteq {\frak K}_{\ge \lambda}
[{\frak s}]$
\sn
\item "{$(c)$}"  $\le_{{\frak K}[{\frak s}']} = \le_{\frak K} \restriction
K[{\frak s}']$
\sn
\item "{$(d)$}"  $K_\lambda[{\frak s}']$ is categorical.
\endroster
\endproclaim
\bigskip

\demo{Proof}  Straight by \scite{600-0.31}, \scite{600-0.34}, \scite{600-1.14}.
\hfill$\square_{\scite{600-1.16}}$
\enddemo
\bigskip

------------------------------------------------------------
\newpage

\head {\S3 Examples} \endhead  \resetall \sectno=3
 \spuriousreset
\bigskip

We show here that the context from \S2 occurs in earlier investigation: in
\cite{Sh:88} = \chaptercite{88r}, \cite{Sh:576}, \cite{Sh:48} (and \cite{Sh:87a},
\cite{Sh:87b}). Of course, also the class $K$ of models of
a superstable (first order) theory $T$, with $\le_{\frak K} = \prec$ and
${\Cal S}^{\text{bs}}(M)$ being the set of 
regular types (also just ``the set non-algebraic types" works)
with $\nonfork{}{}_{}(M_0,M_1,a,M_3)$ iff $M_0 \le_{\frak K} M_1
\le_{\frak K} M_3$ are in $K_\lambda,a \in M_3$ and \ortp$(a,M_1,M_3) \in
{\Cal S}^{\text{bs}}(M_1)$ does not fork over $M_0$, (in the sense of
\cite[III]{Sh:c}, of course).  The
reader may concentrate on \scite{600-Ex.4} (or \scite{600-Ex.1}) below
for easy life. \nl
Note that \scite{600-Ex.1} (or \scite{600-Ex.1A}) will be used to continue 
\cite{Sh:88} = \chaptercite{88r} and also to give an alternative 
proof to the theorem of \cite{Sh:87a}, \cite{Sh:87b} + (deducing
``there is a model in $\aleph_n$" if there are not too many models in
$\aleph_\ell$ for $\ell < n$) and note that
\scite{600-Ex.1A} will be used to continue \cite{Sh:48}, i.e., on $\psi \in
\Bbb L_{\omega_1,\omega}(\bold Q)$ and \scite{600-Ex.4} will be used
to continue \cite{Sh:576}.  Many of the axioms from \scite{600-1.1} are easy.
\bigskip

\subhead {(A) The superstable prototype} \endsubhead

\proclaim{\stag{600-Ex.0} Claim}   Assume $T$ is a first order
complete theory and $\lambda$ be a cardinal $\ge |T| + \aleph_0$;
let ${\frak K} = {\frak K}_{T,\lambda} = 
(K_{T,\lambda} \le_{{\frak K}_{T,\lambda}})$ be defined by:
\mr
\item "{$(a)$}"  $K_{T,\lambda}$ 
is the class of models of $T$ of cardinality $\ge \lambda$
\sn
\item "{$(b)$}"  $\le_{{\frak K}_{T,\lambda}}$ is ``being elementary submodel".
\ermn
0) ${\frak K}$ is an a.e.c. with {\rm LS}$({\frak K}) = \lambda$.
\nl
1) If $T$ is superstable, stable in $\lambda$, \ub{then} 
${\frak s} = {\frak s}_{T,\lambda}$ is a good $\lambda$-frame
when ${\frak s} = ({\frak K}_{T,\lambda}{\Cal S}^{\text{bs}},\nonfork{}{}_{})$
is defined by:
\mr
\item "{$(c)$}"  $p \in {\Cal S}^{\text{bs}}(M)$ \ub{iff} $p = 
\ortp_{{\frak K}_{t,\lambda}}(a,M,N)$ for some $a,N$ such that
$\sftp_{\Bbb L(\tau_T)}(a,M,N)$, see Definition \scite{600-EX.0.7} is a
non-algebraic complete 1-type over $M$
\sn
\item "{$(d)$}"  $\nonfork{}{}_{}(M_0,M_1,a,M_3)$ iff $M_0 \prec M_1
\prec M_3$ are in $K_{T,\lambda}$ and $a \in M_3$ and
$\sftp_{\Bbb L(\tau_T)}(a,M_1,M_3)$ is a type that does not fork over
$M_0$ in the sense of \cite[III]{Sh:c}.
\ermn 
2) Let $\kappa = \text{\rm cf}(\kappa) \le \lambda$.  The model 
$M$ is a $(\lambda,\kappa)$-brimmed model for ${\frak K}_{T,\lambda}$
iff (i)+(ii) or (i)+(iii) where
\mr
\widestnumber\item{$(iii)$}
\item "{$(i)$}"  $T$ is stable in $\lambda$
\sn
\item "{$(ii)$}"  $\kappa = { \text{\rm cf\/}}(\kappa) \ge \kappa(T)$
and $M$ is a saturated model of $T$ of cardinality $\lambda$
\sn
\item "{$(iii)$}"  $\kappa = { \text{\rm cf\/}}(\kappa) < \kappa(T)$
and there is a $\prec$-increasing continuous sequence $\langle M_i:i
\le \kappa \rangle$ (by $\prec$, equivalently
by $\le_{\frak s}$) such that $M = M_\kappa$ and
$(M_{i+1},c)_{c \in M_i}$ is saturated for $i < \kappa$.
\ermn
2A) So there is a $(\lambda,\kappa)$-brimmed model for 
${\frak K}_{T,\lambda}$ iff $T$ is stable in $\lambda$.
\nl
3) $M$ is $(\lambda,\kappa)$-brimmed over $M_0$ in ${\frak K}_{T,\lambda}$ 
iff $(M,c)_{c \in M_0}$ is $(\lambda,\kappa)$-brimmed.
\nl
4) Assume $T$ is superstable first order complete theory stable in $\lambda$
and we define ${\frak s}^{\text{reg}}_{T,\lambda}$ as above only
${\Cal S}^{\text{bs}}(M)$ is the set of regular types 
$p \in {\Cal S}_{{\frak K}_T}(M)$.  \ub{Then} 
${\frak s}^{\text{reg}}_{T,\lambda}$ is a good $\lambda$-frame.
\nl
5) For $\kappa \le \lambda$ or $\kappa = \aleph_\varepsilon$ (abusing
notation), ${\frak s}^\kappa_{T,\lambda}$ is defined similarly
restricting ourselves to $\bold F^a_\kappa$-saturated models.  (Let
${\frak s}^0_{t,\lambda} = {\frak s}_{T,\lambda}$.)  If $T$ is
superstable, stable in $\lambda$ \ub{then} ${\frak s}^\kappa_{t,\lambda}$ is
a good $\lambda$ frame.
\endproclaim
\bigskip

\remark{Remark}  We can replace (c) of \scite{600-Ex.0} by:
\mr
\item "{$(c)'$}"  $p \in {\Cal S}^{\text{bs}}(M)$ iff $p = 
\ortp_{{\frak K}_{T,\lambda}}(a,M,N)$ for some $a,N$ such that
tp$_{\Bbb L(\tau_T)}(a,M,N)$ is a complete 1-type over $M$
\ermn
except that clause (D)(b) of Definition \scite{600-1.1} fail.
In fact the proofs are easier in this case; of course, the two meaning
of types essentially agree.
\endremark
\bigskip

\demo{Proof}  0),1),2),2A),3)  Obvious (see \cite{Sh:c}).
\nl
4) As in (1),  except density of regular types which holds 
by \cite{HuSh:342}). \nl
5) Also by \cite{Sh:c}.  \hfill$\square_{\scite{600-Ex.0}}$
\enddemo
\bn
Recall
\definition{\stag{600-EX.0.7} Definition}   1) For a logic ${\Cal L}$ and
vocabulary $\tau,{\Cal L}(\tau)$ is the set of ${\Cal L}$-formulas in
this vocabulary.
\nl
2) $\Bbb L = \Bbb L_{\omega,\omega}$ is first order logic.
\nl
3) A theory in ${\Cal L}(\tau)$ is a set of sentences from 
${\Cal L}(\tau)$ which we assume has a model if not said otherwise.
Similarly in a language $L(\subseteq {\Cal L}(\tau))$
\enddefinition
\bn
Very central in \chaptercite{88r} but 
peripheral here (except when in (parts of) \S3 we
continue \chaptercite{88r} in our framework) is:
\definition{\stag{600-0.4} Definition}  Let $T_1$ be a theory in 
$\Bbb L(\tau_1), \tau \subseteq \tau_1$ vocabularies,
$\Gamma$ a set of types in $\Bbb L(\tau_1)$; (i.e.
for some $m$, a set of formulas $\varphi(x_0,\dotsc,x_{m-1}) \in 
\Bbb L(\tau_1)$). \newline
1) EC$(T_1,\Gamma) = \{M:M \text{ a } \tau_1 \text{-model of } T_1
\text{ which omits every } p \in \Gamma\}$. \newline
(So \wilog \, $\tau_1$ is reconstructible from $T_1,\Gamma$) and

$$
\text{PC}_\tau(T_1,\Gamma) = \text{ PC}(T_1,\Gamma,\tau) 
= \{M:M \text{ is a } 
\tau \text{-reduct of some } M_1 \in \text{ EC}(T_1,\Gamma)\}.
$$
\mn
2) We say that ${\frak K}$ is PC$^\mu_\lambda$ or PC$_{\lambda,\mu}$
\ub{if} for some $T_1,T_2,\Gamma_1,\Gamma_2$ and $\tau_1$ and 
$\tau_2$ we have: 
($T_\ell$ a first order theory in the vocabulary $\tau_\ell,
\Gamma_\ell$ a set of types in $\Bbb L(\tau_\ell)$ and) 
$K = \text{ PC}(T_1,\Gamma_1,\tau_{\frak K})$ and
$\{(M,N):M \le_{\frak K} N$ and $M,N \in K\} = \text{ PC}
(T_2,\Gamma_2,\tau')$ where
\newline
$\tau' = \tau_{\frak K} \cup \{P\}$ 
($P$ a new one place predicate and $(M,N)$ means the $\tau'$-model
$N^+$ expanding $N$ where $P^{N^+} = |M|$), $|T_\ell| \le \lambda,
|\Gamma_\ell| \le \mu$ for $\ell = 1,2$.
\nl
3)  If $\mu = \lambda$, we may omit $\mu$.
\enddefinition
\bn
\subhead {(B) The 
abstract elementary class which is PC$_{\aleph_0}$} \endsubhead

\proclaim{\stag{600-Ex.1} Theorem}  Assume 
$2^{\aleph_0} < 2^{\aleph_1}$ and consider the statements
\medskip
\roster
\item "{$(\alpha)$}"  ${\frak K}$ is an abstract elementary class with
${\text{\rm LS\/}}({\frak K}) = \aleph_0$ 
(the last phrase follows by clause $(\beta))$ and $\tau = \tau({\frak
K})$ is countable
\sn
\item "{$(\beta)$}"  ${\frak K}$ is ${\text{\rm PC\/}}_{\aleph_0}$, 
equivalently for some $\psi_1,\psi_2 \in \Bbb L_{\omega_1,\omega}
(\tau_1)$ where
$\tau_1$ is a countable vocabulary extending $\tau$ we have
$$
\gather
K = \{M_1 \restriction \tau:M_1 \text{ a model of } \psi_1\} \\
\{(N,M):M \le_{\frak K} N\} = \{(N_1 \restriction \tau,M_1 \restriction
\tau):(N_1,M_1) \text{ a model of } \psi_2\}
\endgather
$$
\noindent
\item "{$(\gamma)$}"  $1 \le \dot I(\aleph_1,{\frak K}) < 2^{\aleph_1}$ 
\sn
\item "{$(\delta)$}"  ${\frak K}$ is categorical in $\aleph_0$, has the
amalgamation property in $\aleph_0$ and is stable in $\aleph_0$
\sn
\item "{$(\delta)^-$}"  like $(\delta)$ but ``stable in $\aleph_0$" is
weakened to: $M \in {\frak K}_{\aleph_0} \Rightarrow |{\Cal S}(M)| 
\le \aleph_1$
\sn
\item "{$(\varepsilon)$}"  all models of ${\frak K}$ are 
$\Bbb L_{\infty,\omega}$-equivalent and $M \le_{\frak K} N \Rightarrow M 
\prec_{{\Bbb L}_{\infty,\omega}} N$.
\ermn
For $M \in {\frak K}_{\aleph_0}$ we define ${\frak K}'_M$ as 
follows: the class of members is \nl
$\{N \in K:N \equiv_{{\Bbb L}_{\infty,\omega}} M\}$ and $N_1
\le_{{\frak K}'_M} N_2$ iff $N_1 \le_{\frak K} N_2 \and N_1
\prec_{{\Bbb L}_{\infty,\omega}} N_2$.
\medskip
\noindent
1) Assume $(\alpha) + (\beta) + (\gamma)$, \ub{then} for some
$M \in {\frak K}_{\aleph_0}$ the class ${\frak K}'_M$ satisfies
$(\alpha) + (\beta) + (\gamma) + (\delta)^- + (\varepsilon)$ 
(any $M \in {\frak K}_{\aleph_0}$ such that 
$({\frak K}'_M)_{\aleph_1} \ne \emptyset$ will do and there are such
$M \in K_{\aleph_0}$). 
Moreover, if ${\frak K}$ satisfies $(\delta)$ then also ${\frak
K}'_M$ satisfies it; also trivially $K'_M \subseteq K$ and
$\le_{{\frak K}'_M} \subseteq \le_{\frak K}$.
\nl
1A) Also there is ${\frak K}'$ such that: ${\frak K}'$ satisfies
$(\alpha) + (\beta) + (\gamma) + (\delta) + (\varepsilon)$, and for every
$\mu$ we have $K'_\mu \subseteq  K_\mu$. 
In fact, in the notation of \marginbf{!!}{\cprefix{88r}.\scite{88r-5.6}} for every $\alpha <
\omega_1$ we can choose ${\frak K}' = {\frak K}_{\bold D_\alpha}$.
\nl
2) Assume $(\alpha) + (\beta) + (\gamma) + (\delta)$.  \ub{Then}
$({\frak K},\nonfork{}{}_{},{\Cal S}^{\text{bs}})$ 
is a good $\aleph_0$-frame for 
some $\nonfork{}{}_{}$ and ${\Cal S}^{\text{bs}}$.  \nl
3)  In fact, in part (2) we can choose ${\Cal S}^{\text{bs}}(M) = 
\{p \in {\Cal S}(M):p \text{ not algebraic}\}$ and 
$\nonfork{}{}_{}$ is defined by \marginbf{!!}{\cprefix{88r}.\scite{88r-5.11}} 
(the definable extensions).
\endproclaim
\bigskip

\remark{Remark}  1) In \marginbf{!!}{\cprefix{88r}.\scite{88r-5.23}} we use the 
assumption $\dot I(\aleph_2,K) < \mu_{\text{unif}}(\aleph_2)$.  
But this Theorem is not used here!
\nl
2) Note that ${\frak K}'_M$ is related to $K^{[M]}$ from Definition
\scite{600-0.33} but is different.
\nl
3) In the proof we relate the types in the sense of 
${\Cal S}_{\frak s}(M)$, and those in \sectioncite[\S5]{88r}.  Now in
\sectioncite[\S5]{88r} we have lift types, from ${\frak K}_{\aleph_0}$ to
any ${\frak K}_\mu$, i.e., define $\bold D(N)$ for $N \in {\frak
K}_\mu$.  In $\mu > \aleph_0$, in general we do not know how to
relate them to types ${\Cal S}_{{\frak K}_{\frak s}}(N)$.  But when
${\frak s}^+$ is defined (in the ``successful" cases, see \S8 here and
\sectioncite[\S1]{705}) we can get the parallel claim.
\endremark
\bn
\ub{Discussion}:  1) What occurs if we do not pass in \scite{600-Ex.1} to the
case ``$\bold D(N)$ countable for every $N \in K_{\aleph_0}$"?  If we
still assume ``${\frak K}$ categorical in $\aleph_0$" then as $|\bold
D(N_0)| \le \aleph_1$, if we assume ``there is a superlimit model in
${\frak K}_{\aleph_1}$" we can find a good $\aleph_1$-frame ${\frak
s}$.  Not clear if we do not assume categoricity in $\aleph_1$.

\demo{Proof}  1) Note that for any $M \in K_{\aleph_0}$, the class 
${\frak K}'_M$ satisfies $(\alpha), (\beta), (\varepsilon)$ 
and it is categorical in $\aleph_0$ and
$(K'_M)_\mu \subseteq K_\mu$ hence $\dot I(\mu,K'_M) \le \dot I(\mu,K)$. 
By Theorem \marginbf{!!}{\cprefix{88r}.\scite{88r-3.6}}, (note: if you use the original version 
(i.e., \cite{Sh:88}) by its proof or use
it and get a less specified class with the desired properties) for some
$M \in K_{\aleph_0}$ we have $({\frak K}'_M)_{\aleph_1} \ne 
\emptyset$.   By \marginbf{!!}{\cprefix{88r}.\scite{88r-3.5}} (or \cite[1.4]{Sh:576}(page 46)1.6(page 48))
we  get that ${\frak K}'_M$ has 
amalgamation in $\aleph_0$ and by \chaptercite{88r} \ub{almost} we get that in
${\frak K}'_M$ the set ${\Cal S}(M)$ is small; be careful - the types there 
are defined differently than here, but by the
amalgamation (in $\aleph_0$)  
and the omitting types theorem in this case they are the same,
see more in the proof of part (3) below.  So by \marginbf{!!}{\cprefix{88r}.\scite{88r-5.1}},\marginbf{!!}{\cprefix{88r}.\scite{88r-5.2}} we 
have $M \in ({\frak K}'_\mu)_{\aleph_0} \Rightarrow 
|{\Cal S}_{{\frak K}'_\mu}(M)| \le \aleph_1$. 
\nl
Also the second and third sentences in (1) are easy. 
\nl
1A) Use \marginbf{!!}{\cprefix{88r}.\scite{88r-5.18}},\marginbf{!!}{\cprefix{88r}.\scite{88r-5.19}}. 
\nl
In more detail, (but not much point in reading without some understanding
of \sectioncite[\S5]{88r}, however we should not use \marginbf{!!}{\cprefix{88r}.\scite{88r-5.23}} as long as
we do not strengthen our assumptions) 
by part (1) we can assume that clauses $(\delta)^- +(\varepsilon)$ 
hold.  (Looking at the old version \cite{Sh:88} of
\chaptercite{88r} remember that there $\prec$ means $\le_{\frak K}$.)  
We can find $\bold D_* = \bold D^*_\alpha,\alpha < \omega_1$,
which is a good countable diagram (see Definition \marginbf{!!}{\cprefix{88r}.\scite{88r-5.6.1}}
and Fact \marginbf{!!}{\cprefix{88r}.\scite{88r-5.6}} or \marginbf{!!}{\cprefix{88r}.\scite{88r-5.16}}, \marginbf{!!}{\cprefix{88r}.\scite{88r-5.17}}.  
So in particular (give the
non-maximality of models below)  such
that for some countable $M_0 <_{\frak K} M_1 <_{\frak K} M_2$ we have
$M_m$ is $(\bold D^*(M_\ell),\aleph_0)$-homogeneous for $\ell < m \le 2$.
In \marginbf{!!}{\cprefix{88r}.\scite{88r-5.18}} we define $(K_{\bold D_*},\le_{\bold D_*})$.
By \marginbf{!!}{\cprefix{88r}.\scite{88r-5.19}} the pair $(K_{\bold D_*},\le_{\bold D_*})$ 
is an abstract elementary class (the choice of $\bold D_*$ a part, e.g.
transitivity = Axiom II which holds 
by the $M_\ell$'s above and \marginbf{!!}{\cprefix{88r}.\scite{88r-5.16}}) categorical in
$\aleph_0$ and no maximal countable model (by $\le_{\bold D_*}$, see
\marginbf{!!}{\cprefix{88r}.\scite{88r-5.6}}(2).
Now $\aleph_0$-stability holds by \marginbf{!!}{\cprefix{88r}.\scite{88r-5.19}}(2) and the
 equality of the three definitions of types in the proof of parts (2),(3)
 and $K_{\bold D_*} \subseteq K$ so we are done by part 3) below. 
\nl
2),3)  The first part of the proof serves also part (1) of the
theorem so we assume $(\delta)^-$ instead of $(\delta)$.  
We should be careful: the notion of type has three relevant
meanings here.
For $N \in K_{\aleph_0}$ the three definitions for 
${\Cal S}^{< \omega}(N)$ and of $\sftp(\bar a,N,M)$ when $\bar a \in
{}^{\omega >}M,N \le_{\frak K} M \in K_{\aleph_0}$ (of course we can 
use just 1-types) are:
\mr
\item "{$(\alpha)$}"  the one we use here (recall \scite{600-0.12}) which uses 
elementary mappings; for the present proof we 
call them ${\Cal S}^{< \omega}_0(M),\sftp_0(\bar a,M,N)$
\sn
\item "{$(\beta)$}"  ${\bold S}_1(N)$ which is (recall: materialzie is
close to but different from realize)
$$
\align
\bold D(N) = \bigl\{p:&\,p \text{ a complete } 
\Bbb L^0_{\aleph_1,\aleph_0}(N) \text{-type over } N \\
  &\text{ (so in each formula only finitely many parameters from } N 
\text{ appear)} \\
  &\text{ such that for some } M,\bar a \in {}^{\omega >}M,
\bar a \text{ materializes } p \text{ in } (M,N) \bigr\}
\endalign
$$
(``materializing a type" is defined in \marginbf{!!}{\cprefix{88r}.\scite{88r-4.2}}(2)) so
$$
{\bold S}_1(N) = \{\text{\sftp}_1(\bar a,N,M):\bar a \in {}^{\omega >}M
\text{ and } N \le_{\frak K} M \in K_{\aleph_0}\}
$$
where
$$
\text{\sftp}_1(\bar a,N,M) = \{\varphi(\bar x) \in 
\Bbb L^0_{\aleph_1,\aleph_0}(N):
M \Vdash^{\aleph_1}_{\frak K} \varphi(\bar a)\}
$$
(see \marginbf{!!}{\cprefix{88r}.\scite{88r-4.2}}(1) 
on the meaning of this forcing relation).
\sn
\item "{$(\gamma)$}"  $\quad {\bold S}_2(N)$ which is
$$
\align
\bold D^*(N) = \bigl\{p:&\,p \text{ a complete } 
\Bbb L^0_{\aleph_1,\aleph_0}(N;N) 
\text{-type over } N \\
  &\text{ (so in each formula all members of } N \text{ may appear)} \\
  &\text{ such that for some } M \in K_{\aleph_0} \text{ and} \\
  &\bar a \in {}^{\omega >}M \text{ satisfying } N \le_{\frak K} M
\text{ the sequence} \\
  &\bar a \text{ materializes } p \text{ in } (M,N) \bigr\}
\endalign
$$
\mn
so

$$
{\bold S}_2(N) = \{\sftp_2(\bar a,N,M):\bar a \in 
{}^{\omega >}M \text{ and } N \le_{\frak K} M \in K_{\aleph_0}\}
$$

$$
\text{\sftp}_2(\bar a,N,M) = \{\varphi(\bar x) \in \Bbb L^0_{\aleph_1,\aleph_0}
(N,N):M \Vdash^{\aleph_1}_{\frak K} \varphi(\bar a)\}.
$$
\ermn
As we have amalgamation in $K_{\aleph_0}$ it is enough to prove for 
$\ell,m < 3$ that
\mr
\item "{$(*)_{\ell,m}$}"  if $k < \omega,
N \le_{\frak K} M \in K_{\aleph_0}$ and $\bar a,\bar b \in {}^k M$, then \nl
$\sftp_\ell(\bar a,N,M) = \sftp_\ell(\bar b,N,M)
\Rightarrow \sftp_m(\bar a,N,M) = \sftp_m(\bar b,N,M)$.
\ermn
Now $(*)_{2,1}$ holds trivially (more formulas) and $(*)_{1,2}$ holds by
\marginbf{!!}{\cprefix{88r}.\scite{88r-5.5}}.  By amalgamation in ${\frak K}_{\aleph_0}$, 
if $\sftp_0(\bar a,N,M) = \sftp_0(\bar b,N,M)$, then for 
some $M',M \le_{\frak K} M' \in 
K_{\aleph_0}$ there is an automorphism $f$ of $M'$ over $N$ such that
$f(\bar a) = \bar b$, so trivially $(*)_{0,1},(*)_{0,2}$ hold
(we use the facts that $\sftp_\ell(\bar a,N,M)$ is preserved by
isomorphism and by replacing $M$ by $M_1$ if $M \le_{\frak K} M_2 \in
K_{\aleph_0}$ and $N \cup \bar a \subseteq M_1 \le_{\frak K}
M_2$).   Lastly we prove $(*)_{2,0}$. \nl
So $N \le_{\frak K} M \in K_{\aleph_0}$, hence 
$\{\text{\ortp}_2(\bar c,N,M):\bar c \in {}^{\omega >} M\} \subseteq 
\bold D^*(N)$ is countable so by \marginbf{!!}{\cprefix{88r}.\scite{88r-5.6}}(b),(c) for some 
countable $\alpha < \omega_1$
we have $\{\sftp_2(\bar c,N,M):\bar c \in {}^{\omega >} M\} \subseteq
\bold D^*_\alpha(N)$.  Now 
there is $M' \in K_{\aleph_0}$ such that $M \le_{\frak K}
M',M'$ is $(\bold D^*_\alpha,\aleph_0)$-homogeneous (by 
\marginbf{!!}{\cprefix{88r}.\scite{88r-5.6}}(e) see Definition \marginbf{!!}{\cprefix{88r}.\scite{88r-5.7}}) 
hence $M'$ is $(\bold D^*_\alpha(N),\aleph_0)$-
homogeneous (by \marginbf{!!}{\cprefix{88r}.\scite{88r-5.6}}(f)), and $\sftp_2(\bar a,N,M') =
\sftp_2(\bar b,N,M')$ by \marginbf{!!}{\cprefix{88r}.\scite{88r-5.4.1}}(3) 
($N$ here means $N_0$ there, 
that is increasing the model preserve the type). \nl
Lastly by Definition \marginbf{!!}{\cprefix{88r}.\scite{88r-5.7}} there is an automorphism $f$ of $M'$
over $N$ mapping $\bar a$ to $\bar b$, so we have proved $(*)_{2,0}$, so the
three definitions of type are the same.

Now we define for $M \in K_{\aleph_0}$:
\mr
\item "{$(a)$}"  ${\Cal S}^{\text{bs}}(M) = \{p \in 
{\Cal S}_{\frak K}(M):p \text{ not
algebraic}\}$
\sn
\item "{$(b)$}"  for $M_0,M_1,M_3 \in K_{\aleph_0}$ and an element 
$a \in M_3$ we define: 
\nl
$\nonfork{}{}_{}(M_0,M_1,a,M_3)$ \ub{iff} $M_0 \le_{\frak K} M_1
\le_{\frak K} M_3$ and $a \in M_3 \backslash M_1$ and \nl
$\sftp_1(a,M_1,M_3)(= \sftp(a,M_1,M_3)$ in \chaptercite{88r}'s notation) 
\nl
is definable over some finite $\bar b \in {}^{\omega >} M_0$ (equivalently
is preserved by every automorphism of $M_1$ over $\bar b$ (see
\marginbf{!!}{\cprefix{88r}.\scite{88r-5.11}}) \nl
(equivalently gtp$(a,M_1,M_3)$ is the stationarization of gtp$(a,M_0,M_3)$).
\ermn
Now we should check the axioms from Definition \scite{600-1.1}.
\bn
\ub{Clause (A)}:  By clause $(\alpha)$ of the assumption.
\mn
\ub{Clauses (B),(C)}:  By clause $(\delta)$ or $(\delta)^-$ of 
the assumption except ``the superlimit
$M \in K_{\aleph_0}$ is not $\le_{\frak K}$-maximal" which holds by clause
$(\gamma)+(\delta)$ or $(\gamma) + (\delta)^-$.
\mn
\ub{Clause (D)}:  By the definition (note that about clause (d), bs-stability,
that it holds by assumption $(\delta)$, and about clause (c), i.e.,
the density is trivial by the way we have defined ${\Cal S}^{\text{bs}}$).
\mn
\ub{Subclause (E)(a)}:  By the definition.
\mn
\ub{Subclause (E)(b)(monotonicity)}:

Let 
$M_0 \le_{\frak K} M'_0 \le_{\frak K} M'_1 \le_{\frak K} M_1 \le_{\frak K}
M_3 \le M'_3$ be all in ${\frak K}_{\aleph_0}$ and assume
$\nonfork{}{}_{}(M_0,M_1,a,M_3)$.  So
$M'_0 \le_{\frak K} M'_1 \le_{\frak K} M_3 \le M'_3$ and 
$a \in M_3 \backslash M_1
\subseteq M'_3 \backslash M'_1$.  Now by the assumption and the
definition of $\nonfork{}{}_{}$, for some $\bar b \in {}^{\omega >}
(M_0)$, gtp$(a,M_1,M_3)$ is definable over $\bar b$.  So the same 
holds for gtp$(a,M'_1,M_3)$ by \marginbf{!!}{\cprefix{88r}.\scite{88r-5.13}}, in fact
(with the same definition) and hence for gtp$(a,M'_1,M'_3) =
\text{ gtp}(a,M'_1,M_3)$ by \marginbf{!!}{\cprefix{88r}.\scite{88r-5.4.1}}(3), so as $\bar b \in
{}^{\omega >}(M_0) \subseteq {}^{\omega >}(M'_0)$ we have gotten
$\nonfork{}{}_{}(M'_0,M'_1,a,M'_3)$.  

For the additional clause in the monotoncity Axiom, assume 
in addition $M'_1 \cup \{a\} \subseteq
M''_3 \le_{\frak K} M'_3$ again by \marginbf{!!}{\cprefix{88r}.\scite{88r-5.4.1}}(3) 
clearly gtp$(a,M'_1,M''_3) = \text{
gtp}(a,M'_1,M'_3)$, so (recalling the beginning of the proof) we are done.
\mn
\ub{Sublcause (E)(c)(local character)}:

So let $\langle M_i:i \le \delta +1 \rangle$ be $\le_{\frak K}$-increasing
continuous in $K_{\aleph_0}$ and $a \in M_{\delta +1}$ and
\ortp$(a,M_\delta,M_{\delta +1}) \in {\Cal S}^{\text{bs}}
(M_\delta)$, so $a \notin
M_\delta$ and gtp$(a,M_\delta,M_{\delta +1})$ is definable over some $\bar b
\in {}^{\omega >}(M_\delta)$ by \marginbf{!!}{\cprefix{88r}.\scite{88r-5.4}}. \nl
As $\bar b$ is finite, for some $\alpha < \delta$ we have 
$\bar b \subseteq M_\alpha$,
hence as in the proof of (B)(b), we have
(\ortp$(a,M_\beta,M_{\delta +1}) \in {\Cal S}^{\text{bs}}
(M_\beta)$ trivially and)
\ortp$(a,M_\delta,M_{\delta +1})$ does not fork over $M_\beta$.
\mn
\ub{Sublcause (E)(d)(transitivity)}:

By \marginbf{!!}{\cprefix{88r}.\scite{88r-5.13}}(2) or even better \marginbf{!!}{\cprefix{88r}.\scite{88r-5.16}}.
\mn
\ub{Subclause (E)(e)(uniqueness)}:

Holds by the Definition \marginbf{!!}{\cprefix{88r}.\scite{88r-5.11}}.
\mn
\ub{Subclause (E)(f)(symmetry)}:

By \marginbf{!!}{\cprefix{88r}.\scite{88r-5.20}} + uniqueness we get (E)(f).  Actually
\marginbf{!!}{\cprefix{88r}.\scite{88r-5.20}} gives this more directly.
\mn
\ub{Subclause (E)(g)(extension existence)}:

By \marginbf{!!}{\cprefix{88r}.\scite{88r-5.11}} (i.e., by \marginbf{!!}{\cprefix{88r}.\scite{88r-5.4}} + all $M \in
K_{\aleph_0}$ are $\aleph_0$-homogeneous).

Alternatively, see \marginbf{!!}{\cprefix{88r}.\scite{88r-5.15}}.
\mn
\ub{Subclause (E)(h)(continuity)}:

Suppose $\langle M_\alpha:\alpha \le \delta \rangle$ is $\le_{\frak K}$-
increasingly continuous, $M_\alpha \in K_{\aleph_0},\delta < \omega_1,
p \in {\Cal S}(M_\delta)$ and $\alpha < \delta \Rightarrow p 
\restriction M_\alpha$ does 
not fork over $M_0$.  Now we shall use (E)(c)+(E)(d).  
As $p \restriction M_\alpha \in {\Cal S}^{\text{bs}}(M_\alpha)$ clearly
$p \restriction M_\alpha$ is not realized in $M_\alpha$ hence $p$ 
is not realized in $M_\alpha$;
as $M_\delta = \dbcu_{\alpha < \delta} M_\alpha$ necessarily $p$ is not
realized in $M_\delta$, hence $p$ is not algebraic.

So $p \in {\Cal S}^{\text{bs}}(M_\delta)$.  For some finite $\bar b \in
{}^{\omega >}
(M_\delta),p$ is definable over $\bar b$, let $\alpha < \delta$
be such that $\bar b \in {}^{\omega >}(M_\alpha)$, so as in the proof
of (E)(c), (or use it directly) the type $p$ does not fork over
$M_\alpha$.  As $p \restriction M_\alpha$ does not fork over $M_0$, by
(E)(d) we get that $p$ does not fork over $M_0$ as required.  Actually
we can derive (E)(h) by \scite{600-1.16A}.
\mn
\ub{Subclause (E)(i)(non-forking amalgamation)}:

One way is by \marginbf{!!}{\cprefix{88r}.\scite{88r-5.20}}; 
(note that in \marginbf{!!}{\cprefix{88r}.\scite{88r-5.23}} we get more, but
assuming, by our present notation $\dot I(\aleph_2,K) < \mu_{\text{wd}}
(\aleph_2)$); \ub{but} another way is just to use \scite{600-1.15}.  \nl
${{}}$  \hfill$\square_{\scite{600-Ex.1}}$
\enddemo
\bn
\subhead {(C) The uncountable cardinality quantifier case, 
$\Bbb L_{\omega_1,\omega}(\bold Q)$} \endsubhead

Now we turn to $\Bbb L$ sentences in $\Bbb L_{\omega_1,\omega}(\bold Q)$.
\demo{\stag{600-Ex.1A} Conclusion}  Assume $\psi \in 
\Bbb L_{\omega_1,\omega}(\bold Q)$ and $1 \le \dot I(\aleph_1,\psi) 
< 2^{\aleph_1}$ and $2^{\aleph_0} < 2^{\aleph_1}$. \nl
\ub{Then} for 
some abstract elementary classes ${\frak K},{\frak K}^+$
(note $\tau_\psi \subset \tau_{\frak K} = \tau_{{\frak K}^+}$) we have:
\mr
\item "{$(a)$}"  ${\frak K}$ satisfies $(\alpha),(\beta),(\delta),
(\varepsilon)$ from \scite{600-Ex.1} with $\tau_{\frak K} \supseteq
\tau_\psi$ countable (for $(\gamma),(b)$ is a replacement)
\sn
\item "{$(b)$}"  for every $\mu > \aleph_0,
\dot I(\mu,{\frak K}(\aleph_1$-saturated)) $\le \dot I(\mu,\psi)$, where
\footnote{much less than saturation suffice, like ``obeying" $<^{**}$}
``$\aleph_1$-saturated" is well defined as ${\frak K}_{\aleph_0}$ has
amalgamation, see \scite{600-0.19}
\sn
\item "{$(c)$}"  for some $\nonfork{}{}_{},
{\Cal S}^{\text{bs}}$ (and $\lambda =
\aleph_0$), the triple $({\frak K},\nonfork{}{}_{},S^{\text{bs}})$ 
is as in \scite{600-Ex.1}(2) so is a good $\aleph_0$-frame 
\sn
\item "{$(d)$}"  every $\aleph_1$-saturated member of ${\frak K}$
belongs to ${\frak K}^+$ and there is an $\aleph_1$-saturated 
member of ${\frak K}$ (and naturally it is uncountable)
\sn
\item "{$(e)$}"  ${\frak K}^+$ is an a.e.c., has LS number $\aleph_1$ 
and $\{M \restriction \tau_\psi:M \in {\frak K}^+\} \subseteq
\{M:M \models \psi\}$ and every $\tau$-model $M$ of $\psi$ has a
unique expansion in ${\frak K}^+$ hence $\mu \ge \aleph_1 \Rightarrow
\dot I(\mu,\psi) = \dot I(\mu,{\frak K}^+)$ and ${\frak K}^+$ is the
class of models of some complete $\psi \in \Bbb
L_{\omega_1,\omega}(\bold Q)$.
\endroster
\enddemo
\bigskip

\demo{Proof}  Essentially by \cite{Sh:48} and \scite{600-Ex.1}.

I feel that upon reading the proof \cite{Sh:48} the proof should not
be inherently difficult, much more so having 
read \scite{600-Ex.1}, but will give full details.
\nl
Recall Mod$(\psi)$ is the class of $\tau_\psi$-models of $\psi$.  We
can find a countable fragment ${\Cal L}$ of $\Bbb
L_{\omega_1,\omega}(\bold Q)(\tau_\psi)$ to which $\psi$ belongs and a
sentence $\psi_1 \in {\Cal L} \subseteq \Bbb L_{\omega_1,\omega}
(\bold Q)(\tau_\psi)$ such that $\psi_1$ is ``nice" for 
\cite[Definition 3.1,3.2]{Sh:48}, \cite[Lemma 3.1]{Sh:48}
\mr
\item "{$\circledast_1$}"  $(a) \quad \psi_1$ has uncountable models
\sn
\item "{${{}}$}"  $(b) \quad \psi_1 \vdash \psi$, i.e., every model of
$\psi_1$ is a model of $\psi$
\sn
\item "{${{}}$}"  $(c) \quad \psi_1$ is $\Bbb
L_{\omega_1,\omega}(\bold Q)$-complete
\sn
\item "{${{}}$}"  $(d) \quad$ every model $M \models \psi_1$ realizes
just countably many complete \nl

\hskip25pt $\Bbb L_{\omega_1,\omega}(\bold
Q)(\tau_\psi)$-types (of any finite arity, over the empty set), 
\nl

\hskip25pt each isolated by a formula in ${\Cal L}$.
\ermn
The proof of $\circledast_1(d)$ is sketched in Theorem \ub{2.5} of
\cite{Sh:48}.  The reference to Keisler \cite{Ke71} is to the generalization of
theorems 12 and 28 of Keisler's book from $\Bbb L_{\omega_1,\omega}$ to
$\Bbb L_{\omega_1,\omega}(\bold Q)$, see \marginbf{!!}{\cprefix{88r}.\scite{88r-0.1}}.
\nl
Let
\mr
\item "{$\circledast_2$}"  $(i) \quad {\frak K}_0 =
(\text{Mod}(\psi),\prec_{\Cal L})$, 
\sn
\item "{${{}}$}"  $(ii) \quad {\frak K}_1 =
(\text{Mod}(\psi_1),\prec_{\Cal L})$
\sn
\item "{$\circledast_3$}"  ${\frak K}_\ell$ is an a.e.c. with
L.S. number $\aleph_1$ for $\ell=0,1$.
\ermn
Toward defining ${\frak K}$, let $\tau_{\frak K} = \tau_\psi
\cup\{R_{\varphi(\bar x)}:\varphi(\bar x) \in {\Cal
L}\},R_{\varphi(\bar x)}$ a new $\ell g(\bar x)$-predicate and let
$\psi_2 = \psi_1 \wedge \{(\forall \bar y)(R_{\varphi(\bar x)}(\bar y)
= \varphi(\bar y):\varphi(\bar x) \in \Bbb L\}$.  For every
$M \in \text{ Mod}(\psi)$ we define $M^+$ by
\mr
\item "{$\circledast_4$}"  $M^+$ is $M$ expanded to a 
$\tau_{\frak K}$-model by letting $R^{M^+}_{\varphi(\bar x)} = \{\bar
a \in {}^{\ell g(\bar x)} M:M \models \varphi[\bar a]\}$
\sn
\item "{$\circledast_5$}"  $(a) \quad {\frak K}^+_0 = (\{M^+:M \in \text{
Mod}(\psi)\},\prec_{\Bbb L})$ is an a.e.c. with LS$({\frak K}^+_0) =
\aleph_1$
\sn
\item "{${{}}$}"  $(b) \quad {\frak K}^+_1 
= (\{M^+:M \in \text{ Mod}(\psi_1)\},\prec_{\Bbb L})$
is an a.e.c. with LS$({\frak K}^+) = \aleph_1$.
\ermn
Clearly
\mr
\item "{$\circledast_6$}"  if $M \models \psi_1$ then
$M^+$ is an atomic model of the complete first-order theory
$T_{\psi_1}$ where $T_{\psi_1}$ is the set of first order consequences
in $\Bbb L(\tau_{\frak K})$ of $\psi_2$.
\ermn
So it is natural to define ${\frak K}$:
\mr
\item "{$\circledast_7$}"  $(a) \quad N \in {\frak K}$ iff
{\roster
\itemitem{ $(i)$ }  $N$ is a $\tau_{\frak K}$-model which is an atomic
model of $T_{\psi_1}$
\sn
\itemitem{ $(ii)$ }  if $\psi_1 \vdash (\forall \bar x)[\varphi_1(\bar
x) = (\bold Q y) \varphi_2(y,\bar x)]$ and 
$\varphi_1,\varphi_2 \in {\Cal L}$ and
$N \models \neg R_{\varphi_1(\bar x)}[\bar a]$ then $\{b \in N:N
\models R_{\varphi_2(y,\bar x)}(b,\bar a)\}$ is countable
\endroster}
\item "{${{}}$}"  $(b) \quad N_1 \le_{\frak K} N_2$ \ub{iff}
$(N_1,N_2 \in K,N_1 \prec_{\Bbb L} N_2$ equivalently $N_1 \subseteq
N_2$ and) for 
\nl

\hskip25pt $\varphi_1(\bar x),\varphi_2(y,\bar x)$ 
as in subclause (ii) of clause (a) above, 
if $\bar a \in {}^{\ell g(\bar x)}(N_1)$,
\nl

\hskip25pt $N_1 \models \neg R_{\varphi_1(\bar x)}[\bar a]$ and $b \in
N_2 \backslash N_1$ then $N_2 \models \neg R_{\varphi_2(y,\bar x)}[b,\bar a]$.
\ermn
Observe
\mr
\item "{$\circledast_8$}"  $N \in {\frak K}$ iff $N$ is an 
atomic $\tau_{\frak K}$-model of the first order 
$\Bbb L(\tau_{\frak K})$-consequences $\psi_2$ (i.e. of $\psi$ and
every $\tau_{\frak K}$ sentence of the form $\forall \bar x
[R_\varphi(\bar x) \equiv \varphi(\bar x)])$ and clause (ii) of
$\circledast_7(a)$ holds
\sn
\item "{$\circledast_9$}"  ${\frak K}$ is an a.e.c. with
LS$({\frak K}) = \aleph_0$ and is PC$_{\aleph_0},
{\frak K}$ is categorical in $\aleph_0$ (and $\le_{\frak K}$ 
called $\le^*$ in \cite[Definition 3.3]{Sh:48}).
\ermn
Note that ${\frak K}_1,{\frak K}$ has the same number of models, but
 ${\frak K}$ has ``more models" than ${\frak K}^+_1$, in
particular, it has countable members and ${\frak K}_0$ has at least as
many models as ${\frak K}_1$.  For $N \in {\frak K}$ to be in 
${\frak K}^+_1 = \{M^+:M
\in \text{ Mod}(\psi_1)\}$ what is missing is the other implications
in $\circledast_7(a)(ii)$. \nl
This is very close to \scite{600-Ex.1}, but ${\frak K}$ may have many
models in $\aleph_1$ (as $\bold Q$ is not necessarily interpreted as
expected).  However, 
\mr
\item "{$\circledast_{10}$}"  constructing $M \in K_{\aleph_1}$ by the
union as $\le_{\frak K}$-increasing continuous chain $\langle M_i:i <
\omega_1 \rangle$, to make sure $M \in {\frak K}^+_1$ it is enough that for
unboundedly many $\alpha < \omega_1,M_\alpha <^{**} M_{\alpha +1}$
where
\sn
\item "{$\circledast_{11}$}"  for $M,N \in {\frak K},M <^{**} N$ iff
{\roster
\itemitem{ $(i)$ }  $M \le_{\frak K} N$
\sn
\itemitem{ $(ii)$ }  in $\circledast_7(b)$ also the inverse direction
holds.
\endroster}
\ermn
It should be clear by now that we have proved clauses (a),(b),(d),(e)
of \scite{600-Ex.1A} using ${\frak K}$.   
We have to prove clause (c); we cannot quote
\scite{600-Ex.1} as clause $(\gamma)$ there is only almost true.  The proof is
similar to (but simpler than) that of \scite{600-Ex.1} quoting
\cite{Sh:48} instead of \chaptercite{88r}; a marked difference is that in
the present case the number of types over a countable model is
countable (in ${\frak K}$) whereas in \chaptercite{88r} it seemingly could
have been $\aleph_1$, generally \cite{Sh:48} situation is more similar
to the first order logic case.

Does ${\frak K}$ have amalgamation in $\aleph_0$?
Now (\cite[Lemma 3.4]{Sh:48}, almost says this but it assumed
$\diamondsuit_{\aleph_1}$ instead of $2^{\aleph_0} < 2^{\aleph_1}$;
and \marginbf{!!}{\cprefix{88r}.\scite{88r-3.5}} almost says this, but the models are 
from ${\frak K}_{\aleph_1}$ rather than ${\frak K}^+_{\aleph_1}$ 
but \marginbf{!!}{\cprefix{88r}.\scite{88r-3.8.4}} fully says this using the so called
$K^F_{\aleph_1}$).  So
\mr
\item "{$\circledast_{12}$}"   ${\frak K}$ has the amalgamation
property in $\aleph_0$.
\ermn
Recall that all models from ${\frak K}$ are atomic (in the
first order sense) and we shall use below $\sftp_{\Bbb L}$.

As ${\frak K}$ has $\aleph_0$-amalgamation (by $\circledast_{12}$), clearly
\cite[\S4]{Sh:48} applies; now by \cite[Lemma 2.1]{Sh:48}(B) 
+ Definition 3.5, being $(\aleph_0,1)$-stable as defined in
\cite[Definition 3.5]{Sh:48}(A) holds.  Hence all clauses of
\cite[Lemma 4.2]{Sh:48} hold, in particular ($(D)(\beta)$ there and
clause (A), i.e., \cite[Def.3.5]{Sh:48}(B)), so
\mr
\item "{$\circledast_{13}$}"  $(i) \quad$ if 
$M \le_{\frak K} N$ and $\bar a \in N$ then
$\sftp_{\Bbb L}(\bar a,M,N)$ is definable over a finite subset of $M$
\sn
\item "{${{}}$}"  $(ii) \quad$ if $M \in {\frak K}_{\aleph_0}$ then
$\sftp_{\Bbb L}(\bar a,M,N):\bar a \in {}^{\omega >} N$ and
$M \le_{\frak K} N\}$ is countable.
\ermn
By \cite[Lemma 4.4]{Sh:48} it follows that
\mr
\item "{$\circledast_{14}$}"  if $M \le_{\frak K} N$ are countable and
$\bar a \in M$ then $\sftp_{\Bbb L}(\bar a,M,N)$ determine 
$\sftp_{\frak K}(\bar a,M,N)$.
\ermn
Now we define ${\frak s} = ({\frak K}_{\aleph_0},{\Cal
S}^{\text{bs}},\nonfork{}{}_{})$ by
\mr
\item "{$\circledast_{15}$}"  ${\Cal S}^{\text{bs}}(M) 
= \{\text{\ortp}_{\frak K}(\bar a,M,N):M \le_{\frak K} N$ are countable 
and $\bar a \in {}^{\omega >} N$ but $\bar a \notin {}^{\omega >}M\}$
\sn
\item "{$\circledast_{16}$}"  \ortp$_{\frak K}(\bar a,M_1,M_3)$ does not
fork over $M_0$ where $M_0 \le_{\frak K} M_1 \le_{\frak K} M_3 \in
{\frak K}_{\aleph_0}$ iff $\sftp_{\Bbb L}(\bar a,M_1,M_3)$ is definable
over some finite subset of $M_0$.
\ermn
Now we check ``${\frak s}$ is a good frame", i.e., all clauses of
Definition \scite{600-1.1}.   
\enddemo
\bn
\ub{Clause (A)}:  By $\circledast_9$ above.
\bn
\ub{Clause (B)}:  As ${\frak K}$ is categorical in $\aleph_0$, has an
uncountable model and LS$({\frak K}) = \aleph_0$ this should be clear.
\bn
\ub{Clause (C)}:  ${\frak K}_{\aleph_0}$ has amalgamation by
$\circledast_{12}$ and has the JEP by categoricity in $\aleph_0$ and
${\frak K}_{\aleph_0}$ has no maximal model by (categoricity and)
having uncountable models (and LS$({\frak K}) = \aleph_0$).
\bn
\ub{Clause (D)}:  Obvious; stability, i.e., (D)(d) holds by
$\circledast_{13}(ii) + \circledast_{14}$.
\bn
\ub{Subclause (E)(a),(b)}:  By the definition.
\bn
\ub{Subclause (E)(c)}:  (Local character).

If $\langle M_i:i \le \delta +1 \rangle$ is $\le_{\frak K}$-increasing
continuous $M_i \in K_{\aleph_0},\bar a \in {}^{\omega >}
(M_{\delta +1})$ and $\bar a \in {}^{\omega >}(M_\delta)$ 
then for some finite $A \subseteq M_\delta,
\sftp_{\Bbb L}(\bar a,M_\delta,M_{\delta +1})$ 
is definable over $A$, so for some $i <
\delta,A \subseteq M_\delta$ hence $j \in [i,\delta) \Rightarrow
\sftp_{\Bbb L}(\bar a,M_i,M_{\delta +1})$ is definable over $A
\Rightarrow \nonfork{}{}_{}(M_i,M_\delta,\bar a,M_{\delta +i})$.
\bn
\ub{Subclause (E)(d)}:  (Transitivity).

As if $M' \le_{\frak K} M'' \in {\frak K}_{\aleph_0}$, two definitions
in $M'$ of complete types, which give the same result in $M'$ give the
same result in $M''$.
\bn
\ub{Sublause (E)(e)}(uniqueness):  By $\circledast_{14}$ and the
justification of transitivity.  
\bn
\ub{Subclause (E)(f)}(symmetry):  By 
\cite[Theorem 5.4]{Sh:48}, we have the symmetry
property see \cite[Definition 5.2]{Sh:48}.  By \cite[5.5]{Sh:48} + the
uniqueness proved above we can finish easily.
\bn
\ub{Subclause (E)(g)}:  Extension existence.

Easy, included in \cite[5.5]{Sh:48}.
\bn
\ub{Subclause (E)(h)}:  Continuity.

As ${\Cal S}^{\text{bs}}_{\frak s}(M)$ is the set of
non-algebraic types this follows from ``finite character", that is by
\scite{600-1.16A}(3)(4). 
\bn
\ub{Subclause (E)(i)}: non-forking amalgamation

By \scite{600-1.15}.  \hfill$\square_{\scite{600-Ex.1A}}$ 
\bigskip

\remark{\stag{600-Ex.1B} Remark}  So if $\psi \in \Bbb L_{\omega_1,\omega}
(\bold Q)$ and
$1 \le \dot I(\aleph_1,\psi) < 2^{\aleph_1}$, we essentially can apply 
Theorem \scite{600-0.A}, exactly see \scite{600-fc.4}. 
\endremark
\bigskip

\subhead {(D) Starting at $\lambda > \aleph_0$} \endsubhead

The next theorem puts the results of \cite{Sh:576} in our context
hence rely on it heavily.
\proclaim{\stag{600-Ex.4} Theorem}  Assume 
$2^\lambda < 2^{\lambda^+} < 2^{\lambda^{++}}$
and
\medskip
\roster
\item "{$(\alpha)$}"  ${\frak K}$ is an abstract elementary class with
${\text{\rm LS\/}}({\frak K}) \le \lambda$
\sn
\item "{$(\beta)$}"  ${\frak K}$ is categorical in $\lambda$ and in
$\lambda^+$
\sn
\item "{$(\gamma)$}"  ${\frak K}$ has a model in $\lambda^{++}$ 
\sn
\item "{$(\delta)$}"  $\dot I(\lambda^{+2},K) < 2^{\lambda^{++}}$ and
${\text{\rm WDmId\/}}(\lambda^+)$ is not $\lambda^{++}$-saturated \ub{or}
just some \footnote{alternatively we, in clause (d) assume just $\dot
I(\lambda^{++},K) < \mu_{\text{unif}}(\lambda^{++},2^{\lambda^+})$, see
\cite{Sh:838}, see on this cardinal in \marginbf{!!}{\cprefix{88r}.\scite{88r-0.wD}}(3).} 
consequences: density of minimal types (see
\cite{Sh:603}) and $\otimes$ of \cite[6.4,pg.99]{Sh:576} proved by
the conclusion of \cite[Th.6.7]{Sh:576}(pg.101).
(so we can use all the results of \cite[\S8-\S10]{Sh:576})
\ermn
\ub{Then} 1) Letting $\mu = \lambda^+$ we can choose 
$\nonfork{}{}_{\mu},{\Cal S}^{\text{bs}}$ such that $({\frak K}_{\ge \mu},
\nonfork{}{}_{\mu},{\Cal S}^{\text{bs}})$ is a $\mu$-good frame.
\nl
2) Moreover, we can let
\medskip
\roster
\item "{$(a)$}"  ${\Cal S}^{\text{bs}}(M) := 
\{ {\text{\rm \ortp\/}}_{\frak K}(a,M,N):
\text{for some } M,N,a \text{ we have }(M,N,a) \in K^3_{\lambda^+}$ \nl

$\qquad \qquad \qquad \qquad \qquad$ and for some 
$M' \le_{\frak K} M \text{ we have } M' \in K_\lambda$ \nl

$\qquad \qquad \qquad \qquad \qquad$ and 
${\text{\rm \ortp\/}}_{\frak K}(a,M',N) \in {\Cal S}_{\frak K}(M')$ is minimal$\}$
\endroster
\medskip
\noindent
(see Definition \cite[2.3]{Sh:576}(4)(pg.56) and
\cite[2.5]{Sh:576}(1),(13),pg.57-58 
\mr
\item "{$(b)$}"   $\nonfork{}{}_{} = \nonfork{}{}_{\mu}$ be defined 
by: $\nonfork{}{}_{}(M_0,M_1,a,M_3)$ \underbar{iff} 
$M_0 \le_{\frak K} M_1 \le_{\frak K} M_3$ are from $K_\mu,a \in M_3 
\backslash M_1$ and for some $N \le_{\frak K} M_0$ of cardinality $\lambda$,
the type ${\text{\rm \ortp\/}}_{\frak K}(a,N,M_3) \in 
{\Cal S}_{\frak K}(N)$ is minimal.
\endroster
\endproclaim
\bigskip

\demo{Proof}  1), 2)  By clause $(\delta)$ of the assumption, we can
use the ``positive" results of \cite{Sh:576} in particular 
\cite[\S8-\S10]{Sh:576} freely.
Also by \cite{Sh:603}, the one point in \cite{Sh:576} in which we have used 
$\lambda \ne \aleph_0$ is eliminated from this assumption.  Now (see
Definition \scite{600-0.13}(2))
\sn
\mr
\item "{$(*)$}"  if $(M,N,a) \in K^3_{\lambda^+}$ and $M' \le_{\frak K} M,
M' \in K_\lambda$ and $p = \text{ \ortp}_{\frak K}(a,M',N)$ is minimal
(see Definition \scite{600-0.12}(0)) \ub{then}
{\roster
\itemitem{ $(a)$ }  if $q \in {\Cal S}_{\frak K}(M)$ is not algebraic and
$q \restriction M' = p$ then $q = \text{\rm \ortp}_{\frak K}(a,M,N)$
\sn
\itemitem{ $(b)$ }  if $\langle M_\alpha:\alpha < \mu \rangle,\langle
N_\alpha:\alpha < \mu \rangle$ are representations of $M,N$ respectively
\ub{then} for a club of $\delta < \mu$ we have \ortp$_{\frak K}
(a,M_\delta,N_\delta) \in {\Cal S}_{\frak K}
(M_\delta)$ is minimal and reduced
\endroster}
[Why?  For clause (b) let $\alpha^* = \text{ Min}\{\alpha:M' \le_{\frak K}
M_\alpha\}$, so $\alpha^*$ is well defined and as $M$ is saturated
(for ${\frak K}$), for a club of $\delta < \mu = \lambda^+$, the model 
$M_\delta$ is $(\lambda,\text{cf}(\delta))$-brimmed 
over $M'$ hence by \cite[7.5]{Sh:576}(2)(pg.106)
we are done.

For clause (a) let $M^0 = M,M^1 = N$ and $a^1=a$ and $M^2,a^2 = a$ 
be such that $(M^0,M^2,a^2) \in K^3_\mu = K^3_{\lambda^+}$ 
and $q = \text{ \ortp}_{\frak K}(a^2,M^0,M)$.  Now we repeat
the proof of \cite[9.5]{Sh:576}(pg.120)
but instead $f(a^2) \notin M^1$ we require
$f(a^2) = a^1$; we are using \cite[10.5]{Sh:576}(1)(pg.125)
which says $<^*_{\lambda^+} = <_{\frak K} \restriction K_{\lambda^+}$.]
\endroster
\bn
In particular we have used
\mr
\item "{$(**)$}"  if $M_0 \le_{{\frak K}_\lambda} M_1,M_1$ is
$(\lambda,\kappa)$-brimmed over $M_0,p \in {\Cal S}_{\frak K}(M_1)$ is not
algebraic and $p \restriction M_0$ is minimal, \ub{then} $p$ is
minimal and reduced.
\endroster
\bn
\underbar{Clause $(A)$}:

This is by assumption $(\alpha)$. 
\bigskip
\noindent
\underbar{Clause $(B)$}:

As $K$ is categorical in $\mu = \lambda^+$, the existence of superlimit
$M \in K_\mu$ follows; the superlimit is not maximal as LS$({\frak K})
\le \lambda \and K_{\mu^+} = K_{\lambda^{++}} \ne \emptyset$ by
assumption $(\gamma)$.
\bigskip
\noindent
\underbar{Clause $(C)$}:

$K_{\lambda^+}$ has the amalgamation property by \marginbf{!!}{\cprefix{88r}.\scite{88r-3.5}} or
\cite[1.4]{Sh:576}(pg.46),1.6(pg.48) and
${\frak K}_\lambda$ has the JEP in
$\lambda^+$ by categoricity in $\lambda^+$.
\bigskip
\noindent
\underbar{Clause $(D)$}: \newline
\underbar{Subclause $(D)(a),(b)$}:

By the definition of ${\Cal S}^{\text{bs}}(M)$ and of minimal types (in
${\Cal S}_{\frak K}(N),N \in K_\lambda$, \nl
\cite[2.5]{Sh:576}(1)+(3)(pg.57),2.3(4)+(6)](pg.56)), 
this is clear.
\bigskip
\noindent
\underbar{Subclause $(D)(c)$}:

Suppose $M \le_{\frak K} N$ are from $K_\mu$ and $M \ne N$; let
$\langle M_i:i < \lambda^+ \rangle, \langle N_i:i < \lambda^+ \rangle$ be a
$\le_{\frak K}$-representation 
of $M,N$ respectively, choose $b \in N \backslash M$
so $E = \{\delta < \lambda^+:N_\delta \cap M = M_\delta \text{ and } b \in
N_\delta\}$ is a club of $\lambda^+$.  Now for $\delta = \text{ Min}(E)$ we
have $M_\delta \ne N_\delta,M_\delta \le_{\frak K} N_\delta$ and there 
is a minimal inevitable $p \in {\Cal S}_{\frak K}(M_\delta)$ by 
\cite[5.3,pg.94]{Sh:576}
and categoricity of $K$ in $\lambda$; so for some 
$a \in N_\delta \backslash M_\delta$ we have $p = \text{ \ortp}_{\frak K}
(a,M_\delta,N_\delta)$.  So \ortp$_{\frak K}(a,M,N)$ is non-algebraic as 
$a \in M \Rightarrow a \in M \cap
N_\delta = M_\delta$, a contradiction, so \ortp$_{\frak K}(a,M,N) \in 
{\Cal S}^{\text{bs}}(M)$ as required.
\bn
\ub{Subclause $(D)(d)$}:  If $M \in K_\mu$ let $\langle M_i:i < \lambda^+
\rangle$ be a $\le_{\frak K}$-representation of $M$, so by $(*)(a)$ above 
$p \in {\Cal S}^{\text{bs}}(M)$ is
determined by $p \restriction M_\alpha$ if $p \restriction M_\alpha$ is
minimal and reduced.  
But for every such $p$ there is such $\alpha(p) < \lambda^+$ by the
definition of ${\Cal S}^{\text{bs}}(M)$ and for each 
$\alpha < \lambda^+$ there are $\le \lambda$ possible such $p \restriction
M_\alpha$ as ${\frak K}$ is stable in $\lambda$ by
\cite[5.7]{Sh:576}(a)(pg.97), so the  
conclusion follows.  Alternatively,
$M \in K_\mu \Rightarrow |{\Cal S}^{\text{bs}}(M)| \le
\mu$ as by \cite[10.5]{Sh:576}(pg.125), we have $\le^*_{\lambda^+} = 
\le_{\frak K} \restriction K_{\lambda^+}$, so we can apply
\cite[9.7]{Sh:576}(pg.121); or use $(*)$ above.
\bn
\underbar{Clause $(E)$}: \newline
\underbar{Subclause $(E)(a)$}:

Follows by the definition.
\bn
\underbar{Subclause $(E)(b)$}: (Monotonicity)

Obvious properties of minimal types in ${\Cal S}(M)$ for $M \in K_\lambda$.
\bn
\underbar{Subclause $(E)(c)$}: (Local character)

Let $\delta < \mu^+ = 
\lambda^{++}$ and $M_i \in K_\mu$ be $\le_{\frak K}$-increasing continuous 
for $i \le \delta$ and $p \in {\Cal S}^{\text{bs}}
(M_\delta)$, so for some $N \le_{\frak K} M_\delta$ we have $N \in 
K_\lambda$ and $p \restriction N \in {\Cal S}_{\frak K}(N)$ is 
minimal.  Without loss of generality 
$\delta = \text{ cf}(\delta)$ and if $\delta = \lambda^+$, there is 
$i < \delta$ such
that $N \subseteq M_i$ and easily we are done.  
So assume $\delta = \text{ cf}
(\delta) < \lambda^+$. \nl
Let $\langle M^i_\zeta:\zeta < \lambda^+ \rangle$ be a $\le_{\frak K}$-representation of
$M_i$ for $i \le \delta$, hence $E$ is a club of $\lambda^+$ where: 

$$
\align
E := \bigl\{ \zeta < \lambda^+:&\zeta \text{ a limit ordinal and for } 
j < i \le \delta \text{ we have} \\
  &M^i_\zeta \cap M_j = M^j_\zeta \text{ and for }
\xi < \zeta,i \le \delta \text{ we have}: \\
  &M^i_\zeta \text{ is } (\lambda,\text{cf}(\zeta))
\text{-brimmed over } M^i_\xi \text{ and } N \le_{\frak K}
M^\delta_\zeta \bigr\}.
\endalign
$$
\mn
Let $\zeta_i$ be the $i$-th member of $E$ for $i \le \delta$, so
$\langle \zeta_i:i \le \delta \rangle$ is increasing continuous,
$\langle M^i_{\zeta_i}:i \le \delta \rangle$ is $\le_{\frak K}$-increasingly
continuous in $K_\lambda$ and $M^{i+1}_{\zeta_{i+1}}$ is $(\lambda,
\text{cf}(\zeta_{i+1}))$-brimmed over $M^{i+1}_{\zeta_i}$ hence also
over $M^i_{\zeta_i}$.  Also $p \restriction M^\delta_{\zeta_\delta}$
is non-algebraic (as $p$ is) and extends 
$p \restriction N$ (as $N \le_{\frak K} 
M^\delta_{\zeta_\delta}$ as $\zeta_\delta 
\in E$) hence $p \restriction M^\delta_{\zeta_\delta}$ is minimal.

Also $M^\delta_{\zeta_\delta}$ is 
$(\lambda,\text{cf}(\zeta_\delta))$-brimmed
over $M^\delta_{\zeta_0}$ hence over $N$, hence by $(**)$ above we get
that $p \restriction M^\delta_{\zeta_\delta}$ is not only minimal but also
reduced.  Hence by \cite[7.3]{Sh:576}(2)(pg.103) applied to $\langle
M^i_{\zeta_i}:i \le \delta \rangle,p \restriction M^\delta_{\zeta_\delta}$
we know that for some $i < \delta$ the type $p \restriction M^i_{\zeta_i} =
(p \restriction M^\delta_{\zeta_\delta}) \restriction M^i_{\zeta_i}$ is 
minimal and reduced, so it witnesses that 
$p \restriction M_j \in {\Cal S}^{\text{bs}}(M_j)$ for every 
$j \in [i,\delta)$, as required.
\bigskip

\noindent
\underbar{Subclause $(E)(d)$}: (Transitivity)

Easy by the definition of minimal.
\bigskip

\noindent
\underbar{Subclause $(E)(e)$}: (Uniqueness)

By $(*)(a)$ above.
\bn
\underbar{Subclause $(E)(f)$}: (Symmetry)

By the symmetry in the situation assume $M_0 \le_{\frak K} M_1 \le_{\frak K}
M_3$ are from $K_\mu$, \newline
$a_1 \in M_1 \backslash M_0,a_2 \in M_3 \backslash
M_1$ and \ortp$_{\frak K}(a_1,M_0,M_3) \in {\Cal S}^{\text{bs}}(M_0)$ 
and \ortp$_{\frak K}(a_2,M_1,M_3)
\in {\Cal S}^{\text{bs}}(M_1)$ does not fork over $M_0$; hence for $\ell =
1,2$ we have \ortp$_{\frak K}(a_\ell,M_0,M_3) \in {\Cal S}^{\text{bs}}(M_0)$.  
By the existence of disjoint
amalgamation (by \cite[9.11]{Sh:576}(pg.122),10.5(1)(pg.125))
there are $M_2,M'_3,f$ such that
$M_0 \le_{\frak K} M_2 \le_{\frak K} M'_3 \in K_\mu$,
$M_3 \le_{\frak K} M'_3,f$ is an isomorphism from $M_3$ onto $M_2$
over $M_0$, and $M_3 \cap M_2 = M_0$.  By \ortp$_{\frak K}(a_2,M_0,M_3) \in {\Cal
S}^{\text{bs}}(M_1)$ and as $f(a_2) \notin M_1$ being in $M_2 
\backslash M_0 = M_2 \backslash M_3$ and $a_2 \notin M_1$ by
assumption and as $a_2,f(a_2)$ realize the same type from ${\Cal
S}_{\frak K}(M_0)$ clearly by $(*)(a)$
we have \ortp$_{\frak K}(a_2,M_1,M'_3) = \text{ \ortp}_{\frak K}(f(a_2),M_1,M'_3)$.
\medskip

Using amalgamation in ${\frak K}_\mu$ (and equality of types) 
there is $M''_3$ such that: \newline
$M'_3 \le_{\frak K} M''_3 \in K_\mu$, and there is an 
$\le_{\frak K}$-embedding $g$ of $M'_3$ into $M''_3$ such that 
$g \restriction M_1 = \text{ id}_{M_1}$ and $g(f(a_2)) = a_2$.  
Note that as $a_1 \notin g(M_2),M_1 \le_{\frak K} g(M_2) \in K_\mu$
and \ortp$_{\frak K}(a_1,M_1,M''_3)$ is minimal then necessarily
\ortp$_{\frak K}(a_1,g(M_2),M''_3)$ is its non-forking extension.
So $g(M_2),M''_3$ are models as required. 
\bn
\underbar{Subclause $(E)(g)$}:  (Extension existence)

Claims \cite[9.11]{Sh:576}(pg.122),10.5(1)(pg.125) do even more.
\bn
\underbar{Subclause $(E)(h)$}:  (Continuity)

Easy.
\bn
\ub{Subclause $(E)(i)$}:  (Non-forking amalgamation)

Like $(E)(f)$ or use \scite{600-1.15}.  \hfill$\square_{\scite{600-Ex.4}}$
\enddemo
\bigskip
\noindent
\margintag{600-Ex.5}\underbar{\stag{600-Ex.5} Question}:  If ${\frak K}$ is categorical in $\lambda$ 
and in $\mu$ and $\mu > \lambda \ge \text{LS}({\frak K})$, can we 
conclude categoricity in $\chi \in (\mu,\lambda)$?
\bigskip

\demo{\stag{600-Ex.6} Fact}  In \scite{600-Ex.4}: \nl
1)  If $p \in {\Cal S}^{\text{bs}}(M)$
and $M \in K_\mu$, \ub{then} for some $N \le_{\frak K} M,N \in
K_\lambda$ and $p \restriction N$ is minimal and reduced. \nl
2) If $M <_{\frak K} N,M \in K_\mu$ and $p \in {\Cal
S}^{\text{bs}}(M)$, \ub{then} some $a \in N \backslash M$ realizes $p$, 
(i.e., ``a strong version of uni-dimensionality" holds).
\enddemo
\bigskip

\demo{Proof}  The proof is included in the proof of \scite{600-Ex.4}.
\enddemo
\bn
\centerline{$* \qquad * \qquad *$}
\bn
\ub{(E) An Example}:

A trivial example (of an approximation to good $\lambda$-frame) is:
\definition{\stag{600-Ex.7.1} Definition/Claim}  1) Assume that ${\frak K}$
is an a.e.c. and $\lambda \ge \text{\rm LS}({\frak K})$
or ${\frak K}$ is a $\lambda$-a.e.c.
We define ${\frak s} = {\frak s}_\lambda[{\frak K}]$ as the triple
${\frak s} = ({\frak K}_\lambda,{\Cal
S}^{\text{na}},\nonfork{}{}_{\text{na}})$ where:
\mr
\item "{$(a)$}"  ${\Cal S}^{\text{na}}(M) = \{\text{ \ortp}_{\frak
K}(a,M,N),M \le_{\frak K} N$ and $a \in N \backslash M\}$
\sn
\item "{$(b)$}"  $\nonfork{}{}_{\text{na}}(M_0,M_1,a,M_3)$ iff $M_0
\le_{{\frak K}_\lambda} M_1 \le_{{\frak K}_\lambda} M_3$ and $a \in
M_3 \backslash M_1$.
\ermn
2) Then ${\frak s}$ satisfies Definition \scite{600-1.1} of good
$\lambda$-frame except possibly: (B), existence of superlimits, (C)
amalgamation and JEP, (D)(d) stability and (E)(e),(f),(g),(i) uniqueness,
symmetry, extension existence and non-forking amalgamation. 
\enddefinition
\newpage

\head {\S4 Inside the framework} \endhead  \resetall \sectno=4
 \spuriousreset
\bigskip

We investigate good $\lambda$-frames.  We prove stability in $\lambda$
(we have assumed in Definition \scite{600-1.1} only stability for basic
types), hence the existence of a $(\lambda,\sigma)$-brimmed
$\le_{\frak K}$-extension in $K_\lambda$ over $M_0 \in K_\lambda$ (see
\scite{600-4a.1}), and we give a sufficient condition for ``$M_\delta$ is
$(\lambda,\text{cf}(\delta))$-brimmed over $M_0$" (in \scite{600-4a.2}).
We define again $K^{3,\text{bs}}_\lambda$ (like $K^3_\lambda$ from
\scite{600-0.13}(2) but the type is
basic) and the natural order $\le_{\text{bs}}$ on them as well as ``reduced"
(Definition \scite{600-4a.3}), 
and indicate their basic properties (\scite{600-4a.4}). \nl
We may like to construct sometimes pairs $N_i \le_{{\frak K}_\lambda}
M_i$ such that 
$M_i,N_i$ are increasing continuous with $i$ and we would like to guarantee
that $M_\gamma$ is $(\lambda,\text{cf}(\gamma))$-brimmed over
$N_\gamma$, of course we need to carry more inductive assumptions.
Toward this we may give a sufficient condition for building a
$(\lambda,\text{cf}(\gamma))$-brimmed extension over $N_\gamma$
where $\langle N_i:i \le \gamma \rangle$ is 
$\le_{{\frak K}_\lambda}$-increasing continuous, 
by a triangle of extensions of the $N_i$'s, with 
non-forking demands of course (see \scite{600-4a.5}).  We
also give conditions on a rectangle of models to get such pairs in
both directions (\scite{600-4a.9}), for this we use nice extensions of
chains (\scite{600-4a.7}, \scite{600-4a.8}).

Then we can deduce that if ``$M_1$ is $(\lambda,\sigma)$-brimmed
over $M_0$" then the isomorphism type of $M_1$ over $M_0$ does not
depend on $\sigma$ (see \scite{600-4a.6}), so the brimmed $N$ over $M_0$
is unique up to isomorphism (i.e. being $(\lambda,\sigma)$-brimmed
over $M_0$ does not depend on $\sigma$).  We finish giving conclusion
about $K_{\lambda^+},K_{\lambda^{++}}$.
\bigskip

\demo{\stag{600-4a.0} Hypothesis}  ${\frak s} = ({\frak K},
\nonfork{}{}_{},{\Cal S}^{\text{bs}})$ is a good $\lambda$-frame.
\enddemo
\bigskip

\proclaim{\stag{600-4a.1} Claim}  1)  ${\frak K}$ is stable in $\lambda$,
i.e., $M \in {\frak K}_\lambda \Rightarrow |{\Cal S}(M)| \le \lambda$.
 \nl
2) For every $M_0 \in K_\lambda$ and $\sigma \le \lambda$
there is $M_1$ such that $M_0 \le_{\frak K} M_1 \in K_\lambda$ and $M_1$ 
is $(\lambda,\sigma$)-brimmed over $M_0$ (see Definition \scite{600-0.21})
and it is universal
\footnote{in fact, this follows} over $M_0$.
\endproclaim
\bigskip

\demo{Proof}  1)  Let 
$M_0 \in K_\lambda$ and we choose by induction on
$\alpha \in [1,\lambda],M_\alpha \in K_\lambda$ such that:
\mr
\item "{$(i)$}"  $M_\alpha$ is $\le_{\frak K}$-increasing continuous
\sn
\item "{$(ii)$}"   if $p \in {\Cal S}^{\text{bs}}(M_\alpha)$ then 
this type is realized in $M_{\alpha +1}$.
\ermn
No problem to carry this:  for clause (i) use $Ax(A)$, for clause (ii) 
use Axiom $(D)(d)$ and amalgamation in ${\frak K}_\lambda$, i.e., Axiom (C).
If every $q \in {\Cal S}(M_0)$ is realized in
$M_\lambda$ we are done.  So let $q$ be a counterexample, so let 
$M_0 \le_{\frak K} N \in K_\lambda$ be such that $q$ is
realized in $N$.  We now try to choose by induction on $\alpha < \lambda$
a triple $(N_\alpha,f_\alpha,\bar{\bold a}_\alpha)$ such that:
\mr
\item "{$(A)$}"  $N_\alpha \in K_\lambda$ is $\le_{\frak K}$-increasingly
continuous
\sn
\item "{$(B)$}"  $f_\alpha$ is a $\le_{\frak K}$-embedding of
$M_\alpha$ into $N_\alpha$
\sn
\item "{$(C)$}"  $f_\alpha$ is increasing continuous
\sn
\item "{$(D)$}"  $f_0 = \text{ id}_{M_0}$ and $N_0 = N$
\sn
\item "{$(E)$}"  $\bar{\bold a}_\alpha = \langle a_{\alpha,i}:i < \lambda
\rangle$ lists the elements of $N_\alpha$
\sn
\item "{$(F)$}"  if there are $\beta \le \alpha,i < \lambda$ such that
\ortp$(a_{\beta,i},f_\alpha(M_\alpha),N_\alpha) \in {\Cal S}^{\text{bs}}
(f_\alpha(M_\alpha))$ \ub{then} for some such pair $(\beta_\alpha,i_\alpha)$ 
we have:
{\roster
\itemitem{ $(i)$ }  the pair $(\beta_\alpha,i_\alpha)$ is minimal in
an appropriate sense, that is: 
if $(\beta,i)$ is another such pair then $\beta + i >
\beta_\alpha + i_\alpha$ or $\beta +i = \beta_\alpha + i_\alpha \and \beta
> \beta_\alpha$ or $\beta+i = \beta_\alpha + i_\alpha \and \beta =
\beta_\alpha \and i \ge i_\alpha$
\sn
\itemitem{ $(ii)$ }  $a_{\beta_\alpha,i_\alpha} 
\in \text{ Rang}(f_{\alpha+1})$.
\endroster}
\ermn
This is easy: for successor $\alpha$ we use the definition of type and
let $N_\lambda := \cup\{N_\alpha:\alpha < \lambda\}$. 
Clearly $f_\lambda := \cup\{f_\alpha:\alpha < \lambda\}$ is a
$\le_{\frak s}$-embedding of $M_\lambda$ into $N_\lambda$ over $M_0$.

As in $N$, the type $q$ is realized and it is not realized in
$M_\lambda$ necessarily $N \nsubseteq f_\lambda(M_\lambda)$ hence
$N_\lambda  \ne f_\lambda(M_\lambda)$ but easily $f_\lambda(M_\lambda)
\le_{\frak K} N_\lambda$.  So by
Axiom $(D)(c)$ for some $c \in N_\lambda \backslash f_\lambda
(M_\lambda)$ we have $p = \text{ \ortp}(c,f_\lambda(M_\lambda),N_\lambda) 
\in {\Cal S}^{\text{bs}}(f_\lambda(M_\lambda))$.  As $\langle
f_\gamma(M_\gamma):\gamma \le \lambda \rangle$ is $\le_{\frak
K}$-increasing continuous, by 
Axiom (E)(c) for some $\gamma < \lambda$ we have
\ortp$(c,f_\lambda(M_\lambda),N_\lambda)$ does not fork over 
$f_\gamma(M_\gamma)$, also as $c \in N_\lambda = 
\dbcu_{\beta < \lambda} N_\beta$ clearly $c \in N_\beta$ for some
$\beta < \lambda$ and let $i < \lambda$ be such that 
$c = a_{\beta,i}$.  Now if $\alpha \in
[\text{max}\{\gamma,\beta\},\lambda)$ then $(\beta,i)$ is a legitimate
candidate for $(\beta_\alpha,i_\alpha)$ that is
\ortp$(a_{\beta,i},f_\alpha(M_\alpha),N_\alpha) \in {\Cal
S}^{\text{bs}}(f_\alpha(M_\alpha))$ by monotonicity of non-forking,
i.e., Axiom (E)(b).  So $(\beta_\alpha,i_\alpha)$ is well defined for
any such $\alpha$ and
$\beta_\alpha + i_\alpha \le \beta+i$ by clause (F)(i).  But $\alpha_1 <
\alpha_2 \Rightarrow a_{\beta_{\alpha_1},i_{\alpha_1}} \ne
a_{\beta_{\alpha_2},i_{\alpha_2}}$ (as one belongs to $f_{\alpha_1 +1}
(M_{\alpha_1})$ 
and the other not), contradiction by cardinality consideration. \nl
2) So ${\frak K}_\lambda$ is stable in $\lambda$ and has amalgamation, hence
(see \scite{600-0.22}) the conclusion holds; alternatively use
\scite{600-4a.2} below.
\hfill$\square_{\scite{600-4a.1}}$ 
\enddemo
\bigskip

\proclaim{\stag{600-4a.2} Claim}  Assume 
\mr
\item "{$(a)$}"  $\delta < \lambda^+$ is a 
limit ordinal divisible by $\lambda$
\sn
\item "{$(b)$}"  $\bar M = \langle M_\alpha:\alpha \le \delta \rangle$ is
$\le_{\frak K}$-increasing continuous sequence in ${\frak K}_\lambda$
\sn
\item "{$(c)$}"  if $i < \delta$ and 
$p \in {\Cal S}^{\text{bs}}(M_i)$, \ub{then}
for $\lambda$ ordinals $j \in (i,\delta)$ there is $c' \in M_{j+1}$
realizing the non-forking extension of $p$ in ${\Cal S}^{\text{bs}}(M_j)$.
\ermn
\ub{Then} $M_\delta$ is $(\lambda,{\text{\rm cf\/}}
(\delta))$-brimmed over $M_0$ and universal over it.
\endproclaim
\bigskip

\remark{\stag{600-4a.2A} Remark}  1)  See end of proof of \scite{600-nf.17}. \nl
2) Of course, by renaming, $M_\delta$ is $(\lambda,\text{cf}
(\delta))$-brimmed over $M_\alpha$ for any $\alpha < \delta$. \nl
3) Why in clause (c) of \scite{600-4a.2} we ask for ``$\lambda$ ordinals $j \in
(i,\delta)$" rather than ``for unboundedly many $j \in (i,\delta)$"?
For $\lambda$ regular there is no difference but for $\lambda$
singular not so.
Think of ${\frak K}$ the class of $(A,R),R$ an equivalence relation on
$A$; (so it is not categorical) but for some $\lambda$-good frames
${\frak s},{\frak K}_{\frak s} = {\frak K}_\lambda$ and 
exemplifies a problem; some equivalence class of $M_\delta$ may be of
cardinality $< \lambda$. 
\endremark
\bigskip

\demo{Proof}  Like \scite{600-4a.1}, but we give details.

Let $g:\delta \rightarrow \lambda$ be a one to
one and choose by induction on $\alpha \le \delta$ a triple $(N_\alpha,
f_\alpha,\bar{\bold a}_\alpha)$ such that
\mr
\item "{$(A)$}"  $N_\alpha \in K_\lambda$ is $\le_{\frak K}$-increasing
continuous
\sn
\item "{$(B)$}"  $f_\alpha$ is a $\le_{\frak K}$-embedding of
$M_\alpha$ into $N_\alpha$
\sn
\item "{$(C)$}"  $f_\alpha$ is increasing continuous
\sn
\item "{$(D)$}"  $f_0 = \text{ id}_{M_0},N_0 = M_0$
\sn
\item "{$(E)$}"  $\bar{\bold a}_\alpha = \langle a_{\alpha,i}:i < \lambda
\rangle$ list the elements of $N_\alpha$
\sn
\item "{$(F)$}"  $N_{\alpha +1}$ is universal over $N_\alpha$
\sn
\item "{$(G)$}"  if $\alpha < \delta$ and
there is a pair $(\beta,i) = (\beta_\alpha,i_\alpha)$
satisfying the condition $(*)^{\beta,i}_{f_\alpha,N_\alpha}$ stated below
and it is minimal in the sense that \nl
$(*)^{\beta',i'}_{f_\alpha,N_\alpha} \Rightarrow 
(**)^{\beta',i',\beta,i}_g$, see below,
\ub{then} $a_{\beta,i} \in \text{ Rang}(f_{\alpha +1})$, \nl
where
{\roster
\itemitem{ $(*)^{\beta,i}_{f_\alpha,N_\alpha}$ }  $(a) \quad \beta \le
\alpha$ and $i < \lambda$
\sn
\itemitem{ {} } $(b) \quad$ \ortp$(a_{\beta,i},f_\alpha(M_\alpha),N_\alpha) 
\in {\Cal S}^{\text{bs}}(f_\alpha(M_\alpha))$
\sn
\itemitem{ {} } $(c) \quad$ some $c \in M_{\alpha +1}$ realizes
$f^{-1}_\alpha(\text{\ortp}(a_{\beta,i},f_\alpha(M_\alpha),N_\alpha)$, so
by
\nl

\hskip25pt  clause (b) it follows that
$c \in M_{\alpha +1} \backslash M_\alpha$ 
\sn
\itemitem{ $(**)^{\beta',i',\beta,i}_g$ }
$\quad \,\,\,\,[g(\beta) + i < g(\beta')+i'] \vee$ \nl

$\qquad \quad [g(\beta) +i = g(\beta')+i'
\and g(\beta) < g(\beta')] \vee [g(\beta)+i = g(\beta')+ i' \and$ \nl

$\qquad \quad g(\beta) = g(\beta') \and i \le i']$.
\endroster}
\ermn
There is no problem to choose $f_\alpha,N_\alpha$.  Now in the end, by
clauses (A),(F) clearly
$N_\delta$ is $(\lambda,\text{cf}(\delta))$-brimmed over $N_0$, i.e.,
over $M_0$, so it suffices to prove that $f_\delta$ is onto $N_\delta$.
If not, then by Axiom (D)(c), the density, 
there is $d \in N_\delta \backslash f_\delta(M_\delta)$ such that
$p := \text{\ortp}(d,f_\delta(M_\delta),N_\delta) \in {\Cal
S}^{\text{bs}}(f_\delta(M_\delta))$  hence
for some $\beta(*) < \delta$ we have $d \in N_{\beta(*)}$ so for some
$i(*) < \lambda,d = a_{\beta(*),i(*)}$.  Also by Axiom (E)(c), (the
local character) for every $\beta < \delta$ large
enough say $\ge \beta_d$ the type $p$ does not fork over $f_\delta(M_\beta)$,
\wilog \, $\beta_d = \beta(*)$.  Let $q = f^{-1}_\delta(\text{\rm
\ortp}(d,f_\delta(M_\delta),N_\delta)$, so it $\in {\Cal
S}^{\text{bs}}(M_\delta)$. 

Let $u = \{\alpha:\beta(*) \le \alpha < \delta$ and $q \restriction
M_\alpha \in {\Cal S}^{\text{bs}}(M_\alpha)$ (note $\beta(*) \le
\alpha$) is realized in $M_{\alpha +1}\}$.  By clause (c) of the
assumption clearly $|u| = \lambda$.
Also by the definition of $v$ for every $\alpha \in u$ the
condition $(*)^{\beta(*),i(*)}_{N_\alpha,f_\alpha}$ holds, hence in
clause (F) the pair
$(\beta_\alpha,i_\alpha)$ is well defined and is 
``below" $(\beta(*),i(*))$ in the sense of clause (G).
But there are only $\le |g(\beta(*)) \times i(*)| < \lambda$ such pairs hence
for some $\alpha_1 < \alpha_2$ in $u$ we have $(\beta_{\alpha_1},
i_{\alpha_1}) = (\beta_{\alpha_2},i_{\alpha_2})$, a contradiction:
$a_{\beta_{\alpha_1},i_{\alpha_1}} \in \text{ Rang}(f_{\alpha_1+1}) \subseteq
\text{ Rang}(f_{\alpha_2}) = f_{\alpha_2}(M_{\alpha_2})$ 
hence \ortp$(a_{\beta_{\alpha_1},i_{\alpha_1}},
f_{\alpha_2}(M_{\alpha_2}),N_{\alpha_2}) \notin {\Cal S}^{\text{bs}}
(f_{\alpha_2}(M_{\alpha_2}))$, contradiction.  So we are done.
\hfill$\square_{\scite{600-4a.2}}$
\enddemo
\bn
\centerline {$* \qquad * \qquad *$}
\bn
The following is helpful for constructions so that we can amalgamate
disjointly preserving non-forking of a type; we first repeat the definition
of $K^{3,\text{bs}}_\lambda,<_{\text{bs}}$.
\definition{\stag{600-4a.3} Definition}  1) Recall 
$(M,N,a) \in K^{3,\text{bs}}_\lambda$ 
if $M \le_{\frak K} N$ are models from $K_\lambda,a \in N \backslash M$
and \ortp$(a,M,N) \in {\Cal S}^{\text{bs}}(M)$.  
Let $(M_1,N_1,a) \le_{\text{bs}} (M_2,N_2,a)$ or write 
$\le^{\frak s}_{\text{bs}}$, when: both triples are in 
$K^{3,\text{bs}}_\lambda,M_1
\le_{\frak K} M_2,N_1 \le_{\frak K} N_2$ and \ortp$(a,M_2,N_2)$ does not fork
over $M_1$. \nl
2) We say $(M,N,a)$ is bs-reduced \ub{when} if it belongs to 
$K^{3,\text{bs}}_\lambda$ and $(M,N,a) \le_{\text{bs}} 
(M',N',a) \in K^{3,\text{bs}}_\lambda \Rightarrow 
N \cap M' = M$. \nl
3) We say $p \in {\Cal S}^{\text{bs}}(N)$ is a (really the) stationarization
of $q \in {\Cal S}^{\text{bs}}(M)$ if $M \le_{\frak K} N$ and $p$ is an
extension of $q$ which does not fork over $M$.
\enddefinition
\bigskip

\remark{Remark}  1) The definition of $K^{3,\text{bs}}_\lambda$ is 
compatible with the one in \scite{600-1.6} by \scite{600-1.8}(1).
\nl
2) We could have strengthened the definition of bs-reduced
(\scite{600-4a.3}), e.g., add: for no $b \in N' \backslash M'$, do we
have \ortp$(b,M',N') \in {\Cal S}^{\text{bs}}(M')$ and there are
$M'',N''$ such that $(M',N',a) \le_{\text{bs}} (M'',N'',a)$ and
\ortp$(b,M'',N'')$ forks over $M'$.
\endremark
\bigskip

\proclaim{\stag{600-4a.4} Claim}  For parts (3),(4),(5) assume ${\frak s}$
is categorical (in $\lambda$).
\nl
1)  If $\kappa \le \lambda,(M,N,a) \in 
K^{3,\text{bs}}_\lambda$, \ub{then} 
we can find $M',N'$ such that: $(M,N,a) \le_{\text{bs}} 
(M',N',a) \in K^{3,\text{bs}}_\lambda,M'$ is
$(\lambda,\kappa)$-brimmed over
$M,N'$ is $(\lambda,\kappa$)-brimmed over $N$ 
and $(M',N',a)$ is {\rm bs}-reduced. \nl
1A) If $(M,N_\ell,a_\ell) \in K^{3,\text{bs}}_\lambda$ for 
$\ell = 1,2$, \ub{then}
we can find $M^+,f_1,f_2$ such that: $M \le_{\frak K} M^+ \in
K_\lambda$ and for $\ell \in \{1,2\},f_\ell$ is a $\le_{\frak
K}$-embedding of $N_\ell$ into $M^+$ over $M$ and
$(M,f_\ell(N_\ell),f_\ell(a_\ell)) \le_{\text{bs}}
(f_{3-\ell}(N_{3-\ell}),M^+,f_\ell(a_\ell))$, equivalently
${\text{\rm \ortp\/}}(f_\ell(a_\ell),f_{3-\ell}(N_{3-\ell}),M^+)$ does 
not fork over $M$.
\nl
2) If $(M_\alpha,N_\alpha,a) \in K^{3,\text{bs}}_\lambda$ is
$\le_{\text{bs}}$-increasing for $\alpha < \delta$ and 
$\delta < \lambda^+$ is a limit ordinal \ub{then} their 
union $(\dbcu_{\alpha < \delta} M_\alpha,\dbcu_{\alpha < \delta}
N_\alpha,a)$ is a $\le_{\text{bs}}$-lub.  If each $(M_\alpha,N_\alpha,a)$
is bs-reduced then so is their union. \nl
3) Let $\lambda$ divide $\delta,\delta < \lambda^+$.  We can find $\langle
N_j,a_i:j \le \delta,i < \delta \rangle$ such that: $N_j \in
K_\lambda$ is $\le_{\frak K}$-increasing continuous,
$(N_j,N_{j+1},a_j) \in K^{3,\text{bs}}_\lambda$ is bs-reduced and if $i <
\delta,p \in {\Cal S}^{\text{bs}}(N_i)$ \ub{then} for 
$\lambda$ ordinals $j \in (i,i + \lambda)$ the type 
${\text{\rm \ortp\/}}(a_j,N_j,N_{j+1})$ is a non-forking extension of
$p$; so $N_\delta$ is $(\lambda,{\text{\rm cf\/}}(\delta))$-brimmed over
each $N_i,i < \delta$.  We can add ``$N_0$ is brimmed".  
\nl
4) For any $(M_0,M_1,a) \in K^{3,\text{bs}}_\lambda$ and 
$M_2 \in K_\lambda$ such that $M_0
\le_{\frak K} M_2$ there are $N_0,N_1$ such that $(M_0,M_1,a)
\le_{\text{bs}} 
(N_0,N_1,a),M_0 = M_1 \cap N_0$ and $M_2,N_0$ are isomorphic over
$M_0$.  (In fact, if $(M_0,M_2,b) \in K^{3,\text{bs}}_\lambda$ we can add that
for some isomorphism $f$ from $M_2$ onto $N_0$ over $M_0$ we have
$(M_0,N_0,f(a)) \le_{\text{bs}} (M_1,N_1,f(a))$.)
\nl
5) If $M_0 \in K_\lambda$ is brimmed and $M_0 \le_{\frak s} M_\ell$
for $\ell=1,2$ and there is a disjoint $\le_{\frak s}$-amalgamation of
$M_1,M_2$ over $M_0$.
\endproclaim
\bigskip

\demo{Proof}  1) We choose $M_i,N_i,b^\ell_i(\ell=1,2),
\bar{\bold c}_i$ by induction on $i < \delta := \lambda$ such that
\mr
\item "{$(a)$}"  $(M_i,N_i,a) \in K^{3,\text{bs}}_{\frak s}$ is
$\le_{\text{bs}}$-increasing continuous
\sn
\item "{$(b)$}"  $(M_0,N_0) = (M,N)$
\sn
\item "{$(c)_1$}"  $b^1_i \in M_{i+1} \backslash M_i$ and \ortp$_{\frak
s}(b^1_i,M_i,M_{i+1})\in {\Cal S}^{\text{bs}}(M_i)$,
\sn
\item "{$(c)_2$}"  $b^2_i \in N_{i+1} \backslash N_i$ and 
\ortp$(b^2_i,N_i,N_{i+1}) \in {\Cal S}^{\text{bs}}(N_i)$
\sn
\item "{$(d)_1$}"  if $i < \lambda$ and $p \in {\Cal
S}^{\text{bs}}(M_i)$ \ub{then} the set $\{j:i \le j < \lambda$ and
\ortp$_{\frak s}(b^1_j,M_j,M_{j+1})$ is a non-forking extension of $p\}$
has order type $\lambda$
\sn
\item "{$(d)_2$}"  if $i < \lambda$ and $p \in {\Cal
S}^{\text{bs}}(N_i)$ then the set $\{j:i \le j < \lambda$ and
\ortp$_{\frak s}(b^2_j,N_j,N_{j+1})$ is the non-forking extension of
$p\}$ has order type $\lambda$
\sn
\item "{$(e)$}"  $\bar{\bold c}_i = \langle c_{i,j}:j < \lambda
\rangle$ list $N_i$
\sn
\item "{$(f)$}"    if $\alpha < \lambda,i \le \alpha,j <
\lambda,c_{i,j} \notin M_\alpha$ but for some $(M'',N'')$ we have
$(M_{\alpha +1},N_{\alpha +1},a) \le_{\text{bs}} (M'',N'',a)$ and
$c_{i,j} \in M''$ \ub{then} for some $i_1,j_1 \le \text{ max}\{i,j\}$
we have $c_{i_1,j_1} \in M_{\alpha +1} \backslash M_\alpha$.
\ermn
Lastly, let $M' = \cup\{M_i:i < \lambda\},N' = \cup\{N_i:i <
\lambda\}$, by \scite{600-4a.2} $M'$ is
$(\lambda,\text{cf}(\lambda))$-brimmed over $M$ (using $(d)_1$), and
$N'$ is $(\lambda,\text{cf}(\lambda))$-brimmed over $N$ (using
$(d)_2$). \nl
Lastly, bs-reduced is by clauses (e)+(f).
\nl
1A),2) Easy.
\nl
3)  For proving part (3) use part (1) and the ``so" is
by using \scite{600-4a.2}. 
\nl
4)  For proving part (4), \wilog \,  $M_2$ is 
$(\lambda,\text{\rm cf}(\lambda))$-brimmed over $M_0$, as we can replace
$M_2$ by $M'_2$ if $M_2 \le_{\frak K} M'_2 \in K_\lambda$.  By part
(3) there is a sequence $\langle a_i:i < \delta \rangle$ and 
an $\le_{\frak K}$-increasing continuous $\langle N_i:i \le 
\delta \rangle$ with $N_0
= M_0,N_\delta = M_2$ and $(N_i,N_{i+1},a_i) \in K^{3,\text{bs}}_\lambda$
is reduced.  Then use (1A) successively. 
\nl
5) By part (3).  \hfill$\square_{\scite{600-4a.4}}$
\enddemo
\bigskip

\proclaim{\stag{600-4a.5} Claim}  Assume
\mr
\item "{$(a)$}"  $\gamma < \lambda^+$ is a limit ordinal
\sn
\item "{$(b)$}"  $\delta_i < \lambda^+$ is divisible by $\lambda$ for
$i \le \gamma,\langle \delta_i:i \le \gamma \rangle$ is increasing continuous
\sn
\item "{$(c)$}"  $\langle N_i:i < \gamma \rangle$ is
$\le_{\frak K}$-increasing continuous in $K_\lambda$
\sn
\item "{$(d)$}"  $\langle M_i:i < \gamma \rangle$ is
$\le_{\frak K}$-increasing continuous in $K_\lambda$
\sn
\item "{$(e)$}"  $N_i \le_{\frak K} M_i$ for $i < \gamma$
\sn
\item "{$(f)$}"  $\langle M_{i,j}:j \le \delta_i \rangle$ is
$\le_{\frak K}$-increasing continuous in $K_\lambda$ for each $i < \gamma$
\sn
\item "{$(g)$}"  $M_{i,0} = N_i,M_{i,\delta_i} = M_i,
a_j \in M_{i,j+1} \backslash M_{i,j}$ and ${\text{\rm \ortp\/}}
(a_j,M_{i,j},M_{i,j+1}) \in {\Cal S}^{\text{bs}}(M_{i,j})$ when $i
< \gamma,j < \delta_i$ 
\sn
\item "{$(h)$}"  if $j \le \delta_{i(*)},i(*) < \gamma$ then
$\langle M_{i,j}:i \in [i(*),\gamma) \rangle$ is
$\le_{\frak K}$-increasing continuous
\sn
\item "{$(i)$}"   ${\text{\rm \ortp\/}}
(a_j,M_{\beta,j},M_{\beta,j+1})$ does not fork over
$M_{i,j}$ when $i < \gamma,j < \delta_i,i \le \beta < \gamma$
\sn
\item "{$(j)$}"  if $i < \gamma,j < \delta_i,p \in {\Cal S}^{\text{bs}}
(M_{i,j})$ \ub{then} for $\lambda$ ordinals $j_1 \in [j,\delta_i)$ 
we have \text{\rm tp}$(a_{j_1},M_{i,j_1},M_{i,j_1+1}) \in {\Cal
S}^{\text{bs}}(M_{i,j_1})$ is a non-forking extension of $p$ or we can
ask less
\sn
\item "{$(j)^-$}"  if $i < \gamma,j < \delta_i$ and $p \in {\Cal
S}^{\text{bs}}(M_{i,j})$ \ub{then} for $\lambda$ ordinals 
$j_1 \in [j,\delta_\gamma)$ for some $i_1 \in [i,\gamma)$ we have 
${\text{\rm \ortp\/}}(a_{j_1},M_{i_1,j_1},M_{i_1,j_1+1}) 
\in {\Cal S}^{\text{bs}}(M_{i_1,j_1})$ 
is a non-forking extension of $p$.
\ermn
\ub{Then} $M_\gamma := \cup\{M_{i,j}:i < \gamma,j < \delta_i\} =
\{M_i:i < \gamma\}$ is
$(\lambda,{\text{\rm cf\/}}(\gamma))$-brimmed over $N_\gamma :=
\cup\{N_i:i < \gamma\}$.
\endproclaim
\bigskip

\demo{Proof}  For $j < \delta_\gamma$ let $M_{\gamma,j} =
\cup\{M_{i,j}:i < \gamma\}$, and let
$M_{\gamma,\delta_\gamma} = M_\gamma$ be $\cup\{M_{\gamma,j}:j <
\delta_\gamma\}$.  Easily $\langle M_{\gamma,j}:j \le \delta_\gamma
\rangle$ is $\le_{\frak K}$-increasing continuous, $M_{\gamma,j} \in
K_\lambda$ and $i \le \gamma \wedge j < \delta_i \Rightarrow M_{i,j}
\le_{\frak K} M_{\gamma,j}$.  Also if $i < \gamma,j < \delta_i$ then
\ortp$(a_j,M_{\gamma,j},M_{\gamma,j+1}) \in {\Cal S}^{\text{bs}}(M_{\gamma,j})$
does not fork over $M_{i,j}$ by Axiom (E)(h), continuity. \nl
Now if $j < \delta_\gamma$ and $p \in {\Cal
S}^{\text{bs}}(M_{\gamma,j})$ then for some $i < \gamma,p$ does not
fork over $M_{i,j}$ (by Ax(E)(c)) and \wilog \, $j < \delta_i$.

Hence if clause (j) holds we have 
$u := \{\varepsilon:j < \varepsilon < \delta_i$ and
\ortp$(a_\varepsilon,M_{i,\varepsilon},M_{i,\varepsilon+1})$ is a
non-forking extension of $p \restriction M_{i,j}\}$ has $\lambda$
members.  But for $\varepsilon \in u$,
\ortp$(a_\varepsilon,M_{\gamma,\varepsilon},M_{\gamma,\varepsilon +1})$
does not fork over $M_{i,\varepsilon}$ (by clause (i) of the
assumption) hence does not fork over $M_{i,j}$ and by monotonicity it
does not fork over $M_{\gamma,i}$ and by uniqueness it extends $p$.
If clause $(j)^-$ holds the proof is similar.
By \scite{600-4a.2} the model $M_\gamma$ is $(\lambda,\text{cf}(\gamma))$-brimmed
over $N_\gamma$.  \hfill$\square_{\scite{600-4a.5}}$
\enddemo
\bigskip

\proclaim{\stag{600-4a.6} Lemma}  1) If $M \in K_\lambda$ and the models
$M_1,M_2 \in K_\lambda$ are $(\lambda,*)$-brimmed over $M$ 
(see Definition \scite{600-0.21}(2)), \ub{then} $M_1,M_2$ are isomorphic 
over $M$. \nl
2)  If $M_1,M_2 \in K_\lambda$ are $(\lambda,*)$-brimmed \ub{then} 
they are isomorphic.
\endproclaim
\bn
We prove some claims before proving \scite{600-4a.6}; we will not much use 
the lemma, but it is of obvious interest and its proof is crucial in 
one point of \S6.
\proclaim{\stag{600-4a.7} Claim}  1)
\mr
\item "{$(E)(i)^+$}"  \ub{long non-forking amalgamation for 
$\alpha < \lambda^+$}: \nl
if $\langle N_i:i \le \alpha 
\rangle$ is $\le_{\frak K}$-increasing continuous
sequence of members of $K_\lambda,a_i \in N_{i+1} \backslash N_i$ for
$i < \alpha,p_i = { \text{\rm \ortp\/}}
(a_i,N_i,N_{i+1}) \in {\Cal S}^{\text{bs}}(N_i)$
and $q \in {\Cal S}^{\text{bs}}(N_0)$, \ub{then} we can find a
$\le_{\frak K}$-increasing continuous sequence $\langle N'_i:i \le \alpha
\rangle$ of members of $K_\lambda$ such that: $i \le \alpha \Rightarrow
N_i \le_{\frak K} N'_i$; some $b \in N'_0 \backslash N_0$ realizes 
$q,{\text{\rm \ortp\/}}(b,N_\alpha,N'_\alpha)$ does not fork over $N_0$ 
and ${\text{\rm \ortp\/}}(a_i,N'_i,N'_{i+1})$ does not fork over $N_i$
for $i < \alpha$.
\ermn
2) Above assume in addition that there are $M,b^*$ such that 
$N_0 \le_{\frak K} M \in K_\lambda,b^* \in M$ and ${\text{\rm \ortp\/}}
(b^*,N_0,M) = q$.  \ub{Then} we can add: there is a
$\le_{\frak K}$-embedding of $M$ into $N'_0$ over $N_0$ mapping $b^*$ to
$b$.
\endproclaim
\bigskip

\demo{Proof}  Straight (remembering Axiom (E)(i) on non-forking
amalgamation of Definition \scite{600-1.1}).  In details
\nl
1) Let $M_0,b^*$ be such that $N_0 \le_{{\frak K}[{\frak s}]} M_0$ and
$q = \text{ \ortp}(b^*,N_0,M_0)$ and apply part (2).
\nl
2) We choose $(M_i,f_i)$ by induction on $i \le \alpha$ such that
\mr
\item "{$\circledast$}"  $(a) \quad M_i \in {\frak K}_{\frak s}$ is
$\le_{\frak K}$-increasing continuous
\sn
\item "{${{}}$}"  $(b) \quad f_i$ is a $\le_{\frak K}$-embedding of
$N_i$ into $M_i$
\sn
\item "{${{}}$}"  $(c) \quad f_i$ is increasing continuous with $i \le
\alpha$
\sn
\item "{${{}}$}"  $(d) \quad M_0 = M$ and $f_0 = \text{ id}_{N_0}$
\sn
\item "{${{}}$}"  $(e) \quad$ \ortp$(b^*,f_i(N_i),M_i)$ does not fork
over $N_0$
\sn
\item "{${{}}$}"  $(f) \quad$ \ortp$(f_{i+1}(a_i),M_i,M_{i+1})$ does not
fork over $f_i(N_i)$.
\ermn
For $i=0$ there is nothing to do.  For $i$ limit take unions.  Lastly,
for $i=j+1$, we can find $(M'_i,f'_i)$ such that $f_j \subseteq f'_i$
and $f'_i$ is an isomorphism from $N_i$ onto $M$.  Hence $f_j(N_j)
\le_{{\frak K}[{\frak s}]} N'_i$.  Now use Ax(E)(i) for
$f_j(N_j),M'_i,N_i,f'_i(a_j),b^*$.

Having carried the induction, we rename to finish.
\hfill$\square_{\scite{600-4a.7}}$ 
\enddemo
\bn

In the claim below, we are given a 
$\le_{{\frak K}_\lambda}$-increasing 
continuous $\langle M_i:i \le \delta\rangle$
and $u_0,u_1,u_2 \subseteq \delta$ such that: $u_0$ is where we are
already given $a_i \in M_{i+1} \backslash M_i,u_1 \subseteq \delta$ is
where we shall choose $a_i(\in M'_{i+1} \backslash M'_i)$ and $u_2
\subseteq \delta$ is the place which we ``leave for future use"; main
case is $u_1 = \delta;u_0 = u_2 = \emptyset$.
\proclaim{\stag{600-4a.8} Claim}  1)  Assume
\mr
\item "{$(a)$}"  $\delta < \lambda^+$ is divisible by $\lambda$
\sn
\item "{$(b)$}"  $u_0,u_1,u_2$ are disjoint subsets of $\delta$
\sn
\item "{$(c)$}"  $\delta = \sup(u_1)$ and ${\text{\rm otp\/}}(u_1)$ 
is divisible by $\lambda$
\sn
\item "{$(d)$}"  $\langle M_i:i \le \delta \rangle$ is 
$\le_{\frak K}$-increasing continuous in ${\frak K}_\lambda$
\sn
\item "{$(e)$}"  $\bold{\bar a} = \langle a_i:i \in u_0 \rangle,a_i \in
M_{i+1} \backslash M_i,{\text{\rm \ortp\/}}(a_i,M_i,M_{i+1}) 
\in {\Cal S}^{\text{bs}}(M_i)$.
\ermn
\ub{Then} we can find $\bar M' = \langle M'_i:i \le \delta \rangle$ and
$\bold{\bar a}' = \langle a_i:i \in u_1 \rangle$ such that
\mr
\item "{$(\alpha)$}"  $\bar M'$ is $\le_{\frak K}$-increasing continuous
in $K_\lambda$
\sn
\item "{$(\beta)$}"  $M_i \le_{\frak K} M'_i$
\sn
\item "{$(\gamma)$}"  if $i \in u_0$ then ${\text{\rm \ortp\/}}
(a_i,M'_i,M'_{i+1})$ is 
a non-forking extension of ${\text{\rm \ortp\/}}(a_i,M_i,M_{i+1})$
\sn
\item "{$(\delta)$}"  if $i \in u_2$ then $M'_i = M'_{i+1}$
\sn
\item "{$(\varepsilon)$}"  if $i \in u_1$ then ${\text{\rm \ortp\/}}
(a_i,M'_i,M'_{i+1})\in {\Cal S}^{\text{bs}}(M'_i)$
\sn
\item "{$(\zeta)$}"  if $i < \delta,p \in {\Cal S}^{\text{bs}}(M'_i)$
then for $\lambda$ ordinals $j \in u_1 \cap (i,\delta)$ the type 
${\text{\rm \ortp\/}}(a_j,M'_j,M'_{j+1})$ is a non-forking extension of $p$.
\ermn
2) If we add in part (1) the assumption
\mr
\item "{$(g)$}"  $M_0 \le_{\frak K} N \in K_\lambda$
\ermn
\ub{then} we can add to the conclusion
\mr
\item "{$(\eta)$}"  there is an $\le_{\frak K}$-embedding $f$ of $N$ into
$M'_0$ over $M_0$ and moreover $f$ is onto.
\ermn
3) If we add in part (1) the assumption
\mr
\item "{$(h)^+$}"  $M_0 \le_{\frak K} N \in K_\lambda$ and $b \in N
\backslash M_0,{\text{\rm \ortp\/}}(b,M_0,N) \in {\Cal S}^{\text{bs}}(M_0)$
\ermn
\ub{then} we can add to the conclusion
\mr
\item "{$(\eta)^+$}"  as in $(\eta)$ and ${\text{\rm \ortp\/}}
(f(b),M_\delta,M'_\delta)$ does not fork over $M_0$.
\ermn
4) We can strengthen clause $(\zeta)$ in part (1) to
\mr
\item "{$(\zeta)^+$}" if $i < \delta$ and $p \in {\Cal S}^{\text{bs}}
(M'_i)$ then
for $\lambda$ ordinals $j$ we have $j \in [i,\delta) \cap u_1$ and
${\text{\rm \ortp\/}}(a_j,M'_j,M'_{j+1})$ is a non-forking extension of 
$p$ and ${\text{\rm otp\/}}(u_1 \cap j \backslash i) < \lambda$.
\ermn
5) If $i \in u_2 \Rightarrow M_i = M_{i+1}$ then we can add $i \in u_2
\Rightarrow M'_i = M'_{i+1}$.
\endproclaim
\bigskip

\demo{Proof}  Straight.  Note that we can find a sequence $\langle
u_{1,i,\varepsilon}:i < \delta,\varepsilon < \lambda \rangle$ such
that: this is a sequence of pairwise disjoint subsets of $u_1$ each of
cardinality $\lambda$ satisfying $u_{1,i,\varepsilon} 
\subseteq \{j:i < j,j \in u_1 \text{ and }
|u_1 \cap (i,j)| < \lambda\}$ (or we can demand that $i \le i_1 < i_2
\le \delta \wedge |u_1 \cap (i_1,i_2)| = \lambda \Rightarrow
|u_{1,i,\varepsilon} \cap (i_1,i_2)| = \lambda$).  
\hfill$\square_{\scite{600-4a.8}}$  
\enddemo
\bn
Toward building our rectangles of models with sides of difference
lengths (and then we shall use \scite{600-4a.5}) we show
(to understand the aim of the clauses in the conclusion of \scite{600-4a.9}
see the proof of \scite{600-4a.6} below):
\proclaim{\stag{600-4a.9} Claim}  Assume
\mr
\item "{$(a)$}"  $\delta_\ell < \lambda^+$ is divisible by $\lambda$
for $\ell = 1,2$
\sn
\item "{$(b)$}"  $\bar M^\ell = \langle M^\ell_\alpha:\alpha \le \delta_\ell
\rangle$ is $\le_{\frak K}$-increasing continuous for $\ell = 1,2$
\sn
\item "{$(c)$}"  $u^\ell_0,u^\ell_1,u^\ell_2$ are disjoint subsets of
$\delta_\ell,{\text{\rm otp\/}}(u^\ell_1)$ is divisible by
$\lambda$ and $\delta_\ell = \sup(u^\ell_1)$ for $\ell = 1,2$
\sn
\item "{$(d)$}"  $\bold{\bar a}^\ell \equiv \langle a^\ell_\alpha:\alpha
\in u^\ell_0 \rangle$ and ${\text{\rm \ortp\/}}(a^\ell_\alpha,M^\ell_\alpha,
M^\ell_{\alpha +1}) \in {\Cal S}^{\text{bs}}(M^\ell_\alpha)$ for
$\ell = 1,2,\alpha \in u^\ell_0$
\sn
\item "{$(e)$}"  $M^1_0 = M^2_0$
\sn
\item "{$(f)$}"  $\alpha \in u^\ell_1 \cup u^\ell_2 
\Rightarrow M^\ell_\alpha = M^\ell_{\alpha +1}$ for $\ell =1,2$.
\ermn
\ub{Then} we can find $\bar f^\ell = \langle f^\ell_\alpha:\alpha \le
\delta_\ell \rangle,\bold{\bar b}^\ell = \langle b^\ell_\alpha:
\alpha \in u^\ell_0 \cup u^\ell_1 \rangle$ for $\ell = 1,2$ and $\bar M =
\langle M_{\alpha,\beta}:\alpha \le \delta_1,\beta \le \delta_2
\rangle$ and functions $\zeta:u^1_1 \rightarrow \delta_2$ and
$\varepsilon:u^2_1 \rightarrow \delta_1$ such that
\mr
\item "{$(\alpha)_1$}"  for each 
$\alpha \le \delta_1,\langle M_{\alpha,\beta}:\beta \le \delta_2 \rangle$
is $\le_{\frak K}$-increasing continuous
\sn
\item "{$(\alpha)_2$}"  for each $\beta \le \delta_2,
\langle M_{\alpha,\beta}:\alpha \le \delta_1 \rangle$
is $\le_{\frak K}$-increasing continuous
\sn
\item "{$(\beta)_1$}"  for $\alpha \in u^1_0,b^1_\alpha$ belongs to
$M_{\alpha +1,0}$ and ${\text{\rm \ortp\/}}
(b^1_\alpha,M_{\alpha,\delta_2},M_{\alpha +1,
\delta_2}) \in {\Cal S}^{\text{bs}}(M_{\alpha,\delta_2})$ does not fork over
$M_{\alpha,0}$
\sn
\item "{$(\beta)_2$}"  for $\beta \in u^2_0,b^2_\beta$ belongs to
$M_{0,\beta +1}$ and ${\text{\rm \ortp\/}}
(b^2_\beta,M_{\delta_1,\beta},M_{\delta_1,\beta +1})
\in {\Cal S}^{\text{bs}}(M_{\delta_1,\beta})$ does not fork over
$M_{0,\beta}$
\sn
\item "{$(\gamma)_1$}" for $\alpha \in u^1_1,\zeta(\alpha) < \delta_2$
and we have $b^1_\alpha \in
M_{\alpha +1,\zeta(\alpha) + 1}$ and ${\text{\rm \ortp\/}}
(b^1_\alpha,M_{\alpha,\delta_2},
M_{\alpha +1,\delta_2})$ does not fork over $M_{\alpha,\zeta(\alpha)+1}$
\sn
\item "{$(\gamma)_2$}"  for $\beta \in u^2_1,
\varepsilon(\beta) < \delta_1$ and we have 
$b^2_\beta \in M_{\varepsilon(\beta)+1,\beta + 1}$ and 
${\text{\rm \ortp\/}}
(b^2_\beta,M_{\delta_1,\beta},M_{\delta_1,\beta +1})$ does not 
fork over $M_{\varepsilon(\beta)+1,\beta}$
\sn
\item "{$(\delta)_1$}"  if $\alpha < \delta_1,\beta < \delta_2$ and $p \in
{\Cal S}^{\text{bs}}(M_{\alpha,\beta})$ or just $p \in {\Cal
S}^{\text{bs}}(M_{\alpha,\beta +1})$ \ub{then} for $\lambda$ ordinals
\footnote{we can add ``and otp$(\alpha' \cap u^1_1 \backslash
\alpha_2) < \lambda$"}
$\alpha' \in [\alpha,\delta_1) \cap u^1_1$, the type
${\text{\rm \ortp\/}}
(b^1_{\alpha'},M_{\alpha',\beta +1},M_{\alpha +1,\beta +1})$ is a
(well defined) non-forking extension of $p$ and $\beta = \zeta(\alpha')$
\sn
\item "{$(\delta)_2$}"  if $\alpha < \delta_1,\beta < \delta_2$ and $p \in
{\Cal S}^{\text{bs}}(M_{\alpha,\beta})$ or just $p \in {\Cal
S}^{\text{bs}}(M_{\alpha +1,\beta})$ \ub{then} for $\lambda$ ordinals
\footnote{we can add ``and otp$(\beta' \cap u^2_1 \backslash \beta_2)
< \lambda$"}
$\beta' \in [\beta,\delta_2) \cap u^2_1$, the type
${\text{\rm \ortp\/}}(b^2_{\beta'},
M_{\alpha +1,\beta'},M_{\alpha+1,\beta'+1})$ is a non-forking
extension of $p$ and $\alpha = \varepsilon(\beta')$
\sn
\item "{$(\varepsilon)$}"  $M_{0,0} = M^1_0 = M^2_0$
\sn
\item "{$(\zeta)_1$}"  $f^1_\alpha$ is an isomorphism from $M^1_\alpha$ onto
$M_{\alpha,0}$ such that $\alpha \in u^1_0 \Rightarrow f^1_\alpha
(a^1_\alpha) = b^1_\alpha$ \nl
$f^1_0 = { \text{\rm id\/}}_{M^1_0}$ and 
$f^1_\alpha$ is increasing continuous with $\alpha$
\sn
\item "{$(\zeta)_2$}"  $f^2_\beta$ is an isomorphism from $M^2_\beta$ onto
$M_{0,\beta}$ such that $\beta \in u^2_0 \Rightarrow f^2_\beta(a^2_\beta) =
b^2_\beta$ \nl
$f^2_0 = { \text{\rm id\/}}_{M^2_0}$ and 
$f^2_\alpha$ is increasing continuous with $\alpha$
\sn
\item "{$(\eta)_1$}" if $\alpha \in u^1_2$ then $M_{\alpha,\beta} =
M_{\alpha +1,\beta}$ for every $\beta \le \delta_2$
\sn
\item "{$(\eta)_2$}"  if $\beta \in u^2_2$ then $M_{\alpha,\beta} =
M_{\alpha,\beta +1}$ for every $\alpha \le \delta_1$.
\endroster
\endproclaim
\bigskip

\demo{Proof}  Straight, divide $u^\ell_1$ to $\delta_{3 - \ell}$ 
subsets large enough), in fact, we can first choose the function
$\zeta(-),\varepsilon(-)$.  Now choose $\langle M_{\alpha,\beta}:\alpha
\le \delta_1,\beta \le \beta^* \rangle,\langle f^1_\alpha:\alpha \le
\delta_1 \rangle,\langle f^2_\beta:\beta \le \beta^* \rangle$ and
$\langle b^1_\alpha:\zeta(\alpha) \in \beta^* \rangle,\langle
b^2_\beta:\beta < \beta^* \rangle$ by induction on $\beta^*$ using
\scite{600-4a.8}.   \hfill$\square_{\scite{600-4a.9}}$
\enddemo
\bigskip

\demo{Proof of \scite{600-4a.6}}  By \scite{600-0.22}(3), i.e., uniqueness of
the $(\lambda,\theta_\ell)$-brimmed model over $M$, it is 
enough to show for any regular
$\theta_1,\theta_2 \le \lambda$ that there is a model $N \in K_\lambda$ which
is $(\lambda,\theta_\ell)$-brimmed over $M$ for $\ell = 1,2$.  
Let $\delta_1 = \lambda
\times \theta_1,\delta_2 = \lambda \times \theta_2$ (ordinal
multiplication, of course), $M^1_\alpha = M^2_\beta =
M$ for $\alpha \le \delta_1,\beta \le \delta_2,u^1_0 = u^2_0 = \emptyset,
u^1_1 = \delta_1,u^2_1 = \delta_2,u^1_2 = u^2_2 = \emptyset$.  So there
are $\langle M_{\alpha,\beta}:\alpha \le \delta_1,\beta \le \delta_2 \rangle,
\langle b^1_\alpha:\alpha < \delta_1 \rangle,\langle b^2_\beta:\beta <
\delta_2 \rangle$ and $\langle f^1_\alpha:\alpha \le \delta_1
\rangle,\langle f^2_\beta:\beta \le \delta_2 \rangle$
as in Claim \scite{600-4a.9}.  Without loss of generality
$f^1_\alpha = f^2_\alpha = \text{ id}_M$.  Now
\mr
\item "{$(*)_1$}"  $\langle M_{\alpha,\delta_2}:\alpha \le \delta_1
\rangle$ is $\le_{\frak K}$-increasing continuous in $K_\lambda$ (by
clause $(\alpha)_1$, of \scite{600-4a.9}). Also
\sn
\item "{$(*)_2$}"  if $\alpha < \delta_1$ and
$p \in {\Cal S}(M_{\alpha,\delta_2})$
\ub{then} for $\lambda$ ordinals 
$\alpha' \in (\alpha,\delta_1) \cap u^1_1$ the type 
\ortp$(b^1_{\alpha',\delta_2},
M_{\alpha',\delta_2},M_{\alpha' +1,\delta_2})$ is a
non-forking extension of $p$.
\ermn
(Easy, by Axiom (E)(c) for some $\beta < \delta_2,p$ does not fork over
$M_{\alpha,\beta +1}$ and use clause $(\delta)_1$ of \scite{600-4a.9}).

So by \scite{600-4a.5}, $M_{\delta_1,\delta_2}$ is 
$(\lambda,\text{cf}(\delta_1))$-brimmed over $M_{0,\delta_2}$ which is $M$.

Similarly $M_{\delta_1,\delta_2}$ is 
$(\lambda,\text{cf}(\delta_2))$-brimmed over $M_{\delta_1,0}$ which
is $M$; so together we are done. \nl
${{}}$   \hfill$\square_{\scite{600-4a.6}}$
\enddemo
\bigskip

\proclaim{\stag{600-4a.10A} Claim}  1) If $M \in K_{\lambda^+}$ and $p \in {\Cal
S}^{\text{bs}}(M_0),M_0 \le_{\frak K} M$ (so $M_0 \in K_\lambda$), 
\ub{then} we can find
$b,\langle N^0_\alpha:\alpha \le \lambda^+\rangle$ and $\langle
N^1_\alpha:\alpha \le \lambda^+ \rangle$ such that
\mr
\item "{$(a)$}"  $\langle N^0_\alpha:\alpha < \lambda^+ \rangle$ is a
$\le_{\frak K}$-representation of $N^0_{\lambda^+} = M$
\sn
\item "{$(b)$}"  $\langle N^1_\alpha:\alpha < \lambda^+ \rangle$ is a
$\le_{\frak K}$-representation of $N^1_{\lambda^+} \in K_{\lambda^+}$
\sn
\item "{$(c)$}"  $N^1_{\alpha +1}$ is $(\lambda,\lambda)$-brimmed over
$N^1_\alpha$ (hence $N^1_{\lambda^+}$ is saturated over $\lambda$ in
${\frak K}$)
\sn
\item "{$(d)$}"  $M_0 \le N^0_0$ and $N^0_\alpha \le_{\frak K} N^1_\alpha$
\sn
\item "{$(e)$}"  ${\text{\rm \ortp\/}}_{\frak
s}(b,N^0_\alpha,N^1_\alpha)$ is a non-forking extension of $p$ for
every $\alpha < \lambda^+$.
\ermn
2) We can add
\mr
\item "{$(f)$}"  for $\alpha < \beta < \lambda^+,N^1_\beta$ is
$(\lambda,*)$-brimmed over $N^0_\beta \cup N^1_\alpha$.
\endroster
\endproclaim
\bigskip

\demo{Proof}  1) Easy by long non-forking amalgamation \scite{600-4a.7} 
(see \scite{600-4a.10}).
\nl
2) Use \scite{600-4a.5}.  \hfill$\square_{\scite{600-4a.10A}}$ 
\enddemo
\bigskip

\demo{\stag{600-4a.12} Conclusion}  1) $K_{\lambda^{++}} \ne \emptyset$. \nl
2) $K_{\lambda^+} \ne \emptyset$. \nl
3) {\rm No} $M \in K_{\lambda^+}$ is $\le_{\frak K}$-maximal.
\enddemo
\bigskip

\demo{Proof}  1) By (2) + (3). \nl
2) By $(B)$ of \scite{600-1.1}. \nl
3) By \scite{600-4a.10A}. \hfill$\square_{\scite{600-4a.12}}$
\enddemo
\bn
\margintag{600-4a.14F}\ub{\stag{600-4a.14F} Exercise}:  1) Let $M \in K_{\frak s}$ be superlimit
and ${\frak t} = {\frak s}_{[M]}$, so $K_{\frak t}$ is categorical.
If $(M,N,a) \in K^{\text{bs}}_{\frak t}$ is reduced for ${\frak t}$,
\ub{then} it is reduced for ${\frak s}$.
\nl
2) In \scite{600-4a.4}(3),(4),(5), we can omit the assumption ``${\frak s}$
is categorical" if:
\mr
\item "{$(a)$}"  we add in aprt (3), each $N_i$ is superlimit
(equivalently brimmed)
\sn
\item "{$(b)$}"  in parts (4),(5) add the assumption ``$M_0$ is
superlimit".
\ermn
2) Some extra assumption in \scite{600-4a.4}(5) is needed.
\newpage

\head {\S5 Non-structure or some unique amalgamation} \endhead  \resetall \sectno=5
 \spuriousreset
\bn
We shall get from essentially $\dot I(\lambda^{++},K) < 2^{\lambda^{++}}$ or
just $\dot I(\lambda^{++},K(\lambda^+$-saturated)) $< 2^{\lambda^{++}}$, 
many cases of uniqueness of amalgamation 
assuming WDmId$(\lambda^+)$ is not $\lambda^{++}$-saturated.
The proof is similar to \cite{Sh:482}, \cite[\S3]{Sh:576}.

We define $K^{3,\text{bt}}_\lambda$, it is a brimmed relative of 
$K^{3,\text{bs}}_\lambda$
hence the choice of bt; it guarantees much brimness
(see Definition \scite{600-nu.1}) hence it guarantees some uniqueness,
that is, if $(M,N,a) \in K^{3,\text{bt}}_\lambda,M$ is unique
(recalling the uniqueness of the brimmed model) and more crucially, we consider
$K^{3,\text{uq}}_\lambda$, (the family of members of
$K^{3,\text{bs}}_\lambda$ for which we have
uniqueness in relevant extensions).  Having enough such triples is the
main conclusion of this section (in \scite{600-nu.6} under ``not too many
non isomorphic models" assumptions).  In \scite{600-nu.2} we give some
properties of $K^{3,\text{bt}}_\lambda,K^{3,\text{uq}}_\lambda$.

To construct models in $\lambda^{++}$ we use approximations of cardianlity in
$\lambda^+$ with ``obligation" on the further construction, which are
presented as pairs $(\bar M,\bar{\bold a}) \in K^{\text{sq}}_\lambda$ ordered
by $\le_{\text{ct}}$, see Definition \scite{600-nu.3}, Claims \scite{600-nu.3A},
\scite{600-nu.4}.  We need more: the triples 
$(\bar M,\bar{\bold a},\bold f) \in K^{\text{mqr}}_S,K^{\text{nqr}}_S$
in Definition \scite{600-nu.7}, Claim \scite{600-nu.8}.  
All this enables us to quote results of \cite[\S3]{Sh:576}, but apart
from believing the reader do not need to know \cite{Sh:576}. 
\bigskip

\demo{\stag{600-nu.0} Hypothesis}
\mr
\item "{$(a)$}"  ${\frak s} = 
({\frak K},\nonfork{}{}_{},
{\Cal S}^{\text{bs}})$ is a good $\lambda$-frame. 
\endroster
\enddemo
\bigskip

\definition{\stag{600-nu.1} Definition}  1) Let $K^{3,\text{bt}}_\lambda =
K^{3,\text{bt}}_{\frak s}$
be the set of triples $(M,N,a)$ such that for some $\sigma = 
\text{ cf}(\sigma) \le \lambda,
M \le_{\frak K} N$ are both $(\lambda,\sigma)$-brimmed members of 
$K_\lambda,a \in N \backslash M$ and \ortp$(a,M,N) \in 
{\Cal S}^{\text{bs}}(M)$. \nl
2) For $(M_\ell,N_\ell,a_\ell) \in K^{3,\text{bt}}_\lambda$ for $\ell =1,2$
let $(M_1,N_1,a_1) <_{\text{bt}} (M_2,N_2,a_2)$ mean $a_1 = a_2$, 
\ortp$(a_1,M_2,N_2)$ does not fork over $M_1$ and for 
some $\sigma_2 = \text{ cf}(\sigma_2) \le \lambda$, the model $M_2$ is 
$(\lambda,\sigma_2)$-brimmed over $M_1$ and the model $N_2$ is 
$(\lambda,\sigma_2)$-brimmed over $N_1$.  Finally $(M_1,N_1,a_2) 
\le_{\text{bt}}(M_2,N_2,a_2)$ means 
$(M_1,N_1,a_1) <_{\text{bt}} 
(M_2,N_2,a_2)$ or $(M_1,N_1,a_1) = (M_2,N_2,a_2)$.
\enddefinition
\bigskip

\definition{\stag{600-nu.1A} Definition}  1) Let ``$(M_0,M_2,a) \in 
K^{3,\text{uq}}_\lambda$" mean: $(M_0,M_2,a) \in
K^{3,\text{bs}}_\lambda$ and: 
for every $M_1$ satisfying $M_0 \le_{\frak K} M_1 \in K_\lambda$,
the amalgamation $M$ of $M_1,M_2$ over $M_0$, with \ortp$(a,M_1,M)$ not
forking over $M_0$, is unique, that is: 
\mr
\item "{$(*)$}"  \ub{if} for $\ell =1,2$ we have $M_0 \le_{\frak K} M_1 
\le_{\frak K} M^\ell \in K_\lambda$ and $f_\ell$ is a 
$\le_{\frak K}$-embedding of $M_2$ into $M^\ell$ over $M_0$ 
(so $f_1 \restriction M_0 = f_2 \restriction 
M_0 = \text{ id}_{M_0}$) such that \ortp$(f_\ell(a),M_1,M^\ell)$ does not
fork over $M_0$, \underbar{then}
{\roster
\itemitem{ $(a)$ }  [uniqueness]:  \nl
for some $M',g_1,g_2$ we have: 
$M_1 \le_{\frak K} M' \in K_\lambda$ and \newline
$g_\ell$ is a $\le_{\frak K}$-embedding of $M^\ell$ into $M'$ 
over $M_1$ for $\ell =1,2$ such that $g_1 \circ f_1 \restriction M_2 = 
g_2 \circ f_2 \restriction M_2$
\sn
\itemitem{ $(b)$ }  [being reduced] $f_\ell(M_2) \cap M_1 = M_0$ \nl
[this is ``for free" in the proofs; and is not really necessary so the
decision if to include it is not important].
\endroster}
\ermn
2) $K^{3,\text{uq}}_\lambda$ is dense (or ${\frak s}$ has density for
$K^{3,\text{uq}}_\lambda$) when $K^{3,\text{uq}}_\lambda$ is dense
in $(K^{3,\text{bs}}_\lambda,\le_{\text{bs}})$, i.e., for every
$(M_1,M_2,a) \in K^{3,\text{bs}}_\lambda$ there is $(M_1,N_2,a) \in
K^{3,\text{uq}}_\lambda$ such that $(M_1,M_2,a) \le_{\text{bs}}
(N_1,N_2,a) \in K^{3,\text{uq}}_\lambda$.
\nl
3) $K^{3,\text{uq}}_\lambda$ has existence or ${\frak s}$ has
existence for $K^{3,\text{uq}}_\lambda$ when for every $M_0 \in
K_\lambda$ and $p \in {\Cal S}^{\text{bs}}(M_0)$ for some $M_1,a$ we
have $(M_0,M_1,a) \in K^{3,\text{uq}}_\lambda$ and $p = \text{\rm
\ortp}(a,M_0,M_1)$.
\nl
4) $K^{3,\text{uq}}_{\frak s} = K^{3,\text{uq}}_\lambda$.
\enddefinition
\bigskip  

\proclaim{\stag{600-nu.2} Claim}  1) The relation $\le_{\text{bt}}$ is 
a partial order on $K^{3,\text{bt}}_\lambda$ that is transitive 
and reflexive (but not necessarily the parallel of Ax V of
{\rm a.e.c.} see Definition \scite{600-0.2}). \nl
2) If $(M_\alpha,N_\alpha,a) \in K^{3,\text{bt}}_\lambda$ 
is $\le_{\text{bt}}$-increasing
continuous for $\alpha < \delta$ where $\delta$ is a limit ordinal
$< \lambda^+$ \ub{then} $(M,N,a) = (\dbcu_{\alpha < \delta} M_\alpha, 
\dbcu_{\alpha < \delta} N_\alpha,a)$ belongs to
$K^{3,\text{bt}}_\lambda$ and 
$\alpha < \delta \Rightarrow (M_\alpha,N_\alpha,a) \le_{\text{bt}} 
(M,N,a)$ and
$(M,N,a)$ is a $\le_{\text{bt}}$-upper bound of 
$\langle (M_\alpha,N_\alpha,a):\alpha < \delta \rangle$. \nl
3) In $(*)$ of \scite{600-nu.1A}, clause (b) follows from (a).
\endproclaim
\bigskip

\demo{Proof}  Easy, e.g. (3) by the uniqueness (i.e., clause (a)) 
and \scite{600-4a.4}(4).  \hfill$\square_{\scite{600-nu.2}}$
\enddemo
\bn
We now define $K^{\text{sq}}_{\lambda^+}$, a family of $\le_{\frak
K}$-increasing continuous sequences (the reason for sq) in $K_\lambda$
of length $\lambda^+$, will be used to approximate stages in
constructing models in $K_{\lambda^{++}}$.  
\definition{\stag{600-nu.3} Definition}  1) Let $K^{\text{sq}}_{\lambda^+}
= K^{\text{sq}}_{\frak s}$
be the set of pairs $(\bar M,\bold{\bar a})$ such that (sq stands for
sequence):
\mr
\item "{$(a)$}"  $\bar M = \langle M_\alpha:\alpha < \lambda^+ \rangle$
is a $\le_{\frak K}$-increasing continuous sequence of models
from $K_\lambda$
\sn
\item "{$(b)$}"  $\bold{\bar a} = \langle a_\alpha:\alpha \in S \rangle$,
where $S \subseteq \lambda^+$ is stationary in $\lambda^+$ and 
$a_\alpha \in M_{\alpha +1} \backslash
M_\alpha$ 
\sn
\item "{$(c)$}"  for some club $E$ of $\lambda^+$ for every $\alpha
\in S \cap E$ we have 
\ortp$(a_\alpha,M_\alpha,M_{\alpha +1}) \in {\Cal S}^{\text{bs}}(M_\alpha)$
\sn
\item "{$(d)$}"  if $p \in {\Cal S}^{\text{bs}}
(M_\alpha)$ \ub{then} for stationarily
many $\delta \in S$ we have: \ortp$(a_\delta,M_\delta,M_{\delta +1}) \in
{\Cal S}^{\text{bs}}(M_\delta)$ does not fork over $M_\alpha$ and extends $p$.
\ermn
In such cases we let $M = \dbcu_{\alpha < \lambda^+} M_\alpha$. \nl
2) When for $\ell = 1,2$ we are given $(\bar M^\ell,\bold{\bar a}^\ell) \in 
K^{\text{sq}}_{\lambda^+}$ we say $(\bar M^1,\bold{\bar a}^1) 
\le_{\text{ct}} (\bar M^2,
\bold{\bar a}^2)$ \ub{if} for some club $E$ of $\lambda^+$, letting
$\bold{\bar a}^\ell = \langle a^\ell_\delta:\delta \in S^\ell \rangle$ for
$\ell =1,2$, of course, we have
\mr
\item "{$(a)$}"  $S^1 \cap E \subseteq S^2 \cap E$
\sn
\item "{$(b)$}"  if $\delta \in S^1 \cap E$ then
{\roster
\itemitem{ $(\alpha)$ }  $M^1_\delta \le_{\frak K} M^2_\delta$,
\sn
\itemitem{ $(\beta)$ }   $M^1_{\delta +1} \le_{\frak K} M^2_{\delta +1}$
\sn
\itemitem{ $(\gamma)$ }  $a^2_\delta = a^1_\delta$
\sn
\itemitem{ $(\delta)$ }  \ortp$(a^1_\delta,M^2_\delta,M^2_{\delta +1})$ does
not fork over $M^1_\delta$, so in particular $a^1_\delta \notin M^2_\delta$.
\endroster}
\endroster
\enddefinition
\bigskip

\demo{\stag{600-nu.3A} Observation}  1) If $(\bar M,\bold{\bar a}) \in
K^{\text{sq}}_{\lambda^+}$ then $M := \dbcu_{\alpha < \lambda^+} M_\alpha \in
K_{\lambda^+}$ is saturated. \nl
2) $K^{\text{sq}}_{\lambda^+}$ is partially ordered by
$\le_{\text{ct}}$.
\hfill$\square_{\scite{600-nu.3A}}$
\enddemo
\bigskip

\proclaim{\stag{600-nu.4} Claim}  Assume 
$\langle (\bar M^\zeta,\bold{\bar a}^\zeta):\zeta < \zeta^* \rangle$ 
is $\le_{\text{ct}}$-increasing in $K^{\text{sq}}_{\lambda^+}$,
and $\zeta^*$ is a limit ordinal $< \lambda^{++}$, \ub{then} the 
sequence has a $\le_{\text{ct}}$-{\rm lub} $(\bar M,\bold{\bar a})$.
\endproclaim
\bigskip

\demo{Proof}  Let $\bold{\bar a}^\zeta = \langle a^\zeta_\delta:\delta \in
S_\zeta \rangle$ for $\zeta < \zeta^*$ and \wilog \, $\zeta^* = \text{ cf}
(\zeta^*)$ and for $\zeta < \xi < \zeta^*$ let $E_{\zeta,\xi}$ be a club of
$\lambda^+$ witnessing $(\bar M^\zeta,\bold{\bar a}^\zeta) \le_{\text{ct}}
(\bar M^\xi,\bold{\bar a}^\xi)$.
\enddemo
\bn
\ub{Case 1}:  $\zeta^* < \lambda^+$.

Let $E = \cap\{E_{\zeta,\xi}:\zeta < \xi < \zeta^*\}$ and for $\delta \in E$
let $M_\delta = \cup\{M^\zeta_\delta:\zeta < \zeta^*\}$ and
$M_{\delta +1} = \cup\{M^\zeta_{\delta +1}:\zeta < \zeta^*\}$ and for any
other $\alpha,M_\alpha = M_{\text{Min}(E \backslash \alpha)}$.  Let $S =
\dbcu_{\zeta < \zeta^*} S_\zeta \cap E$ and for $\delta \in S$ let 
$a_\delta = a^\zeta_\delta$ for every $\zeta$ for which 
$\delta \in S_\zeta$.  Clearly
$M_\alpha \in K_\lambda$ is $\le_{\frak K}$-increasing continuous
and $\zeta < \zeta^* \wedge 
\delta \in E \Rightarrow M^\zeta_\delta \le_{\frak K} M_\delta \and
M^\zeta_{\delta +1} \le_{\frak K} M_{\delta +1}$.

Now if $\delta \in E \cap S_\zeta$ then $\xi \in [\zeta,\zeta^*)$ implies
\ortp$(a_\delta,M^\xi_\delta,M_{\delta+1}) = \text{ \ortp}(a^\zeta_\delta,
M^\xi_\delta,M^\xi_{\delta+1})$ does not fork over $M^\zeta_\delta$
(and $\langle M^\xi_\delta:\xi \in [\zeta,\delta)\rangle,\langle
M^\xi_{\delta+1}:\xi \in [\zeta,\delta)\rangle$ are $\le_{\frak
K}$-increasing continuous); hence
by Axiom (E)(h) we know that \ortp$(a_\delta,M_\delta,M_{\delta +1})$ does not
fork over $M^\zeta_\delta$ and in particular $\in {\Cal
S}^{\text{bs}}(M_\delta)$.   Also if $N \le_{\frak K} M := \dbcu_{\alpha <
\lambda^+} M_\alpha,N \in K_\lambda$ and $p \in {\Cal S}^{\text{bs}}(N)$ then
for some $\delta(*) \in E,N \le_{\frak K} 
M_{\delta(*)}$, let $p_1 \in {\Cal S}
^{\text{bs}}(M_{\delta(*)})$ be a non-forking 
extension of $p$, so for some
$\zeta < \zeta^*,p$ does not fork over $M^\zeta_{\delta(*)}$ hence for
stationarily many $\delta \in S_\zeta,q^0_\delta = \text{ \ortp}(a_\delta,
M^\zeta_\delta,M^\zeta_{\delta +1})$ is a non-forking extension of $p_1
\restriction M^\zeta_{\delta(*)}$, hence this holds for stationarily
many $\delta \in S \cap E$ and for each such $\delta,
q_\delta = \text{ \ortp}(a_\delta,
M_\delta,M_{\delta +1})$ is a non-forking extension of $p_1 \restriction
M^\zeta_{\delta(*)}$, hence of $p_1$ hence of $p$.  Looking at the
definitions, clearly $(\bar M,\bar{\bold a}) \in
K^{\text{sq}}_{\lambda^+}$ and $\zeta < \zeta^* \Rightarrow (\bar
M^\zeta,\bar{\bold a}^\zeta) \le_{\text{ct}} (\bar M,\bar{\bold a})$.  

Lastly, it is easy to check the $\le_{\text{ct}}$-e.u.b.
\bn
\ub{Case 2}:  $\zeta^* = \lambda^+$.

Similarly using diagonal union, i.e., $E = \{\delta <
\lambda^+:\delta$ is a limit ordinal such that $\zeta < \xi < \delta
\Rightarrow \delta \in E_{\zeta,\varepsilon}\}$ and we choose $M_\alpha =
\cup\{M^\zeta_\alpha:\zeta < \alpha\}$ when $\alpha \in E$ and
$M_\alpha = M_{\text{min}(E \backslash (\alpha +1))}$ otherwise. 
\hfill$\square_{\scite{600-nu.4}}$
\bigskip

\demo{\stag{600-nu.13.1} Observation}  Assume $K^{3,\text{uq}}_\lambda$ is
dense in $K^{3,\text{bs}}_\lambda$, i.e., in
$(K^{3,\text{bs}}_\lambda,\le_{\text{bs}})$ and even in 
$(K^{3,\text{bt}}_\lambda,<_{\text{bt}})$.  \ub{Then}
\mr
\item "{$(a)$}"  if $M \in K_\lambda$ is superlimit and $p \in {\Cal
S}^{\text{bs}}(M)$ then there are $N,a$ such that $(M,N,a) \in 
K^{3,\text{uq}}_\lambda$ and $p = \text{\rm \ortp}(a,M,N)$
\sn
\item "{$(b)$}"  if in addition $K_{\frak s}$ is categorical (in
$\lambda$) \ub{then} ${\frak s}$ has existence for
$K^{3,\text{uq}}_\lambda$ (recall that this means 
that for every $M \in K_{\frak s}$ and 
$p \in {\Cal S}^{\text{bs}}(M)$ for some pair $(N,a)$ we have
$(M,N,a) \in K^{3,\text{uq}}_\lambda$ and $p = \text{\rm \ortp}(a,M,N))$.
\endroster
\enddemo
\bigskip

\demo{Proof}  Should be clear.
\enddemo
\bn
Now the assumption of \scite{600-nu.13.1} are justified by the following
theorem (and the categoricity in (b) is justified by Claim \scite{600-0.34}).
\proclaim{\stag{600-nu.6} First Main Claim}  Assume that
\mr
\item "{$(a)$}"  as in \scite{600-nu.0}
\sn
\item "{$(b)$}"  ${\text{\rm WDmId\/}}
(\lambda^+)$ is not $\lambda^{++}$-saturated and
$2^\lambda < 2^{\lambda^+} < 2^{\lambda^{++}}$ (or the parallel versions
for the definitional weak diamond).
\ermn
If $\dot I(\lambda^{++},K) < 
2^{\lambda^{++}}$, \ub{then} for every $(M,N,a) \in
K^{3,\text{bs}}_\lambda$ there is $(M^*,N^*,a) \in
K^{3,\text{bt}}_\lambda$ such that
$(M,N,a) <_{\text{bt}} (M^*,N^*,a)$ and 
$(M^*,N^*,a) \in K^{3,\text{uq}}_\lambda$.
\endproclaim
\bigskip

\demo{\stag{600-nu.6.2} Explanation}  The reader who agrees to believe in
\scite{600-nu.6} can ignore the rest of this section (though it can still
serve as a good exercise).

Let $\langle S_\alpha:\alpha < \lambda^{++}\rangle$ be a sequence of
subsets of $\lambda^+$ such that $\alpha < \beta \Rightarrow |S_\alpha
\backslash S_\beta| \le \lambda$ and $S_{\alpha +1} \backslash
S_\alpha \ne \emptyset$ mod WDmId$(\lambda^+)$, exists by assumption.

Why having $(M,N,a)$ failing the conclusion of \scite{600-nu.6} helps us to
construct many models in $K_{\lambda^{++}}$?  The point is that we can
choose $(\bar M^\alpha,\bar{\bold a}^\alpha) \in
K^{\text{sq}}_{\lambda^+}$ with Dom$(\bar{\bold a}^\alpha) = S_\alpha$
for $\alpha <\lambda^{++},<_{\text{ct}}$-increasing continuous 
(see \scite{600-nu.4}).

Now for $\alpha=\beta +1$, having $(\bar M^\beta,\bar{\bold
a}^\beta)$, \wilog \, $M^\beta_{i+1}$ is brimmed over $M^\beta_i$ 
 and we shall choose $M^\alpha_i$
by induction on $i < \lambda^+$ (for simplicity we pretend $M^\alpha_i
\cap \cup\{M^\beta_j:j < \lambda^+\} = M^\beta_i$) and $M^\beta_i
\le_{\frak K} M^\alpha_i \in K_\lambda$ and
\ortp$(a_i,M^\alpha_i,M^\alpha_{i+1})$ does not fork over $M^\beta_i$ and
$M^\alpha_{i+1}$ is brimmed over $M^\alpha_i$).

Given $(\bar M^\beta,\bar{\bold a}^\beta),\bar M^\beta = \langle
M^\beta_i:i < \lambda^+\rangle,\bar{\bold a}^\beta_i$ toward building
$(\bar M^\alpha,\bar{\bold a}^\alpha),\alpha_{\beta +1}$. 

We start with choosing $(M^\alpha_0,b)$ such that no member of
$K^{3,\text{bs}}_\lambda$ which is $\le_{\text{bs}}$-above
$(M^\beta_0,M^\alpha_0,b) \in K^{3,\text{bs}}_\lambda$ belongs to
$K^{3,\text{uq}}_\lambda$ and will choose $M^\beta_i$ by induction on
$i$ such that $(M^\beta_i,M^\alpha_i,b) \in K^{3,\text{bs}}_\lambda$
is $\le_{\text{bs}}$-increasing continuous and even
$<_{\text{bt}}$-increasing hence in particular that
\ortp$(b,M^\beta_i,M^\alpha_i)$ does not fork over $M^\alpha_0$.  Now in
each stage $i=j+1$, as $M^\beta_i$ is universal over $M^\beta_j$, and
the choice of $M^\alpha_0,b$ we have some freedom.  So it makes sense
that we will have many possible outcomes, i.e., models $M =
\cup\{M^\alpha_i:\alpha < \lambda^{++},i < \lambda^+\}$ which are
in $K_{\lambda^{++}}$.   The combination of what we have above and
\cite[\S3]{Sh:576} gives
that $2^\lambda < 2^{\lambda^+} < 2^{\lambda^{++}}$ is enough to
materialize this intuition.  If in addition $2^\lambda = \lambda^+$
and moreover $\diamondsuit_{\lambda^+}$ it is considerably easier.  In
the end we still have to define $\bar{\bold a}^\alpha \restriction
(S_\alpha \backslash S_\beta)$ as required in Definition \scite{600-nu.3}.

Alternatively when \cite{Sh:839} materializes, it will do it in more
transparent way.  Another alternative is to force a model in
$\lambda^{++}$.  Now below we replace $K^{3,\text{sq}}_{\lambda^+}$ by
$K^{\text{mqr}}_{\lambda^+},K^{\text{nqr}}_S$ but
actually $K^{3,\text{sq}}_{\lambda^+}$ is enough,
but not in the way \cite{Sh:576} is done.  So we need a somewhat more
complicated relative as elaborated below which anyhow seems to me more natural.
\enddemo
\bigskip

\proclaim{\stag{600-nu.6A} Second Main Claim}  Assume 
$2^\lambda < 2^{\lambda^+} < 2^{\lambda^{++}}$ 
(or the parallel versions for the definitional weak 
diamond). If $\dot I(\lambda^{++},K) < \mu_{\text{wd}}(\lambda^{++},
2^{\lambda^+})$, \ub{then} for 
every $(M,N,a) \in K^{3,\text{bt}}_\lambda$ there is $(M^*,N^*,a) \in 
K^{3,\text{bt}}_\lambda$ such that $(M,N,a) <_{\text{bt}} (M^*,N^*,a)$ and 
$(M^*,N^*,a) \in K^{3,\text{uq}}_\lambda$.
\endproclaim
\bn
We shall not prove here \scite{600-nu.6A} and shall not use it; toward
proving \scite{600-nu.6} (by quoting) let
\definition{\stag{600-nu.7} Definition}  Let $S \subseteq \lambda^+$ be a
stationary subset of $\lambda^+$. \nl
1)  Let $K^{\text{mqr}}_S$ or $K^{\text{mqr}}_{\lambda^+}[S]$ 
be the set of triples 
$(\bar M,\bold {\bar a},\bold f)$ such that:
\mr
\item "{$(a)$}"  $\bar M = \langle M_\alpha:\alpha < \lambda^+ \rangle$ is
$\le_{\frak K}$-increasing continuous, $M_\alpha \in K_\lambda$ \nl
(we denote $\dsize \bigcup_{\alpha < \lambda^+} M_\alpha$ by $M$) and
demand $M \in K_{\lambda^+}$
\sn
\item "{$(b)$}"  $\bold {\bar a} = \langle a_\alpha:\alpha < \lambda \rangle$
with $a_\alpha \in M_{\alpha + 1}$
\sn
\item "{$(c)$}"  $\bold f$ is a function from $\lambda^+$ to $\lambda^+$ 
such that for some club $E$ of $\lambda^+$ for every 
$\delta \in E \cap S$ and ordinal 
$i < \bold f(\delta)$ we have \ortp$(a_{\delta +i},M_{\delta +i},
M_{\delta +i+ 1}) \in {\Cal S}^{\text{bs}}(M_{\delta +i})$
\sn
\item "{$(d)$}"  for every $\alpha < \lambda^+$ and $p \in 
{\Cal S}^{\text{bs}}(M_\alpha)$,  stationarily many 
$\delta \in S$ satisfies: for some $\varepsilon < \bold f(\delta)$ we have
\ortp$(a_{\delta + \varepsilon},M_{\delta + \varepsilon},
M_{\delta + \varepsilon +1})$ is a 
non-forking extension of $p$.
\ermn
1A) $K^{\text{nqr}}_{\lambda^+}[S] = 
K^{\text{nqr}}_S$ is the set of triples $(\bar M,
\bar{\bold a},\bold f) \in K^{\text{mqr}}_S$ such that:
\mr
\item "{$(e)$}"  for a club of $\delta < \lambda^+$, if $\delta \in S$ then
$\bold f(\delta)$ is divisible by $\lambda$ and \footnote{if we have an
a priori bound $\bold f^*:\lambda^+ \rightarrow \lambda^+$ which is
a $<_{{\Cal D}_{\lambda^+}}$-upper bound of the ``first"
$\lambda^{++}$ functions in ${}^{\lambda^+}(\lambda^+)/D$, we can use
bookkeeping for $u_i$'s as in the proof of \scite{600-4a.8}}
for every $i < \bold f(\delta)$ if $q \in {\Cal S}^{\text{bs}}
(M_{\delta + i})$ then for $\lambda$ ordinals $\varepsilon \in [i,
\bold f(\delta))$ the type 
\ortp$(a_{\delta + \varepsilon},M_{\delta + \varepsilon},
M_{\delta + \varepsilon +1}) \in {\Cal S}^{\text{bs}}
(M_{\delta + \varepsilon})$ is a stationarization of $q$ (=
non-forking extension of $q$, see Definition \scite{600-4a.3}).
\ermn
2)  Assume $(\bar M^\ell,\bold {\bar a}^\ell,
\bold f^\ell) \in K^{\text{mqr}}_S$ for 
$\ell =1,2$; we say 
$(\bar M^1,\bold{\bar a}^1,\bold f^1) \le^0_S (\bar M^2,\bold{\bar a}^2,
\bold f^2)$ \underbar{iff} for some club $E$ of $\lambda^+$, for every
$\delta \in E \cap S$ we have:
\medskip
\roster
\item "{$(a)$}"  $M^1_{\delta +i} \le_{\frak K} M^2_{\delta + i}$ for
\footnote{could have used (systematically) $i < \bold f^1(\delta)$}
$i \le \bold f^1(\delta)$
\sn
\item "{$(b)$}"  $\bold f^1(\delta) \le \bold f^2(\delta)$
\sn
\item "{$(c)$}"  for $i < \bold f^1(\delta)$ we have 
$a^1_{\delta+i} = a^2_{\delta+i}$ and \newline
\ortp$(a^1_{\delta+i},M^2_{\delta +i},M^2_{\delta +i+1})$ 
does not fork over $M^1_{\delta +i}$.
\ermn
3) We define the relation $<^1_S$ on $K^{\text{mqr}}_S$ as in part (2) adding
\mr
\item "{$(d)$}"  if $\delta \in E$ and $i < \bold f^1(\delta)$ then
$M^2_{\delta +i+1}$ is $(\lambda,*)$-brimmed over $M^1_{\delta +i+1} \cup
M^2_{\delta +i}$.
\endroster
\enddefinition
\bigskip

\proclaim{\stag{600-nu.8} Claim}  0) If $(\bar M,\bar{\bold a},\bold f) \in
K^{\text{mqr}}_S$ \ub{then} 
$\dbcu_{\alpha < \lambda^+} M_\alpha \in K_{\lambda^+}$ is
saturated. \nl
1) The relation $\le^0_S$ is a quasi-order \footnote{quasi order $\le$
is a transitive relation, so we waive $x \le y \le x \Rightarrow x=y$}
on $K^{\text{mqr}}_\lambda$; also $<^1_S$ is. \nl
2) $K^{\text{mqr}}_S \supseteq K^{\text{nqr}}_S 
\ne \emptyset$ for any stationary $S
\subseteq \lambda^+$. \nl
3) For every $(\bar M,\bar{\bold a},\bold f) \in 
K^{\text{mqr}}_\lambda[S]$ for 
some $(\bar M',\bar{\bold a},\bold f') \in K^{\text{nqr}}_\lambda[S]$ we have
$(\bar M,\bar{\bold a},\bold f) <^1_S (\bar M',\bar{\bold a},\bold f')$. \nl
4) For every $(\bar M^1,\bold{\bar a}^1,\bold f^1) \in K^{\text{mqr}}_S$ and
$q \in {\Cal S}^{\text{bs}}(M^1_\alpha),\alpha < \lambda^+$, 
\ub{there is} $(M^2,\bold{\bar a}^2,\bold f^2) \in K^{\text{mqr}}_S$ 
such that $(\bar M^1,\bold{\bar a}^1,\bold f^1) <^1_S
(\bar M^2,\bold{\bar a}^2,\bold f^2) \in K^{\text{nqr}}_S$ 
and $b \in M^2_\alpha$ realizing $q$ such that for 
every $\beta \in [\alpha,\lambda^+)$ we have 
{\rm \ortp}$(b,M^1_\beta,M^2_\beta) \in {\Cal S}^{\text{bs}}
(M^1_\beta)$ does not fork over $M^1_\alpha$. \nl
5) If $\langle(\bar M^\zeta,\bold{\bar a}^\zeta,\bold f^\zeta):\zeta <
\xi(*) \rangle$ is $\le^0_S$-increasing continuous in $K^{\text{mqr}}_S$ and
$\xi(*) < \lambda^{++}$ a limit ordering, \ub{then} the sequence has 
a $\le^0_S$-{\rm lub}.
\endproclaim
\bigskip

\demo{Proof}  0, 1)  Easy. 
\nl
2) The inclusion $K^{\text{mqr}}_S \supseteq K^{\text{nqr}}_S$ is
obvious, so let us prove $K^{\text{nqr}}_S \ne \emptyset$.
We choose by induction on 
$\alpha < \lambda^+,a_\alpha,M_\alpha,p_\alpha$ such that
\mr
\item "{$(a)$}"  $M_\alpha \in K_\lambda$ is a super limit model,
\sn
\item "{$(b)$}"  $M_\alpha$ is $\le_{\frak K}$-increasingly continuous,
\sn
\item "{$(c)$}"  if $\alpha = \beta +1$, then $a_\beta \in M_\alpha 
\backslash M_\beta$ realizes $p_\beta \in {\Cal S}^{\text{bs}}(M_\beta)$,
\sn
\item "{$(d)$}"  if 
$p \in {\Cal S}^{\text{bs}}(M_\alpha)$, then for some $i <
\lambda$, for every $j \in [i,\lambda)$ for at least one ordinal
$\varepsilon \in [j,j + i),p_{\alpha + \varepsilon} \restriction 
M_\alpha = p$ and $p_{\alpha + \varepsilon}$ does not fork over $M_\alpha$.
\ermn
For $\alpha = 0$ choose $M_0 \in K_\lambda$.  For $\alpha$ limit,
$M_\alpha = \dbcu_{\beta < \alpha} M_\beta$ is as required.  
Then use Axiom(E)(g) to take care of clause (d) (with careful
bookkeeping). 
Lastly, let $\bold f:\lambda^+ \rightarrow \lambda^+$ be constantly
$\lambda,\bar M = \langle M_\alpha:\alpha < \lambda \rangle,\bar{\bold
a} = \langle a_\alpha:\alpha < \lambda \rangle$; now for any
stationary $S \subseteq \lambda^+$, the triple $(\bar M,\bar{\bold a}
\restriction S,\bold f \restriction S)$ belong to $K^{\text{nqr}}_S$.
\nl
3) Let $E$ be a club witnessing $(\bar M^1,\bold{\bar a}^1,\bold f^1) \in
K^{\text{mqr}}_S$ such that 
$\delta \in E \Rightarrow \delta + \bold f^1(\delta) <
\text{ Min}(E \backslash (\delta +1))$.  Choose $\bold f^2:\lambda^+
\rightarrow \lambda^+$ such that $\alpha < \lambda^+$ implies
$\bold f^1(\alpha) < \bold f^2(\alpha) < \lambda^+$ and $\bold
f^2(\alpha)$ is divisible by $\lambda$.  We
choose by induction on $\alpha < \lambda^+,f_\alpha,M^2_\alpha,p_\alpha,
a^2_\alpha$ such that:
\mr
\widestnumber\item{$(a),(b),(c)$}
\item "{$(a),(b),(c)$}"  as in the proof of part (2)
\sn
\item "{$(d)$}"  $f_\alpha$ is a $\le_{\frak K}$-embedding of 
$M^1_\alpha$ into $M^2_\alpha$
\sn
\item "{$(e)$}"  $f_\alpha$ is increasing continuous
\sn
\item "{$(f)$}"  if $\delta \in E \cap S$ and $i < \bold f^1(\delta)$ hence
\ortp$(a^1_{\delta +i},M^1_{\delta +i},M^1_{\delta +i+1}) \in 
{\Cal S}^{\text{bs}}(M^1_{\delta +i})$, \ub{then} $f_{\delta +i+1}
(a^1_{\delta +i}) = 
a^2_{\delta +i}$ and $p_{\varepsilon +i} = \text{ \ortp}(a^2_{\delta +i},
M^2_{\delta +i},M^2_{\delta + i+1}) \in {\Cal S}^{\text{bs}}
(M^2_{\delta +i})$
is a stationarization of \ortp$\bigl(f_{\delta + i+1}(a^1_{\delta +i}),
f_{\delta +i}
(M^1_{\delta+i}),f_{\delta +i+1}(M^1_{\delta +i+1}) \bigr) =$
\nl
$\text{ \ortp}(a^2_{\delta +i},
f_{\delta +i}(M^1_{\delta +i}),M^2_{\delta +i+1})$
\sn
\item "{$(g)$}"  if $\delta \in E$ and $i < \bold f^2(\delta),
q \in {\Cal S}^{\text{bs}}(M^2_{\delta +i})$ then for some 
$\lambda$ ordinals $\varepsilon \in (i,\bold f^2(\delta))$ the type 
$p_{\delta + \varepsilon}$ is a stationarization of $q$
\sn
\item "{$(h)$}"  if $\delta \in E,i < \bold f^2(\delta)$ then
$M_{\delta +i+1}$ is $(\lambda,*)$-brimmed over $M_{\delta +i} \cup
f_{\delta +i+1}(M^1_{\delta +i+1})$.
\ermn
The proof is as in part (2) only the bookkeeping is different.
At the end without loss of generality 
$\dbcu_{\alpha < \lambda^*} f_\alpha$ is the identity and we are done.  \nl
4) Similar proof but in some cases we have to use 
Axiom (E)(i), the non-forking amalgamation  
of Definition \scite{600-1.1}, in the appropriate cases. \nl
5)  Without loss of generality cf$(\xi(*)) = \xi(*)$.  First
assume that $\xi(*) \le \lambda$.  For $\varepsilon < \zeta < \xi(*)$ let
$E_{\varepsilon,\zeta}$ be a club of $\lambda^+$ witnessing 
$\bar M^\varepsilon <^0_S \bar M^\zeta$.  Let \newline
$E^* = \dsize \bigcap_{\varepsilon < \zeta < \xi(*)} E_{\varepsilon,\zeta} 
\cap \{ \delta < \lambda^+:\text{for every } \alpha < \delta \text{ we
have } 
\underset{\varepsilon < \xi(*)} {}\to \sup \bold f^\varepsilon(\alpha) < 
\delta\}$, it is a club of $\lambda^+$.
Let $\bold f^{\xi(*)}:\lambda^+ \rightarrow \lambda^+ \text{ be }
\bold f^{\xi(*)}(i) = \underset{\varepsilon < \xi(*)} {}\to \sup
\bold f^\varepsilon(i)$ now define $M^{\xi(*)}_i$ as follows:
\mn
\ub{Case 1}:  If $\delta \in E^*$ and $\varepsilon < \xi(*)$ and $i \le 
\bold f^\varepsilon(\delta)$ and $i \ge \dbcu_{\zeta < \varepsilon}
\bold f^\zeta(\delta)$ then
\mr
\item "{$(\alpha)$}"  $M^{\xi(*)}_{\delta +i} = \bigcup 
\bigl\{ M^\zeta_{\delta+i}:\zeta \in [\varepsilon,\xi(*)) \bigr\}$
\sn
\item "{$(\beta)$}"  $i < \bold f^\varepsilon(\delta) 
\Rightarrow a^{\xi(*)}_{\delta+i} = a^\varepsilon_{\delta+i}$.
\ermn
(Note: we may define $M^{\xi(*)}_{\delta +i}$ twice if $i = \bold
f^\varepsilon(\delta)$, but the two values are the same).
\sn
\ub{Case 2}:  If $\delta \in E^*,i = \bold f^{\xi(*)}(\delta)$ is a limit 
ordinal let

$$
M^{\xi(*)}_{\delta +i} = \dsize \bigcup_{j<i} M^{\xi(*)}_{\delta+i}.
$$
\mn
\ub{Case 3}:  If $M^{\xi(*)}_i$ has not been defined yet, let it be 
$M^{\xi(*)}_{\text{Min}(E^* \backslash i)}$.
\mn
\ub{Case 4}:  If $a^{\xi(*)}_i$ has not been defined yet, let 
$a^{\xi(*)}_i \in M^{\xi(*)}_{i+1}$ be arbitrary.
\medskip

Note that Case 3,4 deal with the ``unimportant" cases. \nl
Let $\varepsilon < \xi(*)$, why 
$(\bar M^\varepsilon,\bold{\bar a}^\varepsilon,
\bold f^\varepsilon) \le^0_S (\bar M^{\xi(*)},\bold{\bar a}^{\xi(*)},
\bold f^{\xi(*)}) \in K^{\text{mqr}}_S$?  Enough to check
that the club $E^*$ witnesses it. \nl
Why \ortp$(a_{\delta +i},M^{\xi(*)}_{\delta +i},M^{\xi(*)}_{\delta +i+1}) \in
{\Cal S}^{\text{bs}}(M^{\xi(*)}_{\delta +i})$ and when $\delta \in E^*,i < 
\bold f^{\xi(*)}(i)$, and does not fork over $M^\varepsilon_{\delta +i}$ when
$i < \bold f^\varepsilon(\delta)$ ? by Axiom (E)(h) of Definition 
\scite{600-1.1}. \nl
Why clause (e) of Definition \scite{600-nu.7}(1A)?  By Axiom (E)(c), local
character of non-forking.

The case $\xi(*) = \lambda^+$ is similar using diagonal intersections.
\hfill$\square_{\scite{600-nu.8}}$
\enddemo
\bigskip

\remark{Remark}  If we use weaker versions of ``good
$\lambda$-frames", we should 
systematically concentrate on successor $i < \bold f(\delta)$.
\endremark
\bigskip

\demo{Proof of \scite{600-nu.6}}  The use of $\lambda^{++} \notin \text{
WDmId}(\lambda^{++})$ is as in the proof of
\cite[3.19]{Sh:576}(pg.79)=3.12t.   But
now we need to preserve saturation in limit stages $\delta <
\lambda^{++}$ of cofinality $< \lambda^+$, we use $<^1_S$, otherwise
we act as in \cite[\S3]{Sh:576}.  \hfill$\square_{\scite{600-nu.6}}$
\enddemo
\bn
Let us elaborate
\definition{\stag{600-nu.11} Definition}  We define $\bold C = 
({\frak K}^+,\bold S eq,\le^*)$ as follows:
\mr
\item "{$(a)$}"  $\tau^+ = \tau \cup \{P,<\},{\frak K}^+$ is the set of
$(M,P^M,<^M)$ where $M \in {\frak K}_{< \lambda},P^M \subseteq M,<^M$ 
a linear ordering of $P^M$ (but $=^M$ may be as in
\cite[3.1]{Sh:576}(2) and $M_1 \le_{{\frak K}^+} M_2$ 
iff $(M_1 \restriction \tau) \le_{\frak K} (M_2
\restriction \tau)$ and $M_1 \subseteq M_2$
\sn
\item "{$(b)$}"  ${\bold Seq}_\alpha = \{\bar M:\bar M = \langle M_i:i \le
\alpha \rangle$ is an increasing continuous sequence of members of
${\frak K}^+$ and $\langle M_i \restriction \tau:i \le \alpha \rangle$ is
$\le_{\frak K}$-increasing, and for \newline
$i < j < \alpha:P^{M_i}$ is a proper initial segment of 
$(P^{M_j},<^{M_j})$ and there is a first element in the difference$\}$ \nl
we denote 
the $<^{M_{i+1}}$-first element of $P^{M_{i+1}} \backslash P^{M_i}$,
by $a_i[\bar M]$ and we demand \ortp$(a_i(\bar M),M_i \tau
\restriction,M_{i+1} \restriction \tau) \in {\Cal S}^{\text{bs}}(M_i
\restriction \tau)$ and if $\alpha = \lambda,M = \cup\{M_i
\restriction \tau:i < \lambda^+\}$ is saturated
\sn
\item "{$(c)$}"  $\bar M <^*_t \bar N$ \underbar{iff} \newline
$\bar M = \langle M_i:i < \alpha^* \rangle,\bar N = \langle N_i:
i < \alpha^{**} \rangle$ are from ${\bold Seq},t$ is a set of pairwise 
disjoint closed intervals of $\alpha^*$ and for any $[\alpha,\beta] \in t$ 
we have $(\beta < \alpha^*$ and): 
\sn
$\gamma \in [\alpha,\beta) \Rightarrow M_\gamma \le_{\frak K} N_\gamma \and
a_\gamma[\bar M] \notin N_\gamma$, moreover \newline
$a_\gamma[\bar M] = a_\gamma[\bar N]$ and \ortp$(a_j[\bar M],N_\gamma
\restriction \tau,N_{\gamma+1},\tau)$ does not fork over $M_\gamma
\restriction \tau$.
\endroster
\enddefinition
\bigskip

\proclaim{\stag{600-nu.12} Claim}  1) $\bold C$ is a
$\lambda^+$-construction framework (see \cite[3.3]{Sh:576}(pg.68). \nl
2) $\bold C$ is weakly nice (see Definition \cite[3.14]{Sh:576}(2)(pg.76). \nl
4) $\bold C$ has the weakening $\lambda^+$-coding property.
\endproclaim
\bn
\ub{Discussion}:  Is it better to use (see \cite[3.14]{Sh:576}(1)(pg.75))
stronger axiomatization in \cite[\S3]{Sh:576} to cover this? \nl
But at present this will be the only case.
\bigskip

\demo{Proof}  Straight.  \hfill$\square_{\scite{600-nu.12}}$
\enddemo
\bn
Now \scite{600-nu.6A} follows by \cite[3.19]{Sh:576}(pg.79).
\newpage

\head {\S6 Non-forking amalgamation in ${\frak K}_\lambda$} \endhead  \resetall \sectno=6
 \spuriousreset
\bigskip

We deal in this section only with ${\frak K}_\lambda$. \newline
We would like to, at least, approximate ``non-forking amalgamation of
models" using as a starting point the conclusion of \scite{600-nu.6}, i.e., 
$K^{3,\text{uq}}_\lambda$ is dense.
We use what looks like a stronger hypothesis: the existence for
$K^{3,\text{uq}}_\lambda$ (also called ``weakly successful"); 
but in our application we can assume
categoricity in $\lambda$; the point being that as we have a
superlimit $M \in K_\lambda$ this assumption is reasonable when we restrict
ourselves to ${\frak K}^{[M]}$, recalling that we believe in first
analyzing the saturated enough models; see \scite{600-nu.13.1}.
By \scite{600-4a.6}, the ``$(\lambda,\text{cf}
(\delta))$-brimmed over" is the same for all limit ordinals $\delta < 
\lambda^+$, (but not for $\delta = 1$); nevertheless for possible
generalizations we do not use this.
\bigskip

It may help the reader to note, that if there is a 4-place relation
NF$_\lambda(M_0,M_1,M_2,M_3)$ on $K_\lambda$, satisfying the expected
properties of ``$M_1,M_2$ are amalgamated in a non-forking = free way over
$M_0$ inside $M_3$", i.e., is a ${\frak K}_\lambda$-non-forking
relation from Definition \scite{600-nf.0X} 
below then Definition \scite{600-nf.2} below (of NF$_\lambda$)
gives it (provably!).  So we have ``a definition" of
NF$_\lambda$ satisfying that: if desirable non-forking relation
exists, our definition gives it (assuming the hypothesis \scite{600-nf.0}).  So
during this section we are trying to get better and better
approximations to the desirable properties; have the feeling of going
up on a spiral, as usual.

For the readers who know on non-forking in stable first order theory
we note that in such context NF$_\lambda(M_0,M_1,M_2,M_3)$ says that
\ortp$(M_2,M_1,M_3)$, the type of $M_2$ over $M_1$ inside $M_3$, does not
fork over $M_0$.  It is natural to say that there are $\langle
N_{1,\alpha},N_{2,\alpha}:\alpha \le \alpha^*\rangle,N_{\ell,\alpha}$
is increasing continuous.  $N_{1,0} = M_0,N_{2,0} = M_2,M_1 \subseteq
M_{1,\alpha},M_3 \subseteq M'_3,N_{2,\alpha} \subseteq
M'_3,N_{\ell,\alpha +2}$ is prime over $N_{\ell,\alpha} + a_\alpha$
for $\ell=1,2$ and \ortp$(a_\alpha,N_{2,\alpha})$ does not fork over
$N_{1,\alpha}$ but this is not available.  The
$K^{3,\text{uq}}_\lambda$ is a substitute.
\bigskip

\definition{\stag{600-nf.0X} Definition}  1) Assume that ${\frak K} =
{\frak K}_\lambda$ is a $\lambda$-a.e.c.  
We say NF is a non-forking relation on
${}^4({\frak K}_\lambda)$ or just a ${\frak K}_\lambda$-non-forking
relation \ub{when}:
\mr
\item "{$\boxtimes_{\text{NF}}(a)$}"  NF is a 4-place relation on
$K_\lambda$ and NF is preserved under isomorphisms  
\sn
\item "{$(b)$}"  NF$(M_0,M_1,M_2,M_3)$ implies $M_0 \le_{\frak K}
M_\ell \le_{\frak K} M_3$ for $\ell = 1,2$ 
\sn
\item "{$(c)_1$}"  (monotonicity):  if NF$(M_0,M_1,M_2,M_3)$ and $M_0
\le_{\frak K} M'_\ell \le_{\frak K} M_\ell$ 
for $\ell = 1,2$ then NF$(M_0,M'_1,M'_2,M_3)$
\sn
\item "{$(c)_2$}"  (monotonicity):  if NF$(M_0,M_1,M_2,M_3)$ and
$M_3 \le_{\frak K} M'_3 \in K_\lambda,M_1 \cup M_2 \subseteq M''_3
\le_{\frak K} M'_3$ \ub{then} NF$(M_0,M_1,M_2,M''_3)$
\sn
\item "{$(d)$}"  (symmetry)   NF$(M_0,M_1,M_2,M_3)$ iff
NF$(M_0,M_2,M_1,M_3)$
\sn
\item "{$(e)$}"  (transitivity)  if NF$(M_i,N_i,M_{i+1},N_{i+1})$ for
$i < \alpha,\langle M_i:i \le \alpha \rangle$ is $\le_{\frak
K}$-increasing continuous and $\langle N_i:i \le \alpha \rangle$ is
$\le_{\frak K}$-increasing continuous \ub{then} 
\nl
NF$(M_0,N_0,M_\alpha,N_\alpha)$
\sn
\item "{$(f)$}"  (existence) if $M_0 \le_{\frak K} M_\ell$ for $\ell
=1,2$ (all in $K_\lambda)$ \ub{then} for some $M_3 \in K_\lambda,f_1,
f_2$ we have
$M_0 \le_{\frak K} M_3,f_\ell$ is a $\le_{\frak K}$-embedding of
$M_\ell$ into $M_3$ over $M_0$ for $\ell = 1,2$ and 
NF$(M_0,f_1(M_1),f_2(M_2),M_3)$
\sn
\item "{$(g)$}"  (uniqueness) if
NF$(M^\ell_0,M^\ell_1,M^\ell_2,M^\ell_3)$ and for $\ell =
1,2$ and $f_i$ is an isomorphism from $M^1_i$ onto $M^2_i$ for $i =
0,1,2$ and $f_0 \subseteq f_1,f_0 \subseteq f_2$ \ub{then} $f_1 \cup
f_2$ can be extended to an embedding $f_3$ of $M^1_3$ into some
$M^2_4,M^2_3 \le_{{\frak K}_\lambda} M^2_4$.
\ermn
2) We say that NF is a weak non-forking relation on ${}^4(K_\lambda)$
or a weak ${\frak K}_\lambda$-non-forking relation
if clauses (a)-(f) of $\boxtimes_{\text{NF}}$ above holds but not
necessarily clause (g).
\nl
3) Assume ${\frak s}$ is a good $\lambda$-frame and NF is a
non-forking relation on ${\frak K}$ or just a formal one.  
We say that NF respects ${\frak s}$ or NF is an ${\frak
s}$-non-forking relation \ub{when}:
\mr
\item "{$(h)$}"  if NF$(M_0,M_1,M_2,M_3)$ and $a \in M_2 \backslash
M_0$, \ortp$_{\frak s}(a,M_0,M_2) \in {\Cal S}^{\text{bs}}(M_0)$ then
\ortp$_{\frak s}(a,M_1,M_3)$ 
does not fork over $M_0$ in the sense of ${\frak s}$.
\endroster
\enddefinition
\bigskip

\demo{\stag{600-nf.0X.2} Observation}  Assume ${\frak K}_\lambda$ is a
$\lambda$-a.e.c. and NF is a non-forking relation on ${}^4({\frak
K}_\lambda)$.
\nl
1) Assume ${\frak K}$ is stable in $\lambda$.  If in clause (g) of
\scite{600-nf.0X}(1) above we assume in addition that $M^\ell_3$ is
$(\lambda,\sigma)$-brimmed over $M^\ell_1 \cup M^\ell_2$, \ub{then} in
the conclusion of (g) we can add $M^2_3 = M^2_4$, i.e., $f_1 \cup f_2$
can be extended to an isomorphism from $M^1_3$ onto $M^2_3$.  This
version of (g) is equivalent to it (assuming stability in $\lambda$;
note that ``${\frak K}_\lambda$ has amalgamation" follows by clause (h)
of Definition \scite{600-nf.0X}).
\nl
2) If $M_0 \le_{\frak K} M_1 \le_{\frak K} M_3$ are from $K_\lambda$
then NF$(M_0,M_0,M_1,M_3)$.
\nl
3) In Definition \scite{600-nf.0X}(1), clause (d), symmetry, it is enough
to demand ``if".
\enddemo
\bigskip

\demo{Proof}  1) Chase arrows.
\nl
2) By clause (f) of $\boxtimes_{\text{NF}}$ of \scite{600-nf.0X}(1) and
clause (c)$_2$, i.e., first apply existence 
with $(M_0,M_0,M_2)$ here standing for $(M_0,M_1,M_2)$
there, then chase arrows and use the monotonicity in (c)$_2$. 
\nl
3) Easy.   \hfill$\square_{\scite{600-nf.0X.2}}$
\enddemo
\bn
The main point of the
following claim shows that there is at most one non-forking
relation respecting ${\frak s}$; so it justifies the definition of
NF$_{\frak s}$ later.  The assumption ``NF respects ${\frak s}$" is
not so strong by \scite{600-nf.0.z.1}.
\proclaim{\stag{600-nf.0Y} Claim}  1) If ${\frak s}$ is a
good $\lambda$-frame and 
{\rm NF} is a non-forking relation on ${}^4({\frak K}_{\frak s})$
respecting ${\frak s}$ and $(M_0,N_0,a) \in K^{3,\text{uq}}_\lambda$ and 
$(M_0,N_0,a) \le_{\text{bs}} (M_1,N_1,a)$ then
${\text{\rm NF\/}}(M_0,N_0,M_1,N_1)$. \nl
2) If ${\frak s}$ is a good $\lambda$-frame, weakly successful (which means
$K^{3,\text{uq}}_{\frak s}$ has 
existence in $K^{3,\text{uq}}_{\frak s}$, i.e.,
${\frak s}$ satisfies hypothesis \scite{600-nf.0} below) and
{\rm NF} is a non-forking relation on ${}^4({\frak K}_{\frak s})$
respecting ${\frak s}$ \ub{then} the relation 
{\rm NF}$_\lambda = \text{\rm NF}_{\frak s}$, i.e., 
$\nonforkin{N_1}{N_2}_{N_0}^{N_3}$ defined in Definition \scite{600-nf.2} below is
equivalent to ${\text{\rm NF\/}}(N_0,N_1,N_2,N_3)$. 
[But see \scite{600-nf.20.7A}] 
\nl
3) If ${\frak s}$ is a weakly successful
good $\lambda$-frame and for $\ell=1,2$, the
relation {\rm NF}$_\ell$ is a non-forking relation on ${}^4({\frak
K}_{\frak s})$ respecting ${\frak s}$, 
\ub{then} ${\text{\rm NF\/}}_1 = { \text{\rm NF\/}}_2$.
\endproclaim
\bigskip

\demo{Proof}  Straightforward but we elaborate.
\nl
1) We can find $(M'_1,N'_1)$ such that NF$(M_0,N_0,M'_1,N'_1)$ and
$M_1,M'_1$ are isomorphic over $M_0$, say $f_1$ is such an
isomorphism from $M_1$ onto $M'_1$ over $M_0$; why such
$(M'_1,N'_1,f_1)$ exists? by
clause (f) of $\boxtimes_{\text{NF}}$ of Definition \scite{600-nf.0X}.

As NF respects ${\frak s}$, see Definition \scite{600-nf.0X}(2), recalling
\ortp$(a,M_0,N_0) \in {\Cal S}^{\text{bs}}(M_0)$ we know
that \ortp$_{\frak s}(a,M'_1,N'_1)$ does not fork over $M_0$, 
so by the definition of
$\le_{\text{bs}}$ we have $(M_0,N_0,a) \le_{\text{bs}}
(M'_1,N'_1,a)$.

As $(M_0,N_0,a) \in K^{3,\text{uq}}_\lambda$, by the definition of
$K^{3,\text{uq}}_\lambda$ (and chasing arrows) we conclude that there
are $N_2,f_2$ such that:
\mr
\item "{$(*)$}"   $N_1 \le_{{\frak K}[{\frak s}]} N_2 \in K_\lambda$ and $f_2$
is a $\le_{\frak K}$-embedding of $N'_1$ into $N_2$ extending
$f^{-1}_1$ and id$_{N_0}$.
\ermn
As NF$(M_0,N_0,M'_1,N'_1)$ and NF is preserved under isomorphisms (see
clause (a) in \scite{600-nf.0X}(1)) it follows that
NF$(M_0,N_0,M_1,f_2(N'_1))$.  By the monotonicity of NF (see clause
$(c)_2$ of Definition \scite{600-nf.0X}) it follows that
NF$(M_0,N_0,M_1,N_2)$.  Again
by the same monotonicity we have NF$(M_0,N_0,M_1,N_1)$, as required.
\nl
2) First we prove that NF$_{\lambda,\bar \delta}(N_0,N_1,N_2,N_3)$,
which is defined in Definition \scite{600-nf.1} 
below implies NF$(N_0,N_1,N_2,N_3)$.  By
definition \scite{600-nf.1}, clause (f) there are $\langle
(N_{1,i},N_{2,i}:i \le \lambda \times \delta_1 \rangle),\langle c_i:i
< \lambda \times \delta_1 \rangle$ as there.  Now we prove by
induction on $j\le \lambda \times \delta_1$ that $i \le j \Rightarrow
\text{ NF}(N_{1,i},N_{2,i},N_{1,j},N_{2,j})$.  For $j=0$ or more
generally when $i=j$ this is trivial by \scite{600-nf.0X.2}(2).  
For $j$ a limit ordinal use the induction hypothesis and
transitivity of NF (see clause (e) of \scite{600-nf.0X}(1)).

Lastly, for $j$ successor by the demands in Definition \scite{600-nf.1} we
know that $N_{1,j-1} \le_{\frak K} N_{1,j} \le_{\frak K}
N_{2,j},N_{1,j-1} \le_{\frak K} N_{2,j-1} \le_{\frak K} N_{2,j}$ are
all in $K_\lambda$, \ortp$(c_{j-1},N_{2,j-1},N_{2,j})$ does not fork over
$N_{1,j-1}$ and $(N_{1,j-1},N_{1,j},c_{j-1}) \in
K^{3,\text{uq}}_\lambda$.  By part (1) of this claim we deduce that
NF$(N_{1,j-1},N_{1,j},N_{2,j-1},N_{2,j})$ hence by symmetry (i.e.,
clause (d) of Definition \scite{600-nf.0X}(1)) we deduce
NF$(N_{1,j-1},N_{2,j-1},N_{1,j},N_{2,j})$.

So we have gotten $i < j \Rightarrow$
NF$(N_{1,i},N_{2,i},N_{1,j},N_{2,j})$.  \nl
[Why?  If $i=j-1$ by the previous sentence and for $i < j-1$ note that
by the induction hypothesis NF$(N_{1,i},N_{2,i},N_{1,j-1},N_{1,j-1})$
so by transitivity (clause (e) of \scite{600-nf.0X}(1) of Definition
\scite{600-nf.0X}) we get NF$(N_{1,i},N_{2,i},N_{1,j},N_{2,j})$].

We have carried the induction so in particular for $i=0,j=\alpha$ we get
 NF$(N_{1,0},N_{2,0},N_{1,\alpha},N_{2,\alpha})$ which means
 NF$(N_0,N_1,N_2,N_3)$  as promised.  So we have proved NF$_{\lambda,\bar
\delta}(N_0,N_1,N_2,N_3) \Rightarrow$ NF$(N_0,N_1,N_2,N_3)$.

Second, if NF$_\lambda(N_0,N_1,N_2,N_3)$ as defined in Definition
\scite{600-nf.2} then there are
$M_0,M_1,M_2,M_3 \in K_\lambda$ such that
NF$_{\lambda,\langle \lambda,\lambda\rangle}
(M_0,M_1,M_2,M_3),N_\ell \le_{\frak K} M_\ell$ for $\ell <
4$ and $N_0 = M_0$.  By what we have proved above we can conclude
NF$(M_0,M_1,M_2,M_3)$.  As $N_0 = M_0 \le_{\frak K} N_\ell \le_{\frak
K} M_\ell$ for $\ell=1,2$ by clause $(c)_1$ of Definition
\scite{600-nf.0X}(1) we get NF$(M_0,N_1,N_2,M_3)$ and by clause $(c)_2$ of
Definition \scite{600-nf.0X}(1) we get NF$(N_0,N_1,N_2,N_3)$.  So we have
proved the implication NF$_\lambda(N_0,N_1,N_2,N_3) \Rightarrow$
NF$(N_0,N_1,N_2,N_3)$.   Using this, the equivalence
follows by the existence, uniqueness and monotonicity.
\nl
3) By the rest of this section, i.e., the main conclusion
\scite{600-nf.20.7}, the relation NF$_\lambda$ defined in \scite{600-nf.2} is
a non-forking relation on ${}^4(K_{\frak s})$ respecting ${\frak s}$.
Hence by part (2) of the present claim we have NF$_1 = \text{ NF}_\lambda =
\text{ NF}_2$.  \hfill$\square_{\scite{600-nf.0Y}}$
\enddemo
\bn
\margintag{600-nf.0Z}\ub{\stag{600-nf.0Z} Example}:  1) Do we need ${\frak s}$ in \scite{600-nf.0Y}(3)?
Yes.

Let ${\frak K}$ be the class of graphs and $M \le_{\frak K} N$ iff $M
\subseteq N$; so ${\frak K}$ is an a.e.c. with LS$({\frak K}) =
\aleph_0$.   For cardinal $\lambda$ and
$\ell=1,2$ we define NF$^\ell =
\{(M_0,M_1,M_2,M_3):M_0\le_{\frak K} M_1 \le_{\frak K} M_3$ and $M_0
\le_{\frak K} M_2 \le_{\frak K} M_3$ and $M_1 \cap M_2 =M_0$ and if $a
\in M_1 \backslash M_0,b \in M_2 \backslash M_0$ then $\{a,b\}$ is an
edge of $M_3$ iff $\ell=2\}$ and NF$^\ell_\lambda := \{(M_0,M_1,M_2,M_3)
\in \text{ NF}:M_0,M_1,M_2,M_3 \in K_\lambda\}$.  
Then NF$^\ell_\lambda$ is a non-forking
relation on ${}^4({\frak K}_\lambda)$ but NF$^1_\lambda \ne$ NF$^2_\lambda$.
\nl
2) So the assumption on ${\frak K}_\lambda$ that for some good
$\lambda$-frame ${\frak s}$ we have ${\frak K}_{\frak s} = 
{\frak K}_\lambda$ is quite a strong demand on ${\frak K}_\lambda$.

However, the assumption ``respect" essentially is not necessary as it can be
deduced when ${\frak s}$ is good enough.
\bn
Below on ``good$^+$" see \sectioncite[\S1]{705} in particular Definition
\marginbf{!!}{\cprefix{705}.\scite{705-stg.1}}.
\proclaim{\stag{600-nf.0.z.1} Claim}  Assume that ${\frak s}$ is a good$^+
\lambda$-frame and {\rm NF} is a non-forking relation on
${}^4({\frak K}_{\frak s})$.  \ub{Then} {\rm NF} respects ${\frak s}$.
\endproclaim
\bigskip

\remark{Remark}  The construction in the proof is similar to the ones
in \scite{600-4a.7}, \scite{600-nf.4}.
\endremark
\bigskip

\demo{Proof}  Assume NF$(M_0,M_1,M_2,M_3)$ and $a \in M_2 \backslash
M_0$, \ortp$_{\frak s}(a,M_0,M_2) \in {\Cal S}^{\text{bs}}_{\frak s}(M_0)$.  We
define $(N_{0,i},N_{1,i},f_i)$ for $i < \lambda^+_{\frak s}$ as
follows:
\mr
\item "{$\otimes_1$}"  $(a) \quad N_{0,i}$ is $\le_{\frak
s}$-increasing continuous and $N_{0,0} = M_0$
\sn
\item "{${{}}$}"  $(b) \quad N_{1,i}$ is $\le_{\frak s}$-increasing
continuous and $N_{1,0} = M_1$
\sn
\item "{${{}}$}"  $(c) \quad$
NF$(N_{0,i},N_{1,i},N_{0,i+1},N_{1,i+1})$
\sn
\item "{${{}}$}"  $(d) \quad f_i$ is a $\le_{\frak K}$-embedding of
$M_2$ into $N_{0,i+1}$ over $M_0 = N_{0,0}$ such that 
\nl

\hskip25pt \ortp$_{\frak s}(f_i(a),N_{0,i},N_{0,i+1})$ 
does not fork over $M_0 = N_{0,0}$.
\ermn
We shall choose $f_i$ together with $N_{0,i+1},N_{1,i+1}$.
\nl
Why can we define?  For $i=0$ there is nothing to do.  For $i$ limit
take unions.  For $i=j+1$ choose $f_j,N_{0,i}$ satisfying clause (d)
and $N_{0,j} \le_{\frak s} N_{0,i}$, this is possible for 
${\frak s}$ as we have the existence of
non-forking extensions of \ortp$_{\frak s}(a,M_0,M_2)$ (and
amalgamation).

Lastly, we take care of the rest 
(mainly clause (c) of $\otimes_1$ by clause (f) of Definition \scite{600-nf.0X}(1),
existence).  Now
\mr
\item "{$\circledast_2$}"  for $i<j<\lambda^+$ we have
NF$(N_{0,i},N_{1,i},N_{0,j},N_{1,j})$
\nl
[why?  by transitivity for NF, i.e., clause (e) of Definition
\scite{600-nf.0X}(1), transitivity]
\sn
\item "{$\circledast_3$}"  for some $i$, \ortp$_{\frak s}(f_i(a),N_{1,i},
N_{1,i+1})$ does not fork over $M_0$
\nl
[why?  by the definition of good$^+$].
\ermn
So for this $i,M_0 \le_{\frak s} f_i(M_2) \le_{\frak s} N_{0,i+1}$ by clause
(d) of $\otimes_1$, hence by clause $(c)_1$ of Definition
\scite{600-nf.0X}, monotonicity
we have NF$(M_0,M_1,f_i(M_2),N_{1,i+1})$.  Now again by the choice of
$i$, i.e., by $\circledast_3$ we have 
\ortp$_{\frak s}(f_i(a),M_1,N_{1,i+1})$ does 
not fork over $M_0$.  By clause (g) of Definition \scite{600-nf.0X}(1), i.e.,
uniqueness of NF (and preservation by isomorphisms)
we get \ortp$_{\frak s}(a,M_1,M_3)$ does not fork over $M_0$ as required.
\hfill$\square_{\scite{600-nf.0.z.1}}$
\enddemo
\bn
We turn to our main task in this section proving that such NF exist;
till \scite{600-nf.20.7} we assume:
\demo{\stag{600-nf.0} Hypothesis}  1) ${\frak s} = ({\frak K},\nonfork{}{}_{},
{\Cal S}^{\text{bs}})$ is a good $\lambda$-frame. \nl
2) ${\frak s}$ is weakly successful which just means that it has
existence for $K^{3,\text{uq}}_\lambda$: for every $M \in
K_\lambda$ and $p \in {\Cal S}^{\text{bs}}(M)$
there are $N,a$ such that $(M,N,a) \in K^{3,\text{uq}}_\lambda$ 
(see Definition \scite{600-nu.1A}) and $p = \text{ \ortp}(a,M,N)$.  (This
follows by $K^{3,\text{uq}}_{\frak s}$ is dense in
$K^{3,\text{bs}}_{\frak s}$; when ${\frak s}$ is categorical, see
\scite{600-nu.13.1}.) 

In this section we deal with models from $K_\lambda$ only.
\enddemo
\bigskip

\proclaim{\stag{600-nf.0A} Claim}  If $M \in K_\lambda$ and $N$ is 
$(\lambda,\kappa)$-brimmed over $M$, \ub{then} we can find 
$\bar M = \langle M_i:i \le \delta \rangle,\le_{\frak K}$-increasing
continuous, $(M_i,M_{i+1},c_i) \in
K^{3,\text{uq}}_\lambda,M_0 = M,M_\delta = N$ and $\delta$ 
any pregiven limit ordinal $< \lambda^+$ of cofinality $\kappa$ 
divisible by $\lambda$.
\endproclaim
\bigskip

\demo{Proof}  Let $\delta$ be given, e.g.,  
$\delta = \lambda \times \kappa$.  By 
\scite{600-nf.0}(2) we can find a $\le_{\frak K}$-increasing sequence $\langle
M_i:i \le \delta \rangle$ of members of $K_\lambda$ and
$\langle a_i:i < \delta \rangle$ such that 
$M_0 = M$ and $i < \delta \Rightarrow (M_i,M_{i+1},a_i) \in
K^{3,\text{uq}}_\lambda$ and for every $i < \delta,p \in {\Cal
S}^{\text{bs}}(M_i)$ for $\lambda$ ordinals $j \in (i,i+ \lambda)$ we have
\ortp$(a_j,M_j,M_{j+1})$ is a non-forking extension of $p$.
So the demands in \scite{600-4a.2} hold hence $M_\delta$
is $(\lambda,\kappa)$-brimmed over $M_0=M$.  Now we are done by the
uniqueness of $N$ being $(\lambda,\kappa)$-brimmed over $M_0$, see
\scite{600-0.22}(3).  \hfill$\square_{\scite{600-nf.0A}}$
\enddemo
\bigskip

\proclaim{\stag{600-nf.0B} Claim}   If 
$M^\ell_0 \le_{\frak K} M^\ell_1 \le_{\frak K} M^\ell_3$ and
$M^\ell_0 \le_{\frak K} M^\ell_2 \le_{\frak K} M^\ell_3,c_\ell \in M^\ell_1$
and $(M^\ell_0,M^\ell_1,c_\ell) \in K^{3,\text{uq}}_\lambda$ 
and ${\text{\rm \ortp\/}}(c_\ell,M^\ell_2,M^\ell_3) \in 
{\Cal S}^{\text{bs}}(M^\ell_2)$ does not fork over $M^\ell_0$ and
$M^\ell_3$ is $(\lambda,\sigma)$-brimmed over $M^\ell_1 \cup
M^\ell_2$ all this for 
$\ell=1,2$ and $f_i$ is an isomorphism from $M^1_i$ onto $M^2_i$ for
$i=0,1,2$ such that $f_0 \subseteq f_1,f_0 \subseteq f_2$ and $f_1(c_1) =
c_2$, \ub{then} $f_1 \cup
f_2$ can be extended to an isomorphism from $M^1_3$ onto $M^2_3$. 
\endproclaim
\bigskip

\demo{Proof}  Chase arrows (and recall definition of
$K^{3,\text{uq}}_\lambda$), that is by \scite{600-nf.0X}(1) and Definition
\scite{600-nf.0X.2}(1) and \scite{600-0.22}(3).  \hfill$\square_{\scite{600-nf.0B}}$
\enddemo
\bigskip

\definition{\stag{600-nf.1} Definition}  Assume $\bar \delta = \langle \delta_1,
\delta_2,\delta_3 \rangle,\delta_1,\delta_2,\delta_3$ 
are ordinals $< \lambda^+$, maybe 1.  We say that 
NF$_{\lambda,\bar \delta}(N_0,N_1,N_2,N_3)$ or, in other words $N_1,N_2$ are 
\ub{brimmedly smoothly amalgamated} in $N_3$ over $N_0$ for 
$\bar \delta$ when:
\mr
\item "{$(a)$}"  $N_\ell \in K_\lambda$ for $\ell \in \{ 0,1,2,3\}$
\sn
\item "{$(b)$}"  $N_0 
\le_{\frak K} N_\ell \le_{\frak K} N_3$ for $\ell = 1,2$
\sn
\item "{$(c)$}"  $N_1 \cap N_2 = N_0$ (i.e. in disjoint amalgamation,
actually follows by clause (f))
\sn
\item "{$(d)$}"  $N_1$ is ($\lambda$,cf$(\delta_1)$)-brimmed over
$N_0$; recall that if 
cf$(\delta_1)=1$ this just means $N_0 \le_{\frak K} N_1$
\sn
\item "{$(e)$}"  $N_2$ is ($\lambda$,cf$(\delta_2)$)-brimmed over
$N_0$; so that if cf$(\delta_2) =1$ this just 
means $N_0 \le_{\frak K} N_2$ and
\sn
\item "{$(f)$}"  there are $N_{1,i},N_{2,i}$ for $i \le \lambda \times
\delta_1$ and $c_i$ for $i < \lambda \times \delta_1$
(called witnesses and $\langle N_{1,i},N_{2,i},c_j:i \le \lambda \times
\delta_1,j < \lambda \times \delta_1 \rangle$ is called 
a witness sequence as well as $\langle
N_{1,i}:i \le \lambda \times \delta_1 \rangle,\langle N_{2,i}:i \le
\lambda \times \delta_1 \rangle$) such that:
{\roster
\itemitem{ $(\alpha)$ }  $N_{1,0} = N_0,N_{1,\lambda \times \delta_1} = N_1$
\sn 
\itemitem{ $(\beta)$ }  $N_{2,0} = N_2$
\sn
\itemitem{ $(\gamma)$ }  $\langle N_{\ell,i}:i \le \lambda \times 
\delta_1 \rangle$ is a $\le_{\frak K}$-increasing continuous sequence of
models for $\ell = 1,2$
\sn
\itemitem{ $(\delta)$ }  $(N_{1,i},N_{1,i+1},c_i) \in 
K^{3,\text{uq}}_\lambda$ 
\sn
\itemitem{ $(\varepsilon)$ }  \ortp$(c_i,N_{2,i},N_{2,i+1}) 
\in {\Cal S}^{\text{bs}}
(N_{2,i})$ does not fork over $N_{1,i}$ and $N_{2,i} \cap N_1 =
N_{1,i}$, for $i < \lambda \times \delta_1$ 
(follows by Definition \scite{600-nu.1A})
\sn
\itemitem{ $(\zeta)$ }  $N_3$ is $(\lambda$,cf$(\delta_3)$)-brimmed
over $N_{2,\lambda \times \delta_1}$; so for cf$(\delta_3) 
= 1$ this means just $N_{2,\lambda \times \delta_1} \le_{\frak K} N_3$
\endroster}
\endroster
\enddefinition
\bigskip

\definition{\stag{600-nf.2} Definition}  1) We say
$\nonforkin{N_1}{N_2}_{N_0}^{N_3}$ (or
$N_1,N_2$ are \ub{smoothly amalgamated} over $N_0$ inside $N_3$ or
NF$_\lambda(N_0,N_1,N_2,N_3)$ or NF$_{\frak s}(N_0,N_1,N_2,N_3)$) 
\ub{when} we can find $M_\ell \in
K_\lambda$ (for $\ell < 4$) such that:
\medskip
\roster
\item "{$(a)$}"  NF$_{\lambda,\langle \lambda,\lambda,\lambda \rangle}
(M_0,M_1,M_2,M_3)$
\sn
\item "{$(b)$}"  $N_\ell \le_{\frak K} M_\ell$ for $\ell < 4$
\sn
\item "{$(c)$}"  $N_0 = M_0$
\sn
\item "{$(d)$}"  $M_1,M_2$ are $(\lambda,\text{cf}(\lambda))$-brimmed 
over $N_0$ (follows by (a) see clauses (d), (e) of \scite{600-nf.1}).
\ermn
2)  We call $(M,N,a)$ \ub{strongly bs-reduced} if 
$(M,N,a) \in K^{3,\text{bs}}_\lambda$
and $(M,N,a) \le_{\text{bs}} 
(M',N',a) \in K^{3,\text{bs}}_\lambda \Rightarrow 
\text{ NF}_\lambda(M,N,M',N')$; not used.
\enddefinition
\bn
Clearly we expect ``strongly bs-reduced" to be 
equivalent to ``$\in K^{3,\text{uq}}_\lambda$", e.g. as this occurs in the
first order case.
We start by proving existence for NF$_{\lambda,\bar \delta}$ from
Definition \scite{600-nf.1}.
\proclaim{\stag{600-nf.3} Claim}  1) Assume 
$\bar \delta = \langle \delta_1,\delta_2,
\delta_3 \rangle,\delta_\ell$ an ordinal $< \lambda^+$ and
$N_\ell \in K_\lambda$ for $\ell < 3$ and $N_1$ is 
$(\lambda,{\text{\rm cf\/}}(\delta_1)$)-brimmed over $N_0$ and $N_2$ is 
$(\lambda,{\text{\rm cf\/}}(\delta_2)$)-brimmed
over $N_0$ and $N_0 \le_{\frak K} N_1$ and $N_0 \le_{\frak K} N_2$ 
and for simplicity
$N_1 \cap N_2 = N_0$.  \underbar{Then} we can find $N_3$ such that
${\text{\rm NF\/}}_{\lambda,\bar \delta}(N_0,N_1,N_2,N_3)$. \newline
2) Moreover, we can choose any $\langle N_{1,i}:i \le \lambda \times
\delta_1 \rangle,\langle c_i:i < \lambda \times \delta_1\rangle$ 
as in \scite{600-nf.1} subclauses (f)$(\alpha),
(\gamma),(\delta)$ as part of the witness.
\nl
3) If {\rm NF}$_\lambda(N_0,N_1,N_2,N_3)$ then $N_1 \cap N_2 = N_0$.
\endproclaim
\bigskip

\demo{Proof}  1) We can find $\langle N_{1,i}:i \le \lambda \times
\delta_1 \rangle$ and $\langle c_i:i < \lambda \times \delta_1 \rangle$
as required in part (2) by Claim \scite{600-nf.0A},
the $(\lambda,\text{cf}(\lambda \times \delta_1))$-brimness holds by
\scite{600-4a.2} and apply part (2). \nl
2) We choose the $N_{2,i}$ (by induction on $i$) by \scite{600-4a.7}
preserving $N_{2,i} \cap N_{1,\lambda \times \delta_2} = N_{1,i}$; in
the successor case use Definition \scite{600-nu.1A} + Claim
\scite{600-nu.2}(3).  We then choose $N_3$ using \scite{600-4a.1}(2).
\nl
3) By the definitions of NF$_\lambda$, NF$_{\lambda,\bar \delta}$.
   \hfill$\square_{\scite{600-nf.3}}$
\enddemo
\bn
The following claim tells us that if we have
``$(\lambda,\text{cf}(\delta_3))$-brimmed" in the end, then we can have
it in all successor stages. 
\proclaim{\stag{600-nf.4} Claim}  In Definition 
\scite{600-nf.1}, if $\delta_3$ is a
limit ordinal and $\kappa = { \text{\rm cf\/}}(\kappa) \ge \aleph_0$, 
\ub{then} without loss of generality
(even without changing $\langle N_{1,i}:i \le \lambda \times
\delta_1\rangle,\langle c_i:i < \lambda \times \delta_1\rangle)$
\mr
\item "{$(g)$}"  $N_{2,i+1}$ is $(\lambda,\kappa)$-brimmed over 
$N_{1,i+1} \cup N_{2,i}$  (which means that it is \newline
$(\lambda,\kappa)$-brimmed over some $N$, 
where $N_{1,i+1} \cup N_{2,i} \subseteq N \le_{\frak K} N_{2,i+1}$).
\endroster
\endproclaim
\bigskip

\demo{Proof}  So assume NF$_{\lambda,\bar \delta}(N_0,N_1,N_2,N_3)$ holds
as being witnessed by $\langle N_{\ell,i}:i \le \lambda \times \delta_1 
\rangle,\langle c_i:i < \lambda \times \delta_1 \rangle$
for $\ell = 1,2$.  Now we choose by induction on $i \le \lambda \times
\delta_1$ a model $M_{2,i} \in K_\lambda$ and $f_i$ such that:
\mr
\widestnumber\item{ (iii) }
\item "{$(i)$}"  $f_i$ is a $\le_{\frak K}$-embedding of
$N_{2,i}$ into $M_{2,i}$
\sn
\item "{$(ii)$}"  $M_{2,0} = f_i(N_2)$
\sn
\item "{$(iii)$}"  $M_{2,i}$ is $\le_{\frak K}$-increasing continuous and
also $f_i$ is increasing continuous
\sn
\item "{$(iv)$}"  $M_{2,j} \cap f_i(N_{1,i}) = 
f_i(N_{1,j})$ for $j \le i$
\sn
\item "{$(v)$}"  $M_{2,i+1}$ is $(\lambda,\kappa)$-brimmed over
$M_{2,i} \cup f_i(N_{2,i+1})$
\sn
\item "{$(vi)$}"  \ortp$(f_{i+1}(c_i),M_{2,i},M_{2,i+1}) \in 
{\Cal S}^{\text{bs}}(M_{2,i})$ does not fork over $f_i(N_{1,i})$.
\ermn
There is no problem to carry the induction.  Using in the successor case
$i=j+1$ the existence Axiom (E)(g) of 
Definition \scite{600-1.1} there is a model $M'_{2,i} \in K_{\frak s}$ 
such that $M_{2,j} \le_{\frak K} M'_{2,i}$ and $f_i
\supseteq f_j$ as required in clauses (i), (iv), (vi) and then use Claim
\scite{600-4a.1} to find a model $M_{2,i} \in K_\lambda$ which is
$(\lambda,\kappa)$-brimmed over $M_{2,j} \cup f_i(N_{2,i})$.

Having carried the induction, \wilog \, $f_i = \text{ id}_{N_{2,i}}$.
Let $M_3$ be such that
$M_{2,\lambda \times \delta_1} \le_{\frak K} M_3 \in K_\lambda$ and 
$M_3$ is $(\lambda$,cf$(\delta_3)$)-brimmed over $M_{2,\lambda \times 
\delta_1}$, it exists by \scite{600-4a.1}(2) but $N_{2,\lambda \times
\delta_1} \le_{\frak K} M_{2,\lambda \times \delta_1}$, hence it
follows that $M_3$ is $(\lambda,\kappa)$-brimmed over $N_{1,\lambda
\times \delta_1}$.
So both $M_3$ and $N_3$ are $(\lambda$,cf$(\delta_3)$)-brimmed over
$N_{2,\lambda \times \delta_1}$, hence they are isomorphic over
$N_{2,\lambda \times \delta_1}$ (by \scite{600-0.22}(1)) so 
let $f$ be an isomorphism from $M_3$
onto $N_3$ which is the identity over $N_{2,\lambda \times \delta_1}$.
\newline
Clearly $\langle N_{1,i}:i \le \lambda \times \delta_1 \rangle,
\langle f(M_{2,i}):i \le \lambda \times \delta_1 \rangle$ are also 
witnesses for \newline
NF$_{\lambda,\bar \delta}(N_0,N_1,N_2,N_3)$ satisfying the extra demand
$(g)$ from \scite{600-nf.4}.  \hfill$\square_{\scite{600-nf.4}}$
\enddemo
\bn
The point of the following claim is that having uniqueness in every
atomic step we have uniqueness in the end (using the same ``ladder"
$N_{1,i}$ for now).
\proclaim{\stag{600-nf.5} Claim}  (Weak Uniqueness).

Assume that for $x \in \{ a,b\}$, we have 
${\text{\rm NF\/}}_{\lambda,\bar \delta^x}
(N^x_0,N^x_1,N^x_2,N^x_3)$ holds as witnessed by
$\langle N^x_{1,i}:i \le \lambda \times \delta^x_1 \rangle,
\langle c^x_i:i < \lambda \times \delta^x_1 \rangle,
\langle N^x_{2,i}:i \le \lambda \times \delta^x_1 \rangle$ 
and $\delta_1 := \delta^a_1 = \delta^b_1,{\text{\rm cf\/}}
(\delta^a_2) = { \text{\rm cf\/}}(\delta^b_2)$ and 
${\text{\rm cf\/}}(\delta^a_3) = { \text{\rm cf\/}}
(\delta^b_3) \ge \aleph_0$. \nl
(Note that ${\text{\rm cf\/}}(\lambda \times \delta^a_1) \ge \aleph_0$ 
by the definition of {\rm NF}).
\medskip

Suppose further that $f_\ell$ is an isomorphism from $N^a_\ell$ onto
$N^b_\ell$ for $\ell = 0,1,2$, moreover: $f_0 \subseteq f_1,f_0 \subseteq
f_2$ and $f_1(N^a_{1,i}) = N^b_{1,i},f_1(c^a_i) = c^b_i$.

\underbar{Then} we can find an isomorphism $f$ from $N^a_3$ onto $N^b_3$
extending $f_1 \cup f_2$.
\endproclaim
\bigskip

\demo{Proof}  Without loss of generality for each $i < \lambda \times
\delta_1$, the model 
$N^x_{2,i+1}$ is $(\lambda,\lambda)$-brimmed over 
$N^x_{1,i+1} \cup N^x_{2,i}$ (by \scite{600-nf.4}, 
note there the statement ``without changing the $N_{1,i}$'s").  Now we 
choose by induction on $i \le \lambda \times \delta_1$ an isomorphism 
$g_i$ from 
$N^a_{2,i}$ onto $N^b_{2,i}$ such that: $g_i$ is increasing with $i$ and 
$g_i$ extends $(f_1 \restriction N^a_{1,i}) \cup f_2$. \newline
For $i = 0$ choose $g_0 = f_2$ and for $i$ limit let $g_i$ be
$\dsize \bigcup_{j < i} g_j$ and for $i = j+1$ it exists by \scite{600-nf.0B},
whose assumptions hold by $(N^x_{1,i},N^x_{1,i+1},c^x_i)
\in K^{3,\text{uq}}_\lambda$ (see \scite{600-nf.1}, clause (f)$(\delta)$) and 
the extra brimness clause from \scite{600-nf.4}.  Now by \scite{600-0.22}(3) 
we can extend 
$g_{\lambda \times \delta_1}$ to an isomorphism from 
$N^a_3$ onto $N^b_3$ as $N^x_3$ is $(\lambda,\text{cf}(\delta_3))$-brimmed 
over $N^x_{2,\lambda \times \delta_1}$ (for $x \in \{ a,b\}$). 
\hfill$\square_{\scite{600-nf.5}}$
\enddemo
\bn
Note that even
knowing \scite{600-nf.5} the choice of $\langle N_{1,i}:i \le \lambda \times 
\delta_1 \rangle,\langle c_i:i < \lambda \times \delta_1 \rangle$ 
still possibly matters.  Now we prove an ``inverted" 
uniqueness, using our ability to construct a ``rectangle" of models
which is a witness for NF$_{\lambda,\bar \delta}$ in two ways.
\proclaim{\stag{600-nf.6} Claim}  Suppose that 
\mr
\item "{$(a)$}"  for $x \in \{ a,b\}$ we have 
${\text{\rm NF\/}}_{\lambda,\bar \delta^x}(N^x_0,N^x_1,N^x_2,N^x_3)$
\sn
\item "{$(b)$}"  $\bar \delta^x = \langle \delta^x_1,\delta^x_2,\delta^x_3
\rangle,\delta^a_1 = \delta^b_2$, 
$\delta^a_2 = \delta^b_1,{\text{\rm cf\/}}(\delta^a_3) = 
{ \text{\rm cf\/}}(\delta^b_3)$, all limit ordinals
\sn
\item "{$(c)$}"  $f_0$ is an isomorphism from $N^a_0$ onto $N^b_0$
\sn
\item "{$(d)$}"  $f_1$ is an isomorphism from $N^a_1$ onto $N^b_2$
\sn
\item "{$(e)$}"  $f_2$ is an isomorphism from $N^a_2$ onto $N^b_1$
\sn
\item "{$(f)$}"  $f_0 \subseteq f_1$ and $f_0 \subseteq f_2$.
\endroster
\medskip
\noindent
\underbar{Then} there is an isomorphism from $N^a_3$ onto $N^b_3$ extending
$f_1 \cup f_2$.
\endproclaim
\bn
Before proving we shall construct a third ``rectangle" of models such
that we shall be able to construct appropriate isomorphisms each of
$N^a_3,N^b_3$ 
\proclaim{\stag{600-nf.6A} Subclaim}  Assume
\mr
\item "{$(a)$}"  $\delta^a_1,\delta^a_2,\delta^a_3 < \lambda^+$ are limit
ordinals
\sn
\item "{$(b)_1$}"  $\bar M^1 = \langle M^1_\alpha:\alpha \le \lambda \times
\delta^a_1 \rangle$ is $\le_{\frak K}$-increasing continuous in $K_\lambda$
\nl
and $(M^1_\alpha,M^1_{\alpha +1},c_\alpha) \in K^{3,\text{bs}}_\lambda$
\sn
\item "{$(b)_2$}"  $\bar M^2 = \langle M^2_\alpha:\alpha \le \lambda \times
\delta^a_2 \rangle$ is $\le_{\frak K}$-increasing continuous in
$K_\lambda$ and $(M^2_\alpha,M^2_{\alpha +1},d_\alpha) 
\in K^{3,\text{bs}}_\lambda$
\sn
\item "{$(c)$}"  $M^1_0 = M^2_0$ 
we call it $M$ and $M^1_\alpha \cap M^2_\beta = M$ for $\alpha \le
\lambda \times \delta^a_1,\beta \le \lambda \times \delta^a_2$.
\ermn
\ub{Then} we can find $M_{i,j}$ (for $i \le \lambda \times
\delta^a_1 \text{ and } j \le \lambda \times \delta^a_2)$ 
and $M_3$ such that:
\mr
\item "{$(A)$}"  $M_{i,j} \in K_\lambda$ and $M_{0,0} = M$ and
$M_{i,0} = M^1_i,M_{0,j} = M^2_j$
\sn
\item "{$(B)$}"  $i_1 \le i_2 \and j_1 \le j_2 \Rightarrow M_{i_1,j_1} 
\le_{\frak K} M_{i_2,j_2}$
\sn
\item "{$(C)$}"  if $i \le \lambda \times \delta^a_1$ is a limit 
ordinal and $j \le \lambda \times \delta^a_2$ \ub{then} $M_{i,j} = 
\dsize \bigcup_{\zeta < i} M_{\zeta,j}$
\sn
\item "{$(D)$}"  if $i \le \lambda \times \delta^a_1$ and $j \le \lambda 
\times \delta^a_2$ is a limit ordinal \ub{then} $M_{i,j} = 
\dsize \bigcup_{\xi < j} M_{i,\xi}$
\sn
\item "{$(E)$}"  $M_{\lambda \times \delta^a_1,j+1}$ is 
$(\lambda,{\text{\rm cf\/}}
(\delta^a_1))$-brimmed over $M^a_{\lambda \times \delta^a_1,j}$ 
for $j < \lambda \times \delta^a_2$
\sn
\item "{$(F)$}"  $M_{i+1,\lambda \times \delta^a_2}$ is 
$(\lambda,{\text{\rm cf\/}}(\delta^a_2))$-brimmed over 
$M_{i,\lambda \times \delta^a_2}$ for $i < \lambda \times \delta^a_1$
\sn
\item "{$(G)$}"  $M_{\lambda \times \delta^a_1,\lambda \times \delta^a_2} 
\le_{\frak K} M_3 \in K_\lambda$ moreover \newline
$M_3$ is $(\lambda,{\text{\rm cf\/}}(\delta^a_3))$-brimmed over 
$M_{\lambda \times \delta^a_1,\lambda \times \delta^a_2}$
\sn
\item "{$(H)$}"  for $i < \lambda \times \delta^a_1,j \le \lambda
\times \delta^a_2$ we have
${\text{\rm \ortp\/}}(c_i,M_{i,j},M_{i+1,j})$ does not fork over $M_{i,0}$
\sn
\item "{$(I)$}"  for $j < \lambda \times \delta^a_2,i \le \lambda
\times \delta^a_1$ we have
${\text{\rm \ortp\/}}(d_j,M_{i,j},M_{i,j+1})$ does not fork over $M_{0,j}$.
\ermn
We can add
\mr
\item "{$(J)$}"  for $i < \lambda \times \delta^a_1,j < \lambda \times
\delta^b_2$ the model
$M_{i+1,j+1}$ is $(\lambda,*)$-brimmed over $M_{i,j+1} \cup
M_{i+1,j}$.
\endroster
\endproclaim
\bigskip

\remark{Remark}  1) We can replace in \scite{600-nf.6A} the ordinals $\lambda
\times \delta^a_\ell \, (\ell =1,2,3)$ by any ordinal $\alpha^a_\ell <
\lambda^+$ (for $\ell = 1,2,3$) we use the 
present notation just to conform with its use in
the proof of \scite{600-nf.6}.
\nl
2) Why do we need $u^\ell_1$ in the proof below?  This is used to get
the brimmness demands in \scite{600-nf.6A}.
\endremark
\bigskip

\demo{Proof}   We first change our towers, repeating models to give
space for bookkeeping.  
That is we define ${}^*M^1_\alpha$ for $\alpha \le \lambda \times \lambda
\times \delta^a_1$ as follows:
\sn

if $\lambda \times \beta < \alpha \le \lambda \times 
\beta + \lambda$ and $\beta < \lambda \times \delta^a_1$ then 
${}^* M^1_\alpha = M^1_{\beta + 1}$
\sn

if $\alpha = \lambda \times \beta$, then ${}^*M^1_\alpha = M^1_\beta$.
\mn
Let $u^1_0 = \{\lambda \beta:\beta < \delta^a_1\},u^1_1 = \lambda
\times \lambda \times
\delta^a_1 \backslash u^1_0,u^1_2 = \emptyset$ and for $\alpha = \lambda
\beta \in u^1_0$ let $a^1_\alpha = c_\beta$.

Similarly let us define ${}^*M^2_\alpha$ (for $\alpha \le \lambda
\times \lambda \times \delta^a_2$),$u^2_0,u^2_1,u^2_2$ and 
$\langle a^2_\alpha:\alpha \in u^2_0
\rangle$.

Now apply \scite{600-4a.9} (check) and get ${}^*M_{i,j},
(i \le \lambda \times \lambda \times
\delta^a_1,j \le \lambda \times 
\lambda \times \delta^a_2)$.  Lastly, for $i \le
\delta^a_1,j \le \delta^a_2$ let $M_{i,j} = {}^*M_{\lambda \times
i,\lambda \times j}$.
By \scite{600-4a.2} clearly ${}^* M_{\lambda \times i + \lambda,
\lambda \times j+ \lambda}$ is
$(\lambda,\text{cf}(\lambda))$-brimmed over 
${}^*M_{\lambda \times i+1,\lambda \times
j+1}$ hence $M_{i+1,j+1}$ is $(\lambda,\text{cf}(\lambda))$-brimmed over
$M_{i+1,j} \cup M_{i,j+1}$.
And, by \scite{600-4a.1}(1) choose $M_3 \in K_\lambda$ which is 
$(\lambda,\text{cf}(\delta^a_3))$-brimmed over 
$M_{\lambda \times \delta^a_1,\lambda \times \delta^a_2}$. 
  \hfill$\square_{\scite{600-nf.6A}}$
\enddemo
\bigskip

\demo{Proof of \scite{600-nf.6}}  We shall let $M_{i,j},M_3$ be as in 
\scite{600-nf.6A} for $\bar \delta^a$ and $\bar M^1,\bar M^2$ determined below.  
For $x \in \{ a,b\}$ as 
NF$_{\lambda,\bar \delta^x}(N^x_0,N^x_1,N^x_2,N^x_3)$, we know that 
there are witnesses $\langle N^x_{1,i}:i \le \lambda \times 
\delta^x_1 \rangle,\langle c^x_i:i < \lambda \times \delta^x_1\rangle,
\langle N^x_{2,i}:i \le \lambda \times \delta^x_1 
\rangle$ for this.  So $\langle N^x_{1,i}:i \le \lambda \times 
\delta^x_1 \rangle$ is $\le_{\frak K}$-increasing continuous
and $(N^x_{1,i},N^x_{1,i+1},c^x_i) \in K^{3,\text{uq}}_\lambda$ for $i
< \lambda \times \delta^x_1$.  Hence by the freedom we 
have in choosing $\bar M^1$ and $\langle c_i:i < \lambda \times
\delta_1\rangle$ without loss of generality there is an 
isomorphism $g_1$ from $N^a_{1,\lambda \times \delta^a_1}$ onto 
$M_{\lambda \times \delta^a_1}$ mapping $N^a_{1,i}$ onto $M^1_i 
= M_{i,0}$ and $c^a_i$ to $c_i$; remember that
$N^a_{1,\lambda \times \delta^a_1} = N^a_1$.  
Let $g_0 = g_1 \restriction N^a_0 = g_1 \restriction N^a_{1,0}$ so 
$g_0 \circ f^{-1}_0$ is an isomorphism from $N^b_0$ onto $M_{0,0}$.
\medskip

Similarly as $\delta^b_1 = \delta^a_2$, 
and using the freedom we have in choosing $\bar M^2$ and $\langle
d_i:i < \lambda \times \delta^b_1\rangle$
 without loss of generality there is an isomorphism $g_2$ from
$N^b_{1,\lambda \times \delta^a_2}$ onto $M^2_j = M_{0,\lambda \times 
\delta^a_2}$ mapping $N^b_{1,j}$ onto $M_{0,j}$ 
(for $j \le \lambda \times \delta^a_2)$ and mapping $c^b_i$ to $d_i$
and $g_2$ extends $g_0 \circ f^{-1}_0$. \newline
Now would like to use the weak uniqueness \scite{600-nf.5} and for this note:
\mr
\item "{$(\alpha)$}"  NF$_{\lambda,\bar \delta^a}(N^a_0,N^a_1,N^a_2,N^a_3)$
is witnessed by the sequences 
$\langle N^a_{1,i}:i \le \lambda \times \delta^a_1 \rangle$, and
$\langle N^a_{2,i}:i \le \lambda \times \delta^a_1 \rangle$ \newline
[why?  an assumption]
\sn
\item "{$(\beta)$}"  NF$_{\lambda,\bar \delta^a}(M_{0,0},M_{\lambda \times 
\delta^a_1,0},M_{0,\lambda \times \delta^a_2},M_3)$ 
is witnessed by the sequences \newline
$\langle M_{i,0}:i \le \lambda \times \delta^a_1 \rangle,\langle
M_{i,\lambda \times \delta^a_2}:i \le \lambda \times \delta^a_1 
\rangle$ \newline
[why? check]
\sn
\item "{$(\gamma)$}"  $g_0$ is an isomorphism from $N^a_0$ onto $M_{0,0}$
\newline
[why?  see its choice]
\sn
\item "{$(\delta)$}"  $g_1$ is an isomorphism from $N^a_1$ onto
$M_{\lambda \times \delta^a_1,0}$ mapping $N^a_{1,i}$ onto 
$M_{i,0}$ for $i < \lambda \times \delta^a_1$ and $c^a_i$ to $c_i$ 
for $i < \lambda \times \delta^a_1$ and extending $g_0$ \newline
[why?  see the choice of $g_1$ and of $g_0$]
\sn
\item "{$(\varepsilon)$}"  $g_2 \circ f_2$ is an isomorphism from
$N^a_2$ onto $M_{0,\lambda \times \delta^a_2}$ extending $g_0$ \newline
[why?  $f_2$ is an isomorphism from $N^a_2$ onto $N^b_1$ and $g_2$ is an
isomorphism from $N^b_1$ onto $M_{0,\lambda \times \delta^a_1}$ 
extending $g_0 \circ f^{-1}_0$ and $f_0 \subseteq f_2$].
\endroster
\medskip

So there is by \scite{600-nf.5} an isomorphism $g^a_3$ from $N^a_3$ onto 
$M_3$ extending both $g_1$ and $g_2 \circ f_2$.
\medskip

We next would like to apply \scite{600-nf.5} to the $N^b_i$'s; so note:
\mr
\item "{$(\alpha)'$}"  NF$_{\lambda,\bar \delta^b}(N^b_0,N^b_1,N^b_2,N^b_3)$
is witnessed by the sequences 
$\langle N^b_{1,i}:i \le \lambda \times \delta^a_2 
\rangle$, \newline
$\langle N^b_{2,i}:i \le \lambda \times \delta^a_2 \rangle$
\sn
\item "{$(\beta)'$}"  NF$_{\lambda,\bar \delta^b}(M_{0,0},M_{0,\lambda 
\times \delta^a_2},M_{\lambda \times \delta^a_1,0},M_3)$ is witnessed by the
sequences \newline 
$\langle M_{0,j}:j \le \lambda \times \delta^a_2 \rangle,
\langle 
M_{\lambda \times \delta^a_1,j}:j \le \lambda \times \delta^a_2 \rangle$
\sn
\item "{$(\gamma)'$}"  $g_0 \circ (f_0)^{-1}$ is an 
isomorphism from $N^b_0$ onto $M_{0,0}$ \newline
[why?  Check.]
\sn
\item "{$(\delta)'$}"  $g_2$ is an isomorphism from $N^b_1$ onto
$M_{0,\lambda \times \delta^a_2}$ mapping $N^b_{1,j}$ onto $M_{0,j}$
and $c^a_j$ to $d_j$ for
$j \le \lambda \times \delta^a_2$ and extending $g_0 \circ (f_2)^{-1}$
\newline
[why?  see the choice of $g_2$: it maps $N^b_{1,j}$ onto $M_{0,j}$]
\sn
\item "{$(\varepsilon)'$}"  $g_1 \circ (f_1)^{-1}$ is an isomorphism from
$N^b_2$ onto $M_{\lambda \times \delta^a_0}$ extending $g_0$ \newline
[why?  remember $f_1$ is an isomorphism from $N^a_1$ onto $N^b_2$ extending
$f_0$ and the choice of $g_1$: it maps $N^a_1$ 
onto $M_{\lambda \times \delta^a_1,0}$].
\endroster
\medskip

\noindent
So there is an isomorphism $g^b_3$ form $N^b_3$ onto $M_3$ extending
$g_2$ and $g_1 \circ (f_1)^{-1}$. \newline
Lastly $(g^b_3)^{-1} \circ g^a_3$ is an isomorphism from $N^a_3$ onto
$N^b_3$ (chase arrows). Also

$$
\align
((g^b_3)^{-1} \circ g^a_3) \restriction N^a_1 &= (g^b_3)^{-1}(g^a_3 
\restriction N^a_1) \\
  &= (g^b_3)^{-1} g_1 = ((g^b_3)^{-1} \restriction M_{\lambda \times
\delta^a_1,0}) \circ g_1 \\
  &= (g^b_3 \restriction N^b_2)^{-1} \circ g_1 = ((g_1 \circ (f_1)^{-1})^{-1})
\circ g_1 \\
  &= (f_1 \circ (g_1)^{-1}) \circ g_1 = f_1.
\endalign
$$
\medskip

\noindent
Similarly  $((g^b_3)^{-1} \circ g^a_3) \restriction N^a_2 = f_2$.
\newline
So we have finished. \hfill$\square_{\scite{600-nf.6}}$
\enddemo
\bn
But if we invert twice we get straight; so
\proclaim{\stag{600-nf.7} Claim}  [Uniqueness].  Assume for 
$x \in \{ a,b\}$ we have \newline
${\text{\rm NF\/}}_{\lambda,\bar \delta^x}(N^x_0,N^x_1,N^x_2,N^x_3)$ and
${\text{\rm cf\/}}(\delta^a_1) = 
{ \text{\rm cf\/}}(\delta^b_1),{\text{\rm cf\/}}(\delta^a_2) 
= { \text{\rm cf\/}}(\delta^b_2),{\text{\rm cf\/}}(\delta^a_3) 
= { \text{\rm cf\/}}(\delta^b_3)$, all 
$\delta^x_\ell$ limit ordinals $< \lambda^+$.

If $f_\ell$ is an isomorphism from $N^a_\ell$ onto $N^b_\ell$ for $\ell < 3$
and $f_0 \subseteq f_1,f_0 \subseteq f_2$ \underbar{then} there is an
isomorphism $f$ from $N^a_3$ onto $N^b_3$ extending $f_1,f_2$.
\endproclaim
\bigskip

\demo{Proof}  Let $\bar \delta^c = \langle \delta^c_1,\delta^c_2,
\delta^c_3 \rangle = \langle \delta^a_2,\delta^a_1,\delta^a_3 \rangle$; 
by \scite{600-nf.3}(1) there are
$N^c_\ell$ (for $\ell \le 3$) such that NF$_{\lambda,\bar \delta^c}(N^c_0,
N^c_1,N^c_2,N^c_3)$ and $N^c_0 \cong N^a_0$.  
There is for $x \in \{ a,b\}$ an isomorphism $g^x_0$
from $N^x_0$ onto $N^c_0$ and without loss of generality 
$g^a_0 = g^b_0 \circ f_0$.  Similarly
for $x \in \{a,b\}$
there is an isomorphism $g^x_1$ from 
$N^x_1$ onto $N^c_2$ extending $g^x_0$ (as $N^x_1$ is
$(\lambda,\text{cf}(\delta^x_1))$-brimmed over $N^x_0$ and also
$N^c_2$ is $(\lambda,\text{cf}(\delta^c_2))$-brimmed over $N^c_0$ and
$\text{cf}(\delta^c_2) = \text{cf}(\delta^a_1) = \text{cf}(\delta^x_1))$ 
and without loss of generality $g^b_1 = g^a_1 \circ f_1$.  Similarly for
$x \in \{a,b\}$ 
there is an isomorphism $g^x_2$ from $N^x_2$ onto $N^c_1$ extending $g^x_0$
(as $N^x_2$ is
$(\lambda,\text{cf}(\delta^x_2))$-brimmed over $N^x_0$ and also $N^c_1$
is $(\lambda,\text{cf}(\delta^c_1))$-brimmed over $N^c_0$ and
$\text{cf}(\delta^c_1) = \text{cf}(\delta^a_2) = \text{cf}(\delta^x_2))$
and without loss of generality $g^a_2 = g^b_2 \circ f_2$.
\newline
So by \scite{600-nf.6} for $x \in \{a,b\}$ there is an isomorphism 
$g^x_3$ from $N^x_3$ onto $N^c_3$
extending $g^x_1$ and $g^x_2$.  
Now $(g^b_3)^{-1} \circ g^a_3$ is an isomorphism
from $N^a_3$ onto $N^b_3$ extending $f_1,f_2$ as required.
\hfill$\square_{\scite{600-nf.7}}$
\enddemo
\bn
So we have proved the uniqueness for NF$_{\lambda,\bar \delta}$ when all
$\delta_\ell$ are limit ordinals; this means that the arbitrary choice of
$\langle N_{1,i}:i \le \lambda \times \delta_1 \rangle$ and
$\langle c_i:i < \lambda \times \delta_1 \rangle$ is immaterial; it
figures in the definition and, e.g. existence proof 
but does not influence the net
result.  The 
power of this result is illustrated in the following conclusion.
\demo{\stag{600-nf.8} Conclusion}  [Symmetry]. \newline

If 
NF$_{\lambda,\langle \delta_1,\delta_2,\delta_3 \rangle}(N_0,N_1,N_2,N_3)$
where $\delta_1,\delta_2,\delta_3$ are limit ordinals $< \lambda^+$
\underbar{then} NF$_{\lambda,\langle \delta_2,\delta_1,
\delta_3 \rangle}(N_0,N_2,N_1,N_3)$.
\enddemo
\bigskip

\demo{Proof}  By \scite{600-nf.6A} we can find $N'_\ell (\ell \le 3)$
such that: $N'_0 = N_0,N'_1$ is
$(\lambda,\text{cf}(\delta_1))$-brimmed over $N'_0,N'_2$ is
$(\lambda,\text{cf}(\delta_2))$-brimmed over $N'_0$ and $N'_3$ is
$(\lambda,\text{cf}(\delta_3))$-brimmed over $N'_1 \cup N'_2$ and
NF$_{\lambda,\langle \delta_1,\delta_2,\delta_3
\rangle}(N'_0,N'_1,N'_2,N'_3)$ and NF$_{\lambda,\langle
\delta_2,\delta_1,\delta_3 \rangle}(N'_0,N'_2,N'_1,N'_3)$.  Let
$f_1,f_2$ be an isomorphism from $N_1,N_2$ onto $N'_1,N'_2$ over
$N_0$, respectively.  By \scite{600-nf.7} (or \scite{600-nf.6}) 
there is an isomorphism $f'_3$
form $N_3$ onto $N'_3$ extending $f_1 \cup f_2$.  As isomorphisms
preserve NF we are done.  \hfill$\square_{\scite{600-nf.8}}$
\enddemo
\bigskip

Now we turn to smooth amalgamation (not necessarily brimmed, see
Definition \scite{600-nf.2}).  If we use Lemma \scite{600-4a.6}, of course, we
do not really need \scite{600-nf.10}.
\proclaim{\stag{600-nf.10} Claim}  1) If 
${\text{\rm NF\/}}_{\lambda,\bar \delta}(N_0,N_1,N_2,N_3)$ and
$\delta_1,\delta_2,\delta_3$ are limit ordinals, 
\underbar{then} ${\text{\rm NF\/}}_\lambda
(N_0,N_1,N_2,N_3)$ (see Definition \scite{600-nf.2}). \newline
2) In Definition \scite{600-nf.2}(1) we can add:
\mr
\item "{$(d)^+$}"  $M_\ell$ is $(\lambda,{\text{\rm cf\/}}(\lambda))$-brimmed 
over $N_0$ and moreover over $N_\ell$,
\sn
\item "{$(e)$}"    $M_3$ is $(\lambda,{\text{\rm cf\/}}(\lambda))$-brimmed
over $M_1 \cup M_2$ (actually this is given by clause $(f)(\zeta)$ of 
Definition \scite{600-nf.1}). 
\ermn
3) If $N_0 \le_{\frak K} N_\ell$ for $\ell=1,2$ and $N_1 \cap N_2 =
N_0$,  \ub{then} we can find
$N_3$ such that {\rm NF}$_\lambda(N_0,N_1,N_2,N_3)$.
\endproclaim
\bigskip

\demo{Proof}  1) Note that even if every $\delta_\ell$ is limit and we
waive the ``moreover" in clause $(d)^+$,  
the problem is in the case that e.g. 
$(\text{cf}(\delta^a),\text{cf}(\delta^b),\text{cf}(\delta^c)) \ne 
(\text{cf}(\lambda),\text{cf}(\lambda),\text{cf}(\lambda))$.
For $\ell=1,2$ we can find $\bar M^\ell = \langle M^\ell_i:i \le
\lambda \times (\delta_\ell + \lambda) \rangle$ and $\langle
c^\ell_i:i < \lambda \times (\delta_i + \lambda)\rangle$ such that
$M^\ell_0 = N_0,\bar M^1$ is $\le_{\frak K}$-increasing continuous
$(M^\ell_i,M^\ell_{i+1},c_i) \in K^{3,\text{uq}}_{\frak s}$ and if $p
\in {\Cal S}^{\text{bs}}(M^\ell_i)$ and $i < \lambda \times (\delta_\ell +
\lambda)$ then for $\lambda$ ordinals $j < \lambda$,
\ortp$(c_i,M^\ell_{i+j},M^\ell_{i+j+1})$ is a non-forking extension of
$p$.  So $M^\ell_{\lambda \times \delta_\ell}$ is
$(\lambda,\text{cf}(\delta_\ell))$-brimmed over $M^\ell_0 = N_0$ and
$M^\ell_{\lambda \times (\delta_\ell + \lambda)}$ is
$(\lambda,\text{cf}(\lambda))$-brimmed over $M^\ell_{\lambda \times
\delta_\ell}$; so \wilog \, $M^\ell_{\lambda \times \delta_\ell} =
N_\ell$ for $\ell=1,2$. \nl
By \scite{600-nf.6A} we can find $M_{i,j}$ for 
$i \le \lambda \times (\delta_1 + \lambda),j \le \lambda \times 
(\delta_2 + \lambda)$ for $\bar \delta' := \langle \delta_1 + \lambda,
\delta_2 + \lambda,\delta_3 \rangle$ such that they are as in
\scite{600-nf.6A} for $\bar M^1,\bar M^2$ so 
$M_{0,0} = N_0$; then choose $M'_3 \in K_\lambda$ which 
is $(\lambda,\text{cf}(\delta_3))$-brimmed over $M_{\lambda \times 
\delta_1,\lambda \times \delta_2}$.  So 
NF$_{\lambda,\bar \delta}(M_{0,0},M_{\lambda \times
\delta_1,0},M_{0,\lambda \times \delta_2},M'_3)$, hence by \scite{600-nf.7} 
without loss of generality $M_{0,0} = N_0,M_{\lambda \times \delta_1,0}
 = N_1,M_{0,\lambda \times \delta_2} = N_2$, and $N_3 = M'_3$.  
Lastly, let $M_3$
be $(\lambda,\text{cf}(\lambda))$-brimmed over $M'_3$.
Now clearly also \newline
NF$_{\lambda,\langle \delta_1 + \lambda,\delta_2 + \lambda,\delta_3 +
\lambda \rangle}(M_{0,0},M_{\lambda \times (\delta_1 + \lambda),0},
M_{0,\lambda \times (\delta_2 + \lambda)},M_3)$ and \newline
$N_0 = M_{0,0},N_1 = M_{\lambda \times \delta_2,0} \le_{\frak K} 
M_{\lambda \times (\delta_2 + \lambda),0},
N_2 = M_{0,\lambda \times \delta_2} \le_{\frak K} 
M_{0,\lambda \times (\delta_2 + \lambda)}$ \newline
and $M_{\lambda \times (\delta_1 + \lambda),0}$ is
$(\lambda,\text{cf}(\lambda))$-brimmed over 
$M_{\lambda \times \delta_1,0}$ and
$M_{0,\lambda \times (\delta_2 + \lambda)}$ is \newline
$(\lambda,\text{cf}(\lambda))$-brimmed over 
$M_{0,\lambda \times \delta_2}$ 
and $N_3 = M'_3 \le_{\frak K} M_3$.  So
we get all the requirements for 
NF$_\lambda(N_0,N_1,N_2,N_3)$ (as witnessed
by $\langle M_{0,0},M_{\lambda \times (\delta_1 + \lambda),0},
M_{0,\lambda \times (\delta_2 + \lambda)},M_3 \rangle$). 
2) Similar proof. \nl
3) By \scite{600-nf.3} and the proof above.  \hfill$\square_{\scite{600-nf.10}}$
\enddemo
\bn
Now we turn to NF$_\lambda$; existence is easy.
\proclaim{\stag{600-nf.10.3} Claim}  {\rm NF}$_\lambda$ has existence,
i.e., clause (f) of \scite{600-nf.0X}(1).
\endproclaim
\bigskip

\demo{Proof}  By \scite{600-nf.10}(3).  \hfill$\square_{\scite{600-nf.10.3}}$
\enddemo
\bn
Next we deal with real uniqueness
\proclaim{\stag{600-nf.11} Claim}  [Uniqueness of smooth amalgamation]: \nl
1) If
${\text{\rm NF\/}}_\lambda(N^x_0,N^x_1,
N^x_2,N^x_3)$ for $x \in \{ a,b\},f_\ell$ an
isomorphism from $N^a_\ell$ onto $N^b_\ell$ for $\ell < 3$ and
$f_0 
\subseteq f_1,f_0 \subseteq f_2$ \underbar{then} $f_1 \cup f_2$ can be
extended to a $\le_{\frak K}$-embedding of $N^a_3$ into some 
$\le_{\frak K}$-extension of $N^b_3$. \nl
2) So if above $N^x_3$ is $(\lambda,\kappa)$-brimmed 
over $N^x_1 \cup N^x_2$ for $x = a,b$, we can 
extend $f_1 \cup f_2$ to an isomorphism from $N^a_3$ onto $N^b_3$.
\endproclaim
\bigskip

\demo{Proof}  1) For $x \in \{ a,b\}$ let the sequence $\langle M^x_\ell:
\ell < 4 \rangle$ be a witness to \newline
NF$_\lambda(N^x_0,N^x_1,N^x_2,N^x_3)$ as in \scite{600-nf.2}, \scite{600-nf.10}(2),
so in particular \newline
NF$_{\lambda,\langle \lambda,\lambda,\lambda \rangle}
(M^x_0,M^x_1,M^x_2,M^x_3)$.  By chasing arrows (disjointness) and 
uniqueness, i.e. 
\scite{600-nf.7} without loss of generality $M^a_\ell = M^b_\ell$ for $\ell < 4$ 
and $f_0 = \text{ id}_{N^a_0}$.  As $M^a_1$ is 
$(\lambda,\text{cf}(\lambda))$-brimmed over $N^a_1$ and also
over $N^b_1$ (by clause $(d)^+$ of \scite{600-nf.10}(2)) 
and $f_1$ is an isomorphism from $N^a_1$ onto $N^b_1$, clearly
by \scite{600-0.22} there is an automorphism $g_1$ of $M^a_1$ such that $f_1 \subseteq g_1$,
hence also $\text{id}_{N^a_0} = f_0 \subseteq f_1 \subseteq g_1$.  
Similarly there is an
automorphism $g_2$ of $M^a_2$ extending $f_2$ hence $f_0$.  So $g_\ell \in
\text{ AUT}(M^a_\ell)$ for $\ell = 1,2$ and
$g_1 \restriction M^a_0 = f_0 =
g_2 \restriction M^a_0$.  By the uniqueness of NF$_{\lambda,\langle
\lambda,\lambda,\lambda \rangle}$ (i.e. Claim \scite{600-nf.7}) there is an 
automorphism $g_3$ of $M^a_3$ 
extending $g_1 \cup g_2$.  This proves the desired conclusion.
\nl
2) Should be clear.  \hfill$\square_{\scite{600-nf.11}}$
\enddemo
\bn
We now show that in the cases the two notions of non-forking
amalgamations are meaningful then they coincide, one implication
already is a case of \scite{600-nf.10}.
\proclaim{\stag{600-nf.12} Claim}  Assume
\mr
\item "{$(a)$}"  $\bar \delta = \langle \delta_1,\delta_2,\delta_3 \rangle,
\delta_\ell < \lambda^+$ is a limit ordinal for $\ell = 1,2,3$; \newline
$N_0 \le_{\frak K} N_\ell \le_{\frak K} N_3$ are in $K_\lambda$ for 
$\ell = 1,2$ 
\sn
\item "{$(b)$}"  $N_\ell$ is $(\lambda,{\text{\rm cf\/}}
(\delta_\ell))$-brimmed over $N_0$ for $\ell =1,2$
\sn
\item "{$(c)$}"  $N_3$ is ${\text{\rm cf\/}}(\delta_3)$-brimmed 
over $N_1 \cup N_2$.
\ermn
\ub{Then}  ${\text{\rm NF\/}}_\lambda(N_0,N_1,N_2,N_3)$ iff 
{\rm NF}$_{\lambda,\bar \delta}(N_0,N_1,N_2,N_3)$.
\endproclaim
\bigskip

\demo{Proof}  The ``if" direction holds by \scite{600-nf.10}(1).  
As for the ``only if" direction, basically it follows from the existence for
NF$_{\lambda,\bar \delta}$ and uniqueness for NF$_\lambda$; in details
by the proof of \scite{600-nf.10}(1) (and Definition \scite{600-nf.1}, 
\scite{600-nf.2}) we can find $M_\ell(\ell \le 3)$ such that $M_0 = N_0$ and
NF$_{\lambda,\bar \delta}(M_0,M_1,M_2,M_3)$ and clauses (b), (c), (d) of
Definition \scite{600-nf.2} and $(d)^+$ of \scite{600-nf.10}(2) hold so by 
\scite{600-nf.10} also
NF$_\lambda(M_0,M_1,M_2,M_3)$.  Easily there are for $\ell < 3$, isomorphisms
$f_\ell$ from $M_\ell$ onto $N_\ell$ such that $f_0 = f_\ell
\restriction M_\ell$ where $f_0 = \text{ id}_{N_0}$.  
By the uniqueness of smooth amalgamations
(i.e., \scite{600-nf.11}(2)) we can find an isomorphism $f_3$ from $M_3$ 
onto $N_3$ extending $f_1 \cup f_2$.  So as
NF$_{\lambda,\bar \delta}(M_0,M_1,M_2,M_3)$ holds also NF$_{\lambda,
\bar \delta},(f_0(M_0),f_3(M_1),f_3(M_2),f_3(M_3))$; that is
NF$_{\lambda,\bar \delta}(N_0,N_1,N_2,N_3)$ is as required. 
\hfill$\square_{\scite{600-nf.12}}$
\enddemo
\bigskip

\proclaim{\stag{600-nf.13} Claim} [Monotonicity]:  If 
${\text{\rm NF\/}}_\lambda(N_0,N_1,N_2,N_3)$ and
$N_0 \le_{\frak K} N'_1 \le_{\frak K} N_1$ and $N_0 \le_{\frak K} 
N'_2 \le_{\frak K} N_2$ and $N'_1 \cup N'_2
\subseteq N'_3 \le_{\frak K} N''_3,N_3 \le_{\frak K} N''_3$
\underbar{then} ${\text{\rm NF\/}}_\lambda(N_0,N'_1,N'_2,N'_3)$.
\endproclaim
\bigskip

\demo{Proof}  Read Definition \scite{600-nf.2}(1).  \hfill$\square_{\scite{600-nf.13}}$
\enddemo
\bigskip

\proclaim{\stag{600-nf.14} Claim}  [Symmetry]:  
${\text{\rm NF\/}}_\lambda(N_0,N_1,N_2,N_3)$ 
holds \underbar{if and only if} \newline
{\rm NF}$_\lambda(N_0,N_2,N_1,N_3)$ holds.
\endproclaim
\bigskip

\demo{Proof}  By Claim \scite{600-nf.8} (and Definition \scite{600-nf.2}).
\hfill$\square_{\scite{600-nf.14}}$
\enddemo
\bn
We observe
\demo{\stag{600-nf.14.1} Conclusion}  If NF$_\lambda(N_0,N_1,N_2,N_3),N_3$
is $(\lambda,\sigma)$-brimmed over $N_1 \cup N_2$ and $\lambda \ge
\sigma,\kappa \ge \aleph_0$, \ub{then} there is $N^+_2$ such that
\mr
\item "{$(a)$}"  NF$_\lambda(N_0,N_1,N^+_2,N_3)$
\sn
\item "{$(b)$}"  $N_2 \le_{\frak K} N^+_2$
\sn
\item "{$(c)$}"  $N^+_2$ is $(\lambda,\kappa)$-brimmed over $N_0$ and
even over $N_2$
\sn
\item "{$(d)$}"  $N_3$ is $(\lambda,\sigma)$-brimmed over $N_1 \cup N^+_2$.
\endroster
\enddemo
\bigskip

\demo{Proof}  Let $N^+_2$ be $(\lambda,\kappa)$-brimmed over $N_2$ be
such that $N^+_2 \cap N_3 = N_2$.  So by existence \scite{600-nf.10.3}
there is $N^+_3$ such that NF$_\lambda(N_0,N_1,N^+_2,N^+_3)$ and
$N^+_3$ is $(\lambda,\sigma)$-brimmed over $N_1 \cup N^+_2$.  By
monotonicity \scite{600-nf.13} we have NF$_\lambda(N_0,N_1,N_2,N^+_3)$.
So by uniqueness (i.e., \scite{600-nf.11}(2)) \wilog \, $N_3 = N^+_3$, so
we are done.  \hfill$\square_{\scite{600-nf.14.1}}$
\enddemo
\bn
The following claim is a step toward proving transitivity for
NF$_\lambda$; so we first deal with NF$_{\lambda,\bar \delta}$.  Note
below: if we ignore $N^c_i$ we have problem showing NF$_{\lambda,\bar
\delta}(N^a_0,N^a_\alpha,N^b_0,N^b_\alpha)$.  Note that it is not clear
at this stage whether, e.g. $N^b_\omega$ is even universal over
$N^a_\omega$, but $N^c_\omega$ is; note that the $N^c_i$ are $\le_{\frak
K}$-increasing with $i$ but not necessarily continuous.
However once we finish proving that NF$_\lambda$ is a non-forking
relation on ${\frak K}_{\frak s}$ respecting ${\frak s}$ this claim
will lose its relevance.
\proclaim{\stag{600-nf.15} Claim}  Assume $\alpha < \lambda^+$ is an ordinal and
for $x \in \{a,b,c\}$ the sequence $\bar N^x = \langle N^x_i:i \le \alpha 
\rangle$ is a $\le_{\frak K}$-increasing sequence of members of
$K_\lambda$, and for $x = a,b$ the sequence $\bar N^x$ is 
$\le_{\frak K}$-increasing 
continuous, $N^b_i \cap N^a_\alpha = N^a_i,N^c_i \cap N^a_\alpha = N^a_i,
N^a_i \le_{\frak K} N^b_i \le_{\frak K} N^c_i$ and $N^b_0$ is
$(\lambda,\delta_2)$-brimmed over $N^a_0$ and
${\text{\rm NF\/}}_{\lambda,\bar \delta^i}(N^a_i,N^a_{i+1},N^c_i,N^b_{i+1})$
(so necessarily $i < \alpha \Rightarrow N^c_i \le_{\frak K}
N^b_{i+1}$) where \newline
$\bar \delta^i = \langle \delta^i_1,\delta^i_2,\delta^i_3 \rangle$
with $\delta^i_1,\delta^i_2,\delta^i_3$ are ordinals $< \lambda^+$ 
 and $\delta_3 < \lambda^+$ is limit, $N^c_\alpha$ is 
$(\lambda,{\text{\rm cf\/}}(\delta_3))$-brimmed over $N^b_\alpha,\delta_1 =
\dsize \sum_{\beta < \alpha} \delta^\beta_1$ and $\delta_3 =
\delta^\alpha_3$  and $\delta_2 = \delta^0_2,
\bar \delta = \langle \delta_1,\delta_2,\delta_3 \rangle$.  \nl
\underbar{Then} ${\text{\rm NF\/}}_{\lambda,\bar \delta}
(N^a_0,N^a_\alpha,N^b_0,N^c_\alpha)$.
\endproclaim
\bigskip

\demo{Proof}  For $i < \alpha$ let 
$\langle N^i_{1,\varepsilon},N^i_{2,\varepsilon},d^i_\zeta:\varepsilon
\le \lambda \times \delta^i_1,\zeta < \lambda \times \delta^i_1 
\rangle$ be a witness to
NF$_{\lambda,\bar \delta^i}(N^a_i,N^a_{i+1},N^c_i,N^b_{i+1})$.  Now we
define a sequence $\langle N_{1,\varepsilon},N_{2,\varepsilon},
d^i_\zeta:\varepsilon \le \lambda \times \delta_1$ and $\zeta <
\lambda \times  \delta_1 \rangle$ where
\mr
\item "{$(a)$}"  $N_{1,0} = N^a_0,N_{2,0} = N^b_0$ and
\sn
\item "{$(b)$}"   if $\lambda \times (\dsize \sum_{j<i}
\delta^j_1) < \zeta \le \lambda \times (\dsize \sum_{j \le i}
\delta^j_1)$ then we let $N_{1,\zeta} = 
N^i_{1,\varepsilon_\zeta},N_{2,\zeta} = 
N^i_{2,\varepsilon_\zeta}$ where
$\varepsilon_\zeta = \zeta - \lambda \times ( \dsize \sum_{j < i}
\delta^j_1)$ and
\sn
\item "{$(c)$}"   if $0 < \zeta = \lambda \times \dsize \sum_{j
< \alpha} \delta^j_1$ we let $N_{1,\zeta} = N^a_i,N_{2,\zeta} = N^b_i
= \alpha$ (if $i$ is non-limit we should note that this is compatible
with clause (b), note that by this if $i = \alpha$ then
$N_{1,\zeta} = N^a_\alpha,N_{2,\zeta} =
\cup\{N^i_{2,\lambda \times \delta_1}:i < \alpha\}$
\sn
\item "{$(d)$}" if $\lambda \times (\dsize \sum_{j<i} \delta^j_1) \le
\zeta < \lambda \times (\dsize \sum_{j \le i} \delta^j_1)$ then we let
$d_\zeta = d^i_{\varepsilon_\zeta}$ where $\varepsilon_\zeta = \zeta -
\lambda \times (\dsize \sum_{j<i} \delta^j_j) =
\cup\{N^*_{2,\zeta}:\zeta < \lambda \times (\dsize \sum_{j<\alpha}
\delta^j_1)$. 
\ermn 
Clearly $\langle
N_{1,\zeta}:\zeta \le \lambda \times \delta_1 \rangle$ is $\le_{\frak
K}$-increasing continuous, and also $\langle N_{2,\zeta}:\zeta \le \lambda
\times \delta_1 \rangle$ is.  Obviously $(N_{1,\zeta},N_{1,\zeta
+1},d_\zeta) \in K^{3,\text{uq}}_\lambda$ as this just means
$(N^i_{1,\varepsilon_\zeta},N^i_{1,\varepsilon_\zeta +1},d^i_\zeta)
\in K^{3,\text{uq}}_\lambda$ when $\lambda \times \dsize \sum_{j<i}
\delta^j_1:j \le \zeta < \lambda \times \dsize \sum_{j \le i}
\delta^j_1$ and $\varepsilon_\zeta$ as above.

Why \ortp$(d_\zeta,N_{2,\zeta},N_{2,\zeta +1})$ does not fork over
$N_{1,\zeta}$ for $\zeta,i$ such that $\lambda \times (\dsize
\sum_{j<i} \delta^j_1)  \zeta < \lambda \times (\dsize \sum_{j \le i}
\delta^j_j)$?  If $\lambda \times \dsize \sum_{j<i}
\delta^j_1 < \zeta$ this holds as it means
\ortp$(d^i_{\varepsilon_\zeta},N^i_{2,\varepsilon_\zeta},N^i_{2,\varepsilon_\zeta
+1})$ does not fork over $N^i_{1,\zeta}$.  If $\lambda \times \dsize
\sum_{j<i} \delta^j_1 = \zeta$ this is not the case but $N^i_{1,0} =
N_{1,\zeta} \le_{\frak K} N_{2,\zeta} \le_{\frak K} N^c_i = N^i_{2,0}$
and we know that \ortp$(d_\zeta,N^i_{2,0},N^i_{2,1})$ does not fork over
$N^i_{1,0} = N_{1,\zeta}$ hence by monotonicity of non-forking 
\ortp$(d_\zeta,N_{2,\zeta},N_{2,\zeta +1})$ does not fork over
$N_{1,\zeta}$ is as required.
\nl
Note that we have not demanded or used 
``$\bar N^c$ continuous"; the $N^c_i$ is really needed for $i$ limit
as we do not know that $N^b_i$ is brimmed over $N^a_i$. 
\hfill$\square_{\scite{600-nf.15}}$
\enddemo
\bigskip

\proclaim{\stag{600-nf.16} Claim}  [transitivity]  1) Assume 
that $\alpha < \lambda^+$ and for $x \in
\{ a,b\}$ we have $\langle N^x_i:i \le \alpha \rangle$ is a
$\le_{\frak K}$-increasing continuous sequence of members 
of $K_\lambda$. \newline 
If {\rm NF}$_\lambda(N^a_i,N^a_{i+1},N^b_i,N^b_{i+1})$ for each $i < \alpha$
\underbar{then} {\rm NF}$_\lambda(N^a_0,N^a_\alpha,N^b_0,N^b_\alpha)$. \nl
2)  Assume that $\alpha_1 < \lambda^+,\alpha_2 < \lambda^+$ and 
$M_{i,j} \in K_\lambda$ (for $i \le \alpha_1,j \le \alpha_2$) 
satisfy clauses (B), (C), (D), 
from \scite{600-nf.6A}, and for each $i < \alpha_1,j < \alpha_2$ we have:

$$
\nonforkin{M_{i,j+1}}{M_{i+1,j}}_{M_{i,j}}^{M_{i+1,j+1}}.
$$
\medskip

$$
\text{\underbar{Then} }
\nonforkin{M_{i,0}}{M_{0,j}}_{M_{0,0}}^{M_{\alpha_1,\alpha_2}} \text{ for }
i \le \alpha_1, j \le \alpha_2.
$$
\endproclaim
\bigskip

\demo{Proof}  1) We first prove special cases and use them to prove more
general cases. \newline
\smallskip

\noindent
\underbar{Case A}:  $N^a_{i+1}$ is $(\lambda,\kappa_i)$-brimmed
over $N^a_i$ and $N^b_{i+1}$ is $(\lambda,\sigma_i)$-brimmed 
over $N^a_{i+1} \cup N^b_i$ for $i < \alpha$ ($\sigma_i$ infinite, of
course). 

In essence the problem is that we do not know ``$N^b_i$ is brimmed
over $N^a_i$" ($i$ limit) so we shall use \scite{600-nf.15}; for this we
introduce appropriate $N^c_i$.

Let $\delta^i_1 = \kappa_i,\delta^i_2 = \kappa_i,\delta^i_3 =
\sigma_i$ where we stipulate $\sigma_\alpha = \lambda$.
For $i \le \alpha$ we can choose 
$N^c_i \in K_\lambda$ such that
\mr
\item "{$(a)$}"  $N^b_i \le_{\frak K} N^c_i \le_{\frak K} N^b_{i+1},
N^c_i$ is $(\lambda,\kappa_i)$-brimmed over $N^b_i$, and \newline
${\text{\rm NF\/}}_{\lambda,\langle \delta^i_1,\delta^i_2,\delta^i_3 
\rangle}(N^a_i,N^a_{i+1},N^c_i,N^b_{i+1})$
\sn
\item "{$(b)$}"  $N^c_\alpha \in K_\lambda$ is 
$(\lambda,\delta^\alpha_3)$-brimmed over $N^b_\alpha$
\sn
\item "{$(c)$}"  $\langle N^c_i:i < \alpha \rangle$ is
$\le_{\frak K}$-increasing (in fact follows)
\ermn
(Possible by \scite{600-nf.14.1}).  Now we can use \scite{600-nf.15}.
\bn
\underbar{Case B}:  
For each $i < \alpha$ we have: $N^a_{i+1}$ is $(\lambda,\kappa_i)$-brimmed
over $N^a_i$.  

In essence our problem is that we do not know anything about brimmness of
the $N^b_i$, so we shall ``correct it".

Let $\bar \delta^i = (\kappa_i,\lambda,\lambda)$. \newline
We can find a $\le_{\frak K}$-increasing sequence
$\langle M^x_i:i \le \alpha \rangle$ of models in $K_\lambda$ 
for $x \in \{ a,b,c\}$, continuous for $x=a,b$ such that 
$i < \alpha \Rightarrow M^a_i \le_{\frak K} M^b_i \le_{\frak K} 
M^c_i \le_{\frak K} M^b_{i+1}$ and $M^b_\alpha \le_{\frak K} M^c_\alpha$ and
$M^c_i$ is $(\lambda,\kappa_i)$-brimmed over $M^b_i$ (hence over
$M^a_i$) and
NF$_{\lambda,\bar \delta^i}(M^a_i,M^a_{i+1},M^c_i,M^b_{i+1})$ by
choosing $M^a_i,M^b_i,M^c_i$ by induction on 
$i,M^a_0 = N^a_0$ and $M^b_0$ is universal over $M^a_0$ recalling that
the NF$_{\lambda,\bar \delta^i}$ implies some brimness condition,
e.g. $M^b_{i+1}$ is $(\lambda,\text{cf}(\delta^i_3))$-brimmed over
$M^a_{i+1} \cup M^b_i$.
By Case A we know that
NF$_\lambda(M^a_0,M^a_\alpha,M^b_0,M^c_\alpha)$ holds.
\medskip

We can now choose an isomorphism $f^a_0$ from $N^a_0$ onto $M^a_0$, as
the identity (exists as $M^a_0 = N^a_0$) and then a $\le_{\frak K}$-embedding $f^b_0$
of $N^b_0$ into $M^b_0$ extending $f^a_0$.  Next we choose by induction on
$i \le \alpha,f^a_i$ an isomorphism from $N^a_i$ onto $M^a_i$ such that:
$j < i \Rightarrow f^a_j \subseteq f^a_i$, possible by ``uniqueness of the
$(\lambda,\kappa_i)$-brimmed model over $M^a_i$" so here we are
using the assumption of this case.
\medskip

Now we choose by induction on $i \le \alpha$, a $\le_{\frak K}$-embedding 
$f^b_i$ of
$N^b_i$ into $M^b_i$ extending $f^a_i$ and $f^b_j$ for $j < i$.  For
$i = 0$ we have done it, for $i$ limit use $\dsize \bigcup_{j < i} f^b_j$, 
lastly for $i$ a successor ordinal let $i = j+1$, now we have
\mr
\item "{$(*)_2$}"  NF$_\lambda(M^a_j,M^a_{j+1},f^b_j(N^b_j),M^b_{j+1})$ 
\newline
[why?  because NF$_{\lambda,\bar \delta^j}(M^a_j,M^a_{j+1},M^c_j,M^b_{j+1})$
by the choice of the \newline
$M^x_\zeta$'s hence by \scite{600-nf.12} we have
NF$_\lambda(M^a_j,M^a_{j+1},M^c_j,M^b_{j+1})$ and as \newline
$M^a_j = f^a_j(N^a_j) \le_{\frak K} f^b_j(N^b_j) \le M^b_j 
\le_{\frak K} M^c_j$ by \scite{600-nf.13} we get $(*)_2$.]
\ermn
By $(*)_2$ and the uniqueness of smooth amalgamation \scite{600-nf.11} and as
$M^b_{j+1}$ is $(\lambda,\text{cf}(\delta^3_j))$-brimmed over
$M^a_{j+1} \cup M^b_j$ hence over $M^a_{j+1} \cup f^b_j(N^b_j)$ clearly 
there is $f^b_i$ as required. \newline
So without loss of generality $f^a_\alpha$ is the identity, so we have
$N^a_0 = M^a_0,N^a_\alpha = M^a_\alpha,N^b_0 \le_{\frak K} M^b_0,N^b_\alpha
\le_{\frak K} M^b_\alpha$; also as said above NF$_\lambda(M^a_0,M^a_\alpha,
M^b_0,M^b_\alpha)$ holds (using Case A) so by 
monotonicity, i.e., \scite{600-nf.13} we get
NF$_\lambda(N^a_0,N^a_\alpha,N^b_0,N^b_\alpha)$ as required.
\mn
\underbar{Case C}:  General case.

We can find $M^\ell_i$ for $\ell < 3,i \le \alpha$ such that (note
that $M^1_0 = M^0_0$):
\medskip
\roster
\item "{$(a)$}"  $M^\ell_i \in K_\lambda$
\sn
\item "{$(b)$}"  for each $\ell < 3,M^\ell_i$ is $\le_{\frak
K}$-increasing in $i$ (but for $\ell=1,2$ they are not required to be
continuous) 
\sn
\item "{$(c)$}"  $M^0_i = N^a_i$
\sn
\item "{$(d)$}"  $M^{\ell +1}_{i+1}$ is $(\lambda,\lambda)$-brimmed over
$M^\ell_{i+1} \cup M^{\ell + 1}_i$ for $\ell < 2,i < \alpha$
\sn
\item "{$(e)$}"  NF$_\lambda(M^\ell_i,M^\ell_{i+1},M^{\ell + 1}_i,
M^{\ell + 1}_{i+1})$ for $\ell < 2,i < \alpha$
\sn
\item "{$(f)$}"  $M^1_0 = M^0_0$ and 
$M^2_0$ is $(\lambda,\text{cf}(\lambda))$-brimmed over $M^1_0$ 
\sn
\item "{$(g)$}"  for $\ell < 2$ and $i < \alpha$ limit we have
$$
M^{\ell +1}_i \text{ is } (\lambda,\lambda) \text{-brimmed over }
\dsize \bigcup_{j < i} M^{\ell +1}_j \cup M^\ell_i
$$
\item "{$(h)$}"  for $i < \alpha$ limit we have
$$
{\text{\rm NF\/}}_\lambda 
(\dsize \bigcup_{j < i} M^1_j,M^1_i,\dsize \bigcup_{j < i} M^2_j,
M^2_i).
$$
\noindent
[How?  As in the proof of \scite{600-nf.6A} or just do by hand.]
\ermn
Now note:
\mr
\item "{$(*)_3$}"  $M^{\ell +1}_i$ is $(\lambda,\text{cf}(\lambda \times
(1 + i)))$-brimmed over $M^\ell_i$ if $\ell = 1 \vee i \ne 0$
\newline
[why?  If $i=0$ by clause $(f)$, if $i$ a successor ordinal by clause $(d)$
and if $i$ is a limit ordinal then by clause (g)]
\sn
\item "{$(*)_4$}"  for $i < \alpha,{\text{\rm NF\/}}_\lambda
(M^0_i,M^0_{i+1},M^2_i,M^2
_{i+1})$. \newline
[Why?  If $i=0$ by clause (e) for $\ell=1,i=0$ we get
NF$_\lambda(M^1_0,M^1_1,M^2_0,M^2_1)$ so by clause (f) (i.e., $M^1_0 =
M^0_0$) and monotonicity (i.e., Claim \scite{600-nf.13}) we 
have NF$_\lambda(M^0_0,M^1_0,M^2_0,M^2_1)$ as
required.  If $i>0$ we use 
Case B for $\alpha = 2$ with $M^0_i,M^0_{i+1},M^1_i,M^1_{i+1}, 
M^2_i,M^2_{i+1}$ here standing for $N^a_0,N^b_0,N^a_1,N^b_1,
N^a_2,N^b_2$ there (and symmetry).]
\ermn
Let us define $N^\ell_i$ for $\ell < 3,i \le \alpha$ by: $N^\ell_i$ is
$M^\ell_i$ if $i$ is non-limit and $N^\ell_i = \cup\{N^\ell_j:j < i\}$
if $i$ is limit.
\mr
\widestnumber\item{$(*)_5(iii)$}
\item "{$(*)_5(i)$}"  $\langle N^\ell_i:i \le \alpha \rangle$ is
$\le_{\frak K}$-increasing continuous, $N^0_i = N^a_i$ and $N^\ell_i
\le_{\frak K} M^\ell_i$
\sn
\item "{$(ii)$}"  for $i < \alpha$, NF$_\lambda(N^0_i,N^0_{i+1},
N^2_i,N^2_{i+1})$ \nl
[why? by $(*)_4 +$ monotonicity of NF$_\lambda$]
\sn
\item "{$(iii)$}"  for $i < \alpha,N^2_{i+1}$ is
$(\lambda,\text{cf}(\lambda))$-brimmed over $N^0_{i+1} \cup N^2_i$ and
even over $N^1_{i+1} \cup N^2_i$
\nl
[why?  by clause (d)]
\sn
\item "{$(*)_6$}"  NF$_{\lambda,\langle \lambda,\lambda,1
\rangle}(N^1_0,N^1_\alpha,N^2_0,N^2_\alpha)$. \nl
[Why?  As we have proved case A (or, if you prefer, by \scite{600-nf.15};
easily the assumption there holds).]  
\ermn
Choose $f^a_i = \text{ id}_{N^a_i}$ for $i \le \alpha$ and let $f^b_0$
be a $\le_{\frak K}$-embedding of $N^b_0$ into $N^2_0$.

Now we continue as in Case B defining by induction on $i$ 
a $\le_{\frak K}$-embedding $f^b_i$ of $N^b_i$ into $N^2_i$, 
the successor case is possible by $(*)_5(ii) + (*)_5(iii)$.
In the end by $(*)_6$ and monotonicity of NF$_\lambda$ (i.e., Claim
\scite{600-nf.13}) we are done.  \newline
2) Apply for each $i < \alpha_2$ part (1) to the sequences $\langle M_{\beta,i}:\beta 
\le \alpha_1 \rangle,\langle M_{\beta,i+1}:\beta \le \alpha_1 \rangle$ so we
get $\nonforkin{M_{\alpha_1,i}}{M_{0,i+1}}_{M_{0,i}}^{M_{\alpha_1,i+1}}$ 
hence by symmetry (i.e., \scite{600-nf.11}) 
we have $\nonforkin{M_{0,i+1}}{M_{\alpha_1,i}}
_{M_{0,i}}^{M_{\alpha_1,i+1}}$.  Applying part (1) to the sequences
$\langle M_{0,j}:j \le \alpha_2 \rangle,\langle M_{\alpha_1,j}:j \le \alpha_2
\rangle$ we get $\nonforkin{M_{0,\alpha_2}}{M_{\alpha_1,0}}_{M_{0,0}}
^{M_{\alpha_1,\alpha_2}}$ hence by symmetry (i.e. \scite{600-nf.11}) we have
$\nonforkin{M_{\alpha_1,0}}{M_{0,\alpha_2}}_{M_{0,0}}
^{M_{\alpha_1,\alpha_2}}$; so we get 
the desired conclusion. \hfill$\square_{\scite{600-nf.16}}$
\enddemo
\bigskip

\proclaim{\stag{600-nf.17} Claim}   Assume $\alpha < \lambda^+,
\langle N^\ell_i:i \le \alpha \rangle$
is $\le_{\frak K}$-increasing continuous sequence of models
for $\ell = 0,1$ where $N^\ell_i \in K_\lambda$ and $N^1_{i+1}$ is 
$(\lambda,\kappa_i)$-brimmed over $N^0_{i+1} \cup N^1_i$ and
${\text{\rm NF\/}}_\lambda(N^0_i,N^1_i,N^0_{i+1},N^1_{i+1})$.

\ub{Then} $N^1_\alpha$ is $(\lambda,{\text{\rm cf\/}}(\dsize 
\sum_{i < \alpha} \kappa_i))$-brimmed over $N^0_\alpha \cup N^1_0$.  
\endproclaim
\bigskip

\remark{\stag{600-nf.17A} Remark}  1) If our framework is uni-dimensional 
(see \sectioncite[\S2]{705}; as for example 
when it comes from \cite{Sh:576}) we can simplify the proof. \nl
2) Assuming only ``$N^1_{i+1}$ is 
universal over $N^0_{i+1} \cup N^1_i$" suffices when $\alpha$ is a
limit ordinal, i.e., we get $N^1_\alpha$ is $(\lambda,\text{\rm
cf}(\alpha))$-brimmed over $N^0_\alpha$.   Why?  We choose $N^2_j$ for
$j \le i$ such that $N^2_j = N^1_j$ if $j=0$ or $j$ a limit ordinal
and $N^2_j$ is a model $\le_{\frak s} N^1_j$ and
$(\lambda,\kappa_1)$-brimmed over $N^0_j \cup N^1_i$ when $j=i+1$.
Now $\langle N^2_j:j \le \alpha\rangle$ satisfies all the requirements
in $\langle N^1_j:j \le \alpha\rangle$ in \scite{600-nf.17}.
\nl
3) We could have proved this earlier and used it, e.g. in \scite{600-nf.16}.
\endremark
\bigskip

\demo{Proof}  The case $\alpha$ not a limit ordinal is trivial so 
assume $\alpha$ is a
limit ordinal.  We choose by induction on $i \le \alpha$, an ordinal
$\varepsilon(i)$ and a sequence
$\langle M_{i,\varepsilon}:\varepsilon \le \varepsilon(i) \rangle$ and
$\langle c_\varepsilon:\varepsilon < \varepsilon(i) \text{ non-limit} 
\rangle$ such that:
\mr
\item "{$(a)$}"  $\langle M_{i,\varepsilon}:\varepsilon \le \varepsilon(i)
\rangle$ is (strictly) $<_{\frak K}$-increasing continuous in $K_\lambda$
\sn
\item "{$(b)$}"  $N^0_i \le_{\frak K} M_{i,\varepsilon} \le_{\frak K}
N^1_i$
\sn
\item "{$(c)$}"  $N^0_i = M_{i,0}$ and $N^1_i = M_{i,\varepsilon(i)}$
\sn
\item "{$(d)$}"  $\varepsilon(i)$ is (strictly) increasing continuous in
$i$ and $\varepsilon(i)$ is divisible by $\lambda$
\sn
\item "{$(e)$}"  $j < i \and \varepsilon \le \varepsilon(j) \Rightarrow
M_{i,\varepsilon} \cap N^1_j = M_{j,\varepsilon}$
\sn
\item "{$(f)$}"  for $j<i$ and $\varepsilon \le \varepsilon(j+1)$, the 
sequence $\langle M_{\beta,\varepsilon}:\beta \in
(j,i] \rangle$ is $\le_{\frak K}$-increasing continuous
\sn
\item "{$(g)$}"  for $j<i,\varepsilon < \varepsilon(j)$ non-limit; the type
\ortp$(c_\varepsilon,M_{i,\varepsilon},M_{i,\varepsilon +1}) \in 
{\Cal S}^{\text{bs}}(M_{i,\varepsilon})$ does not fork over 
$M_{j,\varepsilon}$ (actually, here allowing all $\varepsilon$ is O.K., too)
\sn
\item "{$(h)$}"   $M_{i+1,\varepsilon +1}$ is $(\lambda,
\text{cf}(\lambda))$-brimmed over $M_{i+1,\varepsilon} \cup 
M_{i,\varepsilon +1}$
\sn
\item "{$(i)$}"  if $\varepsilon < \varepsilon(i)$ and $p \in 
{\Cal S}^{\text{bs}}(M_{i,\varepsilon})$ then for $\lambda$ successor
ordinals $\xi \in [\varepsilon,\varepsilon(i))$ the type 
\ortp$(c_\xi,M_{i,\xi},M_{i,\xi +1})$ is a non-forking extension of $p$.
\ermn
If we succeed, then $\langle M_{\alpha,\varepsilon}:\varepsilon \le 
\varepsilon(\alpha) \rangle$ is a (strictly)
$<_{\frak K}$-increasing continuous sequence of models from $K_\lambda,
M_{\alpha,0} = N^0_\alpha$,
and $M_{\alpha,\varepsilon(\alpha)} = N^1_\alpha$.
We can apply \scite{600-4a.2}
and we conclude that $N^1_\alpha = M_{\alpha,\varepsilon(\alpha)}$ is
$(\lambda,\text{cf}(\alpha))$-brimmed over $M_{\alpha,\varepsilon(j)}$
hence over $N^0_\alpha \cup N^1_0$ (both $\le_{\frak K} M_{\alpha,1}$).

Carrying the induction is easy.  For $i=0$, 
there is not much to do.  For $i$ successor we use ``$N^j_{i+1}$ is
brimmed over $N^0_{i+1} \cup N^1_i$" the existence of 
non-forking amalgamations and \scite{600-4a.1}, bookkeeping and the extension
property $(E)(g)$.  For $i$ limit we have no problem.
\hfill$\square_{\scite{600-nf.17}}$
\enddemo
\bigskip

\demo{\stag{600-nf.18} Conclusion}  1) If NF$_\lambda(N_0,N_1,N_2,N_3)$ and
$\langle M_{0,\varepsilon}:\varepsilon \le \varepsilon(*) \rangle$ is an
$\le_{\frak K}$-increasing continuous sequence of models from $K_\lambda$, 
$N_0 \le_{\frak K} M_{0,\varepsilon} \le_{\frak K} N_2$ \underbar{then} 
we can find $\langle M_{1,\varepsilon}:
\varepsilon \le \varepsilon(*) \rangle$ and $N'_3$ such that:
\medskip
\roster
\item "{$(a)$}"  $N_3 \le_{\frak K} N'_3 \in K_\lambda$
\sn
\item "{$(b)$}"  $\langle M_{1,\varepsilon}:\varepsilon \le \varepsilon(*)
\rangle$ is $\le_{\frak K}$-increasing continuous 
\sn
\item "{$(c)$}"  $M_{1,\varepsilon} \cap N_2 = M_{0,\varepsilon}$
\sn
\item "{$(d)$}"  $N_1 \le_{\frak K} M_{1,\varepsilon} \le_{\frak K} N'_3$
\sn
\item "{$(e)$}"  if $M_{0,0} = N_0$ then $M_{1,0} = N_1$
\sn
\item "{$(f)$}"
NF$_\lambda(M_{0,\varepsilon},M_{1,\varepsilon},N_2,N'_3)$, for every
$\varepsilon \le \varepsilon(*)$.
\ermn
2) If $N_3$ is universal over $N_1 \cup N_2$, then without loss of generality
$N'_3 = N_3$.
\nl
3) In part (1) we can add
\mr
\item "{$(g)$}"  $M_{1,\varepsilon +1}$ is brimmed over
$M_{0,\varepsilon +1} \cup M_{1,\varepsilon}$.
\endroster
\enddemo
\bigskip

\demo{Proof}  1) Define $M'_{0,i}$ for $i \le \varepsilon^* := 1 +
\varepsilon(*) +1$ by $M'_{0,0} = N_0,M'_{0,1 + \varepsilon} =
M_{0,\varepsilon}$ for $\varepsilon \le \varepsilon(*)$ and $M'_{0,1 +
\varepsilon(*)+1} = N_2$.  By existence (\scite{600-nf.10.3}) 
we can find an $\le_{\frak K}$-increasing continuous sequence $\langle
M'_{1,\varepsilon}:\varepsilon \le \varepsilon^* \rangle$ with
$M'_{1,0} = N_1$ and $\le_{\frak K}$-embedding $f$ of $N_2$ into 
$M'_{1,\varepsilon^*}$ such that $\varepsilon < \varepsilon^*
\Rightarrow \text{\rm NF}_\lambda(f(M'_{0,\varepsilon}),M'_{1,0},
f(M'_{0,\varepsilon +1}),M'_{1,\varepsilon +1})$.  By transitivity we have
NF$_\lambda(f(M'_{0,0}),M'_{1,0},
f(M'_{0,\varepsilon^*}),M'_{1,\varepsilon^*})$.
By disjointness (i.e., $f(M'_{0,\varepsilon^*}) \cap M'_{1,0} =
M'_{0,0}$, see \scite{600-nf.3}(3))
\wilog \, $f$ is the identity.  By uniqueness for NF there are
$N'_3,N_3 \le_{\frak K} N'_3 \in K_\lambda$ and $\le_{\frak
K}$-embedding of $M'_{1,\varepsilon^*}$ onto $N'_3$ over $N_1 \cup N_2
= M'_{0,\varepsilon^*} \cup M'_{1,0}$ so we are done. 
\nl
2) Follows by (1).
\nl
3) Similar to (1).  \hfill$\square_{\scite{600-nf.18}}$
\enddemo
\bigskip

\proclaim{\stag{600-nf.19} Claim}  {\rm NF}$_\lambda$ respects ${\frak s}$;
that is assume 
${\text{\rm NF\/}}_\lambda(M_0,M_1,M_2,M_3)$ and
$a \in M_1 \backslash M_0$ satisfies {\rm \ortp}$(a,M_0,M_3) 
\in {\Cal S}^{\text{bs}}(M_0)$, 
\ub{then} ${\text{\rm \ortp\/}}(a,M_2,M_3) \in {\Cal S}^{\text{bs}}(M_2)$
does not fork over $M_0$.
\endproclaim
\bigskip

\demo{Proof}  Without loss of generality $M_1$ is $(\lambda,*)$-brimmed
over $M_0$.  [Why?  By the existence we can find $M^+_1$ which is a
$(\lambda,*)$-brimmed extension of $M_1$.  By the existence for
NF$_\lambda$ without loss of generality we can find $M^+_3$ such that
NF$_\lambda(M_1,M^+_1,M_3,M^+_3)$, hence by transitivity for
NF$_\lambda$ we have NF$_\lambda(M_0,M^+_1,M_2,M^+_3)$.]  
By the hypothesis of the section there are $M'_1,a'$ such that
$M_0 \cup \{a'\} \subseteq M'_1$ and $\ortp(a',M_0,M'_1) =
\text{\rm \ortp}(a,M_0,M_1)$  and $(M_0,M'_1,a) \in
K^{3,\text{uq}}_\lambda$; as $M^+_1$ is
$(\lambda,*)$-brimmed over $M_0$ without loss of generality 
$M' \le_{\frak K} M^+_1$
and $a' = a$ and $M_1$ is $(\lambda,*)$-brimmed over
$M'_1$.  We can apply \scite{600-nf.0A} to $M'_1,M^+_1$ getting $\langle
M^*_i,a_i:i \le \delta < \lambda^+ \rangle$ as there.  Let $M'_i$ be:
$M_0$ if $i=0,M^*_j$ if $1+j = i$ so $M'_1 = M^*_0 = M'_1$ and let $a_i$
be $a$ if $i=0,a_j$ if $1+j=i$.  So we can find $M'_3$ and $f$ such
that $M_2 \le_{\frak K} M'_3,f$ is a $\le_{\frak K}$-embedding of
$M^+_1$ into $M'_3$ extending id$_{M_0}$ such that NF$_{\lambda,\langle
\delta,\lambda,\lambda \rangle}
(M_0,f(M^+_1),M_2,M'_3)$ and $M'_3$, this is witnessed by
$\langle f(M'_i):i \le \delta \rangle,\langle M''_i:i \le \delta \rangle,
\langle f(a_i):i < \delta \rangle$ and $M''_0 = M_2$; this is
possible by \scite{600-nf.3}(2).  Hence NF$_\lambda(M_0,f(M^+_1),M_2,N) =
\text{\rm NF}_\lambda(f(M'_0),f(M'_\delta),M''_0,N)$
hence by the uniqueness for NF$_\lambda$ \wilog \, 
$f = \text{ id}_{M^+_1}$ and $M_3 \le_{\frak K} N$.  
By the choice of $f,N$ we have \ortp$(a,M_2,M_3) = \text{\ortp}(a_0,M_2,N) = 
\text{ \ortp}(a_0,M''_0,M'_1) \in {\Cal S}^{\text{bs}}(M''_0) 
= {\Cal S}^{\text{bs}}(M_2)$ does not fork over $M'_0 = M_0$ as required.
\hfill$\square_{\scite{600-nf.19}}$
\enddemo
\bigskip

\demo{\stag{600-nf.20} Conclusion}  If  $M_0 \le_{\frak K} M_\ell 
\le_{\frak K} M_3$ for $\ell = 1,2$ and $(M_0,M_1,a) \in 
K^{3,\text{uq}}_\lambda$ and ${\text{\rm \ortp\/}}(a,M_2,M_3) \in 
{\Cal S}^{\text{bs}}(M_2)$ does not fork over $M_0$ \ub{then}
${\text{\rm NF\/}}(M_0,M_1,M_2,M_3)$.
\enddemo
\bigskip

\demo{Proof}  By the definition of $K^{3,\text{uq}}_\lambda$ and
existence for NF$_\lambda$ and \scite{600-nf.19} (or use \scite{600-nf.0Y} +
\scite{600-nf.20.7}.   \hfill$\square_{\scite{600-nf.20}}$
\enddemo
\bn
We can sum up our work by
\demo{\stag{600-nf.20.7} Main Conclusion}  NF$_\lambda$ is a non-forking
relation on ${}^4({\frak K}_\lambda)$ which respects ${\frak s}$.
\enddemo
\bigskip

\demo{Proof}  We have to check clauses (a)-(g)+(h) from \scite{600-nf.0X}.
Clauses (a),(b) hold by the Definition \scite{600-nf.2} of NF$_\lambda$.
Clauses $(c)_1,(c)_2$, i.e.,  monotonicity hold by \scite{600-nf.13}.
Clause (d), i.e.,  symmetry holds by \scite{600-nf.14}.
Clause (e), i.e.,  transitivity holds by \scite{600-nf.16}.
Clause (f), i.e.,  existence hold by \scite{600-nf.10.3}.
Clause (g), i.e.,  uniqueness holds by \scite{600-nf.11}.

Lastly, clause (h), i.e., NF$_\lambda$ respecting ${\frak s}$ by 
\scite{600-nf.19}.  \hfill$\square_{\scite{600-nf.20.7}}$
\enddemo
\bn
The following definition is not needed for now but is 
natural (of course, we can omit ``there is
superlimit" from the assumption and the conclusion).  For the rest of
the section we stop assuming Hypothesis \scite{600-nf.0}.
\definition{\stag{600-nf.20.7A} Definition}  1) A good $\lambda$-frame 
${\frak s}$ is type-full when for $M \in {\frak K}_{\frak s},{\Cal
S}^{\text{bs}}(M) = {\Cal S}^{\text{na}}_{{\frak K}_\lambda}(M)$.
\nl
2) Assume ${\frak K}_\lambda$ is a $\lambda$-a.e.c. and NF is a
4-place relation on $K_\lambda$.  We define ${\frak t} = {\frak
t}_{{\frak K}_\lambda,\text{NF}} = (K_{\frak t},\nonfork{}{}_{\frak
t},{\Cal S}^{\text{bs}}_{\frak t})$ as follows:
\mr
\item "{$(a)$}"  ${\frak K}_{\frak t}$ is the $\lambda$-a.e.c. ${\frak
K}_\lambda$
\sn
\item "{$(b)$}"  ${\Cal S}^{\text{bs}}_{\frak t}(M)$ is ${\Cal
S}^{\text{na}}_{{\frak K}_\lambda}(M)$ for $M \in {\frak K}_\lambda$
\sn
\item "{$(c)$}"  $\nonfork{}{}_{\frak t}$ is defined by:
$(M_0,M_1,a,M_3) \in \nonfork{}{}_{\frak t}$ when we can find
$M_2,M'_3$ such that $M_0 \le_{{\frak K}_\lambda} M_2 \le_{{\frak
K}_\lambda} M'_3,M_3 \le_{{\frak K}_\lambda} M'_3,a \in M_2 \backslash M_0$ and
NF$(M_0,M_1,M_2,M'_3)$.
\endroster
\enddefinition
\bigskip

\proclaim{\stag{600-nf.20.8} Claim}  1) Assume that
\mr
\item "{$(a)$}"  ${\frak K}_\lambda$ is a $\lambda$-a.e.c. with
amalgamation and a superlimit model
\sn
\item "{$(b)$}"  ${\frak K}_\lambda$ is stable
\sn
\item "{$(c)$}"  {\rm NF} is a ${\frak K}_\lambda$-non-forking
relation, see Definition \scite{600-nf.0X}(1).
\ermn
\ub{Then} ${\frak t} = {\frak t}_{{\frak K}_\lambda,\text{NF}}$ is a 
type-full good $\lambda$-frame.
\nl
2) Assume that ${\frak s}$ is a good $\lambda$-frame which has
existence for $K^{3,\text{uq}}_\lambda$ (see \scite{600-nf.0}(2)) and {\rm
NF} = {\rm NF}$_\lambda$.  \ub{Then} ${\frak t}$ is very close to
${\frak s}$, i.e.:
\mr
\item "{$(a)$}"  ${\frak K}_{\frak s} = {\frak K}_{\frak t}$
\sn
\item "{$(b)$}"  if $p \in {\Cal S}^{\text{bs}}_{\frak s}(M_1)$ and
$M_0 \le_{{\frak K}_\lambda} M_1$ \ub{then} $p \in {\Cal
S}^{\text{bs}}_{\frak t}(M_1)$ and $p$ forks over $M_0$ for ${\frak
s}$ iff $p$ forks over $M_0$ for ${\frak t}$.
\endroster
\endproclaim
\bigskip

\demo{Proof}  For the time being, left to the reader (but before it is
really used it is proved in \marginbf{!!}{\cprefix{705}.\scite{705-9.11A}}).
\enddemo
\bigskip

\remark{Remark}  Note that this actually says that from now on we
could have used type-full ${\frak s}$, but it is not necessary for a
long time.
\endremark
\bigskip

\definition{\stag{600-nf.20.9} Definition}  1) Let ${\frak s}$ be a good
$\lambda$-frame.  We say that NF is a weak
${\frak s}$-non-forking relation where
\mr
\item "{$(a)$}"  NF is  weak ${\frak K}_{\frak s}$-non-forking
relation, see Definition \scite{600-nf.0X}(2), i.e., uniqueness is omitted
\sn
\item "{$(b)$}"  NF respects ${\frak s}$, see Definition
\scite{600-nf.0X}(3)
\sn
\item "{$(c)$}"  NF satisfies \scite{600-nf.18}, (NF-lifting of an
$\le_{\frak K}$-increasing sequence).
\ermn
2) We say ${\frak s}$ is pseudo-successful if some NF is a weak
${\frak s}$-non-forking relation witnesses it.
\enddefinition
\bigskip

\demo{\stag{600-nf.20.10} Observation}  1) If ${\frak s}$ is a good
$\lambda$-frame which is weakly successful (i.e., has 
existence for $K^{3,\text{uq}}_\lambda$, i.e., \scite{600-nf.0}) \ub{then}
{\rm NF}$_\lambda = \text{\rm NF}_{\frak s}$ is a weak ${\frak
s}$-non-forking relation.
\nl
2) If ${\frak s}$ is a good $\lambda$-frame and {\rm NF} is a weak
${\frak s}$-non-forking relation then \scite{600-nf.20} holds.
\nl
3) If ${\frak s}$ is a good $\lambda$-frame and NF is an 
${\frak s}$-forking relation \ub{then} NF is a weak 
${\frak s}$-non-forking relation.
\enddemo
\bigskip

\demo{Proof}  Straight.
\nl
1) Follows by \scite{600-nf.20.7}, NF$_\lambda$ satisfies clauses (a)+(b)
and by \scite{600-nf.18} it satisfies also clause (c) of Definition
\scite{600-nf.0X}(1).
\nl
2) Also easy.
\nl
3) We have just to check the proof of \scite{600-nf.18} still works.
\enddemo
\bigskip

\remark{\stag{600-nf.21.3} Remark}  1) In \chaptercite{705} ,\S1 -\S11 we can
use ``${\frak s}$ is pseudo successful as witnessed by NF + lifting of
decompositions" instead of ``${\frak s}$ is weakly successful".  We
shall return to this elsewhere, see \cite{Sh:838}, \cite{Sh:842}.
\endremark
\newpage

\head {\S7 Nice extensions in $K_{\lambda^+}$} \endhead  \resetall \sectno=7
 \spuriousreset
\bigskip

\demo{\stag{600-ne.0} Hypothesis}  Assume the hypothesis \scite{600-nf.0}. \nl
\enddemo
\bn
So by \S6 we have reasonable control on \underbar{smooth} amalgamation in
$K_\lambda$.  We use this to define ``nice" extensions in
$K_{\lambda^+}$ and prove some basic properties.
This will be treated again in \S8. 
\bigskip

\definition{\stag{600-ne.1} Definition}  1) $K^{\text{nice}}_{\lambda^+}$
is the class of saturated $M \in K_{\lambda^+}$. \nl
2) Let $M_0 \le^*_{\lambda^+} M_1$ mean:
\mr
\item "{{}}"  $M_0 \le_{\frak K} M_1$ and they are from
$K_{\lambda^+}$ and we can find $\bar M^\ell = \langle M^\ell_i:i < \lambda^+
\rangle$, a $\le_{\frak K}$-representation of $M_\ell$ for $\ell = 0,1$ 
such that: \newline
NF$_\lambda(M^0_i,M^0_{i+1},M^1_i,M^1_{i+1}$) for $i < \lambda^+$.
\endroster
\medskip
\noindent
3)  Let $M_0 <^+_{\lambda^+,\kappa} M_1$ mean 
\footnote{Note that $M_0 <^+_{\lambda^+,\kappa} M_1$ implies $M_1 \in
K^{\text{nice}}_{\lambda^+}$ but in general $M_0 \in
K^{\text{nice}}_{\lambda^+}$ does not follow.} 
that
$(M_0,M_1 \in K_{\lambda^+}$ and) $M_0 \le^*_{\lambda^+} M_1$ by
some witnesses $M^\ell_i$ (for $i < \lambda^+,\ell < 2$) such that
NF$_{\lambda,\langle 1,1,\kappa \rangle}(M^0_i,M^0_{i+1},
M^1_i,M^1_{i+1})$ for $i < \lambda^+$; 
of course $M_0 \le_{\frak K} M_1$ in this case.
Let $M_0 \le^+_{\lambda^+,\kappa} M_1$ mean $(M_0 = M_1 \in
K_{\lambda^+}) \vee (M_0 <^+_{\lambda^+,\kappa} M_1)$.  
If $\kappa = \lambda$, we may omit it. \nl
4) Let $K^{3,\text{bs}}_{\lambda^+} = \{(M,N,a):M \le^*_{\lambda^+} N$ 
are from $K_{\lambda^+}$ and $a \in N \backslash M$ and for some $M_0 
\le_{\frak K} M,M_0 \in K_\lambda$ we have $[M_0 \le_{\frak K} M_1 
\le_{\frak K} M \and M_1 \in K_\lambda$ \ub{implies} $\text{ \ortp}(a,M_1,N) \in
{\Cal S}^{\text{bs}}(M_1)$ and does not fork over $M_0]\}$. 
We call $M_0$ or \ortp$(a,M_0,N)$ a witness for $(M,N,a) \in 
K^{3,\text{bs}}_{\lambda^+}$.  (In fact this definition 
on $K^{3,\text{bs}}_{\lambda^+}$ is compatible with the
definition in \S2 for triples such that $M \le^*_{\lambda^+} N$ but we
do not know now whether even
$(K^{\text{nice}}_{\lambda^+},\le^*_{\lambda^+})$ is a $\lambda^+$-a.e.c..)
\enddefinition
\bigskip

\proclaim{\stag{600-ne.2} Claim}  0) $K^{\text{nice}}_{\lambda^+}$ has one
and only one model up to isomorphisms and
$M \in K^{\text{nice}}_{\lambda^+}$ implies
$M \le^*_{\lambda^+} M$ and $M \le^+_{\lambda^+} M$; 
moreover, $M \in K_{\lambda^+} \Rightarrow
M \le^*_{\lambda^+} M$.  Also $\le^*_{\lambda^+}$ is a partial order
and if $M_\ell \in K_{\lambda^+}$ for $\ell=0,1,2$ and
$M_0 \le_{\frak K} M_1 \le_{\frak K} M_2$ and $M_0 \le^*_{\lambda^+}
M_2$ then $M_0 \le^*_{\lambda^+} M_1$. 
\nl
1) If $M_0 \le^*_{\lambda^+} M_1$ and 
$\bar M^\ell = \langle M^\ell_i:
i < \lambda^+ \rangle$ is a representation of $M_\ell$ for $\ell=0,1$ \ub{then}
\mr
\item "{$(*)$}"  for some club $E$ of $\lambda^+$,
{\roster
\itemitem{ $(a)$ }  for every $\alpha < \beta$ from $E$ we have
{\rm NF}$_\lambda(M^0_\alpha,M^0_\beta,M^1_\alpha,M^1_\beta)$
\sn
\itemitem{ $(b)$ }  if $\ell < 2$ and $M_\ell \in
K^{\text{nice}}_{\lambda^+}$ then for $\alpha < \beta$ from $E$ the
model $M^\ell_\beta$ is $(\lambda,*)$-brimmed over $M^\ell_\alpha$.
\endroster}
\ermn
2) Similarly for $<^+_{\lambda^+,\kappa}$: if $M_0
<^+_{\lambda^+,\kappa} M_1,\bar M^\ell = \langle \bar M^\ell_i:i <
\lambda^+ \rangle$ a representation of $M_\ell$ for $\ell = 0,1$
\ub{then} for some club $E$ of $\lambda^+$ for every $\alpha < \beta$
from $E$ we have {\rm NF}$_{\lambda,\langle 1,1,\kappa\rangle}
(M^0_\alpha,M^0_\beta,M^1_\alpha,M^1_\beta)$, moreover
NF$_{\lambda,\langle 1,\text{cf}(\lambda \times (1 +
\beta)),\kappa\rangle}(M^0_\alpha,M^0_\beta,M^1_\alpha,M^1_\beta)$ and if
$(M_\alpha,\bar M^0_\beta,M^1_\alpha,M^1_\beta),M_0 \in
K^{\text{nice}}_{\lambda^+}$ then we can add NF$_{\lambda,\langle
\lambda,\text{cf}(\lambda \times(1 + \beta)),\kappa)}(M^0_\alpha,M^0_\beta,M'_\alpha,M'_\beta)$.
\nl
3) The $\kappa$ in Definition \scite{600-ne.1}(3) does not matter.  
\nl
4) If $M_0 <^+_{\lambda^+,\kappa} M_1$, \ub{then} $M_1 \in
K^{\text{nice}}_{\lambda^+}$. \nl
5) If $M \in K_{\lambda^+}$ is saturated, equivalently $M \in
K^{\text{nice}}_{\lambda^+}$ \ub{then} $M$ has a $\le_{\frak K}$-representation
$\bar M = \langle M_\alpha:\alpha < \lambda^+ \rangle$ such that
$M_{i+1}$ is $(\lambda,\lambda)$-brimmed over $M_i$ for $i <
\lambda^+$ and also the inverse is true.
\nl
6) If $M \le^*_{\lambda^+} N$ and $N_0 \le_{\frak K} N,N_0 \in
K_\lambda$ \ub{then} we can find $M_1 \le_{\frak K} N_1$ from
$K_\lambda$ such that $M_1 \le_{\frak K} M,N_0 \le_{\frak K} N_1
\le_{\frak K} N$ and: for every $M_2 \in K_\lambda$ satisfying 
$M_1 \le_{\frak K} M_2 \le_{\frak K} M$ 
there is $N_2 \le_{\frak K} N$ such that
{\rm NF}$_{\frak s}(M_1,M_2,N_1,N_2)$. 
\endproclaim
\bigskip

\demo{Proof}  0) Obvious by now (for the second sentence use part (1)
and NF$_{\frak s}$ being a non-forking relation on ${\frak K}_{\frak
s}$) in particular transitivity and monotonicity.
\nl
1) Straight by \scite{600-nf.16} as any two representations
agree on a club. \newline
2) Up to ``moreover" quite straight.  For the ``moreover" use
\scite{600-nf.17} to show that $M^1_\beta$ is
$(\lambda,\text{cf}(\beta))$-brimmed over $M^0_\beta$.  Lastly, for
the ``we can add" just use part (5), choosing thin enough club $E$ of
$\lambda^+$ then use $\{\alpha \in E:\text{otp}(\alpha \cap E)$ is
divisible by $\lambda\}$.
\nl
3) By \scite{600-nf.17}. \nl
4) By \scite{600-nf.17}. \nl
5) Trivial.
\nl
6) Easy.   \hfill$\square_{\scite{600-ne.2}}$
\enddemo
\bigskip

\proclaim{\stag{600-ne.3} Claim}  0) For every $M_0 \in K_{\lambda^+}$ for some
$M_1 \in K^{\text{nice}}_{\lambda^+}$ we have $M_0 \le_{\frak K} M_1$. \nl
1) For every $M_0 \in K_{\lambda^+}$ and 
$\kappa = { \text{\rm cf\/}}(\kappa) \le \lambda$ 
for some $M_1 \in K_{\lambda^+}$ we have $M_0 <^+_{\lambda^+,\kappa} M_1$ so
$M_1 \in K^{\text{nice}}_{\lambda^+}$. \newline
1A) Moreover, if $N_0 \le_{\frak K} M_0 \in K_{\lambda^+},N_0 \in K_\lambda,
p \in {\Cal S}^{\text{bs}}(N_0)$ \ub{then} in (1) we can add that for 
some $a,(M_0,M_1,a) \in K^{3,\text{bs}}_\lambda$ as witnessed by $p$. \nl
2)  $\le^*_{\lambda^+}$ and $<^+_{\lambda^+,\kappa}$ are transitive. \newline
3)  If $M_0 \le_{\frak K} M_1 \le_{\frak K} M_2$ are in $K_{\lambda^+}$
and $M_0 \le^*_{\lambda^+} M_2$, \ub{then} 
$M_0 \le^*_{\lambda^+} M_1$. \newline
4) If $M_1 <^+_{\lambda^+,\kappa} M_2$, \ub{then} 
$M_1 <^*_{\lambda^+} M_2$. \nl
5) If $M_0 <^*_{\lambda^+} M_1 
<^+_{\lambda,\kappa} M_2$ \ub{then} $M_0 <^+_{\lambda,\kappa} M_2$.
\endproclaim
\bigskip

\demo{Proof}  0) Easy and follows by the proof of part (1) below. \nl
1), 1A)  Let $\langle M^0_i:i < \lambda^+ \rangle$ 
be a $\le_{\frak K}$-representation of $M_0$ with $M^0_i$ brimmed and
brimmed over $M^0_j$ for $j <i$ and for part (1A) we have
$M^0_0 = N_0$, and for part (1) let $p$ be any member of 
${\Cal S}^{\text{bs}}(M^0_0)$.
We choose by induction on $i$ a model $M^1_i \in K_\lambda$ and 
$a \in M^1_0$ such that $M^1_i$ is $(\lambda,\text{cf}(\lambda \times 
(1+i)))$-brimmed over 
$M^0_i,\langle M^1_i:i < \lambda^+ \rangle$ is 
$<_{\frak K}$-increasing continuous, $M^1_i \cap M_0 = M^0_i$ and
\ortp$(a,M^0_0,M^1_0) = p$ and $M^1_{i+1}$ is
$(\lambda,\kappa)$-brimmed over $M^0_{i+1} \cup M^1_i$ and
NF$_{\lambda,\langle 1,\text{cf}(\lambda \times (1+i)),\kappa 
\rangle}(M^0_i,M^0_{i+1},
M^1_i,M^1_{i+1})$ for $i < \lambda^+$.  
Note that for limit $i$, by \scite{600-nf.17}, $M^1_i$ is
$(\lambda,\text{cf}(i))$-brimmed over $M^0_i \cup M^1_j$ for any $j<i$.

Note that for $i < \lambda^+$, the type \ortp$(a,M^0_i,M^1_i)$ does not
fork over $M^0_0 = N_0$ and extends $p$ by \scite{600-nf.19} (saying
NF$_\lambda$ respects ${\frak s}$) 
\scite{600-nf.14} (symmetry) and \scite{600-nf.12}.  So clearly we are done. \nl
2) Concerning $<^+_{\lambda^+,\kappa}$ use \scite{600-ne.2} 
and \scite{600-nf.16} (i.e. transitivity for smooth amalgamations).  
The proof for $<^*_{\lambda^+}$ is the same. \newline
3) By monotonicity for smooth amalgamations in ${\frak K}_\lambda$; 
i.e., \scite{600-nf.13}. \newline
4), 5)  Check. \hfill$\square_{\scite{600-ne.3}}$
\enddemo
\bigskip

\proclaim{\stag{600-ne.2A} Claim}  1) If $(M_0,M_1,a) \in 
K^{3,\text{bs}}_{\lambda^+}$
and $M_1 \le^*_{\lambda^+} M_2 \in K_{\lambda^+}$ \ub{then} \nl
$(M_0,M_2,a) \in K^{3,\text{bs}}_{\lambda^+}$. \nl
2) If $M_0 <^*_{\lambda^+} M_1$, \ub{then} for some $a,(M_0,M_1,a) \in
K^{3,\text{bs}}_{\lambda^+}$. 
\endproclaim
\bigskip

\demo{Proof}  1) By the transitivity of $\le^*_{\lambda^+}$ which
holds by \scite{600-ne.3}(2).
\nl
2) As in the proof of \scite{600-1.11}, in fact it follows from it. \nl
3) Easy (and is included in \scite{600-ne.3}(1A)).
\enddemo
\bigskip

\remark{Remark}  Note that the parallel to \scite{600-ne.3}(1A) is 
problematic in \S2 as, .e.g. locality may fail, i.e., $(M,N_i,a_i) \in
K^{3,\text{bs}}_{\lambda^+}$ and $M' \le_{\frak K} M \wedge M' \in
K_\lambda \Rightarrow \text{ \ortp}_{\frak s}(a_1,M',N_1) = \text{
\ortp}_{\frak s}(a_2,M',N_2)$ but \ortp$_{K^{\frak
s}_{\lambda^+}}(a_1,M,N_1) \ne 
\ortp_{K^{\frak s}_{\lambda^+}}
(\bar a_2,M,N_2)$.   \hfill$\square_{\scite{600-ne.2A}}$
\endremark
\bigskip

\proclaim{\stag{600-ne.4} Claim}  1) [Amalgamation of $\le^*_{\lambda^+}$
and toward extending types]
If $M_0 \le^*_{\lambda^+} M_\ell$ for
$\ell = 1,2,\kappa = { \text{\rm cf\/}}(\kappa) \le 
\lambda$ and $a \in M_2 \backslash
M_0$ is such that $(M_0,M_2,a) \in K^{3,\text{bs}}_{\lambda^+}$ 
is witnessed by $p$, 
\underbar{then} for some $M_3$ and $f$ we have: $M_1 <^+_{\lambda^+,
\kappa}M_3$ and $f$ is an $\le_{\frak K}$-embedding of $M_2$ into 
$M_3$ over $M_0$ with $f(a) \notin M_1$, moreover, $f(M_2) 
\le^*_{\lambda^+}M_3$ and $(M_1,M_3,f(a)) \in K^{3,\text{bs}}_{\lambda^+}$ is
witnessed by $p$. \newline
2) [uniqueness]  Assume $M_0 <^+_{\lambda^+,\kappa}M_\ell$ for $\ell = 1,2$
\underbar{then} there is an isomorphism $f$ from $M_1$ onto $M_2$ over
$M_0$. \nl
3) [locality] Moreover 
\footnote{the meaning of this will be that types over
$M \in K^{\text{nice}}_{\lambda^+}$ 
can be reduced to basic types over a model in
$K_\lambda$, i.e., locality}, in (2) if $a_\ell 
\in M_\ell \backslash M_0$ for $\ell = 1,2$
and $[N \le_{\frak K} M_0 \and N \in K_\lambda \Rightarrow 
{ \text{\rm \ortp\/}}(a_1,N,M_1) = { \text{\rm \ortp\/}}
(a_2,N,M_2)]$, \ub{then} we can demand $f(a_1) = a_2$
(so in particular ${\text{\rm \ortp\/}}(a_1,M_0,M_1) = 
{ \text{\rm \ortp\/}}(a_2,M_0,M_2)$ where the types are as defined in 
${\frak K}_{\lambda^+}$ and even in $(K_{\lambda^+},\le^*_{\lambda^+}$). \nl
4) Moreover in (2), assume further that for $\ell=1,2$, the following hold:
$N_0 \le_{\frak K} N_\ell \le_{\frak K} M_\ell,
N_0 \in K_\lambda,N_0 \le_{\frak K} M_0,
N_\ell \in K_\lambda$ and $(\forall N \in K_\lambda)
[N_0 \le_{\frak K} N \le_{\frak K} M_0 \rightarrow (\exists N' \in K_\lambda)
(N \cup N_\ell \subseteq N' \le_{\frak K} M_\ell \wedge 
{ \text{\rm NF\/}}_\lambda(N_0,N_\ell,N,N')]$.  If $f_0$ 
is an isomorphism from $N_1$ onto $N_2$ over $N_0$
\ub{then} we can add $f \supseteq f_0$.
\endproclaim
\bigskip

\demo{Proof}  We first prove part (2). \newline
2) By \scite{600-ne.2}(1) + (2) there are representations $\bar M^\ell = \langle
M^\ell_i:i < \lambda^+ \rangle$ of $M_\ell$ for $\ell < 3$ 
such that for $\ell =1,2$ we have: $M^\ell_i \cap M_0 = M^\ell_0$
and NF$_{\lambda,\langle 1,1,\kappa \rangle}
(M^0_i,M^0_{i+1},M^\ell_i,M^\ell_{i+1})$ and \wilog \, $M^\ell_0$ is
$(\lambda,\kappa)$-brimmed over $M^0_0$ for $\ell = 1,2$.
\medskip
\noindent
Now we choose by induction on $i < \lambda^+$ an isomorphism $f_i$ from 
$M^1_i$ onto $M^2_i$, increasing with $i$ and being the identity over
$M^0_i$.  For $i=0$ use ``$M^\ell_0$ is $(\lambda,\kappa)$-brimmed over
$M^0_0$ for $\ell = 1,2$" which we assume above.  
For $i$ limit take unions, for $i$ successor ordinal use 
uniqueness (Claim \scite{600-nf.7}).
\enddemo
\bigskip

\demo{Proof of part (1)}  By \scite{600-ne.3}(1) there 
are for $\ell = 1,2$ models
$N^*_\ell \in K_{\lambda^+}$ such that $M_\ell <^+_{\lambda^+,\kappa}
N^*_\ell$.  Now let $\bar M^\ell = \langle M^\ell_i:i < \lambda^+ \rangle$
be a representation of $M_\ell$ for $\ell = 0,1,2$ and let
$\bar N^\ell = \langle N^\ell_i:i < \lambda^+ \rangle$ be a representation of
$N^*_\ell$ for $\ell = 1,2$.  By \scite{600-ne.3}(4) and \scite{600-ne.2}(2) 
without loss of generality $N^\ell_0$ is $(\lambda,\kappa)$-brimmed
over $M^\ell_0$ and
NF$_\lambda(M^0_i,M^0_{i+1},M^\ell_i,M^\ell_{i+1})$ and
NF$_{\lambda,\langle 1,1,\kappa \rangle}
(M^\ell_i,M^\ell_{i+1},N^\ell_i,N^\ell_{i+1})$ respectively for $\ell
= 1,2$.  Let $M^*_0$ be such that $p \in {\Cal S}^{\text{bs}}(M^*_0),M^*_0 \in
K_\lambda,M^*_0 \le_{\frak K} M_0$; \wilog \, $M^*_0 \le_{\frak K}
M^0_0$ and $a \in N^2_0$.
Now $N^\ell_0$ is $(\lambda,\kappa)$-brimmed over 
$M^\ell_0$ hence over
$M^0_0$ (for $\ell = 1,2$) so there is an isomorphism $f_0$ from
$N^2_0$ onto $N^1_0$ extending id$_{M^0_0}$.  There is
$a' \in N^1_0$ such that \ortp$(a',M^1_0,N^1_0)$ is a non-forking
extension of $p$ and \wilog \, $f_0(a) 
= a'$ hence \ortp$(f_0(a),M^1_0,N^1_0) \in
{\Cal S}^{\text{bs}}(M^1_0)$ does not fork over $M^0_0$.
\medskip

We continue as in the proof of part (2).  In the end
$f = \dsize \bigcup_{i < \lambda^+} f_i$ is an isomorphism of $N^*_2$ onto
$N^*_1$ over $M_0$ and as $f_0(a)$ is well defined and in $N^1_0 \backslash
M^1_0$ clearly \ortp$(f(a),M^1_i,N^1_i)$ does not fork over
$M^1_0$ and extends $p$ hence the pair 
$(N^*_1,f \restriction M_2)$ is as required. 
\enddemo
\bigskip

\demo{Proof of part (3), (4)}  Like part (2).  \hfill$\square_{\scite{600-ne.4}}$
\enddemo
\bigskip

\proclaim{\stag{600-ne.5} Claim}  1) If $\delta$ is a limit ordinal 
$< \lambda^{+2}$ and $\langle M_i:i < \delta \rangle$ is a 
$\le^*_{\lambda^+}$-increasing continuous (in $K_{\lambda^+}$) and
$M_\delta = \dbcu_{i < \delta} M_i$ (so $M_\delta \in K_{\lambda^+}$), 
\ub{then} $M_i \le^*_{\lambda^+} M_\delta$ for each $i < \delta$.
\nl
2) If 
$\delta$ is a limit ordinal $< \lambda^{+2}$ and $\langle M_i:i < \delta
\rangle$ is a $\le^*_{\lambda^+}$-increasing sequence, each $M_i$ is in
$K^{\text{nice}}_{\lambda^+}$, \ub{then} $\dbcu_{i < \delta} M_i$ is 
in $K^{\text{nice}}_{\lambda^+}$. \nl
3) If $\delta$ is a limit ordinal $< \lambda^{+2}$ and $\langle M_i:i <
\delta \rangle$ is a $<^+_{\lambda^+}$-increasing continuous 
(or just $<^*_{\lambda^+}$-increasing continuous, and $M_{2i+1}
<^+_{\lambda^+} M_{2i+2}$ for $i < \delta$), \ub{then}
$i < \delta \Rightarrow M_i <^+_{\lambda^+} \dbcu_{j < \delta} M_j$.
\endproclaim
\bigskip

\demo{Proof}  1) We prove it by induction on $\delta$.  
Now if $C$ is a club of $\delta$, (as $\le^*_{\lambda^+}$ is transitive)
then we can replace $\langle M_j:j < \delta \rangle$ by $\langle M_j:j \in
C \rangle$ so without loss of generality $\delta = \text{ cf}(\delta)$, so
$\delta \le \lambda^+$; similarly it is enough to prove $M_0
\le^*_{\lambda^+} M_\delta := 
\dsize \bigcup_{j < \delta} M_j$.  For each $i \le \delta$ let 
$\langle M^i_\zeta:\zeta < \lambda^+ \rangle$ be a 
$<^*_{\frak K}$-representation of $M_i$.  
\enddemo
\bigskip
\noindent
\underbar{Case A}:  $\delta < \lambda^+$. \newline
Without loss of generality (see \scite{600-ne.2}(1)) 
for every $i < j < \delta$ and $\zeta < \lambda^+$ we have: \newline
$M^j_\zeta \cap M_i = M^i_\zeta$ and NF$_\lambda(M^i_\zeta,M^i_{\zeta +1},
M^j_\zeta,M^j_{\zeta +1})$.  Let
$M^\delta_\zeta = \dsize \bigcup_{i < \delta} M^i_\zeta$, so \newline
$\langle M^\delta_\zeta:\zeta < \lambda^+ \rangle$ 
is $\le_{\frak K}$-increasing continuous sequence of members of $K_\lambda$
with limit $M_\delta$, and for $i < \delta,M^\delta_\zeta \cap M_i =
M^i_\zeta$.  By symmetry (see \scite{600-nf.14}) we have NF$_\lambda(M^i_\zeta,
M^{i+1}_\zeta,M^i_{\zeta +1},M^{i+1}_{\zeta + 1})$ so as $\langle M^i_\zeta:
i \le \delta \rangle$,$\langle M^i_{\zeta +1}:i \le \delta \rangle$ are
$\le_{\frak K}$-increasing continuous, by \scite{600-nf.16} the transitivity
of NF$_{\frak s}$, we know NF$_\lambda(M^0_\zeta,M^\delta_\zeta,
M^0_{\zeta +1},M^\delta_{\zeta +1})$ 
hence by symmetry (\scite{600-nf.14}) we have
NF$_\lambda(M^0_\zeta,M^0_{\zeta +1},M^\delta_\zeta,M^\delta_{\zeta +1})$.
\newline
So $\langle M^0_\zeta:\zeta < \lambda^+ \rangle,\langle M^\delta_\zeta:
\zeta < \lambda^+ \rangle$ are witnesses to $M_0 \le^*_{\lambda^+} M_\delta$.
\bn
\underbar{Case B}:  $\delta = \lambda^+$.

By \scite{600-ne.2}(1) (using normality of the club filter, restricting 
to a club of $\lambda^+$ and renaming), without loss
of generality for $i < j \le 1 + \zeta < 1 + \xi < \lambda^+$ we have
$M^j_\zeta \cap M_i = M^i_\zeta$, and NF$_\lambda(M^i_\zeta,M^i_\xi,M^j
_\zeta,M^j_\xi)$.  Let us define $M^{\lambda^+}_\zeta = \dsize \bigcup
_{j < 1 + \zeta} M^j_\zeta$.  
So $\langle M^{\lambda^+}_\zeta:\zeta < \lambda^+ \rangle$ is
a $<_{\frak K}$-representation of $M_{\lambda^+} = M_\delta$ and 
continue as before. \nl
2) Again without loss of generality $\delta = \text{ cf}(\delta)$ call it
$\kappa$.  Let $\langle M^i_\zeta:\zeta < \lambda^+ \rangle$ be a
$<_{\frak K}$-representation of $M_i$ for $i < \delta$.
\bn
\ub{Case A}:  $\delta = \kappa < \lambda^+$.

Easy by now, yet we give details.
So \wilog \, (see \scite{600-ne.2}(1)) for every $i < j < \delta$ and $\zeta <
\xi < \lambda^+$ we have: $M^j_\zeta \cap M_i = M^i_\zeta$,
NF$_\lambda(M^i_\zeta,M^i_\xi,M^j_\zeta,M^j_\xi)$ and $M^i_{\zeta +1}$
is $(\lambda,\lambda)$-brimmed over $M^i_\zeta$.  Let $M^\delta_\zeta =
\dbcu_{\beta < \delta} M^\beta_\zeta$.  Let $\xi < \lambda^+$. 
Now if $p \in {\Cal S}^{\text{bs}}(M^\delta_\xi)$ 
then by  the  local character Axiom (E)(c) + the uniqueness Axiom
(E)(e), for some $i < \delta,p$ does not fork over 
$M^i_\xi$.
As $M_i$ is $\lambda^+$-saturated above $\lambda$, the type $p \restriction
M^i_\xi$ is realized in
$M_i$.  So let $b \in M_i$ realize $p \restriction M^i_\xi$ and 
by Axiom $(E)(h)$, continuity, it 
suffices to prove that for every $j \in (i,\delta),b$ realizes
$p \restriction M^j_\xi$ in $M_j$ which holds by \scite{600-nf.19}
(note that $b \in M_i \le_{\frak K} M_j$ as $j \in [i,\delta)$).
So $p$ is realized in $M_\delta = \dbcu_{i < \delta} M_i$.  As this holds
for every $\xi < \lambda^+$ and $p \in {\Cal S}^{\text{bs}}
(M^\delta_\xi)$, the model $M_\delta$ is saturated.
\bn
\ub{Case B}:  cf$(\delta) = \lambda^+$.

Straight, in fact true for ${\frak K}$ a.e.c. with the
$\lambda$-amalgamation property. \nl
3) Similar.    \hfill$\square_{\scite{600-ne.5}}$
\bigskip

\proclaim{\stag{600-ne.6} Claim}  1) If 
$M_0 \in K_{\lambda^+}$ \ub{then} there is $M_1$
such that $M_0 <^+_{\lambda^+} M_1 \in K^{\text{nice}}_{\lambda^+}$,
and any such $M_1$ is universal over $M_0$ in $(K_{\lambda^+},
\le^*_{\lambda^+})$. 
\nl
2) Assume $\boxtimes_{\bar N_1,\bar N_2,M_1,M_2}$ below holds.  \ub{Then}
$M_1 <^+_{\lambda^+} M_2$ \ub{iff} for every $\alpha < \lambda^+$ for
stationarily many $\beta < \lambda^+$ there is $N$ such that 
$N^1_\beta \cup N^2_\alpha \subseteq N \le_{\frak K} N^2_\beta$ and
$N^2_\beta$ is $(\lambda,*)$-brimmed over $N$ where
\mr
\item "{$\boxtimes_{\bar N_1,\bar N_2,M_1,M_2}$}"  $M_1
\le^*_{\lambda^+} M_2$ is being witnessed by $\bar N_1,\bar N_2$ that
is $\bar N_\ell = \langle N^\ell_\alpha:\alpha < \lambda^+ \rangle$ is
a $\le_{\frak K}$-representation of $M_\ell$ for $\ell =1,2$ and
$\alpha < \lambda^+ \Rightarrow { \text{\rm NF\/}}_\lambda
(N^1_\alpha,N^1_{\alpha +1},N^2_\alpha,N^2_{\alpha +1})$ 
(hence $\alpha \le \beta < \lambda^+
\Rightarrow { \text{\rm NF\/}}_\lambda
(N^1_\alpha,N^1_\beta,N^2_\alpha,N^2_\beta))$.
\endroster
\endproclaim
\bigskip

\demo{Proof}  1) The existence by \scite{600-ne.3}(1).  Why ``any such
$M_1,\ldots$?" if $M_0 \le^*_{\lambda^+}
M_2$ then for some $M^+_2 \in K^{\text{nice}}_{\lambda^+}$ we have $M_2
<^+_\lambda M^+_2 \in K^{\text{nice}}_{\lambda^+}$ so $M_0 \le^*
_{\lambda^+} M_1 <^+_{\lambda^+} M^+_2$ hence by \scite{600-ne.3}(5) we have
$M_0 <^+_\lambda M^+_2$; so by \scite{600-ne.4}(2) the models
$M^+_2,M_1$ are isomorphic over $M_0$, so $M_2$ can be 
$\le^*_{\lambda^+}$-embedded into $M_1$ over $M_0$, so we are
done. \nl
2) Not hard.   \hfill$\square_{\scite{600-ne.6}}$
\enddemo
\newpage

\head {\S8 Is $\le^*_{\lambda^+}$ equal to $\le_{\frak K}$ on
$K^{\text{nice}}_{\lambda^+}$?} \endhead  \resetall 
 \spuriousreset
\bigskip

\demo{\stag{600-rg.0} Hypothesis}  The hypothesis \scite{600-nf.0}.

An important issue is whether
$(K^{\text{nice}}_{\lambda^+},\le^*_{\lambda^+})$ satisfies Ax IV of
a.e.c.  So a model $M \in K_{\lambda^{++}}$ may be the union of a 
$\le^*_{\lambda^+}$-increasing chain of length $\lambda^{++}$, but we
still do not know if there is a continuous such sequence.

E.g. let $\langle M_\alpha:\alpha < \lambda^{++}\rangle$ be
$\le^*_{\lambda^+}$-increasing with union $M \in K_{\lambda^{++}}$ let
$M'_n = M_n,M'_{\omega + \alpha +1} = M_{\omega + \alpha}$ and
$M'_\delta = \cup \{M_\beta:\beta < \delta\}$ for $\delta$ limit.  So
$\langle M'_\alpha:\alpha < \lambda^{++}\rangle$ is $\le_{\frak
K}$-increasing continuous, $\langle M'_{\alpha +1}:\alpha <
\lambda^{++}\rangle$ is $\le^*_{\lambda^+}$-increasing, but we do not
know whether $M'_\delta \le^*_{\lambda^+} M'_{\delta +1}$ for limit
$\delta < \lambda^{++}$.
\enddemo
\bigskip

\definition{\stag{600-sg.1} Definition}  Let $M \in {\frak K}_{\lambda^{++}}$
be the union of an $\le_{\frak K}$-increasing chain from 
$(K^{\text{nice}}_{\lambda^+},\le^*_{\lambda^+})$ or just
$(K_{\lambda^+},\le^*_{\lambda^+}),
\bar M = \langle M_i:i < \lambda^{++} \rangle$ such that 
$\langle M_i:i < \lambda^{++}$ non-limit$\rangle$ is
$\le^*_{\lambda^+}$-increasing. \nl
1) Let $S(\bar M) = \{\delta:M_\delta \nleq^*_{\lambda^+}
M_{\delta+1}$ (see \scite{600-sg.2}(3) below)$\}$, so 
$S(\bar M) \subseteq \lambda^{++}$.
\nl
2) For such $M$ let $S(M)$ be $S(\bar M)/{\Cal D}_{\lambda^{++}}$
where $\bar M$ is a $\le_{\frak K}$-representation of $M$ and
${\Cal D}_{\lambda^{++}}$ is the club filter on $\lambda^{++}$; it
is well defined by \scite{600-sg.2} below. \nl
3) We say $\langle M_i:i < \delta \rangle$ is non-limit
$<^*_{\lambda^+}$-increasing \ub{if} for non-limit $i <j < \delta$ we have
$M_i \le^*_{\lambda^+} M_j$.
\enddefinition
\bigskip

\proclaim{\stag{600-sg.2} Claim}  1) If $\bar M^\ell = \langle M^\ell_i:i < 
\lambda^{++}\rangle$ for 
$\ell \in \{1,2\}$ is $\le_{\frak K}$-increasing continuous and
$i < j < \lambda^{++} \Rightarrow M_0 \le^*_{\lambda^+} M_{i+1}
\le^*_{\lambda^+} M_{j+1}$ and 
$M = \dbcu_{i < \lambda^{++}} M^1_i = \dbcu_{i < \lambda^{++}} M^2_i$
has cardinality $\lambda^{++}$ \ub{then} $S(\bar M^1) = S(\bar M^2)$ 
{\rm mod} ${\Cal D}_{\lambda^{++}}$. \nl
2) If $M,\bar M$ are as in \scite{600-sg.1} hence $M = \dbcu_{i < \lambda^{++}}
M_i$ \ub{then} 
$S(\bar M)/{\Cal D}_{\lambda^{++}}$ depends just on $M/\cong$.  \nl
3) If $\bar M$ is as in \scite{600-sg.1} and $i < j < \lambda^{++}$, 
\ub{then} $M_i \le^*_{\lambda^+} M_{i+1} \Leftrightarrow M_i
\le^*_{\lambda^+} M_j$. \nl
4) If $M \in {\frak K}_{\lambda^{++}}$ is the union of a
$\le^*_{\lambda^+}$-increasing chain from
$(K^{\text{nice}}_{\lambda^+},\le^*_{\lambda^+})$, not necessarily
continuous, \ub{then} there is $\bar M$ as in Defintion \scite{600-sg.1},
that is $\bar M = \langle M_i:i < \lambda^{++}
\rangle$, a $\le_{\frak K}$-representation of $M$ with $M_i
\le^*_{\lambda^+} M_j$ for non-limit $i<j$.
\endproclaim
\bigskip

\demo{Proof}  1) We can find a club $E$ of $\lambda^{++}$ consisting
of limit ordinals such that $i \in E \Rightarrow M^1_i = M^2_i$.  Now
if $\delta_1 < \delta_2$ are from $E$ then $\delta_1 \in S(\bar M^1)
\Leftrightarrow M^1_{\delta_1} \le^*_{\lambda^+} M^1_{\delta_1+1}
\Leftrightarrow M^1_{\delta_1} \le^*_{\lambda^+} M^1_{\delta_2}
\Leftrightarrow M^2_{\delta_1} \le^*_{\lambda^+} M^2_{\delta_2}
\Leftrightarrow M^2_{\delta_1} \le^*_{\lambda^+} M^2_{\delta_1+1}
\Leftrightarrow \delta_1 \in S(M^2)$.
\nl
[Why?  By the definition of $S(\bar M^1)$, by part (3), by
``$\delta_1,\delta_2 \in E$", by part (3), by the definition of
$S(\bar M^2)$, respectively.]  So we are done.
\nl
2) Follows by part (1). \nl
3) The implication $\Leftarrow$ is by \scite{600-ne.3}(3); for the implication
$\Rightarrow$, note that assuming $M_i <^*_{\lambda^+} M_{i+1}$, as
$\le^*_{\lambda^+}$ is a partial order, noting that by the assumption
on $\bar M$ we have $M_{i+1} \le^*_{\lambda^+} M_{j+1}$,
and by \scite{600-ne.3}(3) we are done. \nl
4) Trivial.   \hfill$\square_{\scite{600-sg.2}}$
\enddemo
\bigskip

\proclaim{\stag{600-sg.3} Claim}  If $(*)$ 
below holds \ub{then} for every stationary
$S \subseteq S^{\lambda^{++}}_{\lambda^+} (= \{\delta < \lambda^{++}:
{\text{\rm cf\/}}(\delta) = 
\lambda^+\})$ for some $\lambda^+$-saturated $M \in K_{\lambda^{++}}$
we have $S(M)$ is well defined and equal 
to $S/{\Cal D}_{\lambda^{++}}$, where
\mr
\item "{$(*)$}"  we can find $\langle M_i:i \le \lambda^+ +1 \rangle$
which is $<_{\frak K}$-increasing continuous of members of
$K^{\text{nice}}_{\lambda^+}$ such that 
$i < j \le \lambda^+ +1 \and (i,j) \ne
(\lambda^+,\lambda^+ + 1) \Rightarrow M_i <^+_{\lambda^+} M_j$ but
$\neg(M_{\lambda^+} \le^*_{\lambda^+} M_{\lambda^+ +1})$.
\endroster
\endproclaim
\bigskip

\demo{Proof}  Fix 
$S \subseteq S^{\lambda^{++}}_{\lambda^+}$ and $\langle M_i:i \le
\lambda^+ +1 \rangle$ as in $(*)$. \nl
Without loss of generality $|M_{\lambda^+ +1} \backslash M_{\lambda^+}|
= \lambda^+$.

We choose by induction on $\alpha < \lambda^{+2}$ a model $M^S_\alpha$
such that:
\mr
\item "{$(a)$}"  $M^S_\alpha \in K^{\text{nice}}_{\lambda^+}$ 
has universe an ordinal $< \lambda^{++}$
\sn
\item "{$(b)$}"  for $\beta < \alpha$ we have $M^S_\beta \le_{\frak K}
M^S_\alpha$
\sn
\item "{$(c)$}"  if $\alpha = \beta + 1$, $\beta \notin S$ then $M^S_\beta
<^+_{\lambda^+} M^S_\alpha$
\sn
\item "{$(d)$}"  if $\alpha = \beta + 1,\beta \in S$ then $(M^S_\beta,
M^S_\alpha) \cong (M_{\lambda^+},M_{\lambda^+ +1})$
\sn
\item "{$(e)$}"  if $\beta < \alpha,\beta \notin S$ then $M^S_\beta
\le^+_{\lambda^+} M^S_\alpha$
\sn
\item "{$(f)$}"  if $\alpha$ is a limit ordinal, then $M_\alpha =
\cup\{M_\beta:\beta < \alpha\}$.
\ermn
We use freely the transitivity and continuity of $\le^*_\lambda$ 
and of $<^+_\lambda$. \nl
\sn
For $\alpha = 0$ no problem.
\mn
For $\alpha$ limit no problem; choose an increasing continuous
sequence $\langle \gamma_i:i < \text{\rm cf}(\alpha)\rangle$ of
ordinals with limit $\alpha$ each of cofinality $< \lambda,\gamma_i
\notin S$, and use \scite{600-ne.5}(3) for clause (e).
\mn
For $\alpha = \beta +1,\beta \notin S$ no problem.
\mn
For $\alpha = \beta + 1,\beta \in S$ so cf$(\beta) = \lambda^+$, let
$\langle \gamma_i:i < \lambda^+ \rangle$ be increasing continuous with limit
$\beta$ and cf$(\gamma_i) \le \lambda$, hence $\gamma_i \notin S$ and
each $\gamma_{i+1}$ a successor ordinal.  By
clause (e) above and \scite{600-ne.3}(5) we have 
$M^S_{\gamma_i} <^+_{\lambda^+}
M^S_{\gamma_{i+1}}$, hence $\langle M_{\gamma_i}:i < \lambda^+ \rangle$
is $<^+_{\lambda^+}$-increasing continuous.  Now there is an
isomorphism $f_\beta$ from $M_{\lambda^+}$ onto $M^S_\beta$ 
mapping $M_i$ onto $M^S_{\gamma_i}$ for $i < \lambda$ (why? choose
$f_\beta \restriction M_i$ by induction on $i$, for $i=0$ by
\scite{600-ne.2}(0), for $i$ successor $M^S_{\gamma_i} <^+_\lambda
M^S_{\gamma_{i+1}}$ by \scite{600-ne.3}(3) as $M^S_{\gamma_i}
<^*_{\lambda^+} M^S_{\gamma_{i+1}} <^+_{\lambda^+} M^S_{\gamma_{i+1}}$
so we can use \scite{600-ne.4}(2)).  So we can choose a 
one-to-one function $f_\alpha$ from $M_{\lambda^+ +1}$ onto some
ordinal $< \lambda^{++}$ extending 
$f_\beta$ and let $M_\alpha = f_\alpha(M_{\lambda^+ +1})$.  

Finally having carried the induction, let 
$M_S = \dbcu_{\alpha < \lambda^{+2}} M^S_\alpha$,
it is easy to check that $M_S \in K_{\lambda^{++}}$ is $\lambda^+$-saturated
 and $\bar M = \langle M^S_\alpha:\alpha < \lambda^{++} \rangle$ 
witnesses that
$S(M_S)/{\Cal D}_{\lambda^{++}}$ is well defined and
$S(M_S)/{\Cal D}_{\lambda^{++}} = S(\langle M^S_\alpha:
\alpha < \lambda^{++} \rangle)/{\Cal D}_{\lambda^{++}} =
S/{\Cal D}_{\lambda^{++}}$ as required.  \hfill$\square_{\scite{600-sg.3}}$
\enddemo
\bn
Below we prove that some versions of non-smoothness are equivalent.
\proclaim{\stag{600-sg.10} Claim}  1) We have $(**)_{M^*_1,M^*_2} \Rightarrow
(***)$ (see below). \nl
2) If $(*)$ then $(**)_{M^*_1,M^*_2}$ for some $M^*_1,M^*_2$ and
trivially $(***) \Rightarrow (*)$.
\nl
3) In part (1) we get $\langle M_i:i \le \lambda^+ +1 \rangle$ as in
$(***)$, see below, 
such that $M_{\lambda^+} = M^*_1,M_{\lambda^+ +1} = M^*_2$ if we waive
$i < \lambda^+ \Rightarrow M_i <^+_\lambda M_{\lambda+1}$ or assume
$M^*_1 <_{\frak K} M^* <^+_\lambda M^*_2$ for some $M^*$. 
\nl
4) If $M^*_1 \le^*_{\lambda^+} M^*_2$ and $M^*_1 
\in K^{\text{nice}}_{\lambda^+}$
and $N_1 \le_{\frak K} N_2 \in K_\lambda,N_\ell \le M^*_\ell$ for
$\ell =1,2$ and $p \in 
{\Cal S}^{\text{bs}}(N_2)$ does not fork over $N_1$ \ub{then} some $c \in
M^*_1$ realizes $p$ \nl
\ub{where}
\smallskip

$(*) \quad$   there are limit $\delta < \lambda^{++},N$ and $\bar M =
\langle M_i:i \le \delta \rangle$ a $\le^*_{\lambda^+}$-increasing \nl

\hskip33pt  continuous sequence with
$M_i,N \in K^{\text{nice}}_{\lambda^+}$ such that: $M_i \le^*_{\lambda^+} N 
\Leftrightarrow i < \delta$
\sn

$(**)_{M^*_1,M^*_2}  \quad (i) \quad 
M^*_1 \in K^{\text{nice}}_{\lambda^+},
M^*_2 \in K^{\text{nice}}_{\lambda^+}$
\sn

\hskip58pt $(ii) \quad  M^*_1 \le_{\frak K} M^*_2$
\sn

\hskip58pt  $(iii) \quad  M^*_1 \nleq^*_{\lambda^+} M^*_2$
\sn

\hskip58pt $(iv) \quad$  if 
$N_1 \le_{\frak K} N_2$ are from $K_\lambda$, \nl

\hskip80pt  $N_\ell \le_{\frak K} M^*_\ell$ for $\ell =1,2$ and 
$p \in {\Cal S}^{\text{bs}}(N_2)$ does not \nl

\hskip80pt  fork over $N_1$, 
\ub{then} some $a \in M^*_1$ realizes $p$ in $M^*_2$
\sn

$(***) \quad$  there is $\bar M = \langle M_i:i \le \lambda^+ +1
\rangle,\le_{\frak K}$-increasing continuous, every \nl

\hskip40pt  $M_i \in K^{\text{nice}}_{\lambda^+}$
and $M_{\lambda^+} \nleq^*_{\lambda^+} M_{\lambda^+ +1}$ but \nl

\hskip40pt  $i < j \le \lambda^+ +1 \and i \ne 
\lambda^+ \Rightarrow M_i <^+_{\lambda^+} M_j$.
\endproclaim
\bigskip

\demo{Proof}  1),3) Let $\langle a^\ell_i:i < \lambda^+ \rangle$ list the
elements of $M^*_\ell$ for $\ell = 1,2$. 
Let $\langle N^*_{2,i}:i < \lambda^+ \rangle$ be a 
$\le_{\frak K}$-representation of $M^*_2$. 
\nl
Let $\langle (p_\zeta,N^*_\zeta,\gamma_\zeta):\zeta < \lambda^+ \rangle$
list the triples $(p,N,\gamma)$ such that $\gamma < \lambda^+,p \in
{\Cal S}^{\text{bs}}(N),N \in \{N^*_{2,i}:i < \lambda^+\}$ with each 
such triple appearing $\lambda^+$ times.  By induction on $\alpha 
< \lambda^+$ we choose $\langle N^\alpha_i:i \le \alpha \rangle,
N_\alpha$ such that:
\mr
\item "{$(a)$}"  $N^\alpha_i \in K_\lambda$ and $N^\alpha_i \le_{\frak K}
M^*_1$ \nl
\sn
\item "{$(b)$}"  $N_\alpha \le_{\frak K} M^*_2$
\sn
\item "{$(c)$}"  $\langle N^\alpha_i:i \le \alpha \rangle$ is
$\le_{\frak K}$-increasing continuous
\sn
\item "{$(d)$}"  $N^\alpha_\alpha \le_{\frak K} N_\alpha,N_\alpha \cap
M^*_1 = N^\alpha_\alpha$
\sn
\item "{$(e)$}"  if $i \le \alpha$ then $\langle N^\beta_i:\beta \in [i,\alpha]
\rangle$ is $\le_{\frak K}$-increasing continuous
\sn
\item "{$(f)$}"  $\langle N_\beta:\beta \le \alpha \rangle$ is
$\le_{\frak K}$-increasing continuous
\sn
\item "{$(g)$}"  if 
$\alpha = \beta +1,i \le \beta$ then NF$_\lambda(N^\beta_i,
N_\beta,N^\alpha_i,N_\alpha)$
\sn
\item "{$(h)$}"  if $\alpha = 2\beta +1$ then $a^2_\beta \in N_{\alpha +1}$
\sn
\item "{$(i)$}"  if $\alpha = 2\beta +2$ and $i < \alpha$ then
$N^\alpha_{i+1}$ is brimmed over $N^\alpha_i \cup N^{2 \beta
+1}_{i+1}$ and $N^\alpha_0$ is brimmed over $N^{2\beta}_0$.
\endroster
\bn
\ub{Why is this enough}?

We let $M_{\lambda^+} = M^*_1,M_{\lambda^+ +1} = M^*_2$ and let
$M'_{\lambda^+ +1} \in K^{\text{nice}}_{\lambda^+}$ be such that
$M_{\lambda^+ +1} <^+_{\lambda^+} M'_{\lambda^+ +1}$ and for
$i < \lambda^+$ we let 
$M_i = \cup \{N^\alpha_i:\alpha \in [i,\lambda^+)\}$; now
\mr
\item "{$(\alpha)$}"  $M^*_1 = \dbcu_{\alpha < \lambda^+} N^\alpha_\alpha
= \dbcu_{i < \lambda^+} M_i$ 
and $M^*_2 = \dbcu_{\alpha < \lambda^+} N_\alpha$ \nl
[why?  the second by clause (h) (and (b) of course), 
the first as $N_\alpha \cap M^*_1 = N^\alpha_\alpha$].
\ermn
Now:
\mr
\item "{$(\beta)$}"  $\langle M_i:i \le \lambda^+ +1 \rangle$ is
$\le_{\frak K}$-increasing continuous \nl
[trivial by clauses (c) + (e) if $i < \lambda^+$ and (d) if $i
= \lambda^+$]
\sn
\item "{$(\gamma)$}"  for $i < \lambda^+,M_i$ is saturated, i.e.,
$\in K^{\text{nice}}_{\lambda^+}$. \nl
[Why?  Clearly $\langle N^\alpha_i:\alpha \in (i,\lambda^+)\rangle$ is
a $\le_{\frak K}$-representation of $M_i$ by clause (e) and the choice
of $M_i$.  If $i=0$ this follows by clauses (i) + (e).  
If $i=j+1$ this follows by clauses (e) + (i).  If $i$ is a
limit ordinal use \scite{600-ne.5}(2) and clause (g)]
\sn
\item "{$(\delta)$}"  for $i < \lambda^+,i < j \le \lambda^+ +1$ we
have $M_i \le^*_{\lambda^+} M_j$. \nl
[Why?  Let $N^\alpha_{\lambda^+} :=
N^\alpha_\alpha,N^\alpha_{\lambda^+ +1} = N_\alpha$ for $\alpha <
\lambda^+$ and let $\gamma$ be $i$ if $j = \lambda^+,\lambda^+ +1$ and
be $j$ if $j < \lambda^+$; so in any case $\gamma < \lambda^+$.  Now
as $\langle N^\alpha_i:\alpha \in [\gamma,\lambda^+)\rangle$ is a
$\le_{\frak K}$-representation of $M_i$ and $\langle N^\alpha_j:\alpha
\in [\gamma,\lambda^+) \rangle$ is a $\le_{\frak K}$-representation of
$M_j$ and if $\gamma \le \beta < \lambda^+$ then by clause (g) we have
NF$_\lambda(N^\beta_i,N_\beta,N^{\beta +1}_i,N_{\beta +1})$ hence by
symmetry NF$_\lambda(N^\beta_i,N^{\beta +1}_i,N_\beta,N_{\beta +1})$
hence by monotonicity \nl
NF$_\lambda(N^\beta_i,N^{\beta
+1}_i,N^\beta_j,N^{\beta +1}_j)$; this suffices]
\sn
\item "{$(\varepsilon)$}"  if $i < j \le \lambda^+$ then
$M_i <^+_{\lambda^+} M_j$ \nl
[why?  by \scite{600-ne.5}(1) it suffices to prove this in the cases
$j=i+1$.  Now claim \scite{600-ne.6}(2), clause (i) guaranteed this.]
\ermn
Clearly $\langle M_i:i \le \lambda^+ +1 \rangle$ is as required for
part (3) and $\langle M_i:i \le \lambda^+ \rangle \char 94 \langle
M'_{\lambda^+ +1}\rangle$ is as required in part (1).
\nl
So we are done.
\bn
So let us carry the construction.
\mn
For $\alpha = 0$ trivially.
\mn
For $\alpha$ limit: straightforward.
\mn
For $\alpha = 2 \beta + 1$ we let $N^\alpha_i = N^{2 \beta}_i$ for $i \le
2 \beta$ and $N_\alpha \in K_\lambda$ is chosen such that $N_{2 \beta} \cup
\{a^2_\beta\} \subseteq N_\alpha \le_{\frak K} M^*_2$ and $N_\alpha
\restriction M^*_1 \le_{\frak K} M^*_1$, easy by the properties of abstract
elementary class and we let $N^\alpha_{2 \beta +1} = N_\alpha \restriction
M^*_1$.
\mn
For $\alpha = 2 \beta + 2$ we choose by induction on
$\varepsilon < \lambda^2$, a triple $(N^\oplus_{\alpha,\varepsilon},
N^\otimes_{\alpha,\varepsilon},a_{\alpha,\varepsilon})$ such that:
\mr
\item "{$(A)$}"  $N^\otimes_{\alpha,\varepsilon} \le_{\frak K} M^*_2$
belongs to $K_\lambda$ and is $\le_{\frak K}$-increasing continuous with
$\varepsilon$
\sn
\item "{$(B)$}"  $N^\otimes_{\alpha,0} = N_{2 \beta +1}$ and
$N^\otimes_{\alpha,\varepsilon} \restriction M^*_1 \le^*_{\frak K} M^*_1$
\sn
\item "{$(C)$}"  $N^\oplus_{\alpha,\varepsilon} \le_{\frak K} M^*_1$ belongs
to $K_\lambda$ and is $\le_{\frak K}$-increasing continuous with 
$\varepsilon$
\sn
\item "{$(D)$}"  $N^\oplus_{\alpha,0} = N^{2 \beta +1}_{2 \beta +1}$
\sn
\item "{$(E)$}"  $(N^\oplus_{\alpha,\varepsilon},
N^\oplus_{\alpha,\varepsilon + 1},
a_{\alpha,\varepsilon}) \in K^{3,\text{uq}}_\lambda$
\sn
\item "{$(F)$}"  \ortp$(a_{\alpha,\varepsilon},N^\otimes_{\alpha,\varepsilon},
M^*_2)$ does not fork over $N^\oplus_{\alpha,\varepsilon}$
\sn
\item "{$(G)$}"  $N^\oplus_{\alpha,\varepsilon} \le_{\frak K}
N^\otimes_{\alpha,\varepsilon}$ 
\sn
\item "{$(H)$}"  for every $p \in {\Cal S}^{\text{bs}}
(N^\oplus_{\alpha,\varepsilon})$ for some odd $\zeta \in [\varepsilon,
\varepsilon + \lambda)$ the type \ortp$(a_{\alpha,\zeta},
N^\otimes_{\alpha,\zeta},N^\otimes_{\alpha,\zeta +1})$ is a
non-forking extension of $p$.
\ermn
No problem to carry this.  [Why?  For $\varepsilon =0$ and
$\varepsilon$ limit there are no problems.  In stage $\varepsilon +1$
by bookkeeping gives you a type $p_\varepsilon \in {\Cal
S}^{\text{bs}}(N^\oplus_{\alpha,\varepsilon})$ and let $q_\varepsilon
\in {\Cal S}^{\text{bs}}(N^\otimes_{\alpha,\varepsilon})$ be a
non-forking extension of $p_\varepsilon$.  By assumption (iv) of
$(**)_{M^*_1,M^*_2}$ there is an element $a_{\alpha,\varepsilon} \in
M^*_1$ realizing $q_\varepsilon$.  Now $M^*_1$ is saturated hence
there is a model $N^\oplus_{\alpha,\varepsilon +1} \in K_\lambda$ such
that $N^\oplus_{\alpha,\varepsilon +1} \le_{\frak K} M^*_1$ and
$(N^\oplus_{\alpha,\varepsilon},N^\oplus_{\alpha,\varepsilon
+1},a_{\alpha,\varepsilon}) \in K^{3,\text{uq}}_\lambda$.

Lastly, choose $N^\otimes_{\alpha,\varepsilon +1}$ satisfying
clauses (A),(B),(G) so we have carried the induction on
$\varepsilon$.]  Note that NF$_\lambda(N^\oplus_{\alpha,\varepsilon},
N^\otimes_{\alpha,\varepsilon},N^\oplus_{\alpha,\varepsilon
+1},N^\otimes_{\alpha,\varepsilon +1})$ for each $\varepsilon <
\lambda^2$ by clauses (E),(F) of \scite{600-nf.20}, 
hence NF$(N^{2 \beta +1}_{2\beta +1},N_{2 \beta +1},
\cup\{N^\oplus_{\alpha,\varepsilon}:\varepsilon < \lambda^2\},
\cup\{N^\otimes_{\alpha,\varepsilon}:\varepsilon <
\lambda^2\})$ by \scite{600-nf.16} as
$(N^\oplus_{\alpha,0},N^\otimes_{\alpha,0}) = (N^{2 \beta +1}_{2
\beta+1},N_{2 \beta+1})$ and the sequence $\langle
N^\oplus_{\alpha,\varepsilon}:\varepsilon < \lambda^+\rangle,\langle
N^\otimes_{\alpha,\varepsilon}:\varepsilon < \lambda^+\rangle$ are
increasing continuous. 

Now let $N_\alpha = \bigcup \{N^\otimes_{\alpha,\varepsilon}:
\varepsilon < \lambda^2\},N^\alpha_\alpha = N_\alpha \cap M^*_1$
recalling clauses (A)+(B).

Now $\cup\{N^\oplus_{\alpha,\varepsilon}:\varepsilon < \lambda^2\}
\le_{\frak K} M^*_1$ is $(\lambda,*)$-brimmed over $N^{2 \beta +1}
_{2 \beta +1}$ by \scite{600-4a.2} (and clause (H) above).  Hence 
there is no problem to choose
$N^\alpha_i \le_{\frak K} N^\alpha_\alpha$ 
for $i \le 2 \beta +1$ as required, that is $N^{2 \beta+1}_i
\le_{\frak K} N^\alpha_i,\langle N^\alpha_i:i \le 2 \beta +1 \rangle$
is $\le_{\frak K}$-increasing continuous, NF$_\lambda(N^{2 \beta
+1}_i,N^{2 \beta +1}_{i+1},N^\alpha_i,N^\alpha_{i+1})$ and $N^\alpha_{i+1}$ is
$(\lambda,*)$-brimmed over $N^{2 \beta +1}_{i+1} \cup N^\alpha_i$ and
$N^\alpha_0$ is $(\lambda,*)$-brimmed over $N^{2 \beta +1}_0$.
So we have finished the induction step on $\alpha = 2 \beta +2$.

Having carried the induction we are done.
\nl
2) So assume $(*)$ and let $M_{\delta +1} := N$ from $(*)$.
It is enough to prove that $(**)_{M_\delta,M_{\delta +1}}$ holds.  Clearly 
clauses (i), (ii), (iii) hold, so we should prove (iv).  
Without loss of generality $\delta =
\text{ cf}(\delta)$ so $\delta = \lambda^+$ or $\delta \le \lambda$. 
For $i \le \delta + 1$ let $\langle M_{i,\alpha}:\alpha <
\lambda^+\rangle$ be a $\le_{\frak K}$-representation 
of $M_i$ and for $i < \delta,j \in 
(i,\delta +1]$ let $E_{i,j}$ be a club of $\lambda^+$ witnessing 
$M_i \le^*_{\lambda^+} M_j$ for $\bar M^i,\bar M^j$.
First assume $\delta \le \lambda$.  Let $E = \cap\{E_{i,j}:i < \delta,j \in
(i,\delta +1]\}$, it is a club of $\lambda^+$.  So assume 
$N_2 \le_{\frak K} M_{\delta +1},N_1 \le_{\frak K} N_2,N_1
\le_{\frak K} M_\delta$ and $N_1,N_2 \in K_\lambda$ and 
$p \in {\Cal S}^{\text{bs}}(N_2)$ does not fork over
$N_1$.  We can choose $\zeta \in E$ such that $N_2 \subseteq
M_{\delta +1,\zeta}$, let $p_1 \in {\Cal S}^{\text{bs}}(M_{\delta +1,\zeta})$
be a non-forking extension of $p$, so $p_1$ does not fork over $N_1$
hence (by monotonicity) 
over $M_{\delta,\zeta}$ so $p_2 := p_1 \restriction M_{\delta,\zeta} \in
{\Cal S}^{\text{bs}}(M_{\delta,\zeta})$.  By Axiom $(E)(c)$
for some $\alpha < \delta,p_2$ does not fork over $M_{\alpha,\zeta}$ hence
$p_2 \restriction M_{\alpha,\zeta} \in {\Cal S}^{\text{bs}}
(M_{\alpha,\zeta})$.  As $M_\alpha \in K^{\text{nice}}_{\lambda^+}$,
i.e., $M_\alpha$ is $\lambda^+$-saturated (above $\lambda$), clearly 
for some $\xi \in (\zeta,\lambda^+) \cap E$ 
some $c \in M_{\alpha,\xi}$ realizes $p_2 \restriction M_{\alpha,\zeta}$
but NF$_\lambda(M_{\alpha,\zeta},M_{\delta +1,\zeta},M_{\alpha,\xi},
M_{\delta +1,\xi})$ hence by \scite{600-nf.19} we know that
\ortp$(c,M_{\delta +1,\zeta},
M_{\delta +1,\xi})$ belongs to ${\Cal S}^{\text{bs}}(M_{\delta +1,\zeta})$ and
does not fork over $M_{\alpha,\zeta}$ hence $c$ realizes 
$p_2$ and even $p_1$ hence $p$ and we are done. \nl

Second, assume $\delta = \lambda^+$, then for some $\delta^* < \delta$
we have $N_1 \le_{\frak K} M_{\delta^*}$, and use the proof above for
$\langle M_i:i \le \delta^* \rangle,M_{\delta +1}$ (or use
$M_{\delta^*} \le^*_{\lambda^+} M_{\delta +1}$). \nl
4) Straight, in fact included the proof of \scite{600-ne.5}(2).   
\hfill$\square_{\scite{600-sg.10}}$
\enddemo
\bn
The definition below has affinity to ``blowing ${\frak K}_\lambda$ to
${\frak K}^{\text{up}}_\lambda$" in \S1.
\definition{\stag{600-rg.6} Definition}  0)  $K^{3,\text{cs}}_{\lambda^+} =
\{(M,N,a) \in K^{3,\text{bs}}_{\lambda^+}:M,N$ are from
$K^{\text{nice}}_{\lambda^+}\}$; we say $N' \in K_\lambda$ (or $p'$) 
witness $(M,N,a) \in K^{3,\text{cs}}_{\lambda^+}$ if it witnesses
$(M,N,a) \in K^{3,\text{bs}}_\lambda$.
\nl
1) ${\Cal S}^{\text{cs}}_{\lambda^+} := 
\{\text{\ortp}(a,M,N):M \le^*_{\lambda^+} N$ are in 
$K^{\text{nice}}_{\lambda^+},a \in N$ and $(M,N,a) \in 
K^{3,\text{cs}}_{\lambda^+}\}$, the type being for ${\frak
K}^{\text{nice}}_{\lambda^+} = (K^{\text{nice}}_{\lambda^+},
\le^*_{\lambda^+}$), see below
\footnote{actually to define \ortp$_{{\frak K}_\lambda}(a,M,N)$ where $M
\le_{{\frak K}_\lambda} N,\bar a \in N$ we need less that ``${\frak
K}_\lambda$ is a $\lambda$-a.e.c.", and we know on
$(K^{\text{nice}}_{\lambda^+},\le^*_{\lambda^+})$ more than enough}
so the notation is justified by \scite{600-rg.7}(2). \nl
2) We define ${\frak K}^\otimes = (K^\otimes,\le^\otimes)$ as follows
\mr
\item "{$(a)$}"  $K^\otimes = {\frak K} \restriction
\{M \in K:M$ is the union of a directed family of $\le_{\frak K}$-submodels
each from $K^{\text{nice}}_{\lambda^+}\}$
\sn
\item "{$(b)$}"  Let 
$M_1 \le^\otimes M_2$ if $M_1,M_2 \in K^\otimes,M_1 \le_{\frak K} M_2$ and:
{\roster
\itemitem{ $(*)_{M_1,M_2}$ }  if $N_\ell \in K_\lambda,
N_\ell \le_{\frak K} M_\ell$,
for $\ell = 1,2,p \in {\Cal S}^{\text{bs}}(N_2)$ does not fork over $N_1$ and
$N_1 \le_{\frak K} N_2$ \ub{then} some $a \in M_1$ realizes $p$ in $M_2$
\endroster}
\item "{$(c)$}"    let $\le^\otimes_{\lambda^+} =
\le^\otimes \restriction K^\otimes_{\lambda^+}$. 
\ermn
3) $\nonfork{}{}_{\lambda^+} = \{(M_0,M_1,a,M_3):M_0 \le^*_{\lambda^+} M_1 
\le^*_{\lambda^+} M_3$ are in $K^{\text{nice}}_{\lambda^+}$ and $(M_1,M_3,a)
\in K^{3,\text{cs}}_{\lambda^+}$ as witnessed by some $N \le_{\frak K} M_0$
from $K_\lambda\}$. 
\nl
4) ${\frak K}^{\text{nice}}_{\lambda^+} = (K^{\text{nice}}_{\lambda^+},
\le^*_{\lambda^+})$, that is
$(K^{\text{nice}}_{\lambda^+},\le^*_{\lambda^+} \restriction 
K^{\text{nice}}_{\lambda^+})$.
\nl
5) We say that $M'$ or $p'$ witness $p = 
\text{ \ortp}_{{\frak K}^{\text{nice}}_{\lambda^+}}(a,M,N)$ when $M'
\le_{\frak K} M,M' \in K_\lambda$ and $[M' \le_{{\frak K}_\lambda} M'' 
\le_{\frak K} M \Rightarrow \text{\rm \ortp}_{\frak s}(a,M'',N)$ does not
fork over $M'$ and $p' = \text{\rm \ortp}_{\frak s}(a,M',N)$.
\enddefinition
\bigskip

\demo{\stag{600-rg.7} Conclusion}   Assume \footnote{this is like
$(**)_{M_1,M_2}$ from \scite{600-sg.10}, particularly see clause (iv)
there} (recalling \scite{600-sg.3}):
\mr
\item "{$\boxtimes$}"  not for every 
$S \subseteq S^{\lambda^{++}}_{\lambda^+}$ is there
$\lambda^+$-saturated $M \in K_{\lambda^{++}}$
such that $S(M) = S/{\Cal D}_{\lambda^{++}}$.
\ermn
0) On $K^{\text{nice}}_{\lambda^+}$, the relations $\le^*_{\lambda^+},
\le^\otimes$ agree. \nl
1) ${\frak K}^{\text{nice}}_{\lambda^+} = 
(K^{\text{nice}}_{\lambda^+},\le^*_{\lambda^+})$ is a 
$\lambda^+$-abstract elementary class and is categorical in $\lambda^+$ and
has no maximal member and has amalgamation. 
\nl
2) $K^\otimes$ is the class of $\lambda^+$-saturated
models in ${\frak K}$ so $K^\otimes_{\lambda^+} = K^{\text{nice}}
_{\lambda^+}$. \nl
3) ${\frak K}^\otimes$ is an a.e.c. with LS$(K^\otimes) = \lambda^+$
and is the lifting of ${\frak K}^{\text{nice}}_{\lambda^+}$.
\nl
4) On $K^{\text{nice}}_{\lambda^+},({\Cal S}^{\text{cs}}_{\lambda^+},
\nonfork{}{}_{\lambda^+})$ are equal to $({\Cal S}^{\text{bs}} 
\restriction K^{\text{nice}}_{\lambda^+},
\nonfork{}{}_{< \infty} \restriction K^{\text{nice}}_{\lambda^+})$
where they are defined in \scite{600-1.6}, \scite{600-1.7}. \nl 
5) $({\frak K}^{\text{nice}}_{\lambda^+},{\Cal S}^{\text{cs}}_{\lambda^+},
\nonfork{}{}_{\lambda^+})$ is a good $\lambda^+$-frame.   \nl
6) For $M_1 \le^*_{\lambda^+} M_2$ from $K^\otimes_{\lambda^+}$ 
and $a \in M_2 \backslash M_1$, the
type \ortp$_{K^\otimes}(a,M_1,M_2)$ is determined by
\ortp$_{{\frak K}_\lambda}(a,N_1,M_2)$ for all 
$N_1 \le_{\frak K} M_1,N_1 \in K_\lambda$.
\enddemo
\bigskip

\demo{Proof}  0) By \scite{600-sg.3} and our assumption $\boxtimes$, we have
$M_1,M_2 \in K^{\text{nice}}_{\lambda^+} \and M_1 \le^\otimes M_2
\Rightarrow M_1 \le^*_{\lambda^+} M_2$ (otherwise $(**)_{M_1,M_2}$ of \scite{600-sg.10} holds
hence $(***)$ of \scite{600-sg.10} holds and by \scite{600-sg.3} we get $\neg \boxtimes$,
contradiction).  The other direction is easier just see \scite{600-sg.10}(4).
\nl
1) We check the axioms for being a $\lambda^+$-a.e.c.:
\sn
\ub{Ax 0}: (Preservation under isomorphisms)  Obviously.
\sn
\ub{Ax I}:  Trivially.
\sn
\ub{Ax II}:  By \scite{600-ne.3}(2).
\sn
\ub{Ax III}:  By \scite{600-ne.5}(2) the union belongs to
$K^{\text{nice}}_{\lambda^+}$ and it $\le^*_{\lambda^+}$-extends each
member of the union by \scite{600-ne.5}(1).
\sn
\ub{Ax IV}:  Otherwise $(*)$ of \scite{600-sg.10} holds, hence by \scite{600-sg.10}
also $(***)$ of \scite{600-sg.10} holds.  So 
by \scite{600-sg.3} our assumption $\boxtimes$ fail, contradiction; this
is the only place we use $\boxtimes$ in the proof of (1).
\sn
\ub{Ax V}:  By \scite{600-ne.3}(3) and Ax V for ${\frak K}$.

Also ${\frak K}^{\text{nice}}_{\lambda^+}$ is categorical by the
 uniqueness of the saturated model in $\lambda^+$ for ${\frak K}$
 has no maximal model by
\scite{600-ne.3}(1).  ${\frak K}^{\text{nice}}_{\lambda^+}$ has
amalgamation by \scite{600-ne.4}(1).
\nl
2) Every member of $K^\otimes$ is $\lambda^+$-saturated in ${\frak K}$ by
\scite{600-ne.5}(2) (prove by induction on the cardinality of the directed
family in Definition \scite{600-rg.6}(2), i.e. by the LS-argument it is
 enough to deal with the index family of $\le \lambda^+$ models each
 of cardinality $\lambda^+$, which holds by part (0) + (1)).  
If $M \in K$ is $\lambda^+$-saturated, then as in \S1, it is easy
to find a directed system of saturated $\le_{\frak K}$-submodels from
$K_{\lambda^+}$ using \scite{600-ne.5}(2) (using the stability of ${\frak K}$
in $\lambda$).  
\nl
3),4) Easy by now (or see \S1). \nl
5) We have to check all the clauses in Definition \scite{600-1.1}.  We shall
use parts (0)-(3) freely.
\bn
\ub{Axiom (A)}:  

By part (3) (of \scite{600-rg.7}).
\bn
\ub{Axiom (B)}:  

There is a superlimit model in $K^\otimes_{\lambda^+} =
K^{\text{nice}}_{\lambda^+}$  by part (1) and uniqueness of the 
saturated model.
\bn
\ub{Axiom (C)}:  

By part (1), i.e., \scite{600-ne.4}(1) we have amalgamation; 
JEP holds as $K^{\text{nice}}
_{\lambda^+}$ is categorical in $\lambda^+$.
``No maximal member in ${\frak K}^\otimes_{\lambda^+}$" holds by
\scite{600-ne.3}(1).
\bn
\ub{Axiom (D)(a),(b)}:

By the definition \scite{600-rg.6}(1).
\bn
\ub{Axiom (D)(c)}:

By \scite{600-1.11} (and Definition \scite{600-rg.6}(1)).  Clearly
$K^{3,\text{cs}}_{\lambda^+} = K^{3,\text{bs}} \restriction
K^{\text{nice}}_{\lambda^+}$.
\bn
\ub{Axiom (D)(d)}:

For $M \in {\frak K}^\otimes_{\lambda^+}$ let $\bar M = \langle M_i:
i < \lambda^+ \rangle \, \le_{\frak K}$-represent $M$, so if
$M \le^\otimes N \in K^\otimes_{\lambda^+}$, (hence $M \le^*_{\lambda^+}
N \in K^\otimes_{\lambda^+} = K^{\text{nice}}_{\lambda^+}$) 
and $a \in N$, \ortp$_{{\frak K}^{\text{nice}}_{\lambda^+}}(a,M,N) \in 
{\Cal S}^{\text{cs}}_{\lambda^+}(M)$, we let
$\alpha(a,N,\bar M) = \text{ Min}\{\alpha:\text{\ortp}(a,M_\alpha,N) \in
{\Cal S}^{\text{bs}}(M_\alpha)$ and for every $\beta \in (\alpha,\lambda^+)$,
\ortp$(a,M_\beta,N) \in {\Cal S}^{\text{bs}}(M_\beta)$ 
is a non-forking extension of \ortp$(a,M_\alpha,N)\}$.  
\nl
Now
\mr
\item "{$(a)$}"  $\alpha(a,N,\bar M)$ is well defined for $a,N$ as above \nl
[Why?  By Defintion \scite{600-1.9} + \scite{600-rg.6}(1)]
\sn
\item "{$(b)$}"  if $a_\ell,N_\ell$ are above for $\ell = 1,2$ and
$\alpha(a_1,N_1,\bar M) = \alpha(a_2,N_2,\bar M)$ call it $\alpha$
and \ortp$_{\frak s}(a_1,M_\alpha,N) = \text{ \ortp}_{\frak s}
(a_2,M_\alpha,N_2)$ then
{\roster
\itemitem{ $(*)$ }  for $\beta < \lambda^+$ we have
\ortp$_{\frak s}(a_1,M_\beta,N_1) = \text{ \ortp}_{\frak s}(a_1,M_\beta,N_2) \in 
{\Cal S}^{\text{bs}}(M_\beta)$ \nl
[Why?  By $\nonfork{}{}_{}$-calculus when $\beta \ge \alpha$ by
monotonicity if $\beta \le \alpha$]
\endroster}
\item "{$(c)$}"  if $a_\ell,N_\ell$ are as above for $\ell=1,2$ and
$(*)$ above holds then
{\roster
\itemitem{ $(**)$ }  \ortp$_{{\frak K}^\otimes_{\lambda^+}}
(a_1,M,N_1) = \text{\rm \ortp}_{{\frak K}^\otimes_{\lambda^+}}(a_2,M,N_2)$ \nl
[Why?  Use \scite{600-ne.4}(3) or by part (6) below].
\endroster}
\ermn
As $\alpha < \lambda \Rightarrow 
|{\Cal S}^{\text{bs}}_{\frak s}(M_\alpha)| \le \lambda$ (by the
stability Axiom (D)(d) for ${\frak s}$), clearly
$|{\Cal S}^{\text{cs}}_{\lambda^+}(M)| \le \dsize \sum_{\alpha < \lambda^+}
|{\Cal S}^{\text{bs}}(M_\alpha)| \le \lambda^+ = \|M\|$ as required.

The reader may ask why do we not just quote the parallel result from
\S2: The answer is that the equality of types there is ``a formal, not
the true one".  The crux of the matter is that we prove locality (in
clause (c) above).
\bn
\ub{Axiom (E)(a)}:

By \scite{600-1.6} - \scite{600-1.9}.
\bn
\ub{Axiom (E)(b); monotonicity}:

Follows by Axiom (E)(b) for ${\frak s}$ and the definition.
\bn
\ub{Axiom (E)(c)}; local character:  

By \scite{600-1.13}(5) or direct by translating it to the ${\frak s}$-case.
\bn
\ub{Axiom (E)(d); (transitivity)}:

By \scite{600-1.13}(4).
\bn
\ub{Axiom (E)(e); uniqueness}:

By \scite{600-ne.4}(3) or by part (6) below.
\bn
\ub{Axiom (E)(f); symmetry}:

So assume $M_0 \le^*_{\lambda^+} M_1 \le^*_{\lambda^+} M_2$ are from
$K^\otimes_{\lambda^+}$ and for $\ell=1,2$ we have $a_\ell \in M_\ell$, 
\ortp$_{{\frak K}^{\text{nice}}_{\lambda^+}}(a_\ell,
M_0,M_\ell) \in {\Cal S}^{\text{cs}}_{\lambda^+}(M_0)$ as witnessed by
$p_\ell \in {\Cal S}^{\text{bs}}_{\frak s}(N^*_\ell),N^*_\ell \in 
K_\lambda,N^*_\ell \le_{\frak K} M_0$ and
\ortp$_{{\frak K}^\otimes_{\lambda^+}}(a_2,
M_1,M_2)$ does not fork (in
the sense of $\nonfork{}{}_{\lambda^+}$) over $M_0$ (note that
$M_0,M_1,M_2$ here stand for $M_0,M_1,M'_3$ in (E)(f)(i) from
Definition \scite{600-1.1}).  As we know the monotonicity \wilog \, 
$M_1 <^+_{\lambda^+} M_2$.
We can finish by \scite{600-ne.4}(4) (and Axiom (E)(e) 
for ${\frak K}_\lambda$).

In more details, we can find $N_0,N_1,N_2$ such that: $N_\ell
\le_{\frak K} M_\ell$ and $N_\ell \in K_\lambda$ for $\ell = 0,1,2$
and $N^*_1 \cup N^*_2 \subseteq N_0 \le_{\frak K} N_1 \le_{\frak K}
N_2$ and $a_1 \in N_1,a_2 \in N_2$
and $N_2$ is $(\lambda,*)$-brimmed over $N_1$ hence over $N_0$,
and $(\forall N \in K_\lambda)[N_0 \le_{\frak K} N 
\le_{\frak K} M_0 \rightarrow
(\exists M \in K_\lambda)(M \le_{\frak K} M_2 \and
\text{ NF}_\lambda(N_0,N,N_2,M))]$. \nl
By Axiom (E)(f) for $({\frak K},{\Cal S}^{\text{bs}},
\nonfork{}{}_{\lambda})$ we can
find $N'$ such that $N_0 \le_{\frak K} N' \le_{\frak K} N_2$ and $a_2
\in N'$ and \ortp$_{\frak s}(a_1,N',N_2)$ does not fork over
$N_0$.  Now we can find $f'_0,M'_1$ such that $M_0 \le^+_{\lambda^+}
M'_1,f'_0$ is a $\le_{\frak K}$-embedding of $N'$ into $M'_1$ and
$(\forall N \in K_\lambda)[N_0 \le_{\frak K} N \le_{\frak K} M_0
\rightarrow (\exists M \in K_\lambda)(M \le_{\frak K} M'_1 \and
\text{\rm NF}_\lambda(N_0,N,f'_0(N'),M))]$.  Next we can find $f''_0,M'_2$ such
that $M'_1 <^+_{\lambda^+} M'_2,f''_0 \supseteq f'_0,f''_0$ is a
$\le_{\frak K}$-embedding of $N_2$ into $M'_2$ and $(\forall N \in
K_\lambda)[N_0 \le_{\frak K} N \le_{\frak K} M_0 \rightarrow (\exists
M \in K_\lambda)(M \le_{\frak K} M'_2 \and \text{\rm NF}_\lambda
(N_0,N,f''_0(N_2),M)]$. \nl
Lastly, by \scite{600-ne.4}(4) there is an isomorphism $f$ from $M_2$
onto $M'_2$ over $M_0$ extending $f''_0$.  Now $f^{-1}(M'_1)$ is a
model as required.
\bn
\ub{Axiom (E)(g); extension existence}:  

Assume $M_0 \le^*_{\lambda^+} M_1$ are from 
$K^{\text{nice}}_{\lambda^+},
p \in {\Cal S}^{\text{cs}}_{\lambda^+}(M_0)$, hence
there is $N_0 \le_{\frak K} M_0,N_0 \in K_\lambda$ such that $(\forall N \in
K_\lambda)(N_0 \le_{\frak K} N <_{\frak K} M_0 \rightarrow p \restriction N$
does not fork over $N_0)$.  By \scite{600-ne.3}(1A) there are $M_2 \in K^\otimes
_{\lambda^+}$ and $a \in M_2$ such that $M_1 \le^*_{\lambda^+} M_2$ and
\ortp$_{{\frak K}^{\text{nice}}_{\lambda^+}}(a,M_1,M_2) 
\in {\Cal S}^{\text{cs}}_{\lambda^+}(M_1)$
is witnessed by $p \restriction N_0$ and by part (6) we have
\ortp$_{{\frak K}^{\text{nice}}_{\lambda^+}}(a,M_0,M_2)=p$.
Checking the definition of does not
fork, i.e., $\nonfork{}{}_{\lambda^+}$ we are done.
\bn
\ub{Axiom (E)(h), (continuity)}:

By \scite{600-1.13}(6).
\bn
\ub{Axiom (E)(i)}:

It follows from the rest by \scite{600-1.15}.
\nl
6) So assume $M \le^*_{\lambda^+} M_\ell,a_\ell \in M_\ell \backslash
M$ for $\ell=1,2$ and $N \le_{\frak K} M \wedge N \in K_\lambda
\Rightarrow \text{\rm \ortp}_{\frak K}(a_1,N,M_1) = \text{\rm \ortp}_{\frak
K}(a_2,N,M_2)$.  By \scite{600-ne.3}(1) there are $M^+_1,M^+_2 \in
K^{\text{nice}}_{\lambda^+}$ such that $M_\ell <^+_{\lambda^+}
M^+_\ell$ for $\ell=1,2$.  By \scite{600-ne.4}(2),(3) there is an
isomorphism $f$ from $M^+_1$ onto $M^+_2$ over $M$ which maps $a_1$ to
$a_2$.  This clearly suffices.  \hfill$\square_{\scite{600-rg.7}}$
\enddemo
\newpage

\head {\S9 Final conclusions} \endhead  \resetall \sectno=9
 \spuriousreset
\bigskip

We now show that we have actually solved our specific test questions
about categoricity and few models.  First we deal with good $\lambda$-frames.
\proclaim{\stag{600-fc.1} Main Lemma}  Assume
\mr
\item "{$(a)$}"  $2^\lambda < 2^{\lambda^+} < 2^{\lambda^{++}} < \ldots <
2^{\lambda^{+n}}$, and $n \ge 2$ and ${\text{\rm WDmId\/}}
(\lambda^{+ \ell})$ is not
$\lambda^{+ \ell +1}$-saturated (normal ideal on $\lambda^{+ \ell}$)
for $\ell = 1,\dotsc,n-1$
\sn
\item "{$(b)$}"  ${\frak s} = ({\frak K},{\Cal S}^{\text{bs}},
\nonfork{}{}_{})$ is a good $\lambda$-frame
\sn
\item "{$(c)$}"  $\dot I(\lambda^{+ \ell},{\frak K}(\lambda^+$-saturated)) 
$< 2^{\lambda^{+ \ell}}$ for $\ell = 2,\dotsc,n$. 
\ermn
\ub{Then}
\mr
\item "{$(\alpha)$}"   $K$ has a member of cardinality $\lambda^{+n+1}$
\sn
\item "{$(\beta)$}"  for $\ell < n$ there is a good 
$\lambda^{+ \ell}$-frame ${\frak s}_\ell = 
({\frak K}^\ell,{\Cal S}^{\text{bs}}_{{\frak s}_\ell},
\nonfork{}{}_{{\frak s}_\ell})$ such that
$K^\ell_{\lambda^{+ \ell}} \subseteq K_{\lambda^{+ \ell}}$ and
$\le_{{\frak K}^\ell} \subseteq \le_{\frak K}$
\sn
\item "{$(\gamma)$}"  ${\frak s}_0 = {\frak s}$ and if $\ell < m < n$
then $K^\ell_{\lambda^{+m}} \supseteq K^m_{\lambda^{+ m}} \and
\le_{{\frak K}^\ell} \restriction K^m \supseteq \le_{{\frak K}^m}$.
\endroster
\endproclaim
\bigskip

\demo{Proof}  We prove this by induction on $n$.

For $n=m+1 \ge 2$, by the induction hypothesis for $\ell = 0,\dotsc,m-1$, there
is a frame ${\frak s}_\ell = ({\frak K}^\ell,\nonfork{}{}_{{\frak s}_\ell},
{\Cal S}^{\text{bs}}_{{\frak s}_\ell})$ which is $\lambda^{+ \ell}$-good and 
$K_{{\frak s}_\ell} \subseteq K^{\frak s}_{\lambda^{+ \ell}}$ and
$\le_{{\frak K}^\ell} \subseteq \le_{\frak K} \restriction {\frak K}^\ell$.
By \scite{600-nu.6} and clause (c) of the assumption we know that ${\frak
s}$ has density for $K^{3,\text{uq}}_{\frak s}$.
Now \wilog \, $K^{m-1}$ is categorical in  $\lambda^{+(m-1)}$ 
(by \scite{600-1.16} really necessary only for $\ell = 0$) and by
Observation \scite{600-nu.13.1} we get the assumption \scite{600-nf.0} of \S6
hence the results of \S6, \S7, \S8 apply.
Now apply \scite{600-rg.7} to $({\frak K}^{m-1},
{\Cal S}^{\text{bs}}_{{\frak s}_{m-1}},\nonfork{}{}_{{\frak
s}_{m-1}})$ and get a $\lambda^{+ m}$-frame
${\frak s}_m$ as required in clause $(\beta)$.  By \scite{600-4a.12} we have
$K^m_{\lambda^{+m+1}} \ne \emptyset$ which is clause $(\alpha)$ in the
conclusion.  Clause $(\beta)$ has already been proved and clause
$(\gamma)$ should be clear.     \hfill$\square_{\scite{600-fc.1}}$
\enddemo
\bn
Second (this fulfills the aim of \cite{Sh:576}).
\proclaim{\stag{600-fc.2} Theorem}  1) Assume $2^{\lambda^{+ \ell}} < 
2^{\lambda^{+(\ell+1)}}$ for $\ell = 0,\dotsc,n-1$ and the normal
ideal ${\text{\rm WDmId\/}}(\lambda^{+ \ell})$ is not 
$\lambda^{+\ell+1}$-saturated for $\ell = 1,\dotsc,n-1$.

If ${\frak K}$ is an abstract elementary class with 
${\text{\rm LS\/}}({\frak K}) \le \lambda$  which is 
categorical in $\lambda,\lambda^+$ and $1 \le \dot I
(\lambda^{+2},K)$ and $\dot I(\lambda^{+m},
{\frak K}) < 2^{\lambda^{+m}}$ 
for $m \in [2,n)$ (or just $\dot I(\lambda^{+m},
{\frak K}(\lambda^+$-saturated)) $< 2^{\lambda^{+m}}$), 
\ub{then} ${\frak K}_{\lambda^{+n}} \ne \emptyset$ (and there are
${\frak s}_\ell(\ell < n)$ as in $(\gamma)$ of \scite{600-fc.1}).
\nl
2) We can omit the assumption ``not $\lambda^{+ \ell +1}$-saturated",
if we strengthen $\dot I(\lambda^{+m},{\frak K}) <
2^{\lambda^{+m}}$ to $\dot I(\lambda^{+m},{\frak K}) <
\mu_{\text{unif}}(\lambda^{+m},2^{\lambda^{+(m-1)}})$, see
\marginbf{!!}{\cprefix{88r}.\scite{88r-0.wD}}(3). 
\endproclaim
\bigskip

\demo{Proof}  1) By \scite{600-Ex.4} and \scite{600-fc.1}.
\nl
2) See \cite{Sh:838}.  \hfill$\square_{\scite{600-fc.2}}$
\enddemo
\bn
Next we fulfill an aim of \chaptercite{88r}.
\proclaim{\stag{600-fc.3} Theorem}  1) Assume $2^{\aleph_\ell} < 
2^{\aleph_{(\ell +1)}}$ for $\ell = 0,\dotsc,n-1$ and $n \ge 2$
and {\rm WDmId}$(\lambda^{+ \ell})$ is not $\lambda^{+ \ell +1}$-saturated
for $\ell = 1,\dotsc,n-1$.

If ${\frak K}$ is an abstract 
elementary class which is ${\text{\rm PC\/}}_{\aleph_0}$
and $1 \le \dot I(\aleph_1,{\frak K}) < 2^{\aleph_1}$ and $\dot I
(\aleph_\ell,
{\frak K}) < 2^{\aleph_\ell}$ for $\ell =2,\dotsc,n$, \ub{then} ${\frak K}$
has a model of cardinality $\aleph_{n+1}$ (and there are ${\frak
s}_\ell(\ell < n)$ as in \scite{600-fc.2}.
\nl
2) We can omit the assumption ``not $\lambda^{+\ell+1}$-saturated" if
we strengthen ``$\dot I(\aleph_\ell,{\frak K}) < 2^{\aleph_\ell}$"
to ``$\dot I(\aleph_\ell,K) < 
\mu_{\text{unif}}(\aleph_\ell,2^{\aleph_{\ell -1}})$.
\endproclaim
\bigskip

\remark{Remark}  Compared with Theorem \scite{600-fc.2} our gains are no
assumption on $\dot I(\lambda,K)$ and weaker assumption on $\dot
I(\lambda^+,K)$, i.e., $< 2^{\aleph_1}$ (and $\ge 1$) rather than
$=1$.  The price is $\lambda = \aleph^+_0$ and being PC$_{\aleph_0}$. 
\endremark
\bigskip

\demo{Proof}  1) By \scite{600-Ex.1} and \scite{600-fc.1}.
\nl
2) See \cite{Sh:838}.  \hfill$\square_{\scite{600-fc.3}}$
\enddemo
\bn
Lastly, we fulfill an aim of \cite{Sh:48}.
\proclaim{\stag{600-fc.4} Theorem}  1) Assume $2^{\aleph_\ell} < 
2^{\aleph_{\ell +1}}$ for $\ell \le n -1$
and ${\text{\rm WDmId\/}}
(\lambda^{+ \ell})$ is not $\lambda^{+ \ell +1}$-saturated
for $\ell = 1,\dotsc,n-1,\psi \in \Bbb L_{\omega_1,\omega}(\bold Q),
\dot I(\aleph_1,\psi) \ge 1$ and $\dot I(\aleph_\ell,\psi) < 
2^{\aleph_\ell}$ for
$\ell = 1,\dotsc,n$.  \ub{Then} $\psi$ has a model in $\aleph_{n+1}$
and there are ${\frak s}_1,\dotsc,{\frak s}_{n-1}$ as in \scite{600-fc.3}
for $K = \text{\rm Mod}_\psi$ and appropriate $\le_{\frak K}$.
\nl
2) We can omit the assumption ``not $\lambda^{+ \ell +1}$-saturated"
if we strengthen ``$\dot I(\aleph_\ell,{\frak K}) < 2^{\aleph_\ell}$
for $\ell=2,\dotsc,n$" to ``$\dot I(\aleph_\ell,{\frak K}) <
\mu_{\text{unif}}(\aleph_\ell,2^{\aleph_{\ell-1}})$". 
\endproclaim
\bigskip

\demo{Proof}  1) By \scite{600-Ex.1A} mainly clauses (c)-(d) and
\scite{600-fc.1}.  Note that this time in \scite{600-fc.1} we use the \nl
$\dot I(\lambda^{+ \ell},{\frak K}
(\lambda^+$-saturated)) $< 2^{\lambda^{+ \ell}}$.
\nl
2) See \cite{Sh:838}.   \hfill$\square_{\scite{600-fc.4}}$
\enddemo

\nocite{ignore-this-bibtex-warning} 
\newpage
    
REFERENCES.  
\bibliographystyle{lit-plain}
\bibliography{lista,listb,listx,listf,liste}

\def\germ{\frak} \def\scr{\cal} \ifx\documentclass\undefinedcs
  \def\bf{\fam\bffam\tenbf}\def\rm{\fam0\tenrm}\fi 
  \def\defaultdefine#1#2{\expandafter\ifx\csname#1\endcsname\relax
  \expandafter\def\csname#1\endcsname{#2}\fi} \defaultdefine{Bbb}{\bf}
  \defaultdefine{frak}{\bf} \defaultdefine{=}{\B} 
  \defaultdefine{mathfrak}{\frak} \defaultdefine{mathbb}{\bf}
  \defaultdefine{mathcal}{\cal}
  \defaultdefine{beth}{BETH}\defaultdefine{cal}{\bf} \def\bbfI{{\Bbb I}}
  \def\mbox{\hbox} \def\text{\hbox} \def\om{\omega} \def\Cal#1{{\bf #1}}
  \def\pcf{pcf} \defaultdefine{cf}{cf} \defaultdefine{reals}{{\Bbb R}}
  \defaultdefine{real}{{\Bbb R}} \def\restriction{{|}} \def\club{CLUB}
  \def\w{\omega} \def\exist{\exists} \def\se{{\germ se}} \def\bb{{\bf b}}
  \def\equivalence{\equiv} \let\lt< \let\gt>
  \def\implies{\Rightarrow}\def\mathfrak{\bf}\def\germ{\frak} \def\scr{\cal}
  \ifx\documentclass\undefinedcs
  \def\bf{\fam\bffam\tenbf}\def\rm{\fam0\tenrm}\fi 
  \def\defaultdefine#1#2{\expandafter\ifx\csname#1\endcsname\relax
  \expandafter\def\csname#1\endcsname{#2}\fi} \defaultdefine{Bbb}{\bf}
  \defaultdefine{frak}{\bf} \defaultdefine{=}{\B} 
  \defaultdefine{mathfrak}{\frak} \defaultdefine{mathbb}{\bf}
  \defaultdefine{mathcal}{\cal}
  \defaultdefine{beth}{BETH}\defaultdefine{cal}{\bf} \def\bbfI{{\Bbb I}}
  \def\mbox{\hbox} \def\text{\hbox} \def\om{\omega} \def\Cal#1{{\bf #1}}
  \def\pcf{pcf} \defaultdefine{cf}{cf} \defaultdefine{reals}{{\Bbb R}}
  \defaultdefine{real}{{\Bbb R}} \def\restriction{{|}} \def\club{CLUB}
  \def\w{\omega} \def\exist{\exists} \def\se{{\germ se}} \def\bb{{\bf b}}
  \def\equivalence{\equiv} \let\lt< \let\gt>
\begin{thebibliography}{HuSh 342}
\makeatletter \renewcommand{\@biblabel}[1]{[#1]} \makeatother
\def\eprintfootnotetext{References of the form {\tt math.XX/$\cdots$}
 refer to {\tt arXiv.org} }
\ifx\documentstyle\undefinedcontrolsequence
   \def\anyfootnote{\footnote{*}}
   \else\def\anyfootnote{\footnote}\fi
\def\eprintfn{\ifEprint\anyfootnote{\eprintfootnotetext}\fi\Eprintfalse }
\newif\ifEprint  \Eprinttrue

\bibitem[HuSh 342]{HuSh:342}Ehud Hrushovski and Saharon Shelah.
\newblock {A dichotomy theorem for regular types}.
\newblock {\em {Annals of Pure and Applied Logic}}, {\bf 45}:157--169, 1989.

\bibitem[JrSh 875]{JrSh:875}Adi Jarden and Saharon Shelah.
\newblock {Good frames minus stability}.
\newblock {\em Preprint}.

\bibitem[Ke71]{Ke71}Jerome~H. Keisler.
\newblock {\em {Model theory for infinitary logic. Logic with countable
  conjunctions and finite quantifiers}}, volume~62 of {\em {Studies in Logic
  and the Foundations of Mathematics}}.
\newblock North--Holland Publishing Co., Amsterdam--London, 1971.

\bibitem[Mw85a]{Mw85a}Johann~A. Makowsky.
\newblock {Abstract embedding relations}.
\newblock In J.~Barwise and S.~Feferman, editors, {\em Model-Theoretic Logics},
  pages 747--791. Springer-Verlag, 1985.

\bibitem[Sh 482]{Sh:482}Saharon Shelah.
\newblock {Compactness in ZFC of the Quantifier on ``Complete embedding of
  BA's''}.
\newblock In {\em {Non structure theory, Ch XI}}, accepted. {Oxford University
  Press}.

\bibitem[Sh 838]{Sh:838}Saharon Shelah.
\newblock {Non-structure in $\lambda^{++}$ using instances of WGCH}.
\newblock {\em Preprint}.

\bibitem[Sh 842]{Sh:842}Saharon Shelah.
\newblock {Solvability and Categoricity spectrum of a.e.c. with amalgamation}.
\newblock {\em Preprint}.

\bibitem[Sh 839]{Sh:839}Saharon Shelah.
\newblock {Stable Frames and weight}.
\newblock {\em Preprint}.

\bibitem[Sh 849]{Sh:849}Saharon Shelah.
\newblock {Weak forms of good frames}.
\newblock {\em Preprint}.

\bibitem[Sh 48]{Sh:48}Saharon Shelah.
\newblock {Categoricity in $\aleph _{1}$ of sentences in $L_{\omega
  _{1},\omega}(Q)$}.
\newblock {\em {Israel Journal of Mathematics}}, {\bf 20}:127--148, 1975.

\bibitem[Sh 87a]{Sh:87a}Saharon Shelah.
\newblock {Classification theory for nonelementary classes, I. The number of
  uncountable models of $\psi \in L_{\omega _{1},\omega }$. Part A}.
\newblock {\em {Israel Journal of Mathematics}}, {\bf 46}:212--240, 1983.

\bibitem[Sh 87b]{Sh:87b}Saharon Shelah.
\newblock {Classification theory for nonelementary classes, I. The number of
  uncountable models of $\psi \in L_{\omega _{1},\omega }$. Part B}.
\newblock {\em {Israel Journal of Mathematics}}, {\bf 46}:241--273, 1983.

\bibitem[Sh 88]{Sh:88}Saharon Shelah.
\newblock {Appendix: on stationary sets (in ``Classification of nonelementary
  classes. II. Abstract elementary classes'')}.
\newblock In {\em {Classification theory (Chicago, IL, 1985)}}, volume 1292 of
  {\em {Lecture Notes in Mathematics}}, pages 419--497. {Springer, Berlin},
  1987.
\newblock {Proceedings of the USA--Israel Conference on Classification Theory,
  Chicago, December 1985; ed. Baldwin, J.T.}

\bibitem[Sh 300]{Sh:300}Saharon Shelah.
\newblock {Universal classes}.
\newblock In {\em {Classification theory (Chicago, IL, 1985)}}, volume 1292 of
  {\em {Lecture Notes in Mathematics}}, pages 264--418. {Springer, Berlin},
  1987.
\newblock {Proceedings of the USA--Israel Conference on Classification Theory,
  Chicago, December 1985; ed. Baldwin, J.T.}

\bibitem[Sh:c]{Sh:c}Saharon Shelah.
\newblock {\em {Classification theory and the number of nonisomorphic models}},
  volume~92 of {\em {Studies in Logic and the Foundations of Mathematics}}.
\newblock {North-Holland Publishing Co., Amsterdam, xxxiv+705 pp}, 1990.

\bibitem[Sh 576]{Sh:576}Saharon Shelah.
\newblock {Categoricity of an abstract elementary class in two successive
  cardinals}.
\newblock {\em {Israel Journal of Mathematics}}, {\bf 126}:29--128, 2001.
\newblock math.LO/9805146.

\bibitem[Sh 603]{Sh:603}Saharon Shelah.
\newblock {Few non minimal types and non-structure}.
\newblock In {\em {Proceedings of the 11 International Congress of Logic,
  Methodology and Philosophy of Science, Krakow August'99; In the Scope of
  Logic, Methodology and Philosophy of Science}}, volume~1, pages 29--53.
  Kluwer Academic Publishers, 2002.
\newblock math.LO/9906023.

\bibitem[Sh:E46]{Sh:E46}{Shelah, Saharon}.
\newblock {Categoricity of an abstract elementary class in two successive
  cardinals, revisited}.

\end{thebibliography}

\enddocument